\documentclass[11pt]{article}
\usepackage{epic,latexsym,amssymb}
\usepackage{color}
\usepackage{tikz}

\newtheorem{thm}{Theorem}
\newtheorem{lem}[thm]{Lemma}

\newtheorem{claim}{Claim}

\newtheorem{ob}{Observation}
\newtheorem{prob}{Problem}
\newtheorem{prop}{Proposition}
\newtheorem{conj}{Conjecture}

\newcommand{\qed}{$\Box$}

\newcommand{\defic}{{\rm def}}

\newcommand{\odd}{{\rm odd}}
\newcommand{\even}{{\rm even}}

\newcommand{\gt}{\gamma_t}

\newcommand{\barG}{{\overline{G}}}

\newcommand{\smallqed}{{\tiny ($\Box$)}}

\newcommand{\TR}[1]{\mbox{$\tau(#1)$}}

\def\ZZ{Z\!\!\!Z}

\newcommand{\ClaimX}[1]{{\bf Claim #1:}}
\newcommand{\ProofClaimX}[1]{{\em Proof of Claim #1:}}

\newcommand{\cF}{{\cal F}}

\newcommand{\cH}{{\cal H}}

\newcommand{\oc}{{\rm oc}}

\def \nH {n_{_H}}
\def \mH {m_{_H}}

\let\oldenumerate\enumerate
\renewcommand{\enumerate}{
  \oldenumerate
  \setlength{\itemsep}{0pt}
  \setlength{\parskip}{0pt}
  \setlength{\parsep}{0pt}
}

\def\vertex(#1){\put(#1){\circle*{2}}}
\def\vertexo(#1){\put(#1){\circle{2}}}
\def\vert(#1){\put(#1){\circle*{1.5}}}
\def\verto(#1){\put(#1){\circle{1.5}}}
\def\lab(#1)#2{\put(#1){\makebox(0,0)[c]{#2}}}
\newcommand{\proof}{\noindent\textbf{Proof. }}

\newcommand{\2}{ \vspace{0.2cm} }
\newcommand{\1}{ \vspace{0.1cm} }
\newcommand{\4}{ \vspace{0.05cm} }
\newcommand{\5}{ \vspace{0.05cm} }

\newcommand{\Vx}[1]{V(#1)}

\newcommand{\hyperedgefour}[4]{
	\hyperedgebig #1#2#3;
	\hyperedgebig #2#3#4;
	\hyperedgebig #3#4#1;
	\hyperedgebig #4#1#2;
}

\newcommand{\hyperedgebig}[9]{
	\pgfmathsetmacro\Done{sqrt((#4-#1)^2+(#5-#2)^2)}
	\pgfmathsetmacro\angleone{(#2>#5)*(180+asin((#4-#1)/ \Done)-asin((#1-#4)/ \Done))+asin((#1-#4)/ \Done)+asin((#3-#6)/\Done)}
	\pgfmathsetmacro\angleone{\angleone-360*(\angleone>0)-360*(\angleone>360)}
	\pgfmathsetmacro\Dtwo{sqrt((#7-#4)^2+(#8-#5)^2)}
	\pgfmathsetmacro\angletwo{(#5>#8)*(180+asin((#7-#4)/ \Dtwo)-asin((#4-#7)/ \Dtwo))+asin((#4-#7)/ \Dtwo)+asin((#6-#9)/\Dtwo)}
	\pgfmathsetmacro\angletwo{\angletwo-360*(\angletwo>0)-360*(\angletwo>360)}
	\draw ([shift=(\angletwo:#6)] #4,#5)--([shift=(\angletwo:#9)]#7,#8);
	\draw (#4,#5)+(\angleone:#6) arc(\angleone:\angletwo+360*(\angletwo<\angleone):#6);
}

\newcommand{\hyperedgetwo}[6]{
	\pgfmathsetmacro\Done{sqrt((#4-#1)^2+(#5-#2)^2)}
	\pgfmathsetmacro\angleone{(#2>#5)*(180+asin((#4-#1)/ \Done)-asin((#1-#4)/ \Done))+asin((#1-#4)/ \Done)+asin((#3-#6)/\Done)}
	\pgfmathsetmacro\angleone{\angleone-360*(\angleone>0)-360*(\angleone>360)}
	\draw ([shift=(\angleone:#3)] #1,#2)--([shift=(\angleone:#6)]#4,#5);
	\pgfmathsetmacro\Dtwo{sqrt((#1-#4)^2+(#2-#5)^2)}
	\pgfmathsetmacro\angletwo{(#5>#2)*(180+asin((#1-#4)/ \Dtwo)-asin((#4-#1)/ \Dtwo))+asin((#4-#1)/ \Dtwo)+asin((#6-#3)/\Dtwo)}
	\pgfmathsetmacro\angletwo{\angletwo-360*(\angletwo>0)-360*(\angletwo>360)}
	\draw ([shift=(\angletwo:#6)] #4,#5)--([shift=(\angletwo:#3)]#1,#2);
	\draw (#1,#2)+(\angletwo:#3) arc(\angletwo:\angleone+360*(\angleone<\angletwo):#3);
	\draw (#4,#5)+(\angleone:#6) arc(\angleone:\angletwo+360*(\angletwo<\angleone):#6);
}

\begin{document}

\title{Transversals in Uniform Linear Hypergraphs}

\author{Michael A. Henning${}^{1,}$\thanks{Research
supported in part by the South African
National Research Foundation and the University of Johannesburg} \, and Anders Yeo${}^{1,2}$\\
\\
${}^1$Department of Pure and Applied Mathematics\\
University of Johannesburg \\
Auckland Park, 2006 South Africa \\
\small\tt Email:  mahenning@uj.ac.za \\
\\
${}^2$Department of Mathematics and Computer Science \\
University of Southern Denmark \\
Campusvej 55, 5230 Odense M, Denmark \\
\small\tt Email:  andersyeo@gmail.com
}

\date{}
\maketitle

\begin{abstract}
The transversal number $\tau(H)$ of a hypergraph $H$ is the minimum number of vertices that intersect every edge of $H$. A linear hypergraph is one in which every two distinct edges intersect in at most one vertex. A $k$-uniform hypergraph has all edges of size~$k$. Very few papers give bounds on the transversal number for linear hypergraphs, even though these appear in many applications, as it seems difficult to utilise the linearity in the known techniques. This paper is one of the first that give strong non-trivial bounds on the transversal number for linear hypergraphs, which is better than for non-linear hypergraphs. It is known that $\tau(H) \le (n + m)/(k+1)$ holds for all $k$-uniform, linear hypergraphs $H$ when $k \in \{2,3\}$ or when $k \ge 4$ and the maximum degree of $H$ is at most two. It has been conjectured (at several conference talks) that $\tau(H) \le (n+m)/(k+1)$ holds for all $k$-uniform, linear hypergraphs $H$. We disprove the conjecture for large $k$, and show that the best possible constant $c_k$ in the bound $\tau(H) \le c_k (n+m)$ has order $\ln(k)/k$ for both linear (which we show in this paper) and non-linear hypergraphs. We show that for those $k$ where the conjecture holds, it is tight for a large number of densities if there exists an affine plane $AG(2,k)$ of order~$k \ge 2$. We raise the problem to find the smallest value, $k_{\min}$, of~$k$ for which the conjecture fails. We prove a general result, which when applied to a projective plane of order $331$ shows that $k_{\min} \le 166$. Even though the conjecture fails for large $k$, our main result is that it still holds for $k=4$, implying that $k_{\min} \ge 5$. The case $k=4$ is much more difficult than the cases $k \in \{2,3\}$, as the conjecture does not hold for general (non-linear) hypergraphs when $k=4$. Key to our proof is the completely new technique of the deficiency of a hypergraph introduced in this paper.
\end{abstract}


\section{Introduction}

In this paper we continue the study of transversals in hypergraphs.
Hypergraphs are systems of sets which are conceived as natural
extensions of graphs.  A \emph{hypergraph} $H = (V,E)$ is a finite
set $V = V(H)$ of elements, called \emph{vertices}, together with a
finite multiset $E = E(H)$ of subsets of $V$, called
\emph{hyperedges} or simply \emph{edges}. The \emph{order} of $H$ is
$n(H) = |V|$ and the \emph{size} of $H$ is $m(H) = |E|$.
A $k$-\emph{edge} in $H$ is an edge of size~$k$. The hypergraph $H$
is said to be $k$-\emph{uniform} if every edge of $H$ is a $k$-edge.
Every (simple) graph is a $2$-uniform hypergraph. Thus graphs are
special hypergraphs.  For $i \ge 2$, we denote the number
of edges in $H$ of size $i$ by $e_i(H)$. The \emph{degree} of a
vertex $v$ in $H$, denoted by $d_H(v)$, is the number of edges of $H$ which contain
$v$. A \emph{degree}-$k$ \emph{vertex} is a vertex of degree~$k$. The minimum
and maximum degrees among the vertices of $H$ is denoted by
$\delta(H)$ and $\Delta(H)$, respectively.

Two vertices $x$ and $y$ of $H$ are \emph{adjacent} if there is an
edge $e$ of $H$ such that $\{x,y\}\subseteq e$. The \emph{neighborhood} of a vertex $v$ in $H$, denoted $N_H(v)$ or simply $N(v)$ if $H$ is clear from the context, is the set of all vertices different from $v$ that are adjacent to $v$. A vertex in $N(v)$ is a \emph{neighbor} of $v$.
The \emph{neighborhood} of a set $S$ of vertices of $H$ is the set $N_H(S) = \cup_{v \in S} N_H(v)$, and the \emph{boundary} of $S$ is the set $\partial_H(S) = N_H(S) \setminus S$. Thus, $\partial_H(S)$ consists of all vertices of $H$ not in $S$ that have a neighbor in $S$. If $H$ is clear from context, we simply write $N(S)$ and $\partial(S)$ rather than $N_H(S)$ and $\partial_H(S)$.
Two vertices $x$ and $y$ of $H$ are \emph{connected} if there is a
sequence $x=v_0,v_1,v_2\ldots,v_k=y$ of vertices of $H$ in which $v_{i-1}$ is adjacent to $v_i$ for $i=1,2,\ldots,k$. A \emph{connected hypergraph} is a hypergraph in which every pair of vertices are connected. A maximal connected subhypergraph of $H$ is a \emph{component} of $H$.  Thus, no edge in $H$ contains vertices from different components. A component of $H$ isomorphic to a hypergraph $F$ we call an $F$-\emph{component} of $H$.

A subset $T$ of vertices in a hypergraph $H$ is a \emph{transversal} (also called \emph{vertex cover} or \emph{hitting set} in many papers) if $T$ has a nonempty intersection with every edge of $H$. The \emph{transversal number} $\TR{H}$ of $H$ is the minimum size of a transversal in $H$. A transversal of size~$\TR{H}$ is called a $\TR{H}$-transversal.  Transversals in hypergraphs are well studied in the
literature (see, for
example,~\cite{BuHeTu12,ChMc,CoHeSl79,DoHe14,HeLo12,HeLo14,HeYe08,HeYe10,HeYe13,HeYe13b,KoMuVe14,LaCh90,ThYe07,Tu90}).

A hypergraph $H$ is called an \emph{intersecting hypergraph} if every two distinct edges of $H$ have a non-empty intersection, while $H$ is
called a \emph{linear hypergraph} if every two distinct edges of $H$ intersect in at most one vertex. We say that two edges in $H$ \emph{overlap} if they intersect in at least two vertices. A linear hypergraph therefore has no overlapping edges. Linear hypergraphs are well studied in the literature (see, for example,~\cite{BlGr16,CaFu00,FrFu86,Fu14,KhNa11,Ko10,MeTy97,Naik82,To98}), as are uniform hypergraphs (see, for example,~\cite{BuHeTu12,BuHeTuYe14,CaFu00,Fu14,GuSe17,HeYe13b,HeYe16,LaCh90,LoWa13,MeTy97,Naik82,To98}).

A set $S$ of vertices in a hypergraph $H$ is \emph{independent} (also called strongly independent in the literature) if no two vertices in $S$ belong to a common edge. Independence in hypergraphs is well studied in the literature (see, for example,~\cite{BaMoSa13,HeYe17,KoMuVe14}, for recent papers on this topic).

Given a hypergraph $H$ and subsets $X,Y \subseteq V(H)$ of vertices, we let $H(X,Y)$ denote the hypergraph obtained by deleting all vertices in $X \cup Y$ from $H$ and removing all (hyper)edges containing vertices from $X$ and removing the vertices in $Y$ from any remaining edges. If $Y = \emptyset$, we simply denote $H(X,Y)$ by $H - X$; that is, $H - X$ denotes that hypergraph obtained from $H$ by removing the vertices $X$ from $H$, removing all edges that intersect $X$ and removing all resulting isolated vertices, if any. Further, if $X = \{x\}$, we simply write $H - x$ rather than $H - X$.
When we use the definition $H(X,Y)$ we furthermore assume that
no edges of size zero are created. That is, there is no edge $e \in
E(H)$ such that $\Vx{e} \subseteq Y \setminus X$. In this case we
note that if add $X$ to any $\tau(H(X,Y))$-set, then we get a
transversal of $H$, implying that $\tau(H) \le |X| + \tau(H(X,Y))$.
We will often use this fact throughout the paper.


In geometry, a \emph{finite affine plane} is a system of points and lines that satisfy the following rules:
(R1) Any two distinct points lie on a unique line. (R2) Each line has at least two points. (R3) Given a point and a line, there is a unique line which contains the point and is parallel to the line, where two lines are called parallel if they are equal or disjoint. (R4) There exist three non-collinear points (points not on a single line).
A finite affine plane $AG(2,q)$ of \emph{order~$q \ge 2$} is a collection of $q^2$ points and $q^2 + q$ lines, such that each line contains $q$ points and each point is contained in $q+1$ lines.

A \emph{total dominating set}, also called a TD-set, of a graph $G$ with no isolated vertex is a set $S$ of vertices of $G$ such that every vertex is adjacent to a vertex in $S$. The \emph{total domination number} of $G$, denoted by $\gt(G)$, is the minimum cardinality of a TD-set of $G$. Total domination in graphs is now well studied in graph theory.
The literature on this subject has been surveyed and detailed in a recent book on this topic that can be found in~\cite{MHAYbookTD}. A survey of total domination in graphs can be found in~\cite{He09}.
We use the standard notation $[k] = \{1,2,\ldots,k\}$.

\section{Motivation}

In this paper we study upper bounds on the transversal number of a linear uniform hypergraph in terms of its order and size, and establish connections between linear uniform hypergraphs with large transversal number and affine planes and projective planes. An affine plane of order $k$ exists if and only if a projective plane of order $k$ exists. Jamison~\cite{Ja77} and Brouwer and Schrijver~\cite{BrSc76} proved that the minimum cardinality of a subset of $AG(2,k)$ which intersects all hyperplanes is $2(k-1) + 1$.

\begin{thm}{\rm (\cite{BrSc76,Ja77})}
 \label{t:Affine_plane}
If there exists a finite affine plane $AG(2,k)$ of order~$k \ge 2$, then the minimum cardinality of a subset of $AG(2,k)$ which intersects all hyperplanes is $2(k-1) + 1$.
\end{thm}

One of the most fundamental results in combinatorics is the result due to Bose~\cite{Bo38} that there are $k-1$ mutually orthogonal Latin squares if and only if there is an affine plane $AG(2,k)$. The prime power conjecture for affine and projective planes, states that there is an affine plane of order~$k$ exists only if $k$ is a prime power. Veblen and Wedderburn~\cite{VeWe07} and Lam~\cite{La91} established the existence of affine and projective planes of small orders.

\begin{figure}[htb]
\begin{center}
\begin{tikzpicture}
[thick,scale=0.8,
	vertex/.style={circle,draw,fill=white,inner sep=0pt,minimum size=1.8mm},
	vertexblack/.style={circle,draw,fill=black,inner sep=0pt,minimum size=1.8mm, text=white},
	vertexlabel/.style={circle,draw=white,inner sep=0pt,minimum size=0mm}]
\def \d {1.2}
\def \r {.3}
\coordinate (A) at (-0.5,0);
\draw (A) +(-0.5,-2) node[rectangle]{(a) $F_4=AG(2,2)$} ;
\coordinate (B) at (5,0);
\draw (B) +(-0.5,-2) node[rectangle]{(b) $F_3$} ;
\draw (A)
{
	node [vertexblack] at +(-1*\d,-1*\d) {}
	node [vertexblack] at +(-1*\d,0*\d) {}
	node [vertexblack] at +(0*\d,-1*\d) {}
	node [vertexblack] at +(0*\d,0*\d) {}
};
\foreach \i in {-1,0} {
	\foreach \j in {0} {
		\draw[rotate=\j*90](A) ++(-1*\d,\i*\d) ++(180:\r)--++(0:1*\d+\r*2.25);
	}
}
\foreach \i in {-1,0} {
	\foreach \j in {-1} {
		\draw[rotate=\j*90](A) ++(-1*\d,\i*\d) ++(180:-0.85)--++(0:1*\d+\r*2.25);
	}
}
\foreach \j in {0} {
	\draw[rotate=\j*90](A) ++(-1*\d,-1*\d)++(-135:\r)--++(45:1*\d*2 + \r*0.1);
}
\foreach \j in {1} {
	\draw[rotate=\j*90](A) ++(-1*\d,0*\d)++(-135:\r)--++(45:1*\d*2 + \r*0.15);
}
\draw (B)
{
	node [vertexblack] at +(-1*\d,-1*\d) {}
	node [vertexblack] at +(-1*\d,0*\d) {}
	node [vertexblack] at +(0*\d,-1*\d) {}
};
\foreach \i in {-1} {
	\foreach \j in {0} {
		\draw[rotate=\j*90](B) ++(-1*\d,\i*\d) ++(180:\r)--++(0:1*\d+\r*2.25);
	}
}
\foreach \i in {-1} {
	\foreach \j in {-1} {
		\draw[rotate=\j*90](B) ++(-1*\d,\i*\d) ++(180:-0.85)--++(0:1*\d+\r*2.25);
	}
}
\foreach \j in {1} {
	\draw[rotate=\j*90](B) ++(-1*\d,0*\d)++(-135:\r)--++(45:1*\d*2 + \r*0.15);
}
\end{tikzpicture}
\end{center}
\vskip -0.75 cm
\caption{The affine plane $AG(2,2)$ and the hypergraph $F_3$}
 \label{AG22}
\end{figure}
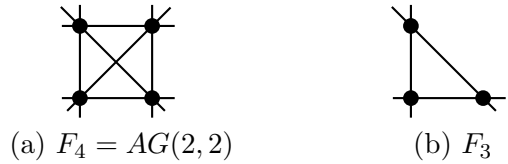

The affine plane $AG(2,2)$ of dimension~$2$ and order~$2$ is illustrated in Figure~\ref{AG22}(a). This is equivalent to a linear $2$-uniform $3$-regular hypergraph $F_4$ (the complete graph $K_4$ on four vertices), where the lines of $AG(2,2)$ correspond to the $2$-edges of $F_4$. Deleting any vertex from $F_4$ yields the $2$-uniform $2$-regular hypergraph $F_3$ (the complete graph $K_3$)  illustrated in  Figure~\ref{AG22}(b), where the lines correspond to the $2$-edges of $F_3$.

\vskip -0.75 cm

\begin{figure}[htb]
\begin{center}
\begin{tikzpicture}
[thick,scale=0.8,
	vertex/.style={circle,draw,fill=white,inner sep=0pt,minimum size=1.8mm},
	vertexblack/.style={circle,draw,fill=black,inner sep=0pt,minimum size=1.8mm, text=white},
	vertexlabel/.style={circle,draw=white,inner sep=0pt,minimum size=0mm}]
\def \d {1.2}
\def \r {.3}
\coordinate (A) at (0,0);
\draw (A) +(0,-3) node[rectangle]{(a) $F_9=AG(2,3)$} ;
\coordinate (B) at (6,0);
\draw (B) +(0,-3) node[rectangle]{(b) $F_8$} ;
\draw (A)
{
	node [vertexblack] at +(-1*\d,-1*\d) {}
	node [vertexblack] at +(-1*\d,0*\d) {}
	node [vertexblack] at +(-1*\d,1*\d) {}
	node [vertexblack] at +(0*\d,-1*\d) {}
	node [vertexblack] at +(0*\d,0*\d) {}
	node [vertexblack] at +(0*\d,1*\d) {}
	node [vertexblack] at +(1*\d,-1*\d) {}
	node [vertexblack] at +(1*\d,0*\d) {}
	node [vertexblack] at +(1*\d,1*\d) {}
};
\foreach \i in {0,1,2,3} {
	\draw[rotate=90*\i](A) ++(0*\d,-1*\d) ++(-45:\r)-- ++(135:1.414214*\d*1.3+\r)
		arc (225:45:1.060660*\d)--++(-45:1.414214*\d*.8+\r);
}
\foreach \i in {-1,0,1} {
	\foreach \j in {0,1} {
		\draw[rotate=\j*90](A) ++(-1*\d,\i*\d) ++(180:\r)--++(0:2*\d+\r*2);
	}
}
\foreach \j in {0,1} {
	\draw[rotate=\j*90](A) ++(-1*\d,-1*\d)++(-135:\r)--++(45:1.414214*\d*2 + \r*2);
}
\draw (B)
{
	node [vertexblack] at +(-1*\d,-1*\d) {}
	node [vertexblack] at +(-1*\d,0*\d) {}
	node [vertexblack] at +(-1*\d,1*\d) {}
	node [vertexblack] at +(0*\d,-1*\d) {}
	node [vertexblack] at +(0*\d,0*\d) {}
	node [vertexblack] at +(0*\d,1*\d) {}
	node [vertexblack] at +(1*\d,-1*\d) {}
	node [vertexblack] at +(1*\d,0*\d) {}
};
\foreach \i in {1,2,3} {
	\draw[rotate=90*\i](B) ++(0*\d,-1*\d) ++(-45:\r)-- ++(135:1.414214*\d*1.3+\r)
		arc (225:45:1.060660*\d)--++(-45:1.414214*\d*.8+\r);
}
\foreach \i in {-1,0} {
	\foreach \j in {0,1} {
		\draw[rotate=\j*-90](B) ++(-1*\d,\i*\d) ++(180:\r)--++(0:2*\d+\r*2);
	}
}
\foreach \j in {1} {
	\draw[rotate=\j*90](B) ++(-1*\d,-1*\d)++(-135:\r)--++(45:1.414214*\d*2 + \r*2);
}
\end{tikzpicture}
\end{center}
\vskip -.6cm \caption{The affine plane $AG(2,3)$ and the hypergraph $F_8$}
 \label{AG23}
\end{figure}
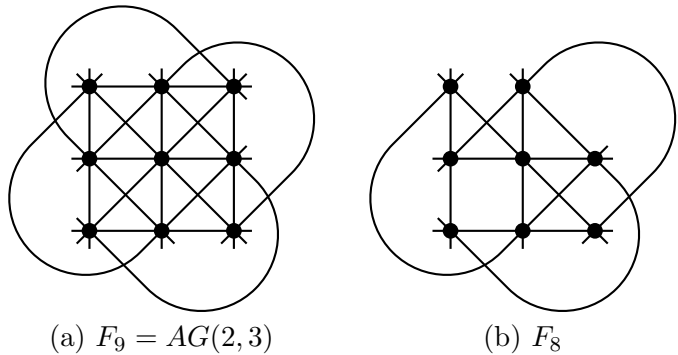

The affine plane $AG(2,3)$ of dimension~$2$ and order~$3$ is illustrated in Figure~\ref{AG23}(a). This is equivalent to a linear $3$-uniform $4$-regular hypergraph $F_9$ on nine vertices, where the lines of $AG(2,3)$ correspond to the $3$-edges of $F_9$. Deleting any vertex from $F_9$ yields the hypergraph $F_8$ (see Figure~\ref{AG23}(b), where the lines correspond to the $3$-edges of $F_8$). Further, deleting any vertex from $F_8$ yields a linear $3$-uniform hypergraph $F_7$. It is known (see, for example,~\cite{HeYe08,HeYe13}) that $\tau(F_7) = 3$, $\tau(F_8) = 4$ and $\tau(F_9) = 5$.

As a consequence of results due to Tuza~\cite{Tu90}, Chv\'{a}tal and McDiarmid~\cite{ChMc}, Henning and Yeo~\cite{HeYe08,HeYe13}, and Dorfling and Henning~\cite{DoHe14}, we have the following result. For $k \ge 2$, let $E_k$ denote the $k$-uniform hypergraph on $k$ vertices with exactly one edge.

\begin{thm}{\rm (\cite{ChMc,DoHe14,HeYe08,Tu90})} \label{t:known23}
For $k \in \{2,3\}$ if $H$ is a $k$-uniform linear connected hypergraph, then
\[
\tau(H) \le \frac{n + m}{k+1}
\]
with equality if and only if $H$ consists of a single edge or $H$ is obtained from the affine plane $AG(2,k)$ of order~$k$ by deleting one or two vertices; that is, $H = E_k$ or $H \in \{F_{k^2-2},F_{k^2-1}\}$.
\end{thm}

\begin{figure}[htb]
\begin{center}
\begin{tikzpicture}[scale=0.40]

\begin{scope}[xshift=0cm,yshift=0cm]
\fill (0,0) circle (0.2cm); \fill (4,0) circle (0.2cm); \fill (8,0)
circle (0.2cm); \fill (2,3.464) circle (0.2cm); \fill (6,3.464)
circle (0.2cm); \fill (4,6.928) circle (0.2cm); \fill (4,2.3093)
circle (0.2cm); \hyperedgetwo{0}{0}{0.5}{8}{0}{0.6};
\hyperedgetwo{4}{6.928}{0.5}{0}{0}{0.6};
\hyperedgetwo{4}{6.928}{0.6}{8}{0}{0.5};
\hyperedgetwo{0}{0}{0.4}{6}{3.464}{0.45};
\hyperedgetwo{4}{6.928}{0.4}{4}{0}{0.45};
\hyperedgetwo{8}{0}{0.4}{2}{3.464}{0.45}; \draw (4,2.3093) circle
(1.9523cm); \draw (4,2.3093) circle (2.6523cm);
\end{scope}

\end{tikzpicture}
\end{center}
\vskip -0.3 cm \caption{The Fano plane $F$} \label{f:Fano}
\end{figure}
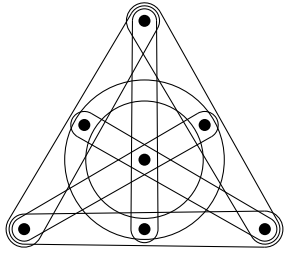

In fact, if $k \in \{2,3\}$ then $\tau(H) \le (n + m)/(k+1)$ holds for all $k$-uniform hypergraphs, even when they are not linear. The case when $k \ge 4$ is more complex as then the bound $\tau(H) \le (n + m)/(k+1)$ does not hold in general.
As remarked in~\cite{DoHe14}, if $k \ge 4$ and $H$ is a $k$-uniform hypergraph that is \textbf{not} linear, then is not always true that $\tau(H) \le (n + m)/(k+1)$ holds, as may be seen for example by taking $k = 4$ and letting $\overline{F}$ be the complement of the Fano plane $F$, where the Fano plane is shown in Figure~\ref{f:Fano} and where its complement $\overline{F}$ is the hypergraph on the same vertex set $V(F)$ and where $e$ is a hyperedge in the complement if and only if $V(F) \setminus e$ is a hyperedge in $F$.

\begin{figure}[htb]
\begin{center}
\tikzstyle{vertexX}=[circle,draw, fill=black!100, minimum size=8pt, scale=0.5, inner sep=0.1pt]
\begin{tikzpicture}[scale=0.25]
 \node (a1) at (1.0,4.0) [vertexX] {};
\node (a2) at (1.0,1.0) [vertexX] {};
\node (a3) at (4.0,1.0) [vertexX] {};
\node (a4) at (2.5,-2.0) {(a) $L_2$};
\draw[color=black!90, thick,rounded corners=4pt] (1.9999433500235204,4.010644094312785) arc (0.6098731973394038:176.4869591045931:1.0); 
\draw[color=black!90, thick,rounded corners=4pt] (0.0018791224802877649,0.9387242800182674) arc (183.5130408954069:359.3901268026606:1.0); 
\draw[color=black!90, thick,rounded corners=4pt] (1.9999433500235204,4.010644094312785) -- (1.9999433500235204,0.989355905687214);
\draw[color=black!90, thick,rounded corners=4pt] (0.0018791224802877649,0.9387242800182674) -- (0.0,2.5) -- (0.0018791224802878759,4.061275719981733);
\draw[color=black!90, thick,rounded corners=4pt] (0.992814465624922,4.899971315151625) arc (90.45745018578945:275.52564338823333:0.9); 
\draw[color=black!90, thick,rounded corners=4pt] (4.899971315151625,4.007185534375078) arc (0.4574501857894546:89.54254981421053:0.9); 
\draw[color=black!90, thick,rounded corners=4pt] (3.1041821183800913,1.0866621195795443) arc (174.47435661176667:359.54254981421053:0.9); 
\draw[color=black!90, thick,rounded corners=4pt] (0.992814465624922,4.899971315151625) -- (4.007185534375078,4.899971315151625);
\draw[color=black!90, thick,rounded corners=4pt] (4.899971315151625,4.007185534375078) -- (4.899971315151625,0.9928144656249211);
\draw[color=black!90, thick,rounded corners=4pt] (3.1041821183800913,1.0866621195795443) -- (3.363603896932107,3.363603896932107) -- (1.086662119579544,3.1041821183800913);
\draw[color=black!90, thick,rounded corners=4pt] (4.004601528639943,0.2000132339006001) arc (-89.6704379697801:87.9060302070261:0.8); 
\draw[color=black!90, thick,rounded corners=4pt] (0.9707691742230907,1.7994657959064915) arc (92.09396979297392:269.6704379697801:0.8); 
\draw[color=black!90, thick,rounded corners=4pt] (4.004601528639943,0.2000132339006001) -- (0.9953984713600572,0.2000132339006001);
\draw[color=black!90, thick,rounded corners=4pt] (0.9707691742230907,1.7994657959064915) -- (2.5,1.8) -- (4.029230825776909,1.7994657959064915);
 \end{tikzpicture}
 \hspace*{2cm}
\tikzstyle{vertexX}=[circle,draw, fill=black!100, minimum size=8pt, scale=0.5, inner sep=0.1pt]
\begin{tikzpicture}[scale=0.25]
 \node (a1) at (1.0,10.0) [vertexX] {};
\node (a2) at (1.0,7.0) [vertexX] {};
\node (a3) at (1.0,4.0) [vertexX] {};
\node (a4) at (1.0,1.0) [vertexX] {};
\node (a5) at (4.0,7.0) [vertexX] {};
\node (a6) at (4.0,4.0) [vertexX] {};
\node (a7) at (4.0,1.0) [vertexX] {};
\node (a8) at (7.0,4.0) [vertexX] {};
\node (a9) at (7.0,1.0) [vertexX] {};
\node (a10) at (10.0,1.0) [vertexX] {};
\node (a11) at (4.5,-2.0) {(b) $L_4$};
\draw[color=black!90, thick,rounded corners=4pt] (1.9999404314142877,4.010914834996838) arc (0.6253863972442666:11.376967106857592:1.0); 
\draw[color=black!90, thick,rounded corners=4pt] (0.007769901048958916,1.1244163604017057) arc (172.85294740996648:359.3903323424423:1.0); 
\draw[color=black!90, thick,rounded corners=4pt] (1.9999984825500685,9.99825790426205) arc (-0.09981478378023212:180.02328137121071:1.0); 
\draw[color=black!90, thick,rounded corners=4pt] (1.9999404314142877,4.010914834996838) -- (1.9999433882011346,0.9893594928299944);
\draw[color=black!90, thick,rounded corners=4pt] (0.007769901048958916,1.1244163604017057) -- (0.0,9.999593663429177);
\draw[color=black!90, thick,rounded corners=4pt] (1.9999984825500685,9.99825790426205) -- (1.9803505536273511,4.197263255581333);
\draw[color=black!90, thick,rounded corners=4pt] (0.9921968514540919,10.919966907487856) arc (90.48597047746559:275.5712207156533:0.92); 
\draw[color=black!90, thick,rounded corners=4pt] (4.919999943253066,10.000323132102855) arc (0.020124028368887714:89.51402952253443:0.92); 
\draw[color=black!90, thick,rounded corners=4pt] (3.0804763018984502,1.0296001457707613) arc (178.15624329998536:359.97987597163115:0.92); 
\draw[color=black!90, thick,rounded corners=4pt] (0.9921968514540919,10.919966907487856) -- (4.007803148545908,10.919966907487856);
\draw[color=black!90, thick,rounded corners=4pt] (4.919999943253066,10.000323132102855) -- (4.919999943253066,0.9996768678971465);
\draw[color=black!90, thick,rounded corners=4pt] (3.0804763018984502,1.0296001457707613) -- (3.35191187357377,9.34702084230403) -- (1.0893163521310583,9.084345813507085);
\draw[color=black!90, thick,rounded corners=4pt] (0.9992511434117376,7.83999966619863) arc (90.05107896148978:272.54430286561285:0.84); 
\draw[color=black!90, thick,rounded corners=4pt] (7.83999966619863,7.000748856588262) arc (0.0510789614897762:89.9489210385102:0.84); 
\draw[color=black!90, thick,rounded corners=4pt] (6.160828076141371,1.0372891701088522) arc (177.45569713438715:359.9489210385102:0.84); 
\draw[color=black!90, thick,rounded corners=4pt] (0.9992511434117376,7.83999966619863) -- (7.000748856588263,7.83999966619863);
\draw[color=black!90, thick,rounded corners=4pt] (7.83999966619863,7.000748856588262) -- (7.83999966619863,0.9992511434117369);
\draw[color=black!90, thick,rounded corners=4pt] (6.160828076141371,1.0372891701088522) -- (6.4060303038033,6.4060303038033) -- (1.0372891701088518,6.160828076141371);
\draw[color=black!90, thick,rounded corners=4pt] (0.9998488636842061,4.759999984972246) arc (90.01139404353056:271.5027608704118:0.76); 
\draw[color=black!90, thick,rounded corners=4pt] (10.759990584449437,4.003783060692465) arc (0.28520303478334164:89.98860595646948:0.76); 
\draw[color=black!90, thick,rounded corners=4pt] (9.24252962708978,1.061956712011301) arc (175.32394489703594:359.71479696521664:0.76); 
\draw[color=black!90, thick,rounded corners=4pt] (0.9998488636842061,4.759999984972246) -- (10.000151136315793,4.759999984972246);
\draw[color=black!90, thick,rounded corners=4pt] (10.759990584449437,4.003783060692465) -- (10.759990584449437,0.9962169393075345);
\draw[color=black!90, thick,rounded corners=4pt] (9.24252962708978,1.061956712011301) -- (9.461451338116294,3.463748809993609) -- (1.0199310897135436,3.2402613925415986);
\draw[color=black!90, thick,rounded corners=4pt] (10.00009704136856,0.32000000692428476) arc (-89.99182343988883:89.97264907846073:0.68); 
\draw[color=black!90, thick,rounded corners=4pt] (0.997529417175017,1.679995511912031) arc (90.20816805916866:269.99182343988883:0.68); 
\draw[color=black!90, thick,rounded corners=4pt] (10.00009704136856,0.32000000692428476) -- (0.9999029586314389,0.32000000692428476);
\draw[color=black!90, thick,rounded corners=4pt] (0.997529417175017,1.679995511912031) -- (4.001070993069664,1.6799991565978922) -- (10.000324607259005,1.679999922522148);
 \end{tikzpicture}
 \hspace*{2cm}
\tikzstyle{vertexX}=[circle,draw, fill=black!100, minimum size=8pt, scale=0.5, inner sep=0.1pt]
\begin{tikzpicture}[scale=0.25]
 \node (a1) at (1.0,16.0) [vertexX] {};
\node (a2) at (1.0,13.0) [vertexX] {};
\node (a3) at (1.0,10.0) [vertexX] {};
\node (a4) at (1.0,7.0) [vertexX] {};
\node (a5) at (1.0,4.0) [vertexX] {};
\node (a6) at (1.0,1.0) [vertexX] {};
\node (a7) at (4.0,13.0) [vertexX] {};
\node (a8) at (4.0,10.0) [vertexX] {};
\node (a9) at (4.0,7.0) [vertexX] {};
\node (a10) at (4.0,4.0) [vertexX] {};
\node (a11) at (4.0,1.0) [vertexX] {};
\node (a12) at (7.0,10.0) [vertexX] {};
\node (a13) at (7.0,7.0) [vertexX] {};
\node (a14) at (7.0,4.0) [vertexX] {};
\node (a15) at (7.0,1.0) [vertexX] {};
\node (a16) at (10.0,7.0) [vertexX] {};
\node (a17) at (10.0,4.0) [vertexX] {};
\node (a18) at (10.0,1.0) [vertexX] {};
\node (a19) at (13.0,4.0) [vertexX] {};
\node (a20) at (13.0,1.0) [vertexX] {};
\node (a21) at (16.0,1.0) [vertexX] {};
\node (a22) at (7.5,-2.0) {(c) $L_6$};
\draw[color=black!90, thick,rounded corners=4pt] (-0.09999999062676213,0.9998563994313584) arc (180.0074797332202:359.98554873389094:1.1); 
\draw[color=black!90, thick,rounded corners=4pt] (2.099894574058543,16.015229115423313) arc (0.7932653786981447:179.99252065999448:1.1); 
\draw[color=black!90, thick,rounded corners=4pt] (-0.09999999062676213,0.9998563994313584) -- (-0.09999999062774756,16.000143593019462);
\draw[color=black!90, thick,rounded corners=4pt] (2.099894574058543,16.015229115423313) -- (2.099973175188146,12.992317951672371) -- (2.099999965011186,0.9997225556107884);
\draw[color=black!90, thick,rounded corners=4pt] (0.9881289814302683,17.029931589435975) arc (90.66036341868451:276.0927890800722:1.03); 
\draw[color=black!90, thick,rounded corners=4pt] (5.029999994074748,16.000110480856918) arc (0.006145715369447746:89.33963658131552:1.03); 
\draw[color=black!90, thick,rounded corners=4pt] (2.9702338531102006,1.0219472713505904) arc (178.7790473823493:359.9938542846306:1.03); 
\draw[color=black!90, thick,rounded corners=4pt] (0.9881289814302683,17.029931589435975) -- (4.011871018569732,17.029931589435975);
\draw[color=black!90, thick,rounded corners=4pt] (5.029999994074748,16.000110480856918) -- (5.029999994074748,0.9998895191430828);
\draw[color=black!90, thick,rounded corners=4pt] (2.9702338531102006,1.0219472713505904) -- (3.27563501887429,15.2677463730928) -- (1.1093230966580256,14.975818150650433);
\draw[color=black!90, thick,rounded corners=4pt] (0.9987325604876904,13.95999916333145) arc (90.07564474577082:272.9252626237068:0.96); 
\draw[color=black!90, thick,rounded corners=4pt] (7.959999986250258,13.000162479244937) arc (0.009697265664371457:89.92435525422917:0.96); 
\draw[color=black!90, thick,rounded corners=4pt] (6.040292490305415,1.0236959033551747) arc (178.58561129673907:359.9903027343356:0.96); 
\draw[color=black!90, thick,rounded corners=4pt] (0.9987325604876904,13.95999916333145) -- (7.00126743951231,13.95999916333145);
\draw[color=black!90, thick,rounded corners=4pt] (7.959999986250258,13.000162479244937) -- (7.959999986250258,0.9998375207550625);
\draw[color=black!90, thick,rounded corners=4pt] (6.040292490305415,1.0236959033551747) -- (6.32151612792989,12.32083902104085) -- (1.0489919549921478,12.041250925243707);
\draw[color=black!90, thick,rounded corners=4pt] (0.9997167498863231,10.889999954926614) arc (90.01823487228567:271.76384343879636:0.89); 
\draw[color=black!90, thick,rounded corners=4pt] (10.889999954926614,10.000283250113677) arc (0.01823487228567444:89.98176512771431:0.89); 
\draw[color=black!90, thick,rounded corners=4pt] (9.11042169700351,1.0273942117587707) arc (178.23615656120364:359.9817651277143:0.89); 
\draw[color=black!90, thick,rounded corners=4pt] (0.9997167498863231,10.889999954926614) -- (10.000283250113677,10.889999954926614);
\draw[color=black!90, thick,rounded corners=4pt] (10.889999954926614,10.000283250113677) -- (10.889999954926614,0.9997167498863229);
\draw[color=black!90, thick,rounded corners=4pt] (9.11042169700351,1.0273942117587707) -- (9.370674964743973,9.370674964743973) -- (1.02739421175877,9.11042169700351);
\draw[color=black!90, thick,rounded corners=4pt] (0.9999133139014178,7.819999995418) arc (90.00605700926887:271.1995742341041:0.82); 
\draw[color=black!90, thick,rounded corners=4pt] (13.819999717326453,7.000680870425277) arc (0.04757439786133322:89.99394299073106:0.82); 
\draw[color=black!90, thick,rounded corners=4pt] (12.180767029993872,1.0354590024526573) arc (177.52160378794775:359.95242560213865:0.82); 
\draw[color=black!90, thick,rounded corners=4pt] (0.9999133139014178,7.819999995418) -- (13.000086686098584,7.819999995418);
\draw[color=black!90, thick,rounded corners=4pt] (13.819999717326453,7.000680870425277) -- (13.819999717326453,0.999319129574723);
\draw[color=black!90, thick,rounded corners=4pt] (12.180767029993872,1.0354590024526573) -- (12.41999110089589,6.420353834688749) -- (1.0171666922065539,6.1801797119620145);
\draw[color=black!90, thick,rounded corners=4pt] (0.999968851494472,4.7499999993531805) arc (90.00237957054055:270.87197684650056:0.75); 
\draw[color=black!90, thick,rounded corners=4pt] (16.74999137943768,4.003595937869617) arc (0.27471047027624174:89.99762042945946:0.75); 
\draw[color=black!90, thick,rounded corners=4pt] (15.252432172122976,1.0603518245227674) arc (175.38446991750504:359.7252895297238:0.75); 
\draw[color=black!90, thick,rounded corners=4pt] (0.999968851494472,4.7499999993531805) -- (16.000031148505528,4.7499999993531805);
\draw[color=black!90, thick,rounded corners=4pt] (16.74999137943768,4.003595937869617) -- (16.74999137943768,0.9964040621303832);
\draw[color=black!90, thick,rounded corners=4pt] (15.252432172122976,1.0603518245227674) -- (15.468546786708977,3.470795425111757) -- (1.0114137096218823,3.25008685354058);
\draw[color=black!90, thick,rounded corners=4pt] (16.000021061175943,0.3200000003261566) arc (-89.99822541692174:89.99654338536715:0.68); 
\draw[color=black!90, thick,rounded corners=4pt] (0.9975298247557323,1.6799955133927449) arc (90.20813371680264:269.9982254299629:0.68); 
\draw[color=black!90, thick,rounded corners=4pt] (16.000021061175943,0.3200000003261566) -- (0.9999789389788338,0.32000000032615195);
\draw[color=black!90, thick,rounded corners=4pt] (0.9975298247557323,1.6799955133927449) -- (4.001219428243573,1.6799989066129142) -- (16.00004102392827,1.6799999987625274);
 \end{tikzpicture}
\end{center}
\vskip -0.75 cm \caption{The hypergraphs $L_2$, $L_4$ and $L_6$.} \label{intersect}
\end{figure}
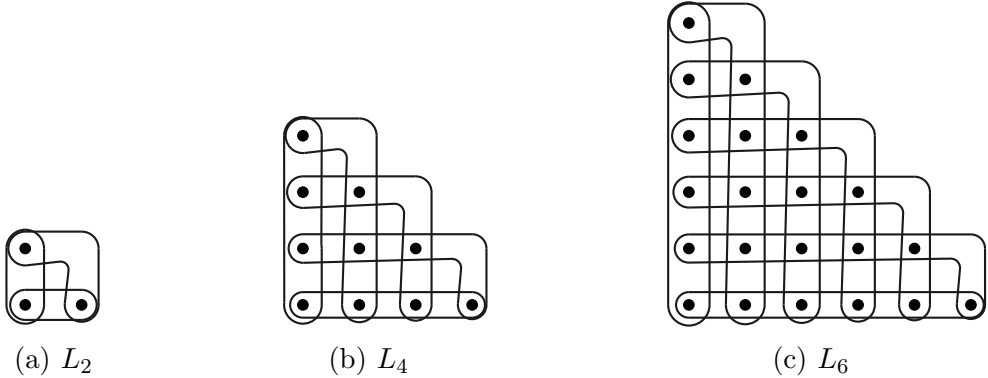

A natural question is whether the upper bound in Theorem~\ref{t:known23} holds for larger values of~$k$ if we impose a linearity constraint. This is indeed the case for all $k \ge 2$ when the maximum degree of a $k$-uniform linear hypergraph is at most~$2$. As shown in~\cite{DoHe14} for $k \ge 2$, there is a unique $k$-uniform, $2$-regular, linear intersecting hypergraph, which we call $L_k$. The hypergraphs $L_2$, $L_4$ and $L_6$, for example, are illustrated in Figure~\ref{intersect}.


\begin{thm}{\rm (\cite{DoHe14})}
 \label{mainthm}
For $k \ge 2$, if $H$ be a $k$-uniform, linear, connected hypergraph on $n$ vertices with $m$ edges satisfying $\Delta(H) \le 2$, then
\[
\tau(H) \le \frac{n + m}{k+1}
\]
with equality if and only if $H$ consists of a single edge or $k$ is even and $H = L_k$.
\end{thm}

The authors conjectured in~\cite{HeYe16} that the bound $\tau(H) \le (n + m)/(k+1)$ also holds for all $k$-uniform linear hypergraphs in the case when $k = 4$ (which we prove in this paper). 

\begin{conj}{\rm (\cite{HeYe16})}
 \label{t:conj4unif}
For $k = 4$ if $H$ is a $k$-uniform linear hypergraph on $n$ vertices with $m$ edges, then
\[
\tau(H) \le \frac{n + m}{k+1}.
\]
\end{conj}

More generally, the following conjecture was posed by the authors at several international conferences in recent years, including in a principal talk at the 8th Slovenian Conference on Graph Theory held in Kranjska Gora in June 2015.

\begin{conj}
 \label{ConjF}
For $k \ge 2$, if $H$ is a $k$-uniform, linear hypergraph on $n$ vertices with $m$ edges, then
\[
\tau(H) \le \frac{n+m}{k+1}.
\]
\end{conj}

\section{Main Results}

We have three immediate aims. Our first aim is to establish the existence of $k$-uniform, linear hypergraph $H$ on $n$ vertices with $m$ edges satisfying $\tau(H) > (n + m)/(k+1)$. More precisely, we shall prove the following general result, a proof of which is given in Section~\ref{S:mainX}. We remark that this result was motivated by an important paper by Alon et al.~\cite{Alon02} on transversal numbers for hypergraphs arising in geometry where each edge in a projective plane is shrunk randomly.

\begin{thm}
\label{mainX}
For $k \ge 2$, let $H$ be an arbitrary $2k$-uniform $2k$-regular hypergraph on $n$ vertices and let
\[
0<c<\frac{1}{\ln(4)} \approx 0.72134
\]

\noindent
be arbitrary. Let $H'$ be the $k$-uniform hypergraph obtained from $H$ by letting $V(H')=V(H)$ and for each edge in $H$ pick $k$ vertices at random from the edge and add these as a hyperedge to $H'$.  We note that $|V(H)|=|E(H)|=|V(H')|=|E(H')|=n$.
If
\[
 5 c \ln(k) \ln(n)  <  k^{1-c\ln(4)},
\]
then
\[
\tau(H') \le \left( \frac{c \, \ln(k)}{k} \right) n
\]

\noindent
occurs with probability strictly less than~$1$.  Therefore, if $5 c \ln(k) \ln(n)  <  k^{1-c\ln(4)}$,
then there exists a $k$-uniform hypergraph, $H^*$ (where the edges of $H^*$ are subsets of the edges of $H$), such that $n=|V(H^*)|=|E(H^*)|$ satisfying
\[
\tau(H^*) >  \left( \frac{c \, \ln(k)}{k} \right) n.
\]
\end{thm}

Applying Theorem~\ref{mainX} when $k=2754$ to a projective plane of order $p=5507$ (which is prime), gives us a counterexample to Conjecture~\ref{ConjF}. We can in fact prove that from a projective plane of order $p=331$ one can construct a $166$-uniform linear hypergraph $H$ on $n$ vertices with $m$ edges satisfying $\tau(H) > (n+m)/(k+1)$. This result, which we discuss in Remark~3 in Section~\ref{S:mainX}, we state formally as follows.

\begin{thm}
 \label{t:thm166}
For $k = 166$, there exist $k$-uniform, linear hypergraph on $n$ vertices with $m$ edges satisfying
\[
\tau(H) > \frac{n+m}{k+1}.
\]
\end{thm}

Our second aim is to show that for those values of $k$ for which Conjecture~\ref{ConjF} holds, the bound is tight for a large number of densities if there exists an affine plane $AG(2,k)$ of order~$k \ge 2$. We shall prove the following result, a proof of which is given in Section~\ref{S:mainY}.

\begin{thm}
\label{t:mainY}
If the affine plane $AG(2,k)$ of order~$k \ge 2$ exists and $\tau(H) \le (n + m)/(k+1)$ holds for all $k$-uniform linear hypergraphs, then this upper bound is tight for a number of linear $k$-uniform hypergraphs $H$ with a wide variety of average degrees.
\end{thm}

Our third aim is to prove Conjecture~\ref{t:conj4unif}. More precisely, we shall prove the following result, a proof of which is given in Section~\ref{S:mainZ}.

\begin{thm}
 \label{t:mainZ}
For $k = 4$ if $H$ is a $k$-uniform linear hypergraph, then
\[
\tau(H) \le \frac{n + m}{k+1}.
\]
\end{thm}

By Theorem~\ref{t:known23}, Conjecture~\ref{ConjF} is true for $k \in \{2,3\}$. Theorem~\ref{t:mainZ} implies that the conjecture is also true for $k = 4$. However as remarked earlier, Conjecture~\ref{ConjF} fails for large $k$. We therefore pose the following new problem that we have yet to settle.

\begin{prob}{\rm (\cite{HeYe13})}
 \label{prob1}
Determine the smallest value, $k_{\min}$, of~$k$ for which there exists a $k$-uniform, linear hypergraph $H$ on $n$ vertices with $m$ edges satisfying
\[
\tau(H) > \frac{n+m}{k+1}.
\]
\end{prob}

As an immediate consequence of Theorem~\ref{t:thm166} and Theorem~\ref{t:mainZ}, we have the following result.

\begin{thm}
 \label{t:main4}
$5 \le k_{\min} \le 166$.
\end{thm}

We proceed as follows. In Section~\ref{S:mainX}, we give a proof of Theorem~\ref{mainX}. In Section~\ref{S:mainY}, we give a proof of Theorem~\ref{t:mainY}.  In Section~\ref{S:applications}, we present several applications of Theorem~\ref{t:mainZ}. In Section~\ref{S:key}, we define fifteen special hypergraphs and we introduce the concept of the deficiency of a hypergraph. Thereafter, we prove a key result, namely Theorem~\ref{t:thm_key}, about the deficiency of a hypergraph that will enable us to deduce our main result. Finally in Section~\ref{S:mainZ}, we present a proof of Theorem~\ref{t:mainZ}.

\section{Proof of Theorem~\ref{mainX}}
\label{S:mainX}

In this section, we give a proof of Theorem~\ref{mainX}. We shall need the following well-known lemma.

\begin{lem}
\label{lemX2}
For all $x > 1$, we have $\left(1 - \frac{1}{x} \right)^{x} < e^{-1} < \left(1 - \frac{1}{x} \right)^{x-1}$.
\end{lem}

We are now in a position to prove Theorem~\ref{mainX}.

\noindent
\textbf{Proof of Theorem~\ref{mainX}.} Let $H$ and $H'$ be defined as in the statement of Theorem~\ref{mainX}. Let $E(H)=\{e_1,e_2,\ldots,e_n\}$ and let
$e_i'$ be the edge in $H'$ obtained by picking $k$ vertices at random from $e_i$ for $i \in [n]$. Let $T$ be a random set of vertices in $H'$ where
\[
|T| = \left\lfloor \frac{c \, \ln(k)}{k} n \right\rfloor
\]
and let $t_i= |T \cap V(e_i)|$ for each $i \in [n]$. The probability, $\Pr(e_i' \mbox{ not covered})$, that the edge $e_i'$ is not covered by $T$ (that is, $V(e_i') \cap T \ne \emptyset$) is given by
\[
\Pr(e_i' \mbox{ not covered}) = \frac{{2k-t_i \choose k}}{{2k \choose k}}.
\]

Given the values $t_i$ for all $i \in [n]$, we note that the probability of an edge $e_i$ being covered by $T$ is independent of an edge $e_j$ being covered by $T$ where $j \in [n] \setminus \{i\}$. Therefore, the probability that $T$ is a transversal of $H'$ is given by
\begin{eqnarray}
\displaystyle{ \Pr(T \mbox{ is a transversal of }H') } & = & \displaystyle{ \prod_{i=1}^{n} \left( 1 - \Pr(e_i' \mbox{ not covered}) \right) \nonumber } \2 \\
& = & \displaystyle{ \prod_{i=1}^{n} \left( 1 -  \frac{{2k-t_i \choose k}}{{2k \choose k}} \right). \label{Eq0} }
\end{eqnarray}
As every vertex in $T$ belongs to $2k$ edges of $H$, we note that
\[
\sum_{i=1}^n t_i = 2k|T|
\]
by double counting. Let $s_1,s_2,\ldots,s_n$ be non-negative integers, such that the expression
\begin{equation}
\label{Eq1n}
\prod_{i=1}^{n} \left( 1 -  \frac{{2k-s_i \choose k}}{{2k \choose k}} \right)
\end{equation}
is maximized, where $\displaystyle{\sum_{i=1}^n s_i = 2k|T|}$. We proceed further with the following series of claims.

\begin{claim}
\label{c:c1}
$|s_i-s_j| \le 1$ for all $i,j \in [n]$.
\end{claim}
\proof For the sake of contradiction, suppose that $s_i \le s_j-2$ for some $i$ and $j$ where $i,j \in [n]$. We will now show that by increasing $s_i$ by one and decreasing $s_j$ by one we increase the value of the maximized expression in (\ref{Eq1n}), a contradiction.  Hence it suffices for us to show that Inequality~(\ref{Eq2n}) below holds.
\begin{equation}
\label{Eq2n}
\left( 1 -  \frac{{2k-(s_i+1) \choose k}}{{2k \choose k}} \right)
\left( 1 -  \frac{{2k-(s_j-1) \choose k}}{{2k \choose k}} \right) >
\left( 1 -  \frac{{2k-s_i \choose k}}{{2k \choose k}} \right)
\left( 1 -  \frac{{2k-s_j \choose k}}{{2k \choose k}} \right)
\end{equation}

We remark that Inequality~(\ref{Eq2n}) is equivalent to the following.

{\small
\[
\left( {2k \choose k} -  {2k-s_i-1 \choose k} \right)
\left( {2k \choose k} -  {2k-s_j+1 \choose k} \right) >
\left( {2k \choose k} -  {2k-s_i \choose k} \right)
\left( {2k \choose k} -  {2k-s_j \choose k} \right).
\]
}

Defining $A_i,A_j,B_i,B_j,C_i,C_j$ as follows,

\[
\begin{array}{rclcrcl} \vspace{0.1cm}
A_i & = & {2k \choose k} -  {2k-s_i \choose k}  & \hspace{1cm} & A_j & = & {2k \choose k} -  {2k-s_j \choose k} \2 \\
B_i & = & {2k-s_i-1 \choose k-1} & &  B_j & = & {2k-s_j+1 \choose k-1} \2 \\
C_i & = & {2k \choose k} -  {2k-s_i-1 \choose k} & &  C_j & = & {2k \choose k} -  {2k-s_j+1 \choose k}, \\
\end{array}
\]
we note that Inequality~(\ref{Eq2n}) is equivalent to
\begin{equation}
\label{Eq3n}
C_iC_j > A_iA_j.
\end{equation}

By Pascal's rule, which states that ${n \choose k} = {n-1 \choose k-1} + {n-1 \choose k}$, we obtain
the following.
\[
\begin{array}{ccccc} \2
C_i & = &   {2k \choose k} -  \left({2k-s_i \choose k} - {2k-s_i-1 \choose k-1} \right) & = & A_i + B_i \1 \\
C_j & = &  {2k \choose k} -  \left({2k-s_j \choose k} + {2k-s_j \choose k-1}  \right)  & = & A_j - B_j. \\
\end{array}
\]

Since $s_i \le s_j-2$, we note that $A_j > C_i = A_i+B_i$ and $B_i \ge B_j$, implying that
\[
\begin{array}{lcl}
C_i C_j & = &  (A_i+B_i)(A_j-B_j) \\
& \ge & (A_i+B_i)(A_j-B_i) \\
& = &  A_iA_j + B_i(A_j-(A_i+B_i)) \\
& > & A_iA_j.
\end{array}
\]
Thus, Inequality~(\ref{Eq3n}) holds, producing the desired contradiction. This completes the proof of Claim~\ref{c:c1}.~\smallqed

\medskip
Let $\Pr(H)$ denote the probability that $\tau(H') \le \left( \frac{c \, \ln(k)}{k} \right) n$. Since $\tau(H')$ is an integer, we note that $\Pr(H) = \Pr(\tau(H') \le |T|)$.

\begin{claim}
\label{c:c2}
$\displaystyle{ \Pr(H) \le {n \choose |T|} \prod_{i=1}^{n} \left( 1 -  \frac{{2k-s_i \choose k}}{{2k \choose k}} \right) }$.
\end{claim}
\proof As there are ${n \choose |T|}$ possible ways of choosing the set $T$, Equation~(\ref{Eq0}) implies that the probability that there exists a transversal of $H'$ of size at most $|T|$ is bounded as follows.
\[
\begin{array}{rcl}
\displaystyle{ \Pr(\tau(H') \le |T|) } & \le & \displaystyle{ {n \choose |T|} \times \Pr(T \mbox{ is a transversal of }H') } \2 \\
& \le & \displaystyle{ n^{|T|} \cdot \prod_{i=1}^{n} \left( 1 -  \frac{{2k-t_i \choose k}}{{2k \choose k}} \right) } \2 \\
& \le & \displaystyle{ n^{|T|} \cdot \prod_{i=1}^{n} \left( 1 -  \frac{{2k-s_i \choose k}}{{2k \choose k}} \right).}  \\
\end{array}
\]
This completes the proof of Claim~\ref{c:c2}.~\smallqed

\medskip
Let $s_{{\rm av}}$ denote the average value of the integers $s_1,s_2,\ldots,s_n$, and let
\[
s^* = 2c \, \ln(k).
\]

\begin{claim}
\label{c:c2A}
$s_{{\rm av}} \le s^*$.
\end{claim}
\proof Since $\displaystyle{\sum_{i=1}^n s_i = 2k|T|}$, the average value of $s_1,s_2,\ldots,s_n$ is bounded as follows.
\[
s_{{\rm av}} = \frac{2k|T|}{n} = \frac{2k\lfloor \frac{c \, \ln(k)}{k} n \rfloor}{n}  \le 2c \, \ln(k) = s^*.
\]
This completed the proof of Claim~\ref{c:c2A}.~\smallqed

\begin{claim}
\label{c:c3}
$\displaystyle{ \frac{{2k-s_i \choose k}}{{2k \choose k}} > \frac{1}{5} \left( \frac{1}{2} \right)^{s_i} }$.
\end{claim}
\proof We note that
\[
\frac{{2k-s_i \choose k}}{{2k \choose k}} = \frac{(2k-s_i)!}{k!(k-s_i)!)} \frac{k! k!}{(2k)!}
= \frac{k(k-1) \cdots (k-s_i+1)}{(2k)(2k-1) \cdots (2k-s_i+1)}
\ge \left( \frac{k-s_i+1}{2k-s_i+1} \right)^{s_i}.
\]
Let
\[
\theta(s_i) = \frac{2k-2s_i+2}{s_i(s_i-1)}.
\]

Since $\theta(s_i) s_i + 1 = \frac{2k-s_i+1}{s_i-1}$, we note by Lemma~\ref{lemX2} that the following holds.

\[
\begin{array}{rcl}
\displaystyle{ \left( \frac{k-s_i+1}{2k-s_i+1} \right)^{s_i} }
& = & \displaystyle{ \left(\frac{1}{2} \times \frac{2k-2s_i+2}{2k-s_i+1} \right)^{s_i} } \2 \\
& = & \displaystyle{ \left(\frac{1}{2} \right)^{s_i} \times \left( 1 - \frac{s_i-1}{2k-s_i+1} \right)^{s_i} } \2 \\
& = & \displaystyle{ \left(\frac{1}{2} \right)^{s_i} \times \left( 1 - \frac{1}{\frac{2k-s_i+1}{s_i-1}} \right)^{s_i} } \2 \\
& = & \displaystyle{ \left(\frac{1}{2} \right)^{s_i} \times \left(\left( 1 - \frac{1}{\theta(s_i) s_i + 1} \right)^{\theta(s_i) s_i} \right)^{\frac{1}{\theta(s_i)}} } \2 \\
& > & \displaystyle{ \left(\frac{1}{2} \right)^{s_i} \times \left( e^{-1} \right)^{\frac{1}{\theta(s_i)}}. }
\end{array}
\]

In order to prove our desired result, it therefore suffices for us to prove that
\[
e^{-\frac{1}{\theta(s_i)}} \ge \frac{1}{5}.
\]

By Claim~\ref{c:c1} and Claim~\ref{c:c2A}, we note that $s_i \le s^*+1$ for all $i \in [n]$. As $\theta(s_i)$ is a decreasing function, the function  $e^{-1/\theta(s_i)}$ is also
a decreasing function in $s_i$. Hence it suffices for us to show that
$e^{-1/\theta(s^*+1)} \ge 1/5$. We note that
\begin{equation}
\label{Eq4n}
e^{-\frac{1}{\theta(s^* + 1)}} = e^{ - \frac{(s^*+1)s^*}{2k-2(s^*+1)+2} }
= e^{ - \frac{(1+2c\ln(k))2c\ln(k)}{2k-4c\ln(k)}.} \end{equation}
Since $c < 1/\ln(4)$,
\[
\frac{(1+2c\ln(k))2c\ln(k)}{2k-4c\ln(k)} \le  \frac{\left(1+2\frac{\ln(k)}{\ln(4)}\right)2\ln(k)}{2\ln(4)k-4\ln(k)}. \]

The maximum value of the function on the right-hand-side of the above inequality is approximately $1.5037$ obtained at $k=3.753$, which is always less than $\ln(5) \approx 1.6094$. Thus,
\[
\frac{(1+2c\ln(k))2c\ln(k)}{2k-4c\ln(k)} < \ln(5),
\]
implying by~(\ref{Eq4n}) that
\[
e^{-\frac{1}{\theta(s^*+1)}} \ge \frac{1}{5}.
\]
This completes the proof of Claim~\ref{c:c3}.~\smallqed

\begin{claim}
\label{c:c4}
$\displaystyle{ \Pr(H) \le {n \choose |T|} \left( \left( 1 -  \frac{1}{5}\left( \frac{1}{2} \right)^{s^*}  \right) \right)^{n}}$.
\end{claim}
\proof For a non-negative real number $s$, let \[f(s) = 1 -  \frac{1}{5}\left( \frac{1}{2} \right)^{s}.\] We show firstly that
\begin{equation}
\label{Eq5n}
f(x-y)f(x+y) \le f(x)^2
\end{equation}

\noindent
holds for all $0 \le y \le x$. This is the case due to the following, where $Y=(1/2)^y$ and $X=(1/2)^x$.

\[
\begin{array}{crcl}
 & 0 & \le & (Y-1)^2 \\
\Updownarrow & & & \\
 & 2Y & \le & Y^2 + 1 \\
\Updownarrow & & & \\
 & 2X & \le & XY + \frac{X}{Y} \\
\Updownarrow & & & \\
 & -\frac{1}{5}XY - \frac{1}{5}\frac{X}{Y} & \le & -\frac{2}{5}X  \\
\Updownarrow & & & \\
 & 1 - \frac{1}{5}\frac{X}{Y} -\frac{1}{5}XY + \frac{1}{25}X^2 & \le & 1 - \frac{2}{5}X + \frac{1}{25}X^2 \\
\Updownarrow & & & \\
 & (1 -\frac{1}{5}\frac{X}{Y})(1-\frac{1}{5}XY) & \le & (1-\frac{1}{5}X)^2 \\
\Updownarrow & & & \\
 & \left(1 - \frac{1}{5}\left(\frac{1}{2}\right)^{x-y}\right)\left(1-\frac{1}{5}\left(\frac{1}{2}\right)^{x+y}\right) 
            & \le & \left(1-\frac{1}{5}\left(\frac{1}{2}\right)^{x}\right)^2 \\
\Updownarrow & & & \\
& f(x-y)f(x+y) & \le & f(x)^2 \\
\end{array}
\]

Therefore the product $f(s_1)f(s_2)\cdots f(s_n)$, where $s_i$ can be reals and not just integers and  $\sum_{i=1}^n s_i = 2k|T|$ is maximized when all $s_i$ have the same value. Furthermore, note that as $f(x)$ is an increasing function the following holds by Lemma~\ref{lemX2}, Claim~\ref{c:c2} and Claim~\ref{c:c3}.

\[
\begin{array}{rcl}  \1
\Pr(H) & \le & \displaystyle{ {n \choose |T|} \prod_{i=1}^{n} \left( 1 -  \frac{{2k-s_i \choose k}}{{2k \choose k}} \right) } \2 \\
 & \le & \displaystyle{ {n \choose |T|} \prod_{i=1}^{n} \left( 1 -  \frac{1}{5}\left( \frac{1}{2} \right)^{s_i} \right) } \2 \\
 & < & \displaystyle{ n^{|T|} \left(  1 -  \frac{1}{5}\left( \frac{1}{2} \right)^{s^{av}} \right)^n } \2 \\
   & \le & \displaystyle{ n^{|T|} \left(  1 -  \frac{1}{5} \left( \frac{1}{2} \right)^{s^{*}}  \right)^n }.
\end{array}
\]
This completes the proof of Claim~\ref{c:c4}.~\smallqed

\medskip
We return to the proof of Theorem~\ref{mainX} one final time. Recall that $|T| = \lfloor \frac{c \, \ln(k)}{k} n \rfloor$ and $s^* = 2c \ln(k)$. Further recall that by assumption, $5 c \ln(k) \ln(n)  <  k^{1-c\ln(4)}$. Thus
\[
5 c \ln(k) \ln(n)  <  k^{1-c\ln(4)} = \frac{k}{k^{2c\ln(2)}} = \frac{k}{2^{\ln(k^{2c})}},
\]
implying that
\[
|T| \ln(n) \le \frac{c \, \ln(k)}{k} n \ln(n) < \frac{n}{5 \cdot 2^{2c\ln(k)}} = \frac{n}{5 \cdot 2^{s^*}}.
\]
Therefore,
\[
|T| \ln(n) - \frac{n}{5 \cdot 2^{s^*}} < 0.
\]
Hence by Claim~\ref{c:c4},
\[
\begin{array}{rcl}  \1
\Pr(H) & \le & \displaystyle{ n^{|T|} \left(  1 -  \frac{1}{5} \left( \frac{1}{2} \right)^{s^{*}}  \right)^n } \1 \\
   & = & \displaystyle{ n^{|T|} \left(\left(  1 -  \frac{1}{5 \cdot 2^{s^*}} \right)^{5 \cdot 2^{s^*}} \right)^{\frac{n}{5 \cdot 2^{s^*}}} } \2 \\
   & < & \displaystyle{ n^{|T|} \left( e^{-1} \right)^{\frac{n}{5 \cdot 2^{s^*}}} } \2 \\
   & = & \displaystyle{ e^{(|T| \ln(n) - \frac{n}{5 \cdot 2^{s^*}})} } \2 \\
   & < & 1.
\end{array}
\]
This completes the proof of Theorem~\ref{mainX}.~\qed

\medskip
We close this section with the following remarks.

\noindent
\textbf{Remark~1.} If we assume that $k$ is large, then the coefficient  $5$ of the expression $c \ln(k) \ln(n)$ given in the statement of Theorem~\ref{mainX} can be reduced. In fact the larger $k$ gets, the closer the number can be made to one. For example, if we assume $k \ge 23$, then the $5$ can be decreased to a $2$, and if $k \ge 54$, then it can be decreased to $1.5$.

\noindent
\textbf{Remark~2.} Let $P$ be a projective plane of order $p$, where $p$ is an odd prime power. Thus, $P$ has $n=p^2+p+1$ points and $m=n$ lines, and each line contains $p+1$ points. Let $H$ be the hypergraph whose vertices are the points of $P$ and whose edges are the lines of $P$. We note that $H$ is a $(p+1)$-regular, $(p+1)$-uniform, linear hypergraph of order and size $|V(H)|=|E(H)|=p^2+p+1$. Letting $k=(p+1)/2$, we note that $n=|V(H)|=|E(H)|=4k^2-2k+1$. Applying the result of Theorem~\ref{mainX} on this hypergraph $H$ with $0<c<\frac{1}{\ln(4)}$, if \[
 5c \ln(k) \ln(4k^2-2k+1)  <  k^{1-c\ln(4)},
\]
then there exists a $k$-uniform hypergraph, $H'$, of order~$n$ satisfying
\[
\tau(H') > \left( \frac{c \, \ln(k)}{k} \right) n
\]
The above holds when $k \ge 2753$ and $c$ is determined such that  $(c \ln(k) / k) n = (n+m)/(k+1)$, which implies that we can use $p=5507$ (which is a prime) and $k=2754$.

\noindent
\textbf{Remark~3.} If we let $H$ be the  hypergraph associated with a projective plane of order $p=331$, then $H$ is a $332$-regular, $332$-uniform, linear hypergraph with $|V(H)|=|E(H)|=p^2+p+1$. Letting $k=166$ and using Claim~\ref{c:c2} in the proof of Theorem~\ref{mainX} we can show that the probability that $H'$ has transversal number less that $(|V(H)|+|E(H)|)/(k+1)$  is strictly less than one. Therefore, there must exist a linear $166$-uniform hypergraph, $H^*$, where $\tau(H^*) > (|V(H^*)|+|E(H^*)|)/(k+1)$. This result is stated formally as Theorem~\ref{t:thm166}.

\begin{prop}
 \label{p:prop1}
For $k = 166$, there exist $k$-uniform, linear hypergraph on $n$ vertices with $m$ edges satisfying
\[
\tau(H) > \frac{n+m}{k+1}.
\]
\end{prop}

\noindent
\textbf{Remark~4.} In Remark~2 and Remark~3, we apply Theorem~\ref{mainX} to projective planes. However, we remark that Theorem~\ref{mainX} can be used on linear hypergraphs $H$ which are not necessarily projective planes as well. Provided the order of $H$ is a polynomial in $k$, the condition of Theorem~\ref{mainX} will always hold for sufficiently large $k$.

\section{Proof of Theorem~\ref{t:mainY}}
\label{S:mainY}

In this section, we give a proof of Theorem~\ref{t:mainY}. For this purpose, we shall prove the following result.

\begin{thm}
\label{t:mainYY}
Let $F_{k^2}$ be the linear, $k$-uniform, $(k+1)$-regular hypergraph of order $k^2$ which is equivalent to the affine plane $AG(2,k)$ of order~$k$ for some $k \ge 2$. Let $e \in E(F_{k^2})$ be an arbitrary edge in $F_{k^2}$ and let $X \subseteq V(e)$ be any non-empty subset of vertices belonging to the edge $e$. If $H = F_{k^2}(X)$ is the linear, $k$-uniform hypergraph obtained from $F_{k^2}$ by deleting the vertex set $X$ and all edges intersecting $X$, then
\[
\tau(H) = \frac{n(H) + m(H)}{k+1}.
\]
\end{thm}
\proof
Let $F_{k^2}$, $e$ and $X$ be defined as in the statement of the theorem. For notational simplicity, as $k$ is fixed, denote $F_{k^2}$ by $F$.  By Theorem~\ref{t:Affine_plane}, the transversal number of $F$ is $2(k-1) + 1$; that is, $\tau(F) =  2(k-1) + 1$. Let $e' \in E(F) \setminus \{e\}$ be an arbitrary edge intersecting $e$.  Since $F$ is equivalent to an affine plane, every edge in $E(F)\setminus \{e,e'\}$ will intersect $e$ or $e'$ (or both $e$ and $e'$). Thus, $V(e) \cup V(e')$ is a transversal in $F$ of size $2(k-1) + 1 = 2k-1$ as $F$ is $k$-uniform.
Let $F(X) = F - X$; that is, $F(X)$ is obtained from $F$ by deleting $X$ and all edges incident with $X$. If $T_X$ is a transversal in $F(X)$, then $T_X \cup X$ is a transversal in $F$, which implies that $\tau(F(X)) \ge 2k-1-|X|$. However, $(V(e) \cup V(e')) \setminus X$ is a transversal in $F(X)$ of size $2k-1-|X|$, and so $\tau(F(X)) \le 2k-1-|X|$. Consequently,
\[
\tau(F(X)) = 2k-1-|X|.
\]

We will now compute the order and size of $F(X)$. The order of $F(X)$ is
\[
n(F(X)) = n(F) - |X| = k^2-|X|.
\]

Let $x$ be an arbitrary vertex in $X$. As $F-x$ is $k$-regular and has order $k^2-1$, the hypergraph $F - x$ contains $k^2-1$ edges. Further since $F$ is linear and $e$ is not an edge in $F-x$, no edge in $F-x$ contains more than one vertex from $X \setminus \{x\}$. Therefore, we remove $k$ edges from $F-x$ for every vertex in $X \setminus \{x\}$ we remove when constructing $F(X)$. Thus, the size of $F(X)$ is
\[
m(F(X)) = k^2-1 - (|X|-1)k = k^2 + k - 1 - k|X|.
\]
Therefore,
\[
\tau(F(X)) = 2k-1-|X| = \frac{(k^2-|X|) + (k^2 + k - 1 - k|X|)}{k+1} = \frac{n(F(X)) + m(F(X))}{k+1}.
\]

This completes the proof of Theorem~\ref{t:mainY}, noting that $F(X)=F_{k^2}(X)$.~\qed

\medskip
In particular, when $k=4$ we know that $AG(2,4)$ exists and therefore, by Theorem~\ref{t:mainY}, we have the following hypergraphs where Theorem~\ref{t:mainZ} is tight.

\begin{center}
\begin{tabular}{|c||c|c|c|c|l|} \hline
\multicolumn{6}{|c|}{$F = F_{16}(X)$} \\ \hline \hline
        &  $\tau(F)$  &  $n(F)$  &  $m(F)$ & $ (n(F)+m(F))/5$ & Average degree  \\ \hline
$|X|=1$ &    6        &   15     &   15    &  6  &  $60/15 = 4$ \\ \hline
$|X|=2$ &    5        &   14     &   11    &  5  &  $44/14 \approx 3.14\ldots$  \\ \hline
$|X|=3$ &    4        &   13     &    7    &  4  &  $28/13 \approx 2.15\ldots$  \\ \hline
$|X|=4$ &    3        &   12     &    3    &  3  &  $12/12  = 1$ \\ \hline
\end{tabular}
\end{center}

\medskip
We close this section with the following remark.

\noindent
\textbf{Remark~5.} We remark that if $\tau(H) \le (n + m)/(k+1)$ holds for all $k$-uniform linear hypergraphs for some $k \ge 2$, then the bound cannot be tight for any hypergraph, $H$, with average degree greater than~$k$ as removing any vertex, $x$, of degree more than~$k$ from $H$ and applying the bound on $\tau(H-x)$ gives us a transversal in $H$ of size at most the following,
\[
\tau(H - x) + |\{x\}| \le \frac{(|V(H)|-1) + (|E(H)|-(k+1))}{k+1} + 1 < \frac{|V(H)|+|E(H)|}{k+1}.
\]

Similarly, if $\tau(H) \le (n + m)/(k+1)$ holds for all $k$-uniform linear hypergraphs for some $k \ge 2$, then the bound cannot be tight if the average degree is less than~$1$ as then there are isolated vertices that can be removed. Therefore, Theorem~\ref{t:mainYY} can be used to show that if $k_{\min} > k$ for some $k \ge 2$ (and the affine plane $AG(2,k)$ exists), then the bound $\tau(H) \le (n + m)/(k+1)$ is tight for average degree $k$ and average degree~$1$ and for a number of average degrees in the interval from $1$ to $k$. This is somewhat surprising as there are no similar kinds of bounds which hold for a wide variety of average degrees if we consider non-linear hypergraphs.

\section{Applications of Theorem~\ref{t:mainZ}}
\label{S:applications}

In this section, we present a few applications to serve as  motivation for the significance of our result given in Theorem~\ref{t:mainZ}.

\subsection{Application~1}

The following conjecture is posed in~\cite{HeYe13}.

\begin{conj}{\rm (\cite{HeYe13})}
 \label{conj1}
If $H$ is a $4$-uniform, linear hypergraph on $n$ vertices with $m$ edges, then
$\tau(H) \le \frac{n}{4} + \frac{m}{6}$.
\end{conj}

We remark that the linearity constraint in Conjecture~\ref{conj1} is essential. Indeed if $H$ is not linear, then Conjecture~\ref{conj1} is not always true, as may be seen, for example, by taking $H$ to be the complement of the Fano plane, $F$, shown in Figure~\ref{f:Fano}.
The second consequence of our main result proves Conjecture~\ref{conj1}.

\begin{thm}
 \label{t:thm2}
Conjecture~\ref{conj1} is true.
\end{thm}
\proof Let $H$ be a $4$-uniform, linear hypergraph on $n$ vertices with $m$ edges. We show that $\tau(H) \le \frac{n}{4} + \frac{m}{6}$. We proceed by induction on~$n$. If $n = 4$, then $H$ consists of a single edge, and $\tau(H) = 1 < \frac{n}{4} + \frac{m}{6}$. Let $n \ge 5$ and suppose that the result holds for all $4$-uniform, linear hypergraphs on fewer than $n$ vertices. Let $H$ be a $4$-uniform, linear hypergraph on $n$ vertices with $m$ edges. Suppose that $\Delta(H) \le 6$. In this case, $2m \le 3n$. By Theorem~\ref{t:mainZ}, $60\tau(H) \le 12n + 12m = 15n + 10m + 2m - 3n \le 15n + 10m$, or, equivalently, $\tau(H) \le \frac{n}{4} + \frac{m}{6}$. Hence, we may assume that $\Delta(H) \ge 7$, for otherwise the desired result follows from Theorem~\ref{t:mainZ}. Let $v$ be a vertex of maximum degree in $H$, and consider the $4$-uniform, linear hypergraph $H' = H - v$ on $n' = n - 1$ vertices with $m'$ edges. We note that $n' = n - 1$ and $m' = m - \Delta(H) \le m - 7$. Every transversal in $H'$ can be extended to a transversal in $H$ by adding to it the vertex~$v$. Hence, applying the inductive hypothesis to $H'$, we have that
$\tau(H) \le 1 + \tau(H') \le 1 + \frac{n'}{4} + \frac{m'}{6} \le 1 + \frac{n-1}{4} + \frac{m-7}{6} < \frac{n}{4} + \frac{m}{6}$.~\qed

\subsection{Application~2}

There has been much interest in determining upper bounds on the
transversal number of a $3$-regular $4$-uniform hypergraph. In
particular, as a consequence of more general results we have the Chv\'{a}tal-McDiarmid bound, the improved Lai-Chang bound, the further improved Thomass\'{e}-Yeo bound, and the recent bound given in~\cite{HeYe16}. These bounds are summarized in Theorem~\ref{t:knownbds}. We observe that $\frac{3}{8} < \frac{8}{21} < \frac{7}{18} < \frac{5}{12}$.

\begin{thm}
 \label{t:knownbds}
Let $H$ be a $3$-regular, $4$-uniform hypergraph on $n$ vertices. Then the following bounds on $\tau(H)$ have been established. \\[-27pt]
\begin{enumerate}
\item  $\tau(H) \le \frac{5}{12} n \approx 0.41667 \, n$ {\rm (Chv\'{a}tal, McDiarmid~\cite{ChMc}).} \1
\item $\tau(H) \le \frac{7}{18}n \approx 0.38888 \, n$ {\rm (Lai, Chang~\cite{LaCh90}).} \1
\item $\tau(H) \le \frac{8}{21}n \hspace*{0.05cm} \approx 0.38095 \, n$ {\rm (Thomass\'{e}, Yeo~\cite{ThYe07}).} \1
\item  $\tau(H) \le \frac{3}{8}n \hspace*{0.15cm} \approx 0.375 \, n$ {\rm (Henning, Yeo~\cite{HeYe16}).}
\end{enumerate}
\end{thm}

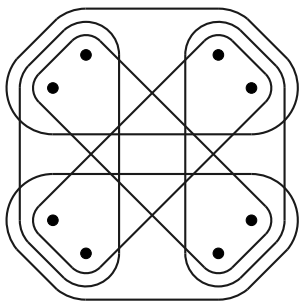
\begin{figure}[htb]
\begin{center}

\tikzstyle{vertexX}=[circle,draw, fill=black!100, minimum size=8pt, scale=0.5, inner sep=0.1pt]
\begin{tikzpicture}[scale=0.44]
 \draw (0,0) node {\mbox{ }};
\draw (4,3) node {\mbox{ }};
\node (a1) at (0.0,1.0) [vertexX] {};
\node (a2) at (1.0,0.0) [vertexX] {};
\node (b1) at (0.0,5.0) [vertexX] {};
\node (b2) at (1.0,6.0) [vertexX] {};
\node (c1) at (5.0,6.0) [vertexX] {};
\node (c2) at (6.0,5.0) [vertexX] {};
\node (d1) at (6.0,1.0) [vertexX] {};
\node (d2) at (5.0,0.0) [vertexX] {};
\draw[color=black!90, thick,rounded corners=4pt] (0.33828113841207375,-0.7497520578289725) arc (228.5689049390958:359.9147461611031:1.0); 
\draw[color=black!90, thick,rounded corners=4pt] (-0.99998850591691,0.9952054161738438) arc (180.27471047027626:221.4310950609044:1.0); 
\draw[color=black!90, thick,rounded corners=4pt] (-0.7497520578289725,5.661718861587926) arc (138.5689049390958:179.72528952972374:1.0); 
\draw[color=black!90, thick,rounded corners=4pt] (1.9999988929874426,6.001487959639652) arc (0.0852538388969144:131.4310950609044:1.0); 
\draw[color=black!90, thick,rounded corners=4pt] (0.33828113841207375,-0.7497520578289725) -- (-0.7497520578289704,0.3382811384120713);
\draw[color=black!90, thick,rounded corners=4pt] (-0.99998850591691,0.9952054161738438) -- (-0.99998850591691,5.004794583826157);
\draw[color=black!90, thick,rounded corners=4pt] (-0.7497520578289725,5.661718861587926) -- (0.3382811384120713,6.74975205782897);
\draw[color=black!90, thick,rounded corners=4pt] (1.9999988929874426,6.001487959639652) -- (1.9999988929874426,-0.0014879596396527753);
\draw[color=black!90, thick,rounded corners=4pt] (5.661718861587927,6.749752057828973) arc (48.56890493909578:179.91474616110307:1.0); 
\draw[color=black!90, thick,rounded corners=4pt] (6.99998850591691,5.004794583826157) arc (0.2747104702762737:41.4310950609044:1.0); 
\draw[color=black!90, thick,rounded corners=4pt] (6.749752057828972,0.3382811384120734) arc (318.5689049390958:359.72528952972374:1.0); 
\draw[color=black!90, thick,rounded corners=4pt] (4.000001107012557,-0.0014879596396528518) arc (180.08525383889693:311.43109506090445:1.0); 
\draw[color=black!90, thick,rounded corners=4pt] (5.661718861587927,6.749752057828973) -- (6.74975205782897,5.661718861587929);
\draw[color=black!90, thick,rounded corners=4pt] (6.99998850591691,5.004794583826157) -- (6.99998850591691,0.9952054161738435);
\draw[color=black!90, thick,rounded corners=4pt] (6.749752057828972,0.3382811384120734) -- (5.6617188615879295,-0.7497520578289695);
\draw[color=black!90, thick,rounded corners=4pt] (4.000001107012557,-0.0014879596396528518) -- (4.000001107012557,6.001487959639653);
\draw[color=black!90, thick,rounded corners=4pt] (-1.090271767160063,5.878241125052609) arc (141.14768868443804:269.77422548308033:1.4); 
\draw[color=black!90, thick,rounded corners=4pt] (0.9830127200795282,7.399896936320993) arc (90.69523094951481:128.85231130303504:1.4); 
\draw[color=black!90, thick,rounded corners=4pt] (5.878241125052609,7.090271767160063) arc (51.14768868443803:89.30476905048522:1.4); 
\draw[color=black!90, thick,rounded corners=4pt] (6.005516697885407,3.600010869312036) arc (-89.77422548308033:38.85231130303504:1.4); 
\draw[color=black!90, thick,rounded corners=4pt] (-1.090271767160063,5.878241125052609) -- (0.12175887518576367,7.090271767352078);
\draw[color=black!90, thick,rounded corners=4pt] (0.9830127200795282,7.399896936320993) -- (5.016987279920471,7.399896936320993);
\draw[color=black!90, thick,rounded corners=4pt] (5.878241125052609,7.090271767160063) -- (7.090271767352078,5.878241124814236);
\draw[color=black!90, thick,rounded corners=4pt] (6.005516697885407,3.600010869312036) -- (-0.005516697885407651,3.600010869312036);
\draw[color=black!90, thick,rounded corners=4pt] (0.1217588749473899,-1.0902717671600624) arc (231.14768868443804:269.3047690504852:1.4); 
\draw[color=black!90, thick,rounded corners=4pt] (-0.005516697885406376,2.399989130687964) arc (90.22577451691966:218.85231130303504:1.4); 
\draw[color=black!90, thick,rounded corners=4pt] (7.090271767160062,0.1217588749473898) arc (-38.85231131556196:89.77422548308027:1.4); 
\draw[color=black!90, thick,rounded corners=4pt] (5.016987279920471,-1.3998969363209934) arc (270.6952309495148:308.85231130303504:1.4); 
\draw[color=black!90, thick,rounded corners=4pt] (0.1217588749473899,-1.0902717671600624) -- (-1.0902717673520785,0.12175887518576367);
\draw[color=black!90, thick,rounded corners=4pt] (-0.005516697885406376,2.399989130687964) -- (6.005516697885408,2.399989130687964);
\draw[color=black!90, thick,rounded corners=4pt] (7.090271767160062,0.1217588749473898) -- (5.878241124814236,-1.0902717673520785);
\draw[color=black!90, thick,rounded corners=4pt] (5.016987279920471,-1.3998969363209934) -- (0.9830127200795284,-1.3998969363209934);
\draw[color=black!90, thick,rounded corners=4pt] (0.5843075616526987,-0.43266591811798066) arc (226.1461890060021:314.98844089708825:0.6); 
\draw[color=black!90, thick,rounded corners=4pt] (-0.42434965297924815,1.42417846717672) arc (135.01155910291172:223.8538109939979:0.6); 
\draw[color=black!90, thick,rounded corners=4pt] (5.4156924383473015,6.432665918117981) arc (46.14618900600209:134.98844089708828:0.6); 
\draw[color=black!90, thick,rounded corners=4pt] (6.424349652979248,4.57582153282328) arc (-44.98844089708825:43.85381099399791:0.6); 
\draw[color=black!90, thick,rounded corners=4pt] (0.5843075616526987,-0.43266591811798066) -- (-0.4326659181179808,0.5843075616526989);
\draw[color=black!90, thick,rounded corners=4pt] (-0.42434965297924815,1.42417846717672) -- (4.57582153282328,6.424349652979248);
\draw[color=black!90, thick,rounded corners=4pt] (5.4156924383473015,6.432665918117981) -- (6.432665918117981,5.4156924383473015);
\draw[color=black!90, thick,rounded corners=4pt] (6.424349652979248,4.57582153282328) -- (1.4241784671767195,-0.4243496529792485);
\draw[color=black!90, thick,rounded corners=4pt] (-0.43266591811798066,5.4156924383473015) arc (136.1461890060021:224.98844089708828:0.6); 
\draw[color=black!90, thick,rounded corners=4pt] (1.4241784671767197,6.424349652979249) arc (45.011559102911725:133.8538109939979:0.6); 
\draw[color=black!90, thick,rounded corners=4pt] (6.432665918117981,0.5843075616526987) arc (-43.85381099399791:44.98844089708825:0.6); 
\draw[color=black!90, thick,rounded corners=4pt] (4.57582153282328,-0.42434965297924815) arc (225.01155910291172:313.8538109939979:0.6); 
\draw[color=black!90, thick,rounded corners=4pt] (-0.43266591811798066,5.4156924383473015) -- (0.5843075616526989,6.432665918117981);
\draw[color=black!90, thick,rounded corners=4pt] (1.4241784671767197,6.424349652979249) -- (6.424349652979248,1.42417846717672);
\draw[color=black!90, thick,rounded corners=4pt] (6.432665918117981,0.5843075616526987) -- (5.4156924383473015,-0.4326659181179808);
\draw[color=black!90, thick,rounded corners=4pt] (4.57582153282328,-0.42434965297924815) -- (-0.4243496529792482,4.57582153282328);
 \end{tikzpicture}

\end{center}
\vskip -0.5 cm
\caption{A $3$-regular $4$-uniform hypergraph, $H_8$, on $n$ vertices with $\tau(H_8) =  \frac{3}{8}n$.}
\label{hyper8}
\end{figure}

The bound in Theorem~\ref{t:knownbds}(d) is best possible, due to the (non-linear) hypergraph, $H_8$, with $n = 8$ vertices and $\tau(H) = 3$ shown in Figure~\ref{hyper8}.

A natural question is whether the upper bound in Theorem~\ref{t:knownbds}(d), namely $\tau(H) \le \frac{3}{8}n$, can be improved if we restrict our attention to linear hypergraphs. We answer this question in the affirmative. If $H$ is a $3$-regular, $4$-uniform, linear hypergraph on $n$ vertices with $m$ edges, then $m = \frac{3}{4}n$, and so, by Theorem~\ref{t:mainZ}, $\tau(H) \le \frac{1}{5}(n+m) = \frac{1}{5}(n + \frac{3}{4}n) = \frac{7}{20}n$. Hence, as an immediate corollary of Theorem~\ref{t:mainZ}, we have the following result.

\begin{thm}
 \label{t:thm3}
If $H$ is a $3$-regular, $4$-uniform, linear hypergraph on $n$ vertices, then $\tau(H) \le \frac{7}{20}n = 0.35 n$.
\end{thm}

\subsection{Application~3}
\label{S:Applic3}

Lai and Chang~\cite{LaCh90} established the following upper bound on the transversal number of a $4$-uniform hypergraph.

\begin{thm}{\rm (\cite{LaCh90})}
 \label{t:LaiChang}
If $H$ is a $4$-uniform hypergraph with $n$ vertices and $m$ edges, then $\tau(H) \le 2(n+m)/9$.
\end{thm}

\vskip -0.55cm
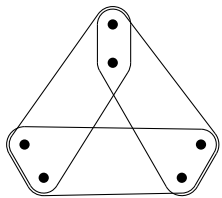
\begin{figure}[htb]
\begin{center}
\begin{tikzpicture}[scale=0.34]

\begin{scope}[xshift=8cm,yshift=0cm]
\fill (-2.7,0) circle (0.2cm); \fill (2.7,0) circle (0.2cm); \fill
(-3.45,1.3) circle (0.2cm); \fill (3.45,1.3) circle (0.2cm); \fill
(0,4.5) circle (0.2cm); \fill (0,6) circle (0.2cm); \hyperedgefour
{{-3.45}{1.3}{0.7}} {{-2.7}{0}{0.7}} {{2.7}{0}{0.6}}
{{3.45}{1.3}{0.6}}; \hyperedgefour {{0}{4.5}{0.7}} {{0}{6}{0.7}}
{{-3.45}{1.3}{0.6}} {{-2.7}{0}{0.6}}; \hyperedgefour {{2.7}{0}{0.7}}
{{3.45}{1.3}{0.7}} {{0}{6}{0.6}} {{0}{4.5}{0.6}};
\end{scope}

\end{tikzpicture}
\end{center}
\vskip -0.6 cm \caption{The hypergraph $T_4$.} \label{f:T4}
\end{figure}

The hypergraph $T_4$, illustrated in Figure~\ref{f:T4}, shows that the Lai-Chang bound is best possible, even if we restrict the maximum degree to be equal to~$2$. Our main result, namely Theorem~\ref{t:mainZ}, improves this upper bound from $\frac{2}{9}(n+m)$ to $\frac{1}{5}(n+m)$ in the case of linear hypergraphs. As an application of the proof of our main result, we show that the $\frac{1}{5}(n+m)$ bound can be further improved to $\frac{3}{16}(n + m) + \frac{1}{16}$ if we exclude the special hypergraph $H_{10}$.

\begin{figure}[htb]
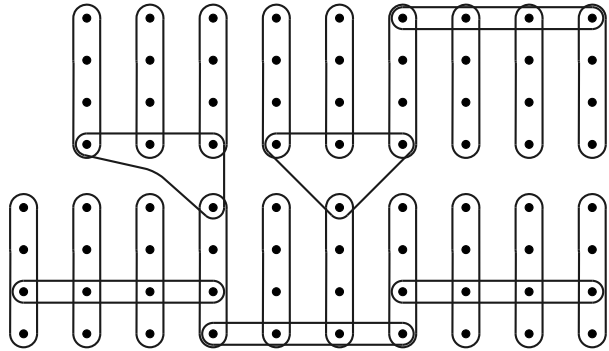

\begin{center}
\tikzstyle{vertexX}=[circle,draw, fill=black!100, minimum size=6pt, scale=0.5, inner sep=0.4pt]
\tikzstyle{vertexY}=[circle,draw, fill=black!40, minimum size=20pt, scale=1.3, inner sep=5.5pt]
\newcommand{\Ss}[1]{{\normalsize {\color{black} #1}}}
\newcommand{\Ssx}[1]{{\footnotesize {\color{blue} #1}}}
\newcommand{\Ssy}[1]{{\footnotesize {\color{green} #1}}}

\end{center}
\vskip -0.6 cm \caption{A hypergraph, $F_6$, in the family~$\cF$.} \label{f:cF}
\end{figure}

For this purpose, we construct a family, $\cF$, of $4$-uniform, connected, linear hypergraphs with maximum degree~$\Delta(H) = 2$ as follows. Let $F_0$ be the hypergraph with one edge (illustrated in Figure~\ref{f:special}(a), but with a different name, $H_4$). For $i \ge 1$, we now build a hypergraph $F_i$ inductively as follows. Let $F_i$ be obtained from $F_{i-1}$ by adding $12$ new vertices, adding three new edges so that each new vertex belongs to exactly one of these added edges, and adding one further edge that contains a vertex in $V(F_{i-1})$ and three additional vertices, one from each of the three newly added edges, in such a way that $\Delta(F_i) = 2$. Let $\cF$ be the family of all such hypergraphs, $F_i$, where $i \ge 0$. A hypergraph, $F_6$, in the family~$\cF$ is illustrated in Figure~\ref{f:cF}.

We are now in a position to state the following result, where $H_{10}$, $H_{14,5}$ and $H_{14,6}$ are the $4$-uniform, linear hypergraphs shown in Figure~\ref{f:special}(b), \ref{f:special}(h) and \ref{f:special}(i), respectively. As observed earlier, $H_4 = F_0$, and so $H_4 \in \cF$. A proof of Theorem~\ref{t:linearDeg2} is given in Section~\ref{S:linearDeg2}.

\begin{thm}
 \label{t:linearDeg2}
Let $H \ne H_{10}$ be a $4$-uniform, connected, linear hypergraph with maximum degree~$\Delta(H) \le 2$ on $n$ vertices with $m$ edges. Then, $\tau(H) \le \frac{3}{16}(n + m) + \frac{1}{16}$, with equality if and only if $H \in \{H_{14,5}, H_{14,6}\}$ or $H \in \cF$.
\end{thm}

\subsection{Application~4}

\begin{figure}[htb]
\tikzstyle{every node}=[circle, draw, fill=black!0, inner sep=0pt,minimum width=.14cm]
\begin{center}
\begin{tikzpicture}[thick,scale=.65]
  \draw(0,0) { 
    +(1.01,3.56) -- +(1.83,3.75)
    +(1.83,3.75) -- +(2.64,3.56)
    +(2.64,3.56) -- +(3.29,3.04)
    +(3.29,3.04) -- +(3.66,2.29)
    +(3.66,2.29) -- +(3.66,1.46)
    +(3.66,1.46) -- +(3.29,0.71)
    +(3.29,0.71) -- +(2.64,0.19)
    +(2.64,0.19) -- +(1.83,0.00)
    +(1.83,0.00) -- +(1.01,0.19)
    +(1.01,0.19) -- +(0.36,0.71)
    +(0.36,0.71) -- +(0.00,1.46)
    +(0.00,1.46) -- +(0.00,2.29)
    +(0.00,2.29) -- +(0.36,3.04)
    +(1.01,3.56) -- +(0.36,3.04)
    +(1.01,3.56) -- +(3.66,1.46)
    +(2.64,3.56) -- +(2.64,0.19)
    +(3.66,2.29) -- +(1.01,0.19)
    +(3.29,0.71) -- +(0.00,1.46)
    +(1.83,0.00) -- +(0.36,3.04)
    +(0.36,0.71) -- +(1.83,3.75)
    +(0.00,2.29) -- +(3.29,3.04)
    +(1.83,3.75) node[fill=black!100]{}
    +(2.64,3.56) node[fill=black!100]{}
    +(3.29,3.04) node[fill=black!100]{}
    +(3.66,2.29) node[fill=black!100]{}
    +(3.66,1.46) node[fill=black!100]{}
    +(3.29,0.71) node[fill=black!100]{}
    +(2.64,0.19) node[fill=black!100]{}
    +(1.83,0.00) node[fill=black!100]{}
    +(1.01,0.19) node[fill=black!100]{}
    +(0.36,0.71) node[fill=black!100]{}
    +(0.00,1.46) node[fill=black!100]{}
    +(0.00,2.29) node[fill=black!100]{}
    +(0.36,3.04) node[fill=black!100]{}
    +(1.01,3.56) node[fill=black!100]{}
      };
\end{tikzpicture}
\end{center}
\vskip -0.6 cm \caption{The Heawood graph.} \label{Ch6:f:Heawood}
\end{figure}
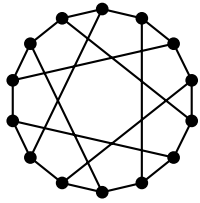

The Heawood graph, shown in Figure~\ref{Ch6:f:Heawood}, is the unique $6$-cage. The bipartite complement of the Heawood graph is the bipartite graph formed by taking the two partite sets of the Heawood graph and joining a vertex from one partite set to a vertex from the other partite set by an edge whenever they are not joined in the Heawood graph.  The bipartite complement of the Heawood graph can also be seen as the incidence bipartite graph of the complement of the Fano plane.

Thomass\'{e} and Yeo~\cite{ThYe07} established the following upper bound on the total domination number of a graph with minimum degree at least~$4$. Recall that $\delta(G)$ denotes the minimum degree of a graph $G$.

\begin{thm}{\rm (\cite{ThYe07})}
 \label{t:tDom}
If $G$ is a graph of order~$n$ with $\delta(G) \ge 4$, then $\gt(G) \le \frac{3}{7}n$.
\end{thm}

The extremal graphs achieving equality in the Thomass\'{e}-Yeo bound of Theorem~\ref{t:tDom} are given by the following result.

\begin{thm}{\rm (\cite{HeYe13,MHAYbookTD})} \label{main_thmD}
If $G$ is a connected graph of order~$n$ with $\delta(G) \ge 4$ that satisfies $\gt(G) = \frac{3}{7}n$, then $G$ is the bipartite complement of the Heawood Graph.
\end{thm}

We remark that every vertex in the bipartite complement of the Heawood Graph belongs to a $4$-cycle. It is therefore a natural question to ask whether the Thomass\'{e}-Yeo upper bound of $\frac{3}{7}n$ can be improved if we restrict $G$ to contain no $4$-cycles.
As a consequence of our main result Theorem~\ref{t:mainZ}, this question can now be answered in the affirmative. For a graph $G$, the \emph{open neighborhood hypergraph}, abbreviated
ONH, of $G$ is the hypergraph $H_G$ with vertex set $V(H_G) = V(G)$
and with edge set $E(H_G) = \{ N_G(x) \mid x \in V \}$ consisting of
the open neighborhoods of vertices in $G$. As first observed in~\cite{ThYe07} (see also~\cite{MHAYbookTD}), the transversal number of the ONH of a graph is precisely the total domination number of the
graph; that is, for a graph $G$, we have $\gt(G) = \tau(H_G)$.

As an application of Theorem~\ref{t:mainZ}, we have the following result, which significantly improves the upper bound of Theorem~\ref{t:tDom} when the graph $G$ contains no $4$-cycle.

\begin{thm}
 \label{t:thm4}
If $G$ is a quadrilateral-free graph of order~$n$ with $\delta(G) \ge 4$, then $\gt(G) \le \frac{2}{5}n$.
\end{thm}
\proof Let $G$ be a quadrilateral-free graph of order~$n$ with $\delta(G) \ge 4$ and let $H_G$ be the ONH of $G$. Then, each edge of $H_G$ has size at least~$4$. Since $G$ is contains no $4$-cycle, the hypergraph $H_G$ contains no overlapping edges and is therefore linear. Let $H$ be
obtained from $H_G$ by shrinking all edges of $H_G$, if necessary, to
edges of size~$4$. Then, $H$ is a $4$-uniform linear hypergraph with $n$
vertices and $n$ edges; that is, $n(H) = m(H) = n(G) = n$.
By Theorem~\ref{t:mainZ} we note that $\tau(H) \le  \frac{1}{5}(n(H)+m(H)) = \frac{2}{5}n$. This completes the proof of the theorem since $\gt(G) = \tau(H_G) \le \tau(H)$.~\qed

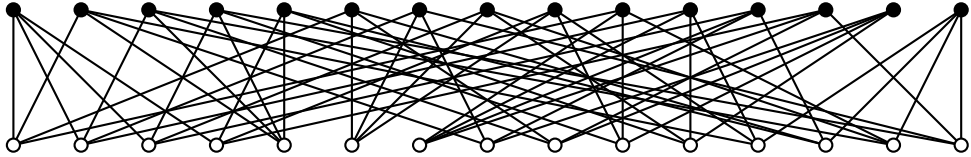
\begin{figure}[htb]
\begin{center}
\begin{tikzpicture}[scale=.9,style=thick,x=1cm,y=1cm]
\def\vr{2.75pt} 
\path (0,0) coordinate (v1);
\path (1,0) coordinate (v2);
\path (2,0) coordinate (v3);
\path (3,0) coordinate (v4);
\path (4,0) coordinate (v5);
\path (5,0) coordinate (v6);
\path (6,0) coordinate (v7);
\path (7,0) coordinate (v8);
\path (8,0) coordinate (v9);
\path (9,0) coordinate (v10);
\path (10,0) coordinate (v11);
\path (11,0) coordinate (v12);
\path (12,0) coordinate (v13);
\path (13,0) coordinate (v14);
\path (14,0) coordinate (v15);
\path (0,2) coordinate (u1);
\path (1,2) coordinate (u2);
\path (2,2) coordinate (u3);
\path (3,2) coordinate (u4);
\path (4,2) coordinate (u5);
\path (5,2) coordinate (u6);
\path (6,2) coordinate (u7);
\path (7,2) coordinate (u8);
\path (8,2) coordinate (u9);
\path (9,2) coordinate (u10);
\path (10,2) coordinate (u11);
\path (11,2) coordinate (u12);
\path (12,2) coordinate (u13);
\path (13,2) coordinate (u14);
\path (14,2) coordinate (u15);
\draw (v1) -- (u1);
\draw (v1) -- (u2);
\draw (v1) -- (u6);
\draw (v1) -- (u10);
\draw (v2) -- (u1);
\draw (v2) -- (u3);
\draw (v2) -- (u7);
\draw (v2) -- (u11);
\draw (v3) -- (u1);
\draw (v3) -- (u4);
\draw (v3) -- (u8);
\draw (v3) -- (u12);
\draw (v4) -- (u1);
\draw (v4) -- (u5);
\draw (v4) -- (u9);
\draw (v4) -- (u13);
\draw (v5) -- (u2);
\draw (v5) -- (u3);
\draw (v5) -- (u4);
\draw (v5) -- (u5);
\draw (v6) -- (u6);
\draw (v6) -- (u7);
\draw (v6) -- (u8);
\draw (v6) -- (u9);
\draw (v7) -- (u10);
\draw (v7) -- (u11);
\draw (v7) -- (u12);
\draw (v7) -- (u13);
\draw (v8) -- (u2);
\draw (v8) -- (u7);
\draw (v8) -- (u12);
\draw (v8) -- (u14);
\draw (v9) -- (u3);
\draw (v9) -- (u6);
\draw (v9) -- (u13);
\draw (v9) -- (u14);
\draw (v10) -- (u4);
\draw (v10) -- (u9);
\draw (v10) -- (u10);
\draw (v10) -- (u14);
\draw (v11) -- (u5);
\draw (v11) -- (u8);
\draw (v11) -- (u11);
\draw (v11) -- (u14);
\draw (v12) -- (u2);
\draw (v12) -- (u9);
\draw (v12) -- (u11);
\draw (v12) -- (u15);
\draw (v13) -- (u5);
\draw (v13) -- (u6);
\draw (v13) -- (u12);
\draw (v13) -- (u15);
\draw (v14) -- (u3);
\draw (v14) -- (u8);
\draw (v14) -- (u10);
\draw (v14) -- (u15);
\draw (v15) -- (u4);
\draw (v15) -- (u7);
\draw (v15) -- (u13);
\draw (v15) -- (u15);
\draw (v1) [fill=white] circle (\vr);
\draw (v2) [fill=white] circle (\vr);
\draw (v3) [fill=white] circle (\vr);
\draw (v4) [fill=white] circle (\vr);
\draw (v5) [fill=white] circle (\vr);
\draw (v6) [fill=white] circle (\vr);
\draw (v7) [fill=white] circle (\vr);
\draw (v8) [fill=white] circle (\vr);
\draw (v9) [fill=white] circle (\vr);
\draw (v10) [fill=white] circle (\vr);
\draw (v11) [fill=white] circle (\vr);
\draw (v12) [fill=white] circle (\vr);
\draw (v13) [fill=white] circle (\vr);
\draw (v14) [fill=white] circle (\vr);
\draw (v15) [fill=white] circle (\vr);
\draw (u1) [fill=black] circle (\vr);
\draw (u2) [fill=black] circle (\vr);
\draw (u3) [fill=black] circle (\vr);
\draw (u4) [fill=black] circle (\vr);
\draw (u5) [fill=black] circle (\vr);
\draw (u6) [fill=black] circle (\vr);
\draw (u7) [fill=black] circle (\vr);
\draw (u8) [fill=black] circle (\vr);
\draw (u9) [fill=black] circle (\vr);
\draw (u10) [fill=black] circle (\vr);
\draw (u11) [fill=black] circle (\vr);
\draw (u12) [fill=black] circle (\vr);
\draw (u13) [fill=black] circle (\vr);
\draw (u14) [fill=black] circle (\vr);
\draw (u15) [fill=black] circle (\vr);
\end{tikzpicture}
\end{center}
\vskip -0.25cm
\caption{A quadrilateral-free $4$-regular graph $G_{30}$ of order~$n = 30$ with $\gamma_t(G_{30}) = \frac{2}{5}n$.} \label{G30}
\end{figure}

\medskip
That the bound in Theorem~\ref{t:thm4} is best possible, may be seen by taking, for example, the $4$-regular bipartite  quadrilateral-free graph $G_{30}$ of order~$n = 30$ illustrated in Figure~\ref{G30} satisfying $\gt(G_{30}) = 12 = \frac{2}{5}n$. We note that the graph $G_{30}$ is the incidence bipartite graph of the linear $4$-uniform hypergraph obtained by removing an arbitrary vertex from the affine plane $AG(2,4)$ of order~$4$.

\section{Key Theorem}
\label{S:key}

In this section, we prove a key theorem that we will need in proving Theorem~\ref{t:mainZ}. First we define a few special hypergraphs. We then introduce the concept of the deficiency of a set of hypergraphs. Thereafter, we prove a key result that will enable us to deduce our main result. Throughout this section, we let $\cH_4$ be the class of all $4$-uniform, linear hypergraphs with maximum degree at most~$3$.

\subsection{Special Hypergraphs}
\label{S:special}

In this section, we define fifteen special hypergraphs, which are shown in Figure~\ref{f:special}.

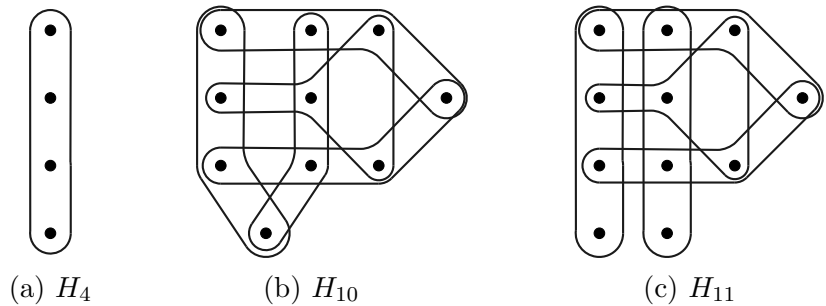
\bigskip

\newcommand{\SsD}[1]{}
\tikzstyle{vertexX}=[circle,draw, fill=black!100, minimum size=8pt, scale=0.5, inner sep=0.1pt]
\newcommand{\Ss}[1]{}
\hspace*{2cm}

 \vskip -0.4cm
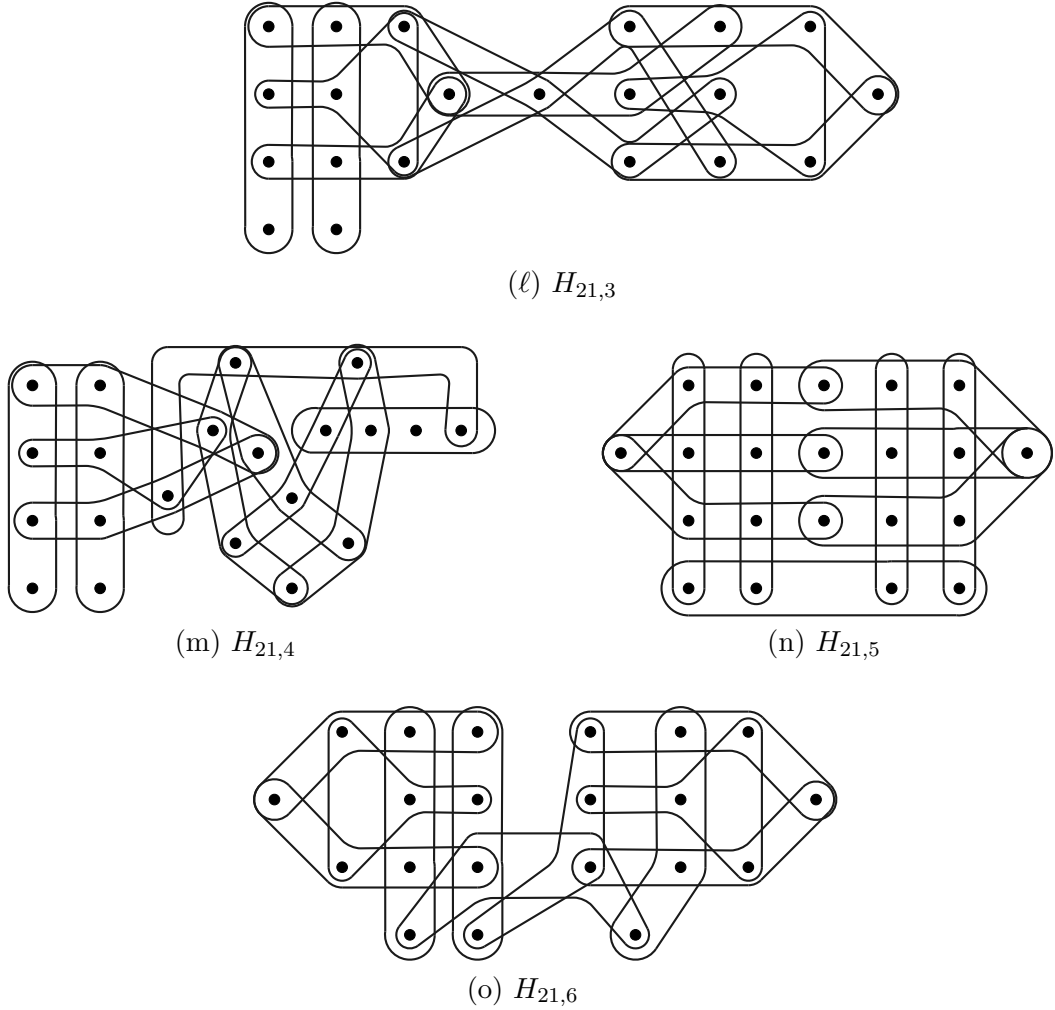
\caption{The fifteen special hypergraphs.} \label{f:special}
\end{figure}


We define
\[
\begin{array}{lcl}
\cH_{14} & = &  \{H_{14,1},H_{14,2},H_{14,3},H_{14,4},H_{14,5},H_{14,6}\} \, \, \mbox{and} \1 \\
 \cH_{21} & =&  \{H_{21,1},H_{21,2},H_{21,3},H_{21,4},H_{21,5},H_{21,6}\}.
 \end{array}
 \]

The following properties of special hypergraphs, which we have verified by computer, will prove to be useful.

\begin{ob} \label{property:special}
Let $H$ be a special hypergraph of order~$\nH$ and size~$\mH$, and let $u$ and $v$ be two arbitrary distinct vertices in $H$.
Then the following hold. \\[-27pt]
\begin{enumerate}
\item[{\rm (a)}]  If $H = H_4$, then $\nH = 4$, $\mH = 1$ and $\tau(H) = 1$. \4
\item[{\rm (b)}] If $H = H_{10}$, then $\nH = 10$, $\mH = 5$ and $\tau(H) = 3$. \4
\item[{\rm (c)}] If $H = H_{11}$, then $\nH = 11$, $\mH = 5$ and $\tau(H) = 3$. \4
\item[{\rm (d)}] If $H \in \cH_{14}$, then $\nH = 14$, $\mH = 7$ and $\tau(H) = 4$. \4
\item[{\rm (e)}] If $H \in \cH_{21}$, then $\nH = 21$, $\mH = 11$ and $\tau(H) = 6$. \4
\item[{\rm (f)}] If $H \in \{H_{10}, H_{14,5},H_{14,6}\}$, then $H$ is $2$-regular. \4
\item[{\rm (g)}] Any given vertex in $H$ belongs to some $\tau(H)$-transversal. \4
\item[{\rm (h)}] If $H \in \{H_{10},H_{14,6}\}$, then there exists a $\tau(H)$-transversal that contains both $u$ and $v$.\4
\item[{\rm (i)}] If $H \ne H_4$, $T \subset V(H)$ and $|T| = 3$, then there exists a $\tau(H)$-transversal that contains at least two vertices of $T$, unless $H = H_{11}$ and $T$ is the set of three vertices in $H$ that have no neighbor of degree~$1$ in $H$. \4
\item[{\rm (j)}] If $H \ne H_4$, $T \subset V(H)$ and $|T| = 4$, then there exists a $\tau(H)$-transversal that contains  at least two vertices of $T$. \4
\item[{\rm (k)}] If $H \ne H_4$ and $T_1$ and $T_2$ are vertex-disjoint subsets of $H$  such that $|T_1| = |T_2| = 2$, then there exists a $\tau(H)$-transversal that contains a vertex from $T_1$ and a vertex from~$T_2$. \4
\item[$(\ell)$] If $T_1$ and $T_2$ are vertex-disjoint subsets of $H$  such that $|T_1| = 3$ and $|T_2| = 1$ and $T_1$ contains two vertices that are not adjacent in $H$, then there exists a $\tau(H)$-transversal that contains a vertex from $T_1$ and a vertex from $T_2$. \4
\item[{\rm (m)}] If $T_1$ and $T_2$ are vertex-disjoint subsets of $H$ such that $|T_1| = 1$ and $|T_2| = 2$ where $T_2$ contains two vertices that are not adjacent in $H$, then there exists a $\tau(H)$-transversal that contains a vertex from $T_1$ and a vertex from $T_2$, except if $H = H_{11}$, and one degree-$1$ vertex in $H$ belongs to $T_1$ and the other degree-$1$ vertex to $T_2$, and the second vertex of $T_2$ is adjacent to the vertex of $T_1$. \4
\item[{\rm (n)}] If $T_1$, $T_2$ and $T_3$ are vertex-disjoint subsets of $H$ such that $|T_1| = |T_2| = 3$ and $ |T_3| \ge 2$, and $T_1$, $T_2$ and $T_3$ are independent sets in $H$, then there exists a $\tau(H)$-transversal that contains a vertex from each of $T_1$, $T_2$ and $T_3$. \4
\item[{\rm (o)}] If $H$ is a special subhypergraph of a $4$-uniform linear hypergraph $F$ where $\Delta(F) \le 3$, and if there are three edges of $F$ each of which intersect $H$ in two vertices, then there is a $\tau(H)$-transversal that covers one of these edges and any specified edge of $H$. \4
\item[{\rm (p)}] If $H = H_{11}$ or if $H \in \cH_{14} \cup \cH_{21}$, and if $v$ is a vertex of degree~$2$ in $H$, then either $H - v$ is connected or $H - v$ is disconnected with exactly two components, one of which consists of an isolated vertex. Further, $H - y$ does not contain two vertex disjoint copies of $H_4$ that are both intersected by a common edge and such that both copies of $H_4$ have three vertices of degree~$1$ and one vertex of degree~$2$.
\end{enumerate}
\end{ob}

\subsection{Known Results and Observations}

We shall need the following theorem of Berge~\cite{Berge} about the matching number of a graph, which is sometimes referred to as the Tutte-Berge formulation for the matching number.

\begin{thm}{\rm (Tutte-Berge Formula)}
For every graph $G$,
\[
\alpha'(G) = \min_{X \subseteq V(G)} \frac{1}{2} \left( |V(G)| +
|X| - \oc(G - X)  \right),
\]
where $\oc(G-X)$ denotes the number of odd components of $G-X$.
\label{Berge}
\end{thm}

We shall also rely heavily on the following well-known theorem due to K\"{o}nig~\cite{k7} and Hall~\cite{Ha35} in 1935.

\begin{thm}{\rm (Hall's Theorem)}
Let $G$ be a bipartite graph with partite sets $X$ and $Y$. Then $X$ can be matched to a subset of $Y$ if and only if $|N(S)| \ge |S|$ for every nonempty subset $S$ of $X$.
\label{Hall}
\end{thm}

\subsection{The Deficiency of a Set}
\label{S:defic}

Let $H$ be a $4$-uniform hypergraph. A set $X$ is a \emph{special $H$-set} if it consists of subhypergraphs of $H$ with the property that every subhypergraph in $X$ is a special hypergraph and further these special hypergraphs are pairwise vertex disjoint. For notational simplicity, we write $V(X)$ and $E(X)$ to denote the set of all vertices and edges, respectively, in $H$ that belong to a subhypergraph $H' \in X$ in the special $H$-set $X$.
Let $X$ be an arbitrary special $H$-set.

A set $T$ of vertices in $V(X)$ is an $X$-\emph{transversal} if $T$ is a minimum set of vertices that intersects every edge from every subhypergraph in $X$.

We define $E_H^*(X)$ to be the set of all edges in $H$ that do not belong to a subhypergraph in $X$ but which intersect at least one subhypergraph in $X$. Hence if $e \in E^*(X)$, then $e \notin E(H')$ for every subhypergraph $H' \in X$ but $V(e) \cap V(H') \ne \emptyset$ for at least one subhypergraph $H' \in X$. If the hypergraph $H$ is clear from context, we simply write $E^*(X)$ rather than $E_H^*(X)$ .

We associate with the set $X$ a bipartite graph, which we denote by $G_X$, with partite sets $X$ and $E^*(X)$, where an edge joins $e \in E^*(X)$ and $H' \in X$ in $G_X$ if and only if the edge $e$ intersects the subhypergraph $H'$ of $X$ in~$H$.

We define a \emph{weak partition} of $X = (X_4,X_{10},X_{11},X_{14},X_{21})$ (where a weak partition is a partition in which some of the sets may be empty) where $X_i \subseteq X$ consists of all subhypergraphs in $X$ of order~$i$, $i \in \{4,10,11,14,21\}$. Thus, $X = X_4 \cup X_{10} \cup X_{11} \cup X_{14} \cup X_{21}$ and $|X| = |X_4| + |X_{10}| + |X_{11}| + |X_{14}| + |X_{21}|$. As an immediate consequence of Observation~\ref{property:special}(a)--(e),
we have the following result.

\begin{ob}
If $X$ is a special $H$-set and $T$ is an $X$-transversal, then
\[
|T| = |X_4| + 3|X_{10}| + 3|X_{11}| + 4|X_{12}| + 6|X_{21}|.
\]
 \label{ob:special1}
\end{ob}

\vskip -1cm
We define the \emph{deficiency} of $X$ in $H$ as
\[
\defic_{H}(X) = 10|X_{10}| + 8|X_4| + 5|X_{14}| + 4|X_{11}| + |X_{21}| - 13|E^*(X)|.
\]
We define the \emph{deficiency} of $H$ by
\[
\defic(H) = \max \defic_H(X)
\]
where the maximum is taken over all special $H$-sets $X$. We note that taking $X = \emptyset$, we have $\defic(H) \ge 0$.

%
%

\subsection{Key Theorem}

Recall that $\cH_4$ is the class of all $4$-uniform linear hypergraphs with maximum degree at most~$3$.
We shall prove the following key result that we will need when proving our main theorem.

\begin{thm}
 \label{t:thm_key}
If $H \in \cH_4$, then $45\tau(H) \le 6n(H) + 13m(H) + \defic(H)$.
\end{thm}
\proof For a $4$-uniform hypergraph $H$, let
\[
\xi(H) = 45\tau(H) - 6n(H) - 13m(H) - \defic(H).
\]
We wish to show that if $H \in \cH_4$, then $\xi(H) \le 0$. Suppose, to the contrary, that the theorem is false and that $H \in \cH_4$ is a counterexample with minimum value of $n(H)+m(H)$. Thus, $\xi(H) > 0$ but every hypergraph $H' \in \cH_4$ with $n(H') + m(H') < n(H) + m(H)$ satisfies $\xi(H') \le 0$. We will now prove a number of claims, where the last claim completes the proof of the theorem.

\ClaimX{A}
The hypergraph $H$ is connected and $\delta(H) \ge 1$.

\ProofClaimX{A}  If $H$ is disconnected, then by the minimality of $H$ we have that the theorem holds for all components of $H$ and therefore
also for $H$, a contradiction. Therefore, $H$ is connected. If $H$ has an isolated vertex, then by the connectivity of $H$ we have that $n(H) = 1$ and $m(H) = 0$, implying that
$45\tau(H) = 0 < 6 = 6n(H) + 13m(H) + \defic(H)$, a contradiction. Hence, $\delta(H) \ge 1$.~\smallqed

\ClaimX{B}
Given a special $H$-set, $X$, there is no $X$-transversal, $T$, such that $|T| = |X_4| + 3|X_{10}| + 3|X_{11}| + 4|X_{14}| + |X_{21}|$ and $T$ intersects every edge in $E^*(X)$.

\ProofClaimX{B} Suppose that there does exist an $X$-transversal, $T$, such that $|T| = |X_4| + 3|X_{10}| + 3|X_{11}| + 4|X_{14}| + |X_{21}|$ and $T$ intersects every edge in $E^*(X)$. Let $H'$ be obtained from $H$ by removing all vertices in $V(X)$ and removing all edges of $H$ that intersect $V(X)$. Since $T$ is an $X$-transversal and $T$ intersects every edge in $E^*(X)$, we remark that $H' = H(T, V(X) \setminus T)$. Then, $H'$ is a $4$-uniform hypergraph with maximum degree $\Delta(H) \le 3$. Further, $n(H') = n(H) - 4|X_4| - 10|X_{10}| - 11|X_{11}| - 14|X_{14}| - 21|X_{21}|$, $m(H') = m(H) - |X_4| - 5|X_{10}| - 5|X_{11}| - 7|X_{14}| - 11|X_{21}| - |E^*(X)|$, and $\tau(H) \le |T| + \tau(H')$. Let
\[
\defic(H') = \defic_{H'}(X')
\]

\noindent for some special $H'$-set, $X'$. Let $X^* = X \cup X'$. Then,
$\defic(H') = \defic_{H}(X^*) - \defic_{H}(X)$.
By the minimality of $n(H) + m(H)$, $H'$ is not a counterexample to our theorem, and so $45\tau(H') \le 6n(H') + 13m(H') + \defic(H')$. Hence,

\[
\begin{array}{lcl}
45\tau(H) & \le & 45\tau(H') + 45|T|
\2 \\
& \le & (6n(H') + 13m(H') + \defic(H')) + \\
& & \hspace*{1cm} 45(|X_4| + 3|X_{10}| + 3|X_{11}| + 4|X_{14}| + 6|X_{21}|) \2 \\
& \le & 6(n(H) - 4|X_4| - 10|X_{10}| - 11|X_{11}| - 14|X_{14}| - 21|X_{21}|) + \\
& & \hspace*{0.5cm} 13(m(H) - |X_4| - 5|X_{10}| - 5|X_{11}| - 7|X_{14}| - 11|X_{21}| - |E^*(X)|) + \\
& & \hspace*{1cm} \defic_{H}(X^*) - \defic_{H}(X) + 45(|X_4| + 3|X_{10}| + 3|X_{11}| + 4|X_{14}| + 6|X_{21}|) \2 \\
 & = & 6n(H) + 13m(H) + 8|X_4| + 10|X_{10}| + 4|X_{11}| + 5|X_{14}|  + |X_{21}| - 13|E^*(X)| +   \\
& & \hspace*{0.5cm} \defic_{H}(X^*) - \defic_{H}(X) \2 \\
& = & 6n(H) + 13m(H) + \defic_{H}(X^*) \1 \\
& \le & 6n(H) + 13m(H) + \defic(H),
\end{array}
\]
contradicting the fact that $H$ is a counterexample.~\smallqed

\medskip
As an immediate consequence of Claim~A and Claim~B, $H$ is not a special hypergraph. We state this formally as follows.

\ClaimX{C} $H$ is not a special hypergraph.

Among all special non-empty $H$-sets, let $X$ be chosen so that

\hspace*{1cm} (1) $|E^*(X)| - |X|$ is minimum.\\
\indent \hspace*{1cm} (2) Subject to (1), $|X|$ is maximum.


\ClaimX{D}
$|E^*(X)| \ge |X| + 1$.

\ProofClaimX{D} Suppose, to the contrary, that $|E^*(X)| \le |X|$. Let $G_X$ be the bipartite graph associated with the special $H$-set $X$. Suppose there exists a matching, $M$, in $G_X$ that matches $E^*(X)$ to a subset of $X$. Then, by Observation~\ref{property:special}(g), there exists a minimum $X$-transversal, $T$, that intersects every edge in $E^*(X)$, contradicting Claim~B. Therefore, no matching in $G_X$ matches $E^*(X)$ to a subset of $X$. By Hall's Theorem, there is a nonempty subset $S \subseteq E^*(X)$ such that $|N_{G_X}(S)|<|S|$. We now consider the special $H$-set, $X' = X \setminus N_{G_X}(S)$. Then, $|X'| = |X| - |N_{G_X}(S)|$ and $|E^*(X')| = |E^*(X)| - |S|$. Thus,
\[
\begin{array}{lcl}
|E^*(X')| - |X'| & =  & (|E^*(X)| - |S|) - (|X| - |N_{G_X}(S)|) \\
& = & (|E^*(X)| - |X|) + (|N_{G_X}(S)| - |S|) \\
& < & |E^*(X)| - |X|,
\end{array}
\]
contradicting our choice of the special $H$-set $X$.~\smallqed

\medskip
As an immediate consequence of Claim~D and by our choice of the special $H$-set $X$, we have the following claim.


\ClaimX{E}
If $X' \ne \emptyset$ is a special $H$-set, then $|E^*(X')| \ge |X'| + 1$.


\ClaimX{F}
$\defic(H) = 0$.

\ProofClaimX{F} Let $X^*$ be a special $H$-set. If $X^* \ne \emptyset$,  then by Claim~E, $|E^*(X^*)| \ge |X^*| + 1$, implying that $\defic_H(X^*) \le 10|X^*| - 13|E^*(X^*)| < 0$. However, taking $X^* = \emptyset$, we note that  $\defic_H(X^*) = 0$. Consequently, $\defic(H) = 0$.~\smallqed

For a hypergraph $H' \in \cH_4$, let
\[
\Phi(H') = \xi(H') - \xi(H).
\]

\ClaimX{G}
If $H' \in \cH_4$ satisfies $n(H') + m(H') < n(H) + m(H)$, then $\Phi(H') < 0$.

\ProofClaimX{G} Since $H$ is a counterexample with minimum value of $n(H)+m(H)$, and since $n(H') + m(H') < n(H) + m(H)$, we note that $\xi(H') \le 0$. Thus,
\[
\begin{array}{lcl}
45\tau(H) & = & \xi(H) + 6n(H) + 13m(H) + \defic(H) \1 \\
& = & \xi(H') - \Phi(H') + 6n(H) + 13m(H) + \defic(H) \1 \\
& \le & - \Phi(H') + 6n(H) + 13m(H) + \defic(H), \\
\end{array}
\]
implying that  $\Phi(H') \le 6n(H) + 13m(H) + \defic(H) - 45\tau(H) = - \xi(H) < 0$.~\smallqed

\ClaimX{H}
$|E^*(X)| \ge |X| + 2$.

\ProofClaimX{H} Suppose, to the contrary, that $|E^*(X)| \le |X| + 1$. Then, by Claim~D, $|E^*(X)| = |X| + 1$. We define an edge $e \in E^*(X)$ to be $X$-\emph{universal} if there exists a special subhypergraph $F \in X$,
such that $e$ intersects $F$ and for every other edge, say $f$, that intersects $F$ there exists a $\tau(F)$-transversal that covers both $e$ and $f$. We proceed further with the following series of subclaims.

\ClaimX{H.1}
If $e \in E^*(X)$, then there exists an $X$-transversal that intersects every edge in $E^*(X) \setminus \{e\}$. Furthermore, no edge in $E^*(X)$ is $X$-universal.


\ProofClaimX{H.1}  Let $G_X$ be the bipartite graph associated with the special $H$-set $X$ and consider the graph $G_X - e$, where $e$ is an arbitrary vertex in $E^*(X)$.
We note that the partite sets of $G_X - e$ are $E^*(X) \setminus \{e\}$ and $X$ and these sets have equal cardinalities.
Suppose that $G_X - e$ does not have a perfect matching.
By Hall's Theorem, there is a nonempty subset $S \subseteq X$ such that in the graph $G_X - e$, we have $|N(S)| < |S|$.
(We remark that here $N(S) \subset E^*(X) \setminus \{e\}$.)
Since the edge $e$ may possibly intersect a subhypergraph of $S$ in $H$, we note that in the graph $G_X$, $|N(S)| \le |S|$.
However, in the graph $G_X$, $|E^*(S)| = |N(S)|$, implying that $|E^*(S)| \le |S|$.
By our choice of the special $H$-set $X$, $|E^*(X)| - |X| \le |E^*(S)| - |S| \le 0$, contradicting Claim~D.
Therefore, $G_X - e$ has a perfect matching.
By Observation~\ref{property:special}(g), this implies the existence of a minimum $X$-transversal, $T$, that intersects every edge in $E^*(X) \setminus \{e\}$.
This proves the first part of Claim~H.1.

For the sake of contradiction, suppose that $e \in E^*(X)$ is a universal edge. Let $F \in X$ be a special subhypergraph of $X$ that is intersected by $e$, such that
for every $f \in E^*(X)$ intersecting $F$ there exists a $\tau(F)$-transversal that covers both $e$ and $f$.
By the above argument there exists a perfect matching, $M$, in $G_X - e$. Let $g$ be the edge that is matched to $F$ in $M$.
By definition there is a $\tau(F)$-transversal intersecting both $e$ and $g$. Due to the matching $M \setminus \{g F\}$ and Observation~\ref{property:special}(g), we therefore obtain
an $X$-transversal, $T$, that intersects every edge in $E^*(X)$, a contradiction to Claim~B. This proves the second part of Claim~H.1.~\smallqed

\ClaimX{H.2}
$|X_{10}| = 0$.

\ProofClaimX{H.2}  Suppose, to the contrary, that $|X_{10}| > 0$.  Let $H'$ be a copy of $H_{10}$ that belongs to the special $H$-set $X$. By Claim~E, at least two edges in $E^*(X)$ intersect $H'$ in $H$. Let $e$ be such an edge of $E^*(X)$ that intersects $H'$. Let $G_X$ be the bipartite graph associated with $X$ and consider the graph $G_X - e$. We note that the partite sets of $G_X - e$ are $E^*(X) \setminus \{e\}$ and $X$ and these sets have equal cardinalities. Suppose that $G_X - e$ does not have a perfect matching. By Hall's Theorem, there is a nonempty subset $S \subseteq E^*(X) \setminus \{e\}$ such that in the graph $G_X - e$ (and in the graph $G_X$), we have $|N(S)| < |S|$. We remark that here $N(S) \subset X$. We now consider the special $H$-set, $X' = X \setminus N(S)$. Then, $X' \ne \emptyset$ and
\[
\begin{array}{lcl}
|E^*(X')| & \le  & |E^*(X)| - |S| \\
& \le & |E^*(X)| - (|N(S)| + 1) \\
& = & |X| + 1 - |N(S)| - 1 \\
& = & |X'|,
\end{array}
\]
contradicting Claim~E. Therefore, $G_X - e$ has a perfect matching, $M$ say. Let $e'$ be the vertex of $E^*(X)$ that is $M$-matched to $H'$ in $G_X - e$. Let $u$ and $v$ be vertices in the subhypergraph $H'$ that belong to the edges $e$ and $e'$, respectively, in $H$. By Observation~\ref{property:special}(h), there exists a $\tau(H')$-transversal that contains both $u$ and $v$.
If $H''$ is a subhypergraph in $X$ different from $H'$ and if $e''$ is the vertex of $E^*(X)$ that is $M$-matched to $H''$ in $G_X - e$, then by Observation~\ref{property:special}(g) there exists a $\tau(H'')$-transversal that contains a vertex of $H''$ that belongs to $e''$ in $H$. This implies the existence of a minimum $X$-transversal, $T$, that intersects every edge in $E^*(X)$, contradicting Claim~B.~\smallqed

\medskip
Recall that the boundary of a set $S$ of vertices in $H$, denoted $\partial_H(S)$ or simply $\partial(S)$ if $H$ is clear from context, is the set $N(S) \setminus S$. Abusing notation, we write $\partial(X)$, rather than $\partial(V(X))$, as the boundary of the set $V(X)$. Let $X'$ be a set of new vertices (not in $H$), where $|X'| = \max(0,4 - |\partial(X)|)$. We note that $|X' \cup \partial(X)| \ge 4$. Further if $|\partial(X)| \ge 4$, then $X' = \emptyset$. Let $H'$ be the hypergraph obtained from $H - V(X)$ by adding the set $X'$ of new vertices and adding a $4$-edge, $e$, containing four vertices in $X' \cup \partial(X)$. Note that $H'$ may not be linear as the edge $e$ may overlap other edges in $H'$.

\ClaimX{H.3}
Either $\defic(H') = 0$ or $\defic(H') = \defic_{H'}(Y)$ where $|Y| = 1$ and $e$ is an edge of the hypergraph in the special $H'$-set $Y$.

\ProofClaimX{H.3} Suppose that $\defic(H') > 0$. Let $Y$ be a special $H'$-set such that $\defic(H') = \defic_{H'}(Y)$. Since $\defic(H') > 0$, we note that in $H'$, $|E^*(Y)| < |Y|$. Suppose that $e \notin E(Y)$. We now consider the special $H$-set $Q = X \cup Y$. Then in $H$, $|E^*(Q)| \le |E^*(X)| + |E^*(Y)| \le (|X| + 1) + (|Y| - 1) = |X| + |Y| = |Q|$, contradicting Claim~E. Therefore, $e \in E(Y)$. Let $e \in E(R)$, where $R \in Y$. Suppose that $|Y| \ge 2$. In this case, we consider the special $H$-set $Q = X \cup (Y \setminus \{R\})$. Then in $H$, $|E^*(Q)| \le |E^*(X)| + |E^*(Y)| \le (|X| + 1) + (|Y| - 1) = |X| + |Y| = |Q| + 1$. By Claim~E, $|E^*(Q)| \ge |Q| + 1$. Consequently, $|E^*(Q)| = |Q| + 1$. However since $|Y| \ge 2$, we note that $|Q| > |X|$, contradicting our choice of the special $H$-set $X$.~\smallqed


\ClaimX{H.4}
$\Phi(H') \ge - 8|X_4| - 5|X_{14}| - 4|X_{11}| - |X_{21}|  - 6|X'| +  13|X|  - \defic(H')$. Furthermore if $H'$ is linear, then both of the following statements hold. \1 \\
\indent \textbf{(a):} $8|X_4| + 5|X_{14}| + 4|X_{11}| + |X_{21}|  > 13|X| - 6|X'| - \defic(H') \ge 13|X| - 6|X'| - 10$. \1 \\
\indent \textbf{(b):} If $X' \ne\emptyset$, then $8|X_4| + 5|X_{14}| + 4|X_{11}| + |X_{21}|  > 13|X| - 6|X'| - 8$.

\ProofClaimX{H.4} By Claim~F, $\defic(H) = 0$. Thus,
\[
\Phi(H') = 45(\tau(H') - \tau(H)) - 6(n(H') - n(H)) - 13(m(H') - m(H)) - \defic(H').
\]

We will first show that $\tau(H) \le \tau(H')+|X_4|+3|X_{10}|+3|X_{11}|+4|X_{14}|+6|X_{21}|$.
Let $T'$ be a minimum transversal in $H'$. If some vertex of $X'$ belongs to $T'$, then we can simply replace such a vertex in $T'$ with a vertex from $V(e) \cap \partial(X)$ since a vertex in $X'$ belongs to the edge $e$ but no other edge of $H'$. We may therefore assume that some vertex in $\partial(X)$ belongs to $T'$, and hence $T'$ covers at least one edge from $E^*(X)$. By Claim~H.1, we can now obtain a transversal of $H$ of size
$\tau(H')+|X_4|+3|X_{10}|+3|X_{11}|+4|X_{14}|+6|X_{21}|$, which implies that  $\tau(H) \le \tau(H')+|X_4|+3|X_{10}|+3|X_{11}|+4|X_{14}|+6|X_{21}|$ as desired.

Let $F \in X$. By Claim~H.2, $X_{10} = \emptyset$.
%
By Observation~\ref{property:special}, if $F \in X_{4}$, then $F$ contributes $4$ to the sum $n(H)  - n(H')$, $1$ to the sum $m(H) - m(H')$, and at most~$1$ to the sum $\tau(H) - \tau(H')$ (due to the above bound on $\tau(H)$), and therefore contributes at least $45 \cdot (-1) - 6 \cdot (-4) - 13 \cdot (-1) = - 8$ to $\Phi(H')$.
Similar arguments show that if $F \in X_{14}$, then $F$ contributes at least $-5$ to $\Phi(H')$. If  $F \in X_{11}$, then $F$  contributes at least $-4$ to $\Phi(H')$, while if $F \in X_{21}$, then $F$ contributes at least $-1$ to $\Phi(H')$. The edges in $E^*(X)$ contribute~$|E^*(X)| = |X| + 1$ to the sum $m(H) - m(H')$ and therefore $13(|X| + 1)$ to $\Phi(H')$, while the added edge $e$ contributes~$1$ to the sum $m(H') - m(H)$ and therefore $-13$ to $\Phi(H')$. If $X' \ne \emptyset$, then each vertex in $X'$ contributes~$1$ to the sum $n(H') - n(H)$, and therefore the vertices in $X'$ contribute~$-6|X'|$ to  $\Phi(H')$. This proves the first part of Claim~H.4.

Suppose, next, that $H'$ is linear. Then,  by Claim~G, $\Phi(H') < 0$. Thus, from our previous inequality established in the first part of Claim~H.4, $0 > \Phi(H') \ge - 8|X_4| - 5|X_{14}| - 4|X_{11}| - |X_{21}|  - 6|X'| +  13|X|  - \defic(H')$. This immediately implies part~(a) as $\defic(H') \le 10$ by Claim~H.3. If  $X' \ne\emptyset$, then the edge $e$ contains a degree-$1$ vertex and is therefore not part of a $H_{10}$ subhypergraph in $H'$.
Part~(b) now again follows from Claim~H.3, noting that in this case $\defic(H') \le 8$.~\smallqed

\ClaimX{H.5}
Suppose the added edge $e$ overlaps with some other edge in $H'$.
\1 \\
\indent \textbf{(a):} If $e$ contains a vertex of degree~$1$ in $H'$, then both of the following statements hold. \\
\hspace*{1cm} \textbf{(i):} $\Phi(H') \le 3 - \defic(H')$. \\
\hspace*{1cm} \textbf{(ii):} $8|X_4| + 5|X_{14}| + 4|X_{11}| + |X_{21}|  \ge  13|X| - 6|X'| - 3$.
\1 \\
\indent \textbf{(b):} If $e$ contains no vertex of degree~$1$ in $H'$, then both of the following statements hold. \\
\hspace*{1cm} \textbf{(i):} $\Phi(H') \le 7 - \defic(H')$. \\
\hspace*{1cm} \textbf{(ii):} $8|X_4| + 5|X_{14}| + 4|X_{11}| + |X_{21}|  \ge  13|X| - 6|X'| - 7$.


\ProofClaimX{H.5} Suppose that the added edge $e$ overlaps with some other edge $e'$ in $H'$. We will first show the following.
\[
 \Phi(H') \le \left\{ \begin{array}{ll}
 3 - \defic(H') & \mbox{if $e$ contains a vertex of degree~$1$ in $H'$} \\
 7 - \defic(H') & \mbox{if $e$ contains no vertex of degree~$1$ in $H'$.}
 \end{array}
 \right.
\]
Suppose, to the contrary, that $\Phi(H') \ge 4 - \defic(H')$ if $e$ contains a vertex of degree~$1$ in $H'$ and that $\Phi(H') \ge 8 - \defic(H')$ if $e$ contains no vertex of degree~$1$ in $H'$. Let $\{x,y\} \subseteq e \cap e'$.

Suppose firstly that $d_{H'}(x) = 3$.
Let $e''$ be the edge incident with $x$ different from $e$ and $e'$, and consider the hypergraph $H^{\star} = H' - x$.
Since $H \in \cH_4$, by construction, $H^{\star} \in \cH_4$.
Then, $m(H^{\star}) = m(H') - 3$ and $\tau(H') \le \tau(H^{\star}) + 1$.
Further, if $e$ contains a vertex of degree~$1$ in $H'$, then $n(H^{\star}) \le n(H') - 2$, while if $e$ contains no vertex of degree~$1$ in $H'$, then $n(H^{\star}) \le n(H') - 1$.
Recall that $\Phi(H') = \xi(H') - \xi(H)$ and $\Phi(H^{\star}) = \xi(H^{\star}) - \xi(H)$. By Claim~G, $\Phi(H^{\star}) < 0$. Thus,

\[
\begin{array}{lcl}
0 & > & \Phi(H^{\star}) \1 \\
& = & \xi(H^{\star}) - \xi(H) \1 \\
& = & \xi(H^{\star}) + \Phi(H') - \xi(H') \1 \\
& = & \Phi(H') - 45(\tau(H') - \tau(H^{\star})) + 6(n(H') - n(H^{\star})) + \\
& & \, \, 13(m(H') - m(H^{\star})) + (\defic(H') - \defic(H^{\star})) \1
\\
& \ge & \Phi(H') - 45 + 6(n(H') - n(H^{\star})) + 13 \cdot 3 + \defic(H') - \defic(H^{\star}) \1
\\
& \ge & \Phi(H') - 6 + 6(n(H') - n(H^{\star})) + \defic(H') - \defic(H^{\star}). 
\end{array}
\]

\ClaimX{H.5.1}
If $e$ contains a vertex of degree~$1$ in $H'$, then $\defic(H^{\star}) \ge 11$.

\ProofClaimX{H.5.1} If $e$ contains a vertex of degree~$1$ in $H'$, then $n(H') - n(H^{\star}) \ge 2$ and, by supposition, $\Phi(H') \ge 4 - \defic(H')$. Thus, in this case, the above inequality chain simplifies to
$0 > \Phi(H^{\star}) \ge (4 - \defic(H')) - 6 + 6 \cdot 2 + \defic(H') - \defic(H^{\star}) = 10 - \defic(H^{\star})$,
and so, $\defic(H^{\star}) > 10$.~\smallqed

\ClaimX{H.5.2}
If $e$ contains no vertex of degree~$1$ in $H'$, then  $\defic(H^{\star}) \ge 9$.

\ProofClaimX{H.5.2} If $e$ contains no vertex of degree~$1$ in $H'$, then $n(H') - n(H^{\star}) \ge 1$ and, by supposition, $\Phi(H') \ge 8 - \defic(H')$. Thus, in this case, our inequality chain simplifies to
$0 > \Phi(H^{\star}) \ge (8 - \defic(H')) - 6 + 6 \cdot 1 + \defic(H') - \defic(H^{\star}) = 8 - \defic(H^{\star})$,
and so, $\defic(H^{\star}) > 8$.~\smallqed

By Claim~H.5.1 and Claim~H.5.2, $\defic(H^{\star}) \ge 9$. Let $Y$ be a special $H^{\star}$-set such that $\defic(H^{\star}) = \defic_{H^{\star}}(Y)$.

\ClaimX{H.5.3}
$|Y| = 1$.

\ProofClaimX{H.5.3} Suppose, to the contrary, that $|Y| \ge 2$. If $|E^*(Y)| \ge |Y| - 1$, then $\defic_{H^{\star}}(Y) \le 10|Y| - 13|E^*(Y)| \le 10|Y| - 13(|Y|-1) = -3|Y| + 13 \le 7$, a contradiction. Hence, $|E^*(Y)| \le |Y| - 2$. We now consider the special $H$-set $Q = X \cup Y$. Then in $H$, $|E^*(Q)| \le |E^*(X)| + |E^*(Y)| + |\{e',e''\}| \le (|X| + 1) + (|Y| - 2) + 2 = |X| + |Y| + 1 = |Q| + 1$. By Claim~E, $|E^*(Q)| \ge |Q| + 1$. Consequently, $|E^*(Q)| = |Q| + 1$. However $|Q| > |X|$, contradicting our choice of the special $H$-set $X$. Therefore, $|Y| = 1$.~\smallqed

By Claim~H.5.3, $|Y| = 1$. Let $Y = \{R\}$. If $R \ne H_{10}$, then $\defic_{H^{\star}}(Y) \le 8|Y| - 13|E^*(Y)| \le 8$, a contradiction. Hence, $R = H_{10}$. If $|E^*(Y)| \ge 1$, then $\defic_{H^{\star}}(Y) = 10 - 13|E^*(Y)| \le -3$, a contradiction. Therefore, $R$ is a component of $H^{\star}$ and $\defic_{H^{\star}}(Y) = 10$. In particular, we note that the edge $e$ therefore contains no vertex of degree~$1$ in $H'$. By Observation~\ref{property:special}(f), $R$ is $2$-regular.

Recall that $\{x,y\} \subseteq e \cap e'$, implying that $d_{H'}(y) \ge 2$. If $d_{H'}(y) = 2$, then both $x$ and $y$ are removed from $H'$ when constructing $H^{\star}$, implying that $n(H') - n(H^{\star}) \ge 2$. In this case, our inequality chain in the proof of Claim~H.5.2 simplifies to $0 > \Phi(H^{\star}) \ge (8 - \defic(H')) - 6 + 6 \cdot 2 + \defic(H') - \defic(H^{\star}) = 14 - \defic(H^{\star})$, and so, $\defic(H^{\star}) > 14$, a contradiction. Therefore, $d_{H'}(y) = 3$, implying that $d_{H^{\star}}(y) = 1$ and $y \notin V(R)$.
%
Let $e'''$ be the edge incident with $y$ different from $e$ and $e'$. If $e'''$ intersects $R$, then since $R$ is a component of $H'$, this would imply $y \in V(R)$, a contradiction. Therefore, $e'''$ does not intersects $R$.

Interchanging the roles of $x$ and $y$, and considering the hypergraph $H^{\star \star} = H' - y$ (instead of $H^{\star} = H' - x$), there exists an $H_{10}$-component, say $R'$, in $H^{\star \star}$. Analogous arguments as employed earlier, show that $d_{H^{\star \star}}(x) = 1$, $x \notin V(R')$  and $e''$ does not intersect $R'$.  Thus, $R$ and $R'$ are $H_{10}$-components in the hypergraph, $H' - \{e,e',e'',e'''\}$, obtained from $H'$ by deleting the edges $e, e', e'', e'''$.

If $R = R'$, then neither $e''$ nor $e'''$ intersect $R$. Hence, $|E^*(X \cup R)| \le |E^*(X)| + |\{e'\}| \le (|X| + 1) + 1 = |X \cup R| + 1$. By Claim~E, $|E^*(X \cup R)| \ge |X \cup R| + 1$. Consequently, $|E^*(X \cup R)| = |X \cup R| + 1$. However $|X \cup R| > |X|$, contradicting our choice of the special $H$-set $X$. Therefore, $R \ne R'$, implying that $R$ and $R'$ are distinct and vertex-disjoint components of $H' - \{e,e',e'',e'''\}$.

As observed earlier, $e'''$ does not intersects $R$. Suppose that $e''$ does not intersect $R$. Then, $|E^*(X \cup R)| \le |E^*(X)| + |\{e'\}| \le (|X| + 1) + 1 = |X \cup R| + 1$, a contradiction. Therefore, $e''$ intersects $R$. Analogously, $e'''$ intersects $R'$.  Let $r'' \in e'' \cap V(R)$ and let $r''' \in e''' \cap V(R')$.

If $R$ and $R'$ contain no vertex in $\partial(X)$, then $|E*(R \cup R')| \le |\{e',e'',e'''\}| = 3 = |R \cup R'| + 1$, implying that $E*(R \cup R') = \{e',e'',e'''\}$ and $|E*(R \cup R')| = |R \cup R'| + 1$. Let $r' \in e' \cap V(R) \cup V(R')$. Without loss of generality, we may assume that $r' \in V(R)$. By Observation~\ref{property:special}(h), there exists a $\tau(R)$-transversal that contains both $r'$ and $r''$, and a $\tau(R')$-transversal that contains $r'''$. Thus, given the special $H$-set, $\{R,R'\}$, we can find an $\{R,R'\}$-transversal that covers every edge in $E^*(R \cup R'\})$, contradicting Claim~B. Therefore, $R$ or $R'$ (or both $R$ and $R'$) contain a vertex in $\partial(X)$. Let $z$ be such a vertex and let $e_z$ be an edge of $E^*(X)$ that contains~$z$. Without loss of generality, we may assume that $r' \in V(R)$. Applying Observation~\ref{property:special}(h) and Claim~H.1, we can find an $(X \cup \{R,R'\})$-transversal that covers every edge in $E^*(X \cup R \cup R')$, contradicting Claim~B. We deduce, therefore, that $d_{H'}(x) = 2$ (recall that in $H'$, $\{x,y\} \subseteq e \cap e'$). Analogously, all vertices in $e \cap e'$ have degree~$2$ in $H'$. In particular, $d_{H'}(y) = 2$.

We once again consider the hypergraph $H^{\star} = H - x$, noting that here both $x$ and $y$ are removed from $H'$ when constructing $H^{\star}$.
%
Since $H^{\star} \in \cH_4$, as before $\defic(H^{\star}) \ge 0$. Further, $m(H^{\star}) = m(H') - 2$ and $\tau(H') \le \tau(H^{\star}) + 1$.
If $e$ contains a vertex of degree~$1$ in $H'$, then $n(H') - n(H^{\star}) \ge 3$ and, by supposition, $\Phi(H') \ge 4 - \defic(H')$. Thus, in this case, our inequality chain simplifies to
$0 > \Phi(H^{\star}) \ge (4 - \defic(H')) - 45 + 6 \cdot 3 + 13 \cdot 2 + \defic(H') - \defic(H^{\star}) = 3 - \defic(H^{\star})$,
and so, $\defic(H^{\star}) \ge 4$.
 If $e$ contains no vertex of degree~$1$ in $H'$, then $n(H') - n(H^{\star}) \ge 2$ and, by supposition, $\Phi(H') \ge 8 - \defic(H')$. Thus, in this case, our inequality chain simplifies to
$0 > \Phi(H^{\star}) \ge (8 - \defic(H')) - 45 + 6 \cdot 2 + 13 \cdot 2 + \defic(H') - \defic(H^{\star}) = 1 - \defic(H^{\star})$,
and so, $\defic(H^{\star}) \ge 2$.  In both cases, $\defic(H^{\star}) > 0$.  Let $Y$ be a special $H^{\star}$-set such that $\defic(H^{\star}) = \defic_{H^{\star}}(Y)$. Since $\defic_{H^{\star}}(Y) > 0$, we note that, $|E^*(Y)| \le |Y| - 1$.
We now consider  the special $H$-set $X \cup Y$. In $H$, $|E^*(X \cup Y)| \le |E^*(X)| + |E^*(Y)| + |\{e'\}| \le (|X| + 1) + (|Y| - 1) + 1 = |X| + |Y| + 1 = |X \cup Y| + 1$, contradicting our choice of the special $H$-set $X$. This completes the proof of the following.
 \[
 \Phi(H') \le \left\{ \begin{array}{ll}
 3 - \defic(H') & \mbox{if $e$ contains a vertex of degree~$1$ in $H'$} \\
 7 - \defic(H') & \mbox{if $e$ contains no vertex of degree~$1$ in $H'$,}
 \end{array}
 \right.
 \]
thereby establishing part~(a)(i) and part~(b)(i). By Claim~H.4, $\Phi(H') \ge - 8|X_4| - 5|X_{14}| - 4|X_{11}| - |X_{21}|  - 6|X'| +  13|X|  - \defic(H')$. Thus, if $e$ contains  a vertex of degree~$1$ in $H'$, then
\[
 3 - \defic(H') \ge \Phi(H') \ge - 8|X_4| - 5|X_{14}| - 4|X_{11}| - |X_{21}|  - 6|X'| +  13|X|  - \defic(H').
\]
This immediately completes the proof of part~(a)(ii), as $\defic(H')$ cancels out. Part~(b) is easily proved analogously.~\smallqed

\ClaimX{H.6}
The case $|X| \ge 2$ and $|\partial(X)| \ge 4$ cannot occur.

\ProofClaimX{H.6}
Note that $8|X_4| + 5|X_{14}| + 4|X_{11}| + |X_{21}| \le 8|X| \le 13|X|-6|X'| -10$, as $|X| \ge 2$ and $|X'|=0$.
If $H'$ is linear, then we are obtain a contradiction to Claim~H.4(a), and if $H'$ is not linear we obtain a contradiction to Claim~H.5.~\smallqed


\ClaimX{H.7}
The case $|X| \ge 3$ and $|\partial(X)|  = 3$ cannot occur.

\ProofClaimX{H.7}
Note that $8|X_4| + 5|X_{14}| + 4|X_{11}| + |X_{21}| \le 8|X| \le 13|X|-6|X'| -9$, as $|X| \ge 3$ and $|X'|=1$.
If $H'$ is linear, then we are obtain a contradiction to Claim~H.4(b), and if $H'$ is not linear we obtain a contradiction to Claim~H.5.~\smallqed


\ClaimX{H.8}
The case $|X| = 2$ and $|\partial(X)|  = 3$ cannot occur.

\ProofClaimX{H.8} Suppose, to the contrary, that $|X|  = 2$ and $|\partial(X)|  = 3$. Then, $|E^*(X)| = 3$ and $|X'| = 1$.
If $H'$ is not linear, then we obtain a contradiction by Claim~H.5(a), as $8|X_4| + 5|X_{14}| + 4|X_{11}| + |X_{21}| \le 16$
and $13|X|-6|X'| - 3 = 17$. Therefore, $H'$ is linear.

\ClaimX{H.8.1} $|X_4| \le 1$.

\ProofClaimX{H.8.1}
Suppose, to the contrary, that $X$ consists of two (vertex-disjoint) copies of $H_4$.
By Claim~H.4(a), $16 = 8|X_4| + 5|X_{14}| + 4|X_{11}| + |X_{21}| > 13|X| - 6|X'| - \defic(H') = 26 - 6 -\defic(H')$, implying that $\defic(H') > 4$. If $e$ is a $H_4$-component in $H'$, then, by the connectivity and linearity of $H$, we note that $H=H_{11}$, contradicting Claim~C.
Hence, $H' \ne H_4$, implying by Claim~H.3 and the fact that $e$ contains a degree-$1$ vertex (and so, $Y \ne H_{10}$) that $Y \in X_{14}$ and $\defic(H') = \defic(Y) = 5$.
Since $e \in E(Y)$ and $e$ contains a vertex of degree~$1$ in $H'$, we note that $H'  \ne H_{14,5}$ and $H'  \ne H_{14,6}$.
Therefore, $H' = H_{14,i}$ for some $i \in [4]$, implying that $H = H_{21,i}$, contradicting Claim~C.~\smallqed

\ClaimX{H.8.2} $H' = H_4$, $|X_4| = 1$ and $|X_{14}| = 1$.

\ProofClaimX{H.8.2}
Recall that $H'$ is linear and $|X'|=1$.  By Claim~H.4(a), $8|X_4| + 5|X_{14}| + 4|X_{11}| + |X_{21}| > 13|X|-6|X'| - \defic(H') = 20-\defic(H')$. By Claim~H.8.1, $|X_4| \le 1$. By Claim~H.3 and the fact that $e$ contains a degree-$1$ vertex,  $\defic(H') \le 8$. Thus, $8 + 5 \ge 8|X_4| + 5|X_{14}| + 4|X_{11}| + |X_{21}| > 20-\defic(H')$, implying that $|X_4|=1$,  $|X_{14}|=1$ and $\defic(H')=8$, as desired.~\smallqed

\medskip
By Claim~H.8.2, $H' = H_4$, $|X_4| = 1$ and $|X_{14}| = 1$. 
Let $E^*(X) = \{e_1,e_2,e_3\}$. Let $F_1$ and $F_2$ be the (vertex-disjoint) hypergraphs that belong to $X$, where $F_1 = H_{14}$ and $F_2 = H_4$. Further, let $E(F_2) = \{f_2\}$. Let $\partial(X) = \{x,y,z\}$.

\ClaimX{H.8.3}
All three edges in $E^*(X)$ intersect $F_1$.

\ProofClaimX{H.8.3} Suppose, to the contrary, that at most two edges in $E^*(X)$ intersect $F_1$. Therefore, by Claim~E, exactly two edges in $E^*(X)$ intersect $F_1$. Renaming edges in $E^*(X)$ if necessary, we may assume that $e_1$ and $e_2$ intersect $F_1$, and that $e_3$ does not intersect $F_1$. By the linearity of $H$, the edge $e_3$ contains one vertex of $F_2$ and all three vertices of $\partial(X)$. This in turn implies by the linearity of $H$, that both edges $e_1$ and $e_2$ intersect $F_1$ in at least two vertices. Let $T_1 = e_1 \cap V(F_1)$ and $T_2 = e_2 \cap V(F_1)$. Thus, $|T_1| \ge 2$ and $|T_2| \ge 2$. If $e_1$ and $e_2$ intersect in a common vertex of $F_1$, then, by Observation~\ref{property:special}(g), we can cover $e_1$ and $e_2$ by a $\tau(F_1)$-transversal and we can cover $e_3$ by a $\tau(F_2)$-transversal, implying that there is an $X$-transversal that intersects every edge in $E^*(X)$, contradicting Claim~B. Hence, $T_1 \cap T_2 = \emptyset$. By Observation~\ref{property:special}(l), there exists a $\tau(F_1)$-transversal that contains one vertex from $T_1$ and one vertex from $T_2$. Once again, this implies that there is an $X$-transversal that intersects every edge in $E^*(X)$, contradicting Claim~B.~\smallqed

\medskip
By Claim~H.8.3, all three edges in $E^*(X)$ intersect $F_1$. If two edges of $E^*(X)$ intersect in a common vertex of $F_1$, then we can cover these two edges by a $\tau(F_1)$-transversal. If no two edges of $E^*(X)$ intersect in a common vertex of $F_1$, then by Observation~\ref{property:special}(i), we can cover two edges of $E^*(X)$ by a $\tau(F_1)$-transversal. In both cases, two edges of $E^*(X)$ can be covered by a $\tau(F_1)$-transversal. Renaming edges in $E^*(X)$ if necessary, we may assume that $e_1$ and $e_2$ can be covered by a $\tau(F_1)$-transversal. If $e_3$ intersects $F_2$, then there is an $X$-transversal that intersects every edge in $E^*(X)$, contradicting Claim~B. Hence, $e_3$ does not intersect $F_2$, implying that both $e_1$ and $e_2$ intersect $F_2$.

If $|e_3 \cap V(F_1)| \ge 3$, then by Observation~\ref{property:special}(l), we can cover $e_1$ and $e_3$ by a $\tau(F_1)$-transversal. Since we can cover $e_2$ by a $\tau(F_2)$-transversal, there is an $X$-transversal that intersects every edge in $E^*(X)$, contradicting Claim~B. Hence, $|e_3 \cap V(F_1)| \le 2$.

If $|(e_1 \cup e_2) \cap V(F_1)| = 1$, then both $e_1$ and $e_2$ contain one vertex from $F_1$, one vertex from $F_2$ and two vertices from $\partial(X)$. However since $|\partial(X)| = 3$, $e_1$ and $e_2$ would then overlap, a contradiction. Hence,  $|(e_1 \cup e_2) \cap V(F_1)| \ge 2$.
Thus, if $|e_3 \cap V(F_1)| = 2$, then by Observation~\ref{property:special}(k), there exists a $\tau(F_1)$-transversal that contains a vertex from $e_3 \cap V(F_1)$ and a vertex from $(e_1 \cup e_2) \cap V(F_1)$, implying that we can cover $e_3$ and one of $e_1$ and $e_2$  by a $\tau(F_1)$-transversal and we can cover the remaining edge in $E^*(X)$ by a $\tau(F_2)$-transversal, once again contradicting Claim~B. Therefore, $|e_3 \cap V(F_1)| = 1$, implying that the edge $e_3$ contains one vertex of $F_1$ and all three vertices in $\partial(X)$, namely $x$, $y$, and $z$. Thus, both $e_1$ and $e_2$ contain one vertex from $F_2$, one vertex from $\partial(X)$ and two vertices from $F_1$. In particular, we note that $e_1 \cap V(F_1)$ is an independent set in $F_1$, as is $e_2 \cap V(F_1)$. Further, since $e_1$ and $e_2$ do not overlap, $|(e_1 \cup e_2) \cap V(F_1)| \ge 3$. Thus, by Observation~\ref{property:special}(l), there exists a $\tau(F_1)$-transversal that contains a vertex from $e_3 \cap V(F_1)$ and a vertex from $(e_1 \cup e_2) \cap V(F_1)$, implying that we can cover $e_3$ and one of $e_1$ and $e_2$  by a $\tau(F_1)$-transversal and we can cover the remaining edge in $E^*(X)$ by a $\tau(F_2)$-transversal, once again contradicting Claim~B. This completes the proof of Claim~H.8.~\smallqed

\ClaimX{H.9}
If $|X| \ge 2$ and $|\partial(X)|  = 2$, then $H' \in \cH_4$.

\ProofClaimX{H.9} Let $|X| \ge 2$ and $|\partial(X)|  = 2$ and suppose, to the contrary, that $H' \notin \cH_4$. Thus, $H'$ is not a linear hypergraph, implying that the edge $e$ overlaps in $H'$ with some other edge, $e'$ say. Since $e$ contains two vertices of degree~$1$, namely the two vertices in $X'$, we note that in $H'$, $e \cap e' = \partial(X)$. We now consider the hypergraph $H'' = H - e'$. Since $H''$ has no overlapping edges, $H'' \in \cH_4$, and so by Claim~G, $\Phi(H'') < 0$.  By Claim~H.2, $|X_{10}| = 0$. An analogous proof to that of Claim~H.4, shows that
$\Phi(H'') \ge - 8|X_4| - 5|X_{14}| - 4|X_{11}| - |X_{21}|  - 6|X'| +  13(|X|+1)  - \defic(H'')$,
noting that the deleted edge $e'$ contributes~$1$ to the sum $m(H) - m(H'')$. Thus, since $|X'| = 2$ and $|X| = 2$,
\[
\begin{array}{lcl}
\Phi(H'') & \ge  & -8|X| - 12 + 13(|X|+1)  - \defic(H'')  \\
& = & 5|X| + 1 - \defic(H'') \\
& = & 11 - \defic(H'').
\end{array}
\]

If $\defic(H'') \le 11$, then $\Phi(H'') \ge 0$, a contradiction. Hence, $\defic(H'') \ge 12$.
Let $Y$ be a special $H''$-set such that $\defic(H'') = \defic_{H''}(Y)$. If $|Y| = 1$, then $\defic_{H''}(Y) \le 10$, a contradiction. Hence, $|Y| \ge 2$. If $|E^*(Y)| \ge |Y| - 1$ in $H''$, then $\defic_{H^{\star}}(Y) \le 10|Y| - 13|E^*(Y)| \le 10|Y| - 13(|Y|-1) = -3|Y| + 13 < 12$, a contradiction. Therefore, $|E^*(Y)| \le |Y| - 2$ in $H''$. We now consider the special $H$-set $X \cup Y$. Suppose that $e \notin E(Y)$. Then in $H$, $|E^*(X \cup Y)| \le |E^*(X)| + |E^*(Y)|  + |\{e'\}| \le (|X| + 1) + (|Y| - 2) + 1 = |X| + |Y|$, contradicting Claim~E. Therefore, $e \in E(Y)$. Let $e \in E(R)$, where $R \in Y$. We consider the special $H$-set $Q =  X \cup (Y \setminus \{R\})$. Then in $H$, $|E^*(Q)| \le |E^*(X)| + |E^*(Y)|  + |\{e'\}| \le (|X| + 1) + (|Y| - 2) + 1 = |X| + |Y| = |Q| + 1$. By Claim~E, $|E^*(Q)| \ge |Q| + 1$. Consequently, $|E^*(Q)| = |Q| + 1$. However since $|Y| \ge 2$, we note that $|Q| > |X|$, contradicting our choice of the special $H$-set $X$.~\smallqed

\ClaimX{H.10}
The case $|X| \ge 4$ and $|\partial(X)|  = 2$ cannot occur.

\ProofClaimX{H.10}
Note that $8|X_4| + 5|X_{14}| + 4|X_{11}| + |X_{21}| \le 8|X| \le 13|X|-6|X'| - 8$, as $|X| \ge 4$ and $|X'|=2$.
As $H'$ is linear by Claim~H.9, this contradicts Claim~H.4(b).~\smallqed

\ClaimX{H.11}
The case $|X| = 3$ and $|\partial(X)|  = 2$ cannot occur.

\ProofClaimX{H.11} Suppose, to the contrary, that $|X| = 3$ and $|\partial(X)|  = 2$, which implies that $|X'| = 2$. By Claim~H.9, $H' \in \cH_4$, and so, by Claim~H.4(a), we have $24 = 8|X| \ge 8|X_4| + 5|X_{14}| + 4|X_{11}| + |X_{21}|  > 13|X| - 6|X'| - \defic(H') = 39 - 12 - \defic(H')$, and so $\defic(H') > 3$. As $e$ contains two degree-$1$ vertices, we must have
$Y=H_4$ and $\defic(H') = \defic(Y)=8$, by Claim~H.3. Therefore, by Claim~H.4(a), $8|X_4| + 5|X_{14}| + 4|X_{11}| + |X_{21}|  > 39 - 12 - \defic(H') = 19$, implying that $|X_4| \ge 2$. By the connectivity of $H$, we note therefore that $V(H) = V(X) \cup \partial(X)$. Let $E^*(X) = \{e_1,e_2,e_3,e_4\}$.
Let $F_1$, $F_2$ and $F_3$ be the (vertex-disjoint) hypergraphs that belong to $X$. By Claim~H.2, $|X_{10}| = 0$.

\ClaimX{H.11.1}
$|X_4| = 2$.

\ProofClaimX{H.11.1}
As observed earlier, $|X_4| \ge 2$. Suppose, to the contrary, that $|X_4| \ge 3$, implying that $|X| = |X_4| = 3$. Thus, every hypergraph in $X$ is a copy of $H_4$. Let $\partial(X) = \{x,y\}$. Thus, $F_1$, $F_2$ and $F_3$ are all copies of $H_4$. Let $E(F_i) = \{f_i\}$ for $i \in [3]$. Further, $V(H) = V(X) \cup \{x,y\}$ and $E(H) = E(X) \cup E^*(X)$. In particular, $n(H) = 14$ and $m(H) = 6$. Renaming $x$ and $y$, if necessary, we may assume that $d_H(x) \ge d_H(y)$.

Suppose that every vertex in $V(X)$ has degree at most~$2$ in $H$. Then, by the linearity of $H$, there are three cases to consider. If $d_H(x) = 3$ and $d_H(y) = 2$, then $H = H_{14,1}$. If  $d_H(x) = 3$ and $d_H(y) = 1$, then $H = H_{14,2}$. If $d_H(x) = 2$, then $d_H(y) = 2$ and $H = H_{14,5}$. In all three cases, we contradict Claim~C.
Therefore, some vertex in $V(X)$ has degree~$3$. We may assume that $F_1$ contains a vertex, $z_1$ say, of degree~$3$ in $H$ and that $e_1$ and $e_4$ contain~$z$. By Claim~E, at least three edges in $E^*(X)$ intersect $F_2 \cup F_3$. Renaming $e_2$ and $e_3$, if necessary, we may assume that $e_2$ intersects $F_2$. If $e_3$ intersects $F_3$, then in this case letting $z_2 \in e_2 \cap f_2$ and $z_3 \in e_3 \cap f_3$, the set $\{z_1,z_2,z_3\}$ is an $X$-transversal that intersects every edge in $E^*(X)$, contradicting Claim~B. Hence, $e_3$ does not intersects $F_3$, implying by the linearity of $H$ that $e_3$ contains both $x$ and $y$, and intersects both $F_1$ and $F_2$. This in turn implies that each of $e_1$, $e_2$ and $e_3$ contain exactly one vertex in $\{x,y\}$ and therefore exactly one vertex of $F_i$ for each $i \in [3]$. In this case, letting $z_2 \in e_3 \cap f_2$ and $z_3 \in e_2 \cap f_3$, the set $\{z_1,z_2,z_3\}$ is an $X$-transversal that intersects every edge in $E^*(X)$, contradicting Claim~B. Therefore, $|X| = 2$.~\smallqed

\medskip By Claim~H.11.1, $|X_4| = 2$. We may assume that $F_1 \ne H_4$. Thus, both $F_2$ and $F_3$ are copies of $H_4$. Let $E(F_2) = \{f_2\}$ and $E(F_3) = \{f_3\}$.

\ClaimX{H.11.2}
All four edges in $E^*(X)$ intersect $F_1$.

\ProofClaimX{H.11.2} Suppose, to the contrary, that at most three edges in $E^*(X)$ intersect $F_1$. If two edges in $E^*(X)$ do not intersect $F_1$, then these two edges would both contain $x$ and $y$, contradicting the linearity of $H$. Hence, at least three edges in $E^*(X)$ intersect $F_1$. Consequently, exactly three edges in $E^*(X)$ intersect $F_1$.  Renaming edges in $E^*(X)$ if necessary, we may assume that $e_1$, $e_2$ and $e_3$ intersect $F_1$. By the linearity of $H$, the edge $e_4$ contains both $x$ and $y$, and intersects both $F_2$ and $F_3$.

Suppose that two of the edges $e_1$, $e_2$ and $e_3$, say $e_1$ and $e_2$, can be covered by a minimum transversal, $T_1$ say, in $F_1$. If $e_3$ intersects $F_2 \cup F_3$, say $e_3$ intersects $F_2$, then we can extend $T_1$ to an $X$-transversal that intersects every edge in $E^*(X)$ by adding to it the vertex in $e_3 \cap f_2$ and the vertex in $e_4 \cap f_3$, contradicting Claim~B. Hence, $e_3$ does not intersect $F_2 \cup F_3$, implying that $e_3$ contains at most one of $x$ and $y$ and at least three vertices of $F_1$. Further this implies that both $e_1$ and $e_2$ intersect $F_2 \cup F_3$ since by Claim~E, at least three edges in $E^*(X)$ intersect $F_2 \cup F_3$. By Observation~\ref{property:special}(l), there exists a $\tau(F_1)$-transversal that covers $e_3$ and one of $e_1$ and $e_2$, say $e_1$. We can then cover $e_2$ and $e_4$ from $F_2 \cup F_3$, once again contradicting Claim~B.

Therefore, we can only cover at most one of the three edges $e_1$, $e_2$ and $e_3$ by a minimum transversal in $F_1$. By Observation~\ref{property:special}(i), this implies that $F_1 = H_{11}$. Further, for $i \in [3]$ the edge $e_i$ intersects $F_1$ in exactly one vertex and this vertex has no neighbor of degree~$1$ in $F_1$. Let $e_i \cap V(F_1) = \{u_i\}$ for $i \in [3]$. Thus, $U = \{u_1,u_2,u_3\}$ is the set of three vertices in $F_1$ all of whose neighbors in $F_1$ have~$2$ in $F_1$. Let $H'$ be the hypergraph obtained from $H$ by deleting the eight vertices in $V(F_1) \setminus U$ (and their incident edges), adding a new vertex $u'$, and then adding the edge $f' = \{u',u_1,u_2,u_3\}$. Then, $H'$ has order~$n(H') = 14$ and size~$m(H') = 7$. We note that $E(H') = \{e_1,e_2,e_3,e_4,f',f_2,f_3\}$. Recall that the edge $e_4$ contains both $x$ and $y$, and intersects both $F_2$ and $F_3$. Therefore, the edge $e_i$ contains exactly one vertex from each of $F_2$ and $F_3$, and exactly one of $x$ and $y$, for each $i \in [3]$.
If every vertex in $F_2 \cup F_3$ has degree at most~$2$ in $H'$, then by the linearity of $H$, $H' = H_{14,1}$, implying that $H = H_{21,1}$, contradicting Claim~C. Therefore, some vertex in $F_2 \cup F_3$ has degree~$3$ in $H'$. We may assume that $z_2$ is such a vertex and that $z_2 \in V(F_2)$. Further, we may assume that $e_1$ and $e_2$ contain the vertex~$z_2$. Let $z_1$ and $z_3$ be the vertices of $e_3$ and $e_4$, respectively, that belongs to $F_1$ and $F_3$. Then, $\{z_1,z_2,z_3\}$ is a $\tau(H')$-transversal that can be extended to a $\tau(H)$-transversal by adding to it two additional vertices of $V(H) \setminus V(H')$. This produces an $X$-transversal that intersects every edge in $E^*(X)$, contradicting Claim~B. Therefore, all four edges in $E^*(X)$ intersect $F_1$.~\smallqed

\medskip
By Claim~H.11.2, all four edges in $E^*(X)$ intersect $F_1$.

\ClaimX{H.11.3}
All four edges in $E^*(X)$ intersect $F_2 \cup F_3$.

\ProofClaimX{H.11.2} Suppose, to the contrary, that at most three edges in $E^*(X)$ intersect $F_2 \cup F_3$. By Claim~E, at least three edges in $E^*(X)$ intersect $F_2 \cup F_3$. Therefore, exactly three edges in $E^*(X)$ intersect $F_2 \cup F_3$. We may assume that $e_1$ does not intersect $F_2 \cup F_3$. Thus, $e_2,e_3,e_4$ all intersect $F_2 \cup F_3$. We note that $|e_1 \cap V(F_1)| \ge 2$. Further since $\Delta(H) \le 3$, we note that $|(e_2 \cup e_3 \cup e_4) \cap V(F_1)| \ge 2$. Thus, by Observation~\ref{property:special}(k), there is a $\tau(F_1)$-transversal, $T_1$ say, that covers $e_1$ and covers one of the edges $e_2, e_3, e_4$, say $e_2$. Renaming $F_2$ and $F_3$, if necessary, we may assume that $e_3$ intersects $F_2$. If $e_4$ intersects $F_3$, then we can extend $T_1$ to an $X$-transversal that intersects every edge in $E^*(X)$ by adding to it the vertex in $e_3 \cap f_2$ and the vertex in $e_4 \cap f_3$, contradicting Claim~B. Therefore, $e_4$ does not intersect $F_3$, implying that $e_4$ intersects $F_2$. If $e_3$ intersects $F_3$, then analogously we can extend $T_1$ to an $X$-transversal that intersects every edge in $E^*(X)$, a contradiction. Therefore, neither $e_3$ nor $e_4$ intersects $F_3$, implying that $e_2$ is the only possible edge in $E^*(X)$ that intersect $F_3$, contradicting Claim~E.~\smallqed

\medskip
By Claim~H.11.3, all four edges in $E^*(X)$ intersect $F_2 \cup F_3$.

\ClaimX{H.11.4}
At least three edges in $E^*(X)$ intersect $F_3$.

\ProofClaimX{H.11.4} Suppose, to the contrary, that at most two edges in $E^*(X)$ intersect $F_3$.  Then, by Claim~E, exactly two edges in $E^*(X)$ intersect $F_3$. We may assume that $e_1$ and $e_2$ intersect $F_3$, and therefore $e_3$ and $e_4$ do not intersect $F_3$. By the linearity of $H$, this implies that at least one of $e_3$ and $e_4$ intersects $F_1$ in at least two vertices. Renaming $e_3$ and $e_4$ if necessary, we may assume that $|e_3 \cap V(F_1)| \ge 2$.

If $|(e_1 \cup e_2) \cap V(F_1)| \ge 2$, then by Observation~\ref{property:special}(k), there is a $\tau(F_1)$-transversal, $T_1$ say, that covers $e_3$ and covers one of the edges $e_1$ and $e_2$, say $e_1$. We can now extend $T_1$ to an $X$-transversal that intersects every edge in $E^*(X)$ by adding to it the vertex in $e_2 \cap f_3$ and the vertex in $e_4 \cap f_2$, contradicting Claim~B. Therefore, $|(e_1 \cup e_2) \cap V(F_1)| = 1$. By the linearity of $H$, this implies that both $e_1$ and $e_2$ intersect $F_2$.

Since $\Delta(H) \le 3$, we note that $|(e_1 \cup e_2 \cup e_4) \cap V(F_1)| \ge 2$. Thus, considering the sets $|e_3 \cap V(F_1)| \ge 2$ and the set $|(e_1 \cup e_2 \cup e_4) \cap V(F_1)| \ge 2$, by Observation~\ref{property:special}(k), there is a $\tau(F_1)$-transversal, $T_1$ say, that covers $e_3$ and covers one of the edges $e_1, e_2, e_4$. If $T_1$ covers $e_1$ (and therefore also $e_2$), then we can extend $T_1$ to an $X$-transversal that intersects every edge in $E^*(X)$ by adding to it the vertex in $e_4 \cap f_2$ and any vertex in $f_3$, contradicting Claim~B. Therefore, $T_1$ covers $e_4$. Since both $e_1$ and $e_2$ intersect both $F_2$ and $F_3$, we can once again extend $T_1$ to an $X$-transversal that intersects every edge in $E^*(X)$, a contradiction.~\smallqed

\medskip
By Claim~H.11.4, at least three edges in $E^*(X)$ intersect $F_3$. Analogously, at least three edges in $E^*(X)$ intersect $F_2$. Recall that by Claim~H.11.2, all four edges in $E^*(X)$ intersect $F_1$. On the one hand, if two edges of $E^*(X)$ intersect $F_1$ in a common vertex, then there is a $\tau(F_1)$-transversal that covers two edges of $E^*(X)$. On the other hand, if no two edges of $E^*(X)$ intersect $F_1$ in a common vertex, then by Observation~\ref{property:special}(k), we can once again find a $\tau(F_1)$-transversal that covers two edges of $E^*(X)$ (by considering, for example, the set $|(e_1 \cup e_2) \cap V(F_1)| \ge 2$ and the set $|(e_1 \cup e_2) \cap V(F_1)| \ge 2$). In both cases, there is a $\tau(F_1)$-transversal, $T_1$,  that covers two edges of $E^*(X)$. Renaming edges in $E^*(X)$ if necessary, we may assume that $T_1$ covers $e_1$ and $e_2$. Recall that by Claim~H.11.3, all four edges in $E^*(X)$ intersect $F_2 \cup F_3$. In particular, both $e_3$ and $e_4$ intersect $F_2 \cup F_3$. Renaming $F_2$ and $F_3$ if necessary, we may assume that $e_3$ intersects $F_2$. If $e_4$ intersects $F_3$, then we can extend $T_1$ to an $X$-transversal that intersects every edge in $E^*(X)$ by adding to it the vertex in $e_3 \cap f_2$ and the vertex in $e_4 \cap f_3$, contradicting Claim~B. Therefore, $e_4$ does not intersect $F_3$, implying that $e_4$ intersects $F_2$. If $e_3$ intersects $F_3$, then analogously we can extend $T_1$ to an $X$-transversal that intersects every edge in $E^*(X)$, a contradiction. Therefore, neither $e_3$ nor $e_4$ intersects $F_3$, implying that $e_1$ and $e_2$ are the only possible edge in $E^*(X)$ that intersect $F_3$, contradicting Claim~H.11.4. This completes the proof of Claim~H.11.~\smallqed

\ClaimX{H.12}
The case $|X| = 2$ and $|\partial(X)|  = 2$ cannot occur.

\ProofClaimX{H.12} Suppose, to the contrary, that $|X| = 2$ and $|\partial(X)|  = 2$. We note that $|E^*(X)| = 3$ and $|X'| = 2$. Let $E^*(X) = \{e_1,e_2,e_3\}$. By Claim~H.9, $H' \in \cH_4$, and so by Claim~G, $\Phi(H') < 0$. Since at most one of the edges in $E^*(X)$ can contain both vertices in $\partial(X)$, we note that
\begin{equation}
\sum_{i=1}^3 |e_i \cap V(F_1)| + \sum_{i=1}^3 |e_i \cap V(F_2)|\ge 8.
\label{Eq:sum}
\end{equation}

Suppose that at most two edges in $E^*(X)$ intersect $F_1$. Then, by Claim~E, exactly two edges in $E^*(X)$ intersect $F_1$. We may assume that both $e_1$ and $e_2$ intersect $F_1$, and therefore that $e_3$ does not intersect $F_1$. Thus, $|e_3 \cap V(F_1)| = 0$ and $|e_3 \cap V(F_2)| \ge 2$
If $|e_1 \cap V(F_1)| + |e_2 \cap V(F_1)| \ge 4$, then by Observation~\ref{property:special}, we can find a $\tau(F_1)$-transversal, $T_1$, that covers both $e_1$ and $e_2$. In this case, we can extend $T_1$ to an $X$-transversal that intersects every edge in $E^*(X)$ by adding to it a $\tau(F_2)$-transversal that covers $e_3$, a contradiction. Therefore, $|e_1 \cap V(F_1)| + |e_2 \cap V(F_1)| \le 3$, implying by Inequality~(\ref{Eq:sum}) that
\begin{equation}
\sum_{i=1}^3 |e_i \cap V(F_2)|\ge 5.
\label{Eq:sum2}
\end{equation}

Therefore, by Inequality~(\ref{Eq:sum2}), we note that
$(|e_1 \cap V(F_2)| + |e_3 \cap V(F_2)|) + (|e_2 \cap V(F_2)| + |e_3 \cap V(F_2)|) \ge 5 + |e_3 \cap V(F_2)| \ge 7$. Renaming $e_1$ and $e_2$ if necessary, we may assume that $|e_1 \cap V(F_2)| + |e_3 \cap V(F_2)| \ge 4$. By Observation~\ref{property:special}, we can find a $\tau(F_2)$-transversal, $T_2$, that covers both $e_1$ and $e_3$. In this case, we can extend $T_3$ to an $X$-transversal that intersects every edge in $E^*(X)$ by adding to it a $\tau(F_1)$-transversal that covers $e_2$, a contradiction. Hence, all three edges in $E^*(X)$ intersect $F_1$. Analogously, all three edges in $E^*(X)$ intersect $F_2$. Renaming $F_1$ and $F_2$ if necessary, we may assume by Inequality~(\ref{Eq:sum}) that
$\sum_{i=1}^3 |e_i \cap V(F_1)| \ge 4$.
%
By Observation~\ref{property:special}, we can find a $\tau(F_1)$-transversal, $T_1$, that covers two of the edges in $E^*(X)$. The third edge in $E^*(X)$ can be covered by a $\tau(F_2)$-transversal noting that all three edges in $E^*(X)$ intersect $F_2$. This produces an $X$-transversal that intersects every edge in $E^*(X)$, a contradiction.~\smallqed

\medskip
We proceed further with some additional notation. We associate with the set $X$ a bipartite multigraph, which we denote by $M_X$, with partite sets $X$ and $E^*(X)$ as follows. If an edge $e \in E^*(X)$ intersects a subhypergraph $H' \in X$ in $H$ in $k$ vertices, then we add $k$ multiple edges joining $e \in E^*(X)$ and $H' \in X$ in $M_X$. Two multiple edges (also called parallel edges in the literature) joining two vertices in $M_X$ we call \emph{double edges}, while three multiple edges joining two vertices in $M_X$ we call \emph{triple edges}. An edge that is not a multiple edge we call a \emph{single edge}. By supposition in our proof of Claim~H, $|E^*(X)| = |X| + 1$.
We say that a pair $(e_1,e_2)$ of edges in $E^*(X)$ form an $E^*(X)$-\emph{pair} if they intersect a common subhypergraph, $F$, of $X$ in $H$, and, further,
there exists a $\tau(F)$-transversal that covers both $e_1$ and $e_2$ in $H$. We call $F$ a subhypergraph of $X$ associated with the $E^*(X)$-pair, $(e_1,e_2)$.

\ClaimX{H.13}
$|E(M_X)| \ge 4|E^*(X)| - 3|\partial(X)| = 4(|X|+1) - 3|\partial(X)|$.
Furthermore, the following holds.  \\
\indent (a) $M_X$ does not contain triple edges. \\
\indent (b) No $F \in X$ has double edges to three distinct vertices of $E^*(X)$ in $M_X$. \\
\indent (c) Every $e \in E^*(X)$ has degree at least $4-|\partial(X)|$ in $M_X$. \\
\indent (d) Every $F \in X$ has degree at least $2$ in $M_X$.

\ProofClaimX{H.13}
Consider the bipartite multigraph, $M_X'$, which is identical to $M_X$, except we add a new vertex $b$ to $M_X$ and for each vertex $e \in E^*(X)$ we add $r$ multiple edges joining $e$ and $b$ in $M_X'$ if $e$ intersects $\partial(X)$ in $r$ vertices in $H$. In $M_X'$, every vertex in $E^*(X)$ has degree~$4$. Further, since $\Delta(H) \le 3$, the vertex $b$ has degree at most $3|\partial(X)|$. Therefore, $M_X = M_X' - b$ contains at least $4|E^*(X)| - 3|\partial(X)|$ edges. This proves the first part of Claim~H.13.

To prove part~(a), for the sake of contradiction, suppose that $M_X$ does contain triple edges that join vertices $e \in E^*(X)$ and $F \in X$. By Observation~\ref{property:special}($\ell$), we note that $e$ is an $X$-universal edge, a contradiction to Claim~H.1. This proves part~(a).

To prove part~(b), for the sake of contradiction, suppose that some $F \in X$ has double edges to three distinct vertices, say $e_1,e_2,e_3$, in $E^*(X)$ in $M_X$. By Observation~\ref{property:special}(o), we note that at least one of $e_1$, $e_2$ or $e_3$ is an $X$-universal edge, a contradiction to Claim~H.1. This proves part~(b).

By the construction of $M_X$, every vertex $e \in E^*(X)$ has degree $4 - |V(e) \cap \partial(X)| \ge 4-|\partial(X)|$ in $M_X$.
By Claim~E (used on $F$), we note that at least two edges in $E^*(X)$ intersect $F$ in $H$, and therefore the degree of $F$ in $M_X$ is at least~$2$. This proves part~(c) and part~(d). ~\smallqed

\ClaimX{H.14}
The case $|X| \ge 5$ and $|\partial(X)|  = 1$ cannot occur.

\ProofClaimX{H.14} Suppose, to the contrary, that $|X| \ge 5$ and $|\partial(X)|  = 1$. We note that in this case $|X'| = 3$ and $H'$ is linear. Let $\partial(X) = \{y\}$.

\ClaimX{H.14.1} $|X| = 5$ and $X = X_4$.

\ProofClaimX{H.14.1}
By Claim~H.4(b), $8|X| \ge 8|X_4| + 5|X_{14}| + 4|X_{11}| + |X_{21}| > 13|X|-6|X'| - 8 = 8|X| + 5|X| - 26$, and so $5|X| < 26$, implying that $|X| \le 5$.
This clearly implies that $|X|=5$. By supposition, $|X| \ge 5$. Consequently, $|X|=5$, implying that $8|X_4| + 5|X_{14}| + 4|X_{11}| + |X_{21}| >  8|X| - 1$, which in turn implies that $X=X_4$.~\smallqed

\medskip
By Claim~H.14.1, $X = X_4$ and $|X_4| = 5$.
Let $X = \{x_1,x_2,\ldots,x_5\}$ and let $E^*(X) = \{e_1,e_2,\ldots,e_6\}$.
We now consider the bipartite multigraph $M_X$ defined earlier. Since $X=X_4$ and $H$ is linear, the multigraph $M_X$ is in this case a graph. By Claim~H.13(c), since here $|\partial(X)|  = 1$, $d_{M_X}(e) \ge 3$ for each $e \in E^*(X)$. By Claim~H.13(d), $d_{M_X}(F) \ge 2$ for each $F \in X$.
Also by Claim~H.13, we note that $\sum_{F \in X} d_{M_X}(F) = |E(M_X)| \ge 4(|X|+1) - 3|\partial(X)| = 24 - 3 = 21$.

\ClaimX{H.14.2} Suppose the edges $e_i$ and $e_j$ form an $E^*(X)$-pair and that $x_{\ell}$ is the vertex in $X$ corresponding to a copy of $H_4$ of $X$ that contains a vertex covering both $e_1$ and~$e_2$. Then there exists a vertex $x_k \in X \setminus \{x_{\ell}\}$ in $M_X$ such that $e_i$ and $e_j$ are both neighbors of $x_k$ in $M_X$ and one of the following holds. \\
\indent {\rm (a)} $d_{M_X}(x_k) = 2$. \\
\indent {\rm (b)} $d_{M_X}(x_k) = 3$ and $N_{M_X}(x_{k'}) \subseteq N_{M_X}(x_k)$ for some vertex $x_{k'} \in X \setminus \{x_k,x_{\ell}\}$.

\ProofClaimX{H.14.2} For notational convenience, we may assume that $i = 1$ and $j = 2$. Thus, $(e_1,e_2)$ is an $E^*(X)$-pair. Further, we may assume that $x_1$ is the vertex in $X$ corresponding to a copy of $H_4$ of $X$ that contains a vertex covering both $e_1$ and $e_2$. We now consider the bipartite graph $M_X' = M_X - \{e_1,e_2,x_1\}$ with partite sets $X' = X \setminus \{x_1\}$ and $E'_X = E^*(X) \setminus \{e_1,e_2\}$. If there exists a matching in $M_X'$ that matches $E'_X$ to $X'$, then there exists a minimum $X$-transversal that intersects every edge in $E^*(X)$, contradicting Claim~B. Therefore, no matching in $M_X'$ matches $E'_X$ to $X'$. By Hall's Theorem, there is a nonempty subset $S \subseteq E'_X$ such that in $M_X'$, $|N(S)|<|S|$. If $S = E'_X$, then let $x_k \in X' \setminus N_{M_X}(S)$. In this case, the vertices $e_1$ and $e_2$ are the only possible neighbors of $x_k$ in $M_X$, implying that $N_{M_X}(x_k) = \{e_1,e_2\}$ and $d_{M_X}(x_k) = 2$. Thus, Part~(a) holds, as desired. Hence, we may assume that $|S| \le 3$ and that no vertex in $X' \setminus \{x_1\}$ has degree~$2$ with $e_1$ and $e_2$ as its neighbors.

Since $|S| \le 3$, we note that $|N_{M_X'}(S)| \le 2$. Since every vertex in $E^*(X)$ has degree at least~$3$ in $M_X$, we note that $|N_{M_X'}(S)| \ge 2$. Consequently, $|N_{M_X'}(S)| = 2$, implying that $|S| = 3$.
Renaming vertices if necessary, we may assume that $S = \{e_3,e_4,e_5\}$ and that $N_{M_X'}(S) = \{x_2,x_3\}$. Since $d_{M_X}(e) \ge 3$ for each $e \in E^*(X)$, this implies that for each vertex $e \in S$, we have $N_{M_X}(e) = \{x_1,x_2,x_3\}$. We note that the only possible neighbors of $x_4$ and $x_5$ in $M_X$ are $e_1$, $e_2$ and $e_6$, and so $d_{M_X}(x_4) \le 3$ and $d_{M_X}(x_5) \le 3$.

We show next that $x_4$ or $x_5$ dominate both $e_1$ and $e_2$ in $M_X$. Suppose, to the contrary, that neither $x_4$ nor $x_5$ dominate both $e_1$ and $e_2$ in $M_X$. Since both $x_4$ and $x_5$ have degree at least~$2$ in $M_X$, this implies that $d_{M_X}(x_4) = d_{M_X}(x_5) = 2$ and that $e_6$ is adjacent to both $x_4$ and $x_5$. Renaming $e_1$ and $e_2$, if necessary, we may assume that $e_2$ is adjacent to $x_5$, and so
$N_{M_X}(x_5) = \{e_2,e_6\}$.
%
Since $\sum_{x \in X} d_{M_X}(x) \ge 21$, the degree sequence of vertices of $X$ in $M_X$ is therefore either $2,2,5,6,6$ or $2,2,6,6,6$. If the degree sequence is given by $2,2,6,6,6$, then $e_6$ would be adjacent to every vertex of $X$ in $M_X$, which is not possible since $d_{M_X}(e) \le 4$ for each $e \in E^*(X)$. Thus, the degree sequence of vertices of $X$ in $M_X$ is $2,2,5,6,6$. In particular, this implies that $d_{M_X}(e_6) = 4$, which in turn implies that $x_2$ and $x_3$ are both adjacent to $e_1$ and $e_2$ in $M_X$. Thus, $M_X - \{e_6,x_4,x_5\} = K_{3,5}$. By assumption, $e_2$ is adjacent to $x_5$. Thus, since $d_{M_X}(e_2) = 4$, we note that $N_{M_X}(e_2) = \{x_1,x_2,x_3,x_5\}$ and therefore $N_{M_X}(x_4) = \{e_1,e_6\}$. The graph $M_X$ is therefore determined.

We note that the two vertices of $X$ of degree~$6$ in $M_X$ both give rise to at least two $E^*(X)$-pairs, while the vertex of $X$ of degree~$5$ in $M_X$ gives rise to at least one $E^*(X)$-pairs. Further since there are no overlapping edges, these five $E^*(X)$-pairs are distinct. Therefore, there are at least five distinct $E^*(X)$-pairs. Recall that $N_{M_X}(x_4) = \{e_1,e_6\}$ and $N_{M_X}(x_5) = \{e_2,e_6\}$. Suppose that there is a vertex in the copy of $H_4$ in $H$ corresponding to the vertex $x_4$ that contains a vertex covering both $e_1$ and $e_6$. Then the structure of $M_X$ implies that there exists a minimum $X$-transversal that intersects every edge in $E^*(X)$, contradicting Claim~B. Analogously, if there is a vertex in the copy of $H_4$ in $H$ corresponding to the vertex $x_5$ that contains a vertex covering both $e_2$ and $e_6$, we produce a contradiction. Therefore, since $e_6$ has only two neighbors in $M_X$ different from $x_4$ and $x_5$, at most two $E^*(X)$-pairs contain $e_6$. Since there are at least five $E^*(X)$-pairs, there exists an $E^*(X)$-pair,
$(e_{i_1},e_{i_2})$ say, that is not the pair $(e_1,e_2)$ and such that $e_6 \notin \{e_{i_1},e_{i_2}\}$.

Renaming $i_1$ and $i_2$, if necessary, we may assume that $e_{i_1} \notin \{e_1,e_2\}$. Thus, $e_{i_1} \in \{e_3,e_4,e_5\}$. Let $x_j$ be the vertex in $X$ corresponding to a copy of $H_4$ of $X$ that contains a vertex covering both $e_{i_1}$ and $e_{i_2}$. Since $e_{i_1}$ is adjacent to neither $x_4$ nor $x_5$, we note that $j \in \{1,2,3\}$. We now proceed analogously as we did with the $E^*(X)$-pair  $(e_1,e_2)$. We consider the bipartite graph $M_X'' = M_X - \{e_{i_1},e_{i_2},x_j\}$  with partite sets $X'' = X \setminus \{x_j\}$ and $E''_X = E^*(X) \setminus \{e_{i_1},e_{i_2}\}$. No matching in $M_X''$ matches $E''_X$ to $X''$. By Hall's Theorem, there is a nonempty subset $R \subseteq E''_X$ such that in $M_X''$, $|N(R)|<|R|$. If $R = E'_X$, then let $x' \in X \setminus N_{M_X}(R)$. The only possible neighbors of $x'$ are the vertices $e_{i_1}$ and $e_{i_2}$. Since $d_{M_X}(x') = 2$, we note that $x' \in \{x_4,x_5\}$. However, $x'$ is not adjacent to $e_{i_1}$, and so $d_{M_X}(x') \le 1$, a contradiction. Therefore, $|R| = 3$ and in $M_X''$, $|N(R)| = 2$. We note that both vertices in $N(R)$ are adjacent to all three vertices in $R$, implying that the two vertices of $X$ not in $R \cup \{x_j\}$ are the two vertices of $X$ of degree~$2$ in $M_X$, namely $x_4$ and $x_5$. Further, $e_6 \notin R$ and $e_6$ is the only common neighbor of $x_4$ and $x_5$. If $i_2 = 1$, then $x_5$ is adjacent to neither $e_{i_1}$ nor $e_{i_2}$, implying that $d_{M_X}(x_5) = 1$. If $i_2 \ne 1$, then $x_4$ is adjacent to neither $e_{i_1}$ nor $e_{i_2}$, implying that $d_{M_X}(x_4) = 1$. Both cases produce a contradiction. Therefore, at least one of $x_4$ and $x_5$ dominate both $e_1$ and $e_2$ in $M_X$. Renaming $x_4$ and $x_5$, if necessary, we may assume that $x_4$ dominate both $e_1$ and $e_2$ in $M_X$.

By our earlier assumption, no vertex in $X' \setminus \{x_1\}$ has degree~$2$ with $e_1$ and $e_2$ as its neighbors. Hence, $N_{M_X}(x_4) = \{e_1,e_2,e_6\}$. As observed earlier, the only possible neighbors of $x_5$ in $M_X$ are $e_1$, $e_2$ and $e_6$, and so $N_{M_X}(x_5) \subseteq  N_{M_X}(x_4)$. Taking $x_k = x_4$ and $x_{k'} = x_5$, Part~(b) holds. This completes the proof of Claim~H.14.2.~\smallqed

\medskip
We now consider the degree sequence of vertices of $X$ in $M_X$. Let this degree sequence, in nondecreasing order, be given by $s \colon d_1,d_2,d_3,d_4,d_5$, and so $2 \le d_1 \le d_2 \le \cdots \le d_5 \le 6$. As observed earlier, $\sum_{i=1}^5 d_i \ge 21$. By Claim~H.14.2, $d_1 = 2$ or $d_1 = d_2 = 3$.

\ClaimX{H.14.3}
There are at most three $E^*(X)$-pairs.

\ProofClaimX{H.14.3}
Suppose, to the contrary, that there are at least four $E^*(X)$-pairs. Renaming vertices of $E^*(X)$, if necessary, we may assume that $(e_1,e_2)$ and $(e_3,e_4)$ are $E^*(X)$-pairs. Let $x_{12}$ (respectively, $x_{34}$) be the vertex in $X$ corresponding to a copy of $H_4$ of $X$ that contains a vertex covering both $e_1$ and~$e_2$ (respectively, $e_3$ and $e_4$). Possibly, $x_{12} = x_{34}$. Let $x_1 \in X \setminus \{x_{12}$ be adjacent to both $e_1$ and $e_2$ in $M_X$ and such that either $d_{M_X}(x_1) = 2$ or $d_{M_X}(x_1) = 3$ and there exists a vertex $x_{1'} \in X \setminus \{x_1,x_{12}\}$ with $N_{M_X}(x_{1'}) \subseteq N_{M_X}(x_1)$. Further, let $x_2 \in X \setminus \{x_{34}$ be adjacent to both $e_3$ and $e_4$ in $M_X$ and such that either $d_{M_X}(x_2) = 2$ or $d_{M_X}(x_2) = 3$ and there exists a vertex $x_{2'} \in X \setminus \{x_2,x_{34}\}$ with $N_{M_X}(x_{2'}) \subseteq N_{M_X}(x_2)$. We note that $x_1$ and $x_2$ exists by Claim~H.14.2. Further, $x_{12} \notin \{x_1,x_2\}$, $x_{34} \notin \{x_1,x_2\}$ and $x_1 \ne x_2$. We may assume that $d_{M_X}(x_1) \le d_{M_X}(x_2)$.

Suppose that $d_{M_X}(x_1) = 2$, and so $N_{M_X}(x_1) = \{e_1,e_2\}$. Then, $d_1 = 2$. Since $d_2 \le 3$ and $\sum_{i=1}^5 d_i \ge 21$, this implies that $d_3 \ge 4$, and so $x_1$ and $x_2$ are the only vertices in $X$ of degree at most~$3$ in $M_X$. If $d_{M_X}(x_2) = 3$, then by Claim~H.14.2, $N_{M_X}(x_1) \subseteq N_{M_X}(x_2)$, a contradiction since at least one of $e_1$ and $e_2$ is not adjacent to $x_2$. Thus, $d_{M_X}(x_2) = 2$, and so $d_2  = 2$ and $N_{M_X}(x_2) = \{e_3,e_4\}$. Let $(e_i,e_j)$ be an $E^*(X)$-pair different from $(e_1,e_2)$ and $(e_3,e_4)$. By Claim~H.14.2, there exists a vertex in $X$ of degree at most~$3$ adjacent to both $e_i$ and $e_j$. The vertices $x_1$ and $x_2$ are the only two vertices in $X$ of degree at most~$3$ in $M_X$. However, neither $x_1$ nor $x_2$ is adjacent to both $e_i$ and $e_j$, a contradiction. Therefore, $d_{M_X}(x_1) = d_{M_X}(x_2) = 3$. Since $\sum_{i=1}^5 d_i \ge 21$, we note that $d_3 \ge 3$. This implies that $N_{M_X}(x_{1'}) = N_{M_X}(x_1)$ and $N_{M_X}(x_{2'}) = N_{M_X}(x_2)$.

Suppose $d_3 = 3$. Then, $s$ is given by $3,3,3,6,6$. Let $x_3$ be the vertex of $X \setminus \{x_1,x_2\}$ of degree~$3$ in $M_X$. Since $N_{M_X}(x_{1}) \ne N_{M_X}(x_2)$, we note that $x_{1'} = x_3$ and $x_{2'} = x_3$. But then $\{e_1,e_2,e_3,e_4\} \subseteq N_{M_X}(x_{3})$, and so $d_{M_X}(x_3) \ge 4$, a contradiction. Therefore, $d_3 \ge 4$. Thus, $x_1$ and $x_2$ are the only vertices in $X$ of degree at most~$3$ in $M_X$. This implies that $x_{1'} = x_2$ (and $x_{2'} = x_1$). But then $\{e_1,e_2,e_3,e_4\} \subseteq N_{M_X}(x_2)$, and so $d_{M_X}(x_2) \ge 4$, a contradiction.~\smallqed

\ClaimX{H.14.4}
The degree sequence $s$ is given by $3,3,5,5,5$. 

\ProofClaimX{H.14.4}
We show firstly that $d_1 = 3$. Suppose, to the contrary, that $d_1 = 2$.
If $d_2 = 2$, then $s$ is given by $2,2,5,6,6$ or $2,2,6,6,6$, implying that there are at least five $E^*(X)$-pairs.
If $d_2 = 3$, then $s$ is given by $2,3,4,6,6$ or $2,3,5,5,6$ or $2,3,5,6,6$ or $2,3,6,6,6$, implying that there are at least four $E^*(X)$-pairs.
In both cases, we contradict Claim~H.14.3. Therefore, $d_1 = d_2 = 3$.

We show next that $d_6 \le 5$. Suppose, to the contrary, that $d_6 = 6$. If $d_3 = 3$, then $s$ is given by $3,3,3,6,6$, implying that there are at least four $E^*(X)$-pairs, a contradiction. Hence, $d_3 \ge 4$. Let $x_1$ and $x_2$ be the two vertices of degree~$3$ in $M_X$. By Claim~H.14.2, $N_{M_X}(x_1) = N_{M_X}(x_2)$. We may assume that $e_1$, $e_2$ and $e_3$ are the three neighbors of $x_1$ and $x_2$. By Claim~H.14.2, if $(e_i,e_j)$ is an $E^*(X)$-pair, then $1 \le i,j \le 3$. We may assume that $x_4$ is a vertex of degree~$6$ in $M_X$. Then, $x_4$ gives rise to two $E^*(X)$-pairs that comprise of four distinct vertices. At least one such pair is distinct from $(e_1,e_2)$, $(e_1,e_3)$ and $(e_2,e_3)$, a contradiction. Therefore, $d_3 \le 5$. Since $\sum_{i=1}^5 d_i \ge 21$, the degree sequence $s$ is given by $3,3,5,5,5$.~\smallqed

\medskip
Renaming vertices in $X$ if necessary, we may assume that $d_{M_X}(x_i) = d_i$. By Claim~H.14.4, $x_1$ and $x_2$ have degree~$3$ in $M_X$, while $x_3$, $x_4$ and $x_5$ have degree ~$3$ in $M_X$. By Claim~H.14.2, $N_{M_X}(x_1) = N_{M_X}(x_2)$. We may assume that $e_1$, $e_2$ and $e_3$ are the three neighbors of $x_1$ and $x_2$. Further, we may assume that $e_1$ is not adjacent to $x_3$, $e_2$ is not adjacent to $x_4$, and $e_3$ is not adjacent to $x_5$. By Claim~H.14.2, the three $E^*(X)$-pairs are $(e_1,e_2)$, $(e_1,e_3)$, and $(e_2,e_3)$. The graph $M_X$ is completely determined.

Let $F_i$ be the copy of $H_4 \in X_4$ associated with the vertex $x_i$ in $M_X$ for $i \in [5]$. We note that $e_4$, $e_5$ and $e_6$ all have degree~$3$ in $M_X$ and are all adjacent in $M_X$ to $x_3$, $x_4$ and $x_5$. Thus, in $H$, each of the edges $e_4$, $e_5$ and $e_6$ contain the vertex~$y$ and one vertex from each copy of $H_4$ associated with $x_3$, $x_4$ and $x_5$. In $H$, each of the edges $e_1$, $e_2$ and $e_3$ contain one vertex from each copy of $H_4$ associated with $x_3$, $x_4$ and $x_5$. The copy of $H_4$ associated with $x_3$ contains a vertex covering both $e_2$ and $e_3$. The copy of $H_4$ associated with $x_4$ contains a vertex covering both $e_1$ and $e_3$. The copy of $H_4$ associated with $x_5$ contains a vertex covering both $e_1$ and $e_2$. By the linearity of $H$, and since $(e_1,e_2)$, $(e_1,e_3)$, and $(e_2,e_3)$ are the only three $E^*(X)$-pairs, the graph $H$ is now completely determined, and $H = H_{21,2}$. By Observation~\ref{property:special}, $\xi(H) = 0$, contradicting the fact that $H$ is a counterexample to the theorem. This completes the proof of Claim~H.14.~\smallqed

\ClaimX{H.15}
The case $|X| = 4$, and $|\partial(X)| = 1$ cannot occur.

\ProofClaimX{H.15}
By Claim~H.13, there are no triple edges in $M_X$. By Claim~H.4(b) we note that $8|X_4| + 5|X_{14}| + 4|X_{11}| + |X_{21}| > 13|X|-6|X'| - 8 = 26$, as $|X| = 4$ and $|X'|=3$. This implies that $|X_4| \ge 3$. Since $|\partial(X)| = 1$, by Claim~H.13(c), each vertex of $E^*(X)$ has degree at least~$3$ in $M_X$. By Claim~H.13(d), each vertex of $X$ has degree at least~$2$ in $M_X$. We now prove the following subclaims.

\ClaimX{H.15.1}
There does not exist a double edge between $e \in E^*(X)$ and $F \in X$ in $M_X$, such that $e$ belongs to a $E^*(X)$-pair associated with $F$.

\ProofClaimX{H.15.1}
For the sake of contradiction, suppose that there is a double edge between $e \in E^*(X)$ and $F \in X$ in $M_X$, such that $e$ belongs to a $E^*(X)$-pair, say $(e,e')$, associated with $F$. As there is a double edge incident with $F$ in $M_X$ we note that $F \notin X_4$, implying that $|X_4| = 3$. We now consider the bipartite multigraph $M_X^*$, where $M_X^*$ is defined as follows. Let $M_X^*$ be obtained from $M_X$ by removing the vertex $e \in E^*(X)$ and removing all edges $(e'',F)$, where $(e,e'')$ is not a $E^*(X)$-pair associated with $F$. We note that $(e',F)$ is an edge in $M_X^*$ and, therefore, $d_{M_X^*}(F) \ge 1$.

First suppose that there is a perfect matching, $M$, in $M_X^*$. Let $(e^*,F)$ be an edge of the matching incident with $F$. By the definition of $M_X^*$, we note that $(e^*,e)$ is a $E^*(X)$-pair associated with $F$. By Observation~\ref{property:special}(g), we can find a $\tau(X)$-transversal covering $E^*(X)$, contradicting Claim~B.  Therefore, there is no perfect matching in $M_X^*$. By Hall's Theorem, there is a nonempty subset $S \subseteq E^*(X)\setminus\{e\}$ such that $|N_{M_X^*}(S)| < |S|$. Let $X^* = X \setminus  N_{M_X^*}(S)$.

We will show that $F \in X^*$. For the sake of contradiction, suppose that $F \notin X^*$, and so $F \in N_{M_X^*}(S)$. Since $F \notin X^*$, the neighbors of a vertex of $X^*$ in $M_X$ remain unchanged in $M_X^*$. We note that no vertex in $X^*$ has a neighbor that belongs to $S$ in $M_X^*$, and therefore also no neighbor that belongs to $S$ in $M_X$. This implies that in $M_X$ the vertices in $X^*$ are adjacent to at most $|E^*(X) \setminus S| = |E^*(X)| - |S| = |X|+1-|S| \le |X^*|$ vertices in $E^*(X)$. Thus, $|E^*(X^*)| \le |X^*|$ in $H$, contradicting Claim~D. Therefore, $F \in X^*$. Since all three subhypergraphs in $X \setminus \{F\}$ belong to $X_4$, all edges in $M_X$ from $S$ to vertices in $N_{M_X^*}(S)$ are therefore single edges.

If $|S|=4$, then the vertices in $X^*$ have no edge to $S = E^*(X) \setminus \{e\}$ in $M_X^*$. In particular, the vertex $F \in X^*$ is isolated in $M_X^*$, contradicting our earlier observation that $d_{M_X^*}(F) \ge 1$. Therefore, $|S| \le 3$.

Suppose that $|S|=3$, and let $S=\{e_1,e_2,e_3\}$. As observed earlier, all edges in $M_X$ from $S$ to vertices in $N_{M_X^*}(S)$ are single edges. Since each vertex of $S$ has degree at least~$3$ in $M_X$ and $|N_{M_X^*}(S)| \le 2$, each vertex of $S$ therefore has $F$ as a neighbor in $M_X$, but not in $M^*_X$ since $F \in X^*$. Since $\Delta(H) \le 3$, this implies that  $|(V(e_1) \cup V(e_2) \cup V(e_3)) \cap V(F)| \ge 2$. By Observation~\ref{property:special}(k),
either $(e,e_1)$ or $(e,e_2)$ or $(e,e_3)$ is a $E^*(X)$-pair associated with $F$. Renaming edges in $E^*(X)$ if necessary,  we may assume that $(e,e_1)$ is a $E^*(X)$-pair associated with $F$. However, since the edge $(e_1,F)$  was removed from $M_X$ when constructing $M_X^*$, this implies that $(e,e_1)$ is not a $E^*(X)$-pair associated with $F$, a contradiction. Therefore, $|S| \le 2$.

Suppose that $|S|=2$, and let $S=\{e_1,e_2\}$. We note that there is no edge in $M_X^*$ joining $e_1$ and $F$. If $e_1$ has a double edge to $F$ in $M_X$, then $(e_1,e)$ would be a $E^*(X)$-pair associated with $F$, by Observation~\ref{property:special}(k), and therefore the edge $(e_1,F)$ would still exist in $M_X^*$, a contradiction. Hence, either $e_1$ has a single edge to $F$ in $M_X$ or is not adjacent to $F$ in $M_X$. This implies that the degree of $e_1$ is at most~$2$ in $M_X$ since it can only be adjacent to the vertex in $N_{M_X^*}(S)$ and to $F$.
This contradicts Claim~H.13. Therefore, $|S| = 1$. However if $|S|=1$, we analogously get a contradiction as the vertex in $S$ has degree at most~$1$ in $M_X$. This completes the proof of Claim~H.15.1.~\smallqed

\ClaimX{H.15.2}
There does not exist a double edge in $M_X$.

\ProofClaimX{H.15.2}
Suppose, to the contrary, that there is a double edge between $e \in E^*(X)$ and $F \in X$ in $M_X$. By Claim~H.15.1, for every edge $(e',F)$ in $M_X$, the pair $(e,e')$ is not a $E^*(X)$-pair associated with $F$. If there is a double edge between some vertex $e' \in E^*(X) \setminus \{e\}$ and $F$ in $M_X$, then $(e,e')$ would be a $E^*(X)$-pair associated with $F$, a contradiction. Therefore, there is no other double edge incident to $F$, except for the double edge that joins it to $e \in E^*(X)$. By Claim~D, the set $N_{M_X}(F)$ contains at least two vertices. If the set $N_{M_X}(F)$ contains four or more vertices, then by  by Observation~\ref{property:special}(k), we would get a $E^*(X)$-pair associated with $F$ containing $e$, a contradiction. Therefore, the set $N_{M_X}(F)$
contains at most three vertices. That is, $|N_{M_X}(F)| \in \{2,3\}$.

As there is a double edge incident with $F$ in $M_X$, we note that $F \notin X_4$, implying that $|X_4| = 3$ and $X \setminus \{F\} = X_4$. By Claim~H.13, $|E(M_X)| \ge 4(|X|+1) - 3|\partial(X)| =  4 \times 5 - 3 = 17$.

We first consider the case when $|N_{M_X}(F)|=2$, and let $N_{M_X}(F) = \{e,e'\}$. That is $d_{M_X}(F) = 3$, since the double edge between $e$ and $F$ counts~$2$ to the degree of $F$ in $M_X$. Since $X \setminus \{F\} = X_4$, we have no double edges incident with a vertex in $X \setminus \{F\}$ in $M_X$. Further, as $|E(M_X)| \ge 17$, the degree-sequence of the four vertices of $X$ in $M_X$ is either $(3,4,5,5)$ or $(3,5,5,5)$. Let $F_1$ and $F_2$ be two vertices of degree~$5$ in $M_X$ that belong to $X$, and let $F_3$ be the remaining vertex in $X \setminus \{F\}$. We note that in $H$, both $F_1$ and $F_2$ belong to $X_4$, and have an associated $E^*(X)$-pair in $M_X$. Since $H$ is linear, there are in fact distinct $E^*(X)$-pairs associated with $F_1$ and $F_2$. Renaming $F_1$ and $F_2$, if necessary, we may assume that the $E^*(X)$-pair, say $(e_1,e_2)$, associated with $F_1$ is distinct from $\{e,e'\}$. Thus, there is a vertex in $F_1$ that covers both $e_1$ and $e_2$, and we can cover a vertex in $\{e,e'\} \setminus \{e_1,e_2\}$ using a $\tau(F)$-transversal. Let $e_3$ be such a vertex covered from $F$. The vertex $F_3$ in $M_X$ has degree at least~$4$ in $M_X$, and can be used to cover a vertex in $E^*(X) \setminus \{e_1,e_2,e_3\}$, say $e_4$. Finally, the vertex $F_2$ in $M_X$ which has degree~$5$ in $M_X$, can be used to cover the vertex in $E^*(X) \setminus \{e_1,e_2,e_3,e_4\}$. We thereby obtain a contradiction to Claim~B.

We next consider the case when $|N_{M_X}(F)|=3$. Let $N_{M_X}(F) = \{e_1,e_2,e_3\}$, where $e = e_1$, and let $E^*(X) \setminus \{e_1,e_2,e_3\} = \{e_4,e_5\}$.  Since there is no $\tau(F)$-transversal covering $e_1$ and a vertex in $V(e_2) \cup V(e_3)$, we note
that $e_2$ and $e_3$ intersect $F$ in the same vertex (by Observation~\ref{property:special}(k)). Therefore, $(e_2,e_3)$ is a $E^*(X)$-pair associated with $F$. Since $d_{M_X}(e_1) \ge 3$, and since there are no triple edges in $M_X$, there is a vertex $F'$ in $X \setminus \{F\}$ adjacent to $e_1$ in $M_X$. Since neither $e_4$ nor $e_5$ has an edge to $F$, but they do have degree at least~$3$ in $M_X$, both $e_4$ and $e_5$ must be adjacent in $M_X$ to all three vertices in $X \setminus \{F\}$. We can therefore cover the edge $e_1$ from $F'$, cover both edges $e_2$ and $e_3$ from $F$, and cover the remaining two edges, $e_4$ and $e_5$, from $X \setminus \{F,F'\}$, thereby obtaining a $\tau(X)$-transversal covering $E^*(X)$, contradicting Claim~B. This completes the proof of Claim~H.15.2.~\smallqed

\medskip
We now return to the proof of Claim~H.15. By Claim~H.15.2, there does not exist a double edge in $M_X$, implying that $M_X$ is therefore a graph. As observed earlier, $|X_4| \ge 3$. By Claim~H.13, $|E(M_X)| \ge 4(|X|+1) - 3|\partial(X)| = 17$. Recall that each vertex of $X$ has degree at least~$2$ in $M_X$. If some vertex of $X$ has degree~$2$ in $M_X$, then the degree sequence of $X$ in $M_X$ is $(2,5,5,5)$, implying that there is a unique vertex of degree~$2$ in $M_X$. Since $|X_4| \ge 3$, there exist two distinct vertices $F,F' \in X_4$ of degree~$5$ in $M_X$. Both $F$ and $F'$ are associated with $E^*(X)$-pairs. Further, since $H$ is linear, the $E^*(X)$-pairs associated with $F$ and $F'$ are distinct. Renaming $F$ and $F'$, if necessary, we may assume that $F$ has an associated $E^*(X)$-pair, say $(e_1,e_2)$, such that the neighborhood of the degree-$2$ vertex in $X$ is not the set $\{e_1,e_2\}$.
If every vertex of $X$ has degree at least~$3$ in $M_X$, then we must still have a vertex, $F$, in $X$ of degree~$5$ in the graph $M_X$, and such a vertex is necessarily associated with an $E^*(X)$-pair, say $(e_1,e_2)$.  In both cases, we have therefore determined a vertex $F \in X$ of degree~$5$ in $M_X$ with an associated $E^*(X)$-pair, $(e_1,e_2)$, such that the neighborhood of the degree-$2$ vertex in $X$, if it exists, is not the set $\{e_1,e_2\}$.

We now consider the bipartite graph $M_X' = M_X - \{e_1,e_2,F\}$ with partite sets $X' = X \setminus \{F\}$ and $E'_X = E^*(X) \setminus \{e_1,e_2\}$. We note that $|E_X'| = 3$, and that there is no matching in $M_X'$ that matches $E'_X$ to $X'$, by Claim~B. By Hall's Theorem, there is a nonempty subset $S \subseteq E'_X$ such that in $M_X'$, $|N(S)| < |S|$. Since every vertex in $E^*(X)$ has degree at least~$3$ in $M_X$, we note that $|N_{M_X'}(S)| \ge 2$. Consequently, $|S| = 3$ and $|N_{M_X'}(S)| = 2$. Thus, $S = E'_X$. Let $F''$ be the vertex in $X$ that is adjacent to no vertex of $S$ in $M_X$. Therefore, $d_{M_X}(F'')=2$ and $N_{M_X}(F'') = \{e_1,e_2\}$. However, this is a contradiction to our choice of the $E^*(X)$-pair $(e_1,e_2)$, which was chosen so that the neighborhood of the degree-$2$ vertex in $X$, if it exists, is not the set $\{e_1,e_2\}$. This completes the proof of Claim~H.15.~\smallqed

\medskip

\ClaimX{H.16}
The case $|X| = 3$ and $|\partial(X)| = 1$ cannot occur.

\ProofClaimX{H.16} Suppose, to the contrary, that $|X| = 3$ and $|\partial(X)| = 1$.
By Claim~H.13, there are no triple edges in $M_X$. Also, by Claim~H.13(b), every vertex $F \in X$ has double edges to at most two vertices in $E^*(X)$ in $M_X$. Since $|\partial(X)| = 1$, by Claim~H.13(c), each vertex of $E^*(X)$ has degree at least~$3$ in $M_X$. By Claim~H.13(d), each vertex of $X$ has degree at least~$2$ in $M_X$.


Suppose that some vertex $F \in X$ has double edges to two distinct vertices, say $\{e_1,e_2\}$, of $E^*(X)$ in $M_X$. Let $E^*(X) \setminus \{e_1,e_2\} = \{e_3,e_4\}$. If there exists a perfect matching in $M_X - \{e_1,e_2,F\}$,
then we obtain a $X$-transversal covering all edges of $E^*(X)$ (by   Observation~\ref{property:special}(g) and \ref{property:special}(k)), contradicting Claim~B. Therefore, there does not exist a perfect matching in $M_X - \{e_1,e_2,F\}$. By Hall's Theorem, there is a nonempty subset $S \subseteq \{e_3,e_4\}$ such that  $|N_{M_X}(S) \setminus \{F\}| < |S|$.
If $|S|=1$, say $S=\{e_3\}$, then $e_3$ has degree at most~$1$ in $M_X$, a contradiction.  Therefore, $|S|=2$ and
$S=\{e_3,e_4\}$. Let $F' \in X \setminus \{F\}$ be defined such that $N_{M_X}(S) \setminus \{F\} = \{F'\}$. Since both $e_3$ and $e_4$ have degree at least~$3$ in $M_X$, we must have an edge from $e_3$ to $F$ and two edges from $e_3$ to $F'$.
Analogously, there is one edge from $e_4$ to $F$ and two edges from $e_4$ to $F'$. Therefore, by  Observation~\ref{property:special}(k),
there exists a $\tau(F')$-transversal covering both $e_3$ and $e_4$. Since $e_1$ and $e_2$ can be covered by a $\tau(F)$-transversal, we can therefore cover all edges of $E^*(X)$ with a
$X$-transversal, a contradiction
to Claim~B. Hence, every $F \in X$ has double edges to at most one vertex of $E^*(X)$ in $M_X$.

By Claim~H.13, $|E(M_X)| \ge 4(|X|+1) - 3|\partial(X)| =  16 - 3 = 13$. By the Pigeonhole Principle, there is therefore a vertex $F$ of $X$ of degree at least~$5$ in $M_X$. By Observation~\ref{property:special}
this implies that there exists a $E^*(X)$-pair, say $(e_1,e_2)$, associated with $F$. Let $E^*(X) \setminus \{e_1,e_2\} = \{e_3,e_4\}$. If there exists a perfect matching in $M_X - \{e_1,e_2,F\}$, then, as before, we obtain a $X$-transversal covering all edges of $E^*(X)$. Therefore, there is no perfect matching in $M_X - \{e_1,e_2,F\}$. By Hall's Theorem, there is a nonempty subset $S \subseteq \{e_1,e_2\}$ such that  $|N_{M_X}(S) \setminus \{F\}| < |S|$. If $|S|=1$, say $S=\{e_3\}$, then $e_3$ has degree at most~$2$ in $M_X$, a contradiction.  Therefore, $|S|=2$ and
$S=\{e_3,e_4\}$. Let $F' \in X \setminus \{F\}$ be defined such that $N_{M_X}(S) \setminus \{F\} = \{F'\}$.

Suppose firstly that the $E^*(X)$-pair, say $(e_1,e_2)$, can be chosen so that either $e_1$ or $e_2$ has a double edge to $F$ in $M_X$. In this case, both $e_3$ and $e_4$ have at most one edge to $F$. Since both $e_3$ and $e_4$ have degree at least~$3$ in $M_X$, this implies that both $e_3$ and $e_4$ have double edges to $F'$, contradicting the fact that every vertex in $X$ has double edges to at most one vertex of $E^*(X)$ in $M_X$. Therefore, all $E^*(X)$-pairs, $(e_1,e_2)$, have only single edges to $F$. Renaming $e_3$ and $e_4$, if necessary, we may assume that $e_3$ has a single edge to $F'$ in $M_X$. Thus, $e_3$ has a double edge to $F$ in $M_X$, which in turn implies that $e_4$ has a single edge to $F$ and a double edge to $F'$ in $M_X$. However, this implies that $F$ has a double edge to $e_3$, and single edges to $e_1,e_2$ and $e_4$. Therefore, by Observation~\ref{property:special}(k), there exists a $E^*(X)$-pair containing $e_3$ and one of $e_1,e_2$ and $e_4$, contradicting the fact that all $E^*(X)$-pairs have only single edges to $F$. This completes the proof of Claim~H.16.~\smallqed

\ClaimX{H.17}
The case $|X| = 2$ and $|\partial(X)| = 1$ cannot occur.

\ProofClaimX{H.17} Suppose, to the contrary, that $|X| = 2$ and $|\partial(X)| = 1$.
By Claim~H.13, there are no triple edges in $M_X$, and $|E(M_X)| \ge 4(|X|+1) - 3|\partial(X)| =  12 - 3 = 9$. By the Pigeonhole Principle, there is therefore a vertex $F$ of $X$ of degree at least~$5$ in $M_X$. By Observation~\ref{property:special}
this implies that there exists a $E^*(X)$-pair, say $(e_1,e_2)$, associated with $F$. Let  $E^*(X)=\{e_1,e_2,e_3\}$ and $X=\{F,F'\}$. If there is an edge from $e_3$ to $F'$ in $M_X$ we obtain a contradiction to Claim~B analogously to the above cases. If there is no edge from $e_3$ to $F'$ in $M_X$ we obtain a contradiction to $e_3$ having degree at least three in $M_X$. This contradiction completes
the proof of Claim~H.17.~\smallqed

\medskip
\ClaimX{H.18}
If $|\partial(X)| = 0$, then $|X| \le 6$ and the following holds. \1 \\
\hspace*{0.65cm} (a) If $|X|=6$, then $X=X_4$. \\
\hspace*{0.65cm} (b) If $|X|=5$, then $|X_4| \ge 3$.\\
\hspace*{0.65cm} (c) All vertices in $E^*(X)$ have degree~$4$ in $M_X$ and $|E(M_X)|=4|X|+4$. \\
\hspace*{0.65cm} (d) $M_X$ does not contain triple edges.\\
\hspace*{0.65cm} (e) No $F \in X$ has double edges to three distinct vertices of $E^*(X)$ in $M_X$.

\ProofClaimX{H.18}
Suppose that $\partial(X) = \emptyset$, and so $V(H) = V(X)$. Further, $H'$ is obtained from $H - V(X)$ by adding the set $X'$ of four new vertices and adding a $4$-edge, $e$, containing these four vertices. Thus, $H' = H_4$ and $|X'| = 4$. By Claim~H.4(b), we have
$8|X| \ge 8|X_4| + 5|X_{14}| + 4|X_{11}| + |X_{21}|  > 13|X| - 6|X'| - 8 = 8|X| + 5|X| - 32$. In particular, $5|X| \le 32$, and so $|X| \le 6$. Furthermore if $|X|=6$, then we must have $X=X_4$, while if $|X|=5$, then $|X_4| \ge 3$. The last three statements in Claim~H.18 follows directly from Claim~H.13.~\smallqed

\medskip
\ClaimX{H.19}
If $|\partial(X)| = 0$, then there does not exist a double edge between $e \in E^*(X)$ and $F \in X$ in $M_X$, such that $e$ belongs to a $E^*(X)$-pair associated with $F$.

\ProofClaimX{H.19} For the sake of contradiction, suppose that there is a double edge between $e \in E^*(X)$ and $F \in X$ in $M_X$, such that $e$ belongs to a $E^*(X)$-pair, say $(e,e')$, associated with $F$. We assume, further, that $e$ and $F$ are chosen so that the number of $E^*(X)$-pairs, $(e,e'')$, associated with $F$, where $e'' \in E^*(X) \setminus \{e,e'\}$, is maximized. As there is a double edge incident with $F$ in $M_X$, we note that $F \notin X_4$.

We now consider the bipartite multigraph $M_X^*$, where $M_X^*$ is defined as follows. Let $M_X^*$ be obtained from $M_X$ by removing the vertex $e \in E^*(X)$ and removing all edges $(e'',F)$, where $(e,e'')$ is not a $E^*(X)$-pair associated with $F$. We note that $(e',F)$ is an edge in $M_X^*$ and, therefore, $d_{M_X^*}(F) \ge 1$. Also due to  Observation~\ref{property:special}(k), we note that at most two edges (and no double edges) have been removed from $F$ to $E^*(X) \setminus \{e\}$ when constructing $M_X^*$.

First suppose that there is a perfect matching, $M$, in $M_X^*$. Let $(e^*,F)$ be an edge of the matching incident with $F$. By the definition of $M_X^*$, we note that $(e^*,e)$ is a $E^*(X)$-pair associated with $F$. By Observation~\ref{property:special}(g), we can find a $\tau(X)$-transversal covering $E^*(X)$, contradicting Claim~B.  Therefore, there is no perfect matching in $M_X^*$. By Hall's Theorem, there is a nonempty subset $S \subseteq E^*(X)\setminus\{e\}$ such that $|N_{M_X^*}(S)| < |S|$. Let $X^* = X \setminus  N_{M_X^*}(S)$.

We will show that $F \in X^*$. For the sake of contradiction, suppose that $F \notin X^*$, and so $F \in N_{M_X^*}(S)$. Since $F \notin X^*$, the neighbors of a vertex of $X^*$ in $M_X$ remain unchanged in $M_X^*$. We note that no vertex in $X^*$ has a neighbor that belongs to $S$ in $M_X^*$, and therefore also no neighbor that belongs to $S$ in $M_X$. This implies that in $M_X$ the vertices in $X^*$ are adjacent to at most $|E^*(X) \setminus S| = |E^*(X)| - |S| = |X|+1-|S| \le |X^*|$ vertices in $E^*(X)$. Thus, $|E^*(X^*)| \le |X^*|$ in $H$, contradicting Claim~D. Therefore, $F \in X^*$.
By Claim~H.18, $|X| \le 6$. As observed earlier, $F \notin X_4$, which implies that $|X| \ne 6$, by Claim~H.18. Thus, $|X| \le 5$.

Suppose that $|X|=5$. If $|S|=5$, then since $F$ is adjacent to $e'$ in $M_X^*$, this would imply that $F \in N_{M_X^*}(S)$, contradicting the fact that $F \in X^*$. Hence, $|S| \le 4$. Suppose that $|S|=4$ and $S=\{e_1,e_2,e_3,e_4\}$. We note that in this case, $E^*(X) \setminus S = \{e,e'\}$ and that $(e,e')$ is the only $E^*(X)$-pair associated with $F$. As observed earlier, at most two edges (and no double edges) were removed from $F$ to $E^*(X) \setminus \{e\}$ when constructing $M_X^*$. Thus, at least two vertices in $S$ did not have any edges to $F$ in $M_X$. Renaming vertices of $S$, if necessary, we may assume that neither $e_1$ nor $e_2$
has an edge to $F$ in $M_X$. By Claim~H.18(c), all vertices in $E^*(X)$ have degree~$4$ in $M_X$, implying that both $e_1$ and $e_2$ have a double edge to a vertex in $N_{M_X^*}(S)$. Since $|X_4| \ge 3$ and $F \notin X_4$, $e_1$ and $e_2$ both have double edges to the same vertex $F' \in N_{M_X^*}(S)$. Thus, $X \setminus \{F,F'\} = X_4$. Furthermore, as $e_3$ and $e_4$ have degree~$4$ in $M_X$, and therefore degree at least~$3$ in $M_X^*$, both $e_3$ and $e_4$ are adjacent to $F'$ in $M_X^*$. Therefore, by Observation~\ref{property:special}(k) and~\ref{property:special}($\ell$), there is a $E^*(X)$-pair containing one vertex from $\{e_1,e_2\}$ and one vertex from $\{e_3,e_4\}$ associated with $F'$. Renaming vertices of $S$, if necessary, we may assume that $(e_1,e_3)$ is a $E^*(X)$-pair associated with $F'$. If we had used this $E^*(X)$-pair instead of $(e,e')$, and $F'$ instead of $F$, then we note that $e_1$ is a double edge to $F'$ and both $(e_1,e_2)$ and $(e_1,e_3)$ are $E^*(X)$-pairs associated with $F'$, while $e$ is a double edge to $F$ but there is only one $E^*(X)$-pair of the form $(e,e'')$, where $e'' \in E^*(X) \setminus \{e,e'\}$, (namely the pair $(e,e')$) that is associated with $F$. This contradicts our choice of $e$ and $F$. Therefore, $|S| \le 3$.
Suppose that $|S|=3$ and $S=\{e_1,e_2,e_3\}$. Since at most two edges (and no double edges) were removed from $F$ to $E^*(X) \setminus \{e\}$ when constructing $M_X^*$, at least one vertex in $S$, say $e_1$, has the same degree in $M_X^*$ as in $M_X$. Thus, $e_1$ has degree~$4$ in $M_X^*$, implying that $e_1$ has double edges to both vertices in $N_{M_X^*}(S)$. However, this implies that no vertex in $N_{M_X^*}(S)$ belongs to $X_4$. As observed earlier, $F \notin X_4$. This implies that $|X_4| \le 2$, a contradiction. Therefore, $|S| \le 2$.
However in this case the vertices in $S$ cannot have degree~$4$ in $M_X$, a contradiction. This completes the case when $|X|=5$.

Suppose that $|X|=4$. If $|S|=4$, then since $F$ is adjacent to $e'$ in $M_X^*$, this would imply that $F \in N_{M_X^*}(S)$, contradicting the fact that $F \in X^*$. Hence, $|S| \le 3$. Suppose that $|S|=3$ and $S=\{e_1,e_2,e_3\}$. We note that in this case, $E^*(X) \setminus S = \{e,e'\}$. Since at most two edges (and no double edges) were removed from $F$ to $E^*(X) \setminus \{e\}$  when constructing $M_X^*$, at least one vertex in $S$, say $e_3$, has the same degree in $M_X^*$ as in $M_X$. Thus, $e_3$ has degree~$4$ in $M_X^*$, implying that $e_3$ has double edges to both vertices in $N_{M_X^*}(S) = \{F_1,F_2\}$.  Since both $e_1$ and $e_2$ have degree~$4$ in $M_X$, and degree at least~$3$ in $M_X^*$, both $e_1$ and $e_2$ have a double edge to $F_1$ or $F_2$. They cannot have a double edge to the same vertex in $\{F_1,F_2\}$, for otherwise such a vertex would then have double edges to all three vertices in $S$, contradicting Claim~H.18(e). Renaming $e_1$ and $e_2$, if necessary, we may assume that $e_1$ has a double edge to $F_1$ and $e_2$ has a double edge to $F_2$ in $M_X^*$ (and in $M_X$). We can now cover $\{e,e'\}$ from $F$, $\{e_1,e_3\}$ from $F_1$, and $e_2$ from $F_2$, contradicting Claim~B.  Therefore, $|S| \le 2$.  However in this case the vertices in $S$ cannot have degree~$4$ in $M_X$, a contradiction. This completes the case when $|X|=4$.

Consider the case when $|X|=3$. If $|S|=3$, then we contradict the fact that $F \in X^*$. Hence, $|S| \le 2$. However in this case the vertices in $S$ cannot have degree~$4$ in $M_X$, a contradiction. This completes the case when $|X|=3$. Analogously, we get a contradiction when $|X|=2$.  This completes the proof of Claim~H.19.~\smallqed

\medskip
\ClaimX{H.20}
The case $|X| \ge 6 $ and $|\partial(X)| = 0$ cannot occur.

\ProofClaimX{H.20}
Suppose, to the contrary, that $|X| \ge 6$ and $|\partial(X)| = 0$. By Claim~H.18,  $|X|=6$, $X=X_4$ and $|E(M_X)| = 4|X| + 4 = 28$. Since $X = X_4$, the multigraph $M_X$ is a graph. By the Pigeonhole Principle, there is a vertex $F$ of $X$ of degree at least~$5$ in $M_X$. By Observation~\ref{property:special}
this implies that there exists a $E^*(X)$-pair, say $(e_1,e_2)$, associated with $F$. Among all such $E^*(X)$-pairs, we choose the pair $(e_1,e_2)$ so that the number of neighbors of $e_1$ and $e_2$ in $M_X$ of degree less than~$4$ is a minimum.
Since $F \in X_4$, we note that in $H$, the edges $e_1$ and $e_2$ intersects $F$ in the same vertex.

We now consider the bipartite multigraph $M_X' = M_X - \{e_1,e_2,F\}$ with partite sets $X' = X \setminus \{F\}$ and $E'_X = E^*(X) \setminus \{e_1,e_2\}$. If there exists a perfect matching in $M_X'$, then we can find a $\tau(X)$-transversal covering $E^*(X)$, contradicting Claim~B. Therefore, there is no perfect matching in $M_X'$. By Hall's Theorem, there is a nonempty subset $S \subseteq E'_X$ such that in $M_X'$, $|N(S)|<|S|$. By Claim~H.18(c), all vertices in $E^*(X)$ have degree~$4$ in $M_X$. Thus, each vertex of $E'_X$ has degree at least~$3$ in $M_X'$, implying that in $M_X'$, $|S| > |N(S)| \ge 3$. Hence, $|S| \in \{4,5\}$.

Suppose that $|S|=5$. If $|N_{M_X}(S)| \le 3$, then let $X^* = X \setminus N_{M_X}(S)$ and note that $|E^*(X^*)| \le |\{e_1,e_2\}| \le |X^*|$, a contradiction to
Claim~D. Therefore, $|N_{M_X}(S)|=4$ and let $F(e_1,e_2)$ denote the vertex in $X \setminus (N_{M_X}(S) \cup \{F\})$. The degree of every vertex of $X$ is at least~$2$ in $M_X$, implying that the vertex $F(e_1,e_2)$ is adjacent in $M_X$ to both $e_1$ and $e_2$ (and to no other vertex of $E^*(X)$), and therefore has degree~$2$ in $M_X$. Hence, in this case when $|S|=5$, the vertices $e_1$ and $e_2$ have a neighbor of degree~$2$. Before completing our proof of this case when $|S|=5$, we first consider the case when $|S| = 4$.

Suppose that $|S|=4$. Since every vertex of $E'_X$ has degree at least~$3$ in $M_X'$,  we note that in $M_X'$, $|N(S)| = 3$. Further, in $M_X'$, every vertex of $S$ is adjacent to every vertex of $N(S)$. This implies that in $M_X$, the vertex $F$ is adjacent to all four vertices of $S$, and therefore $F$ has degree at least~$6$ in $M_X$. Let $S = \{e_3,e_4,e_5,e_6\}$ and let $N_{M_X'}(S) = \{F_1,F_2,F_3\}$. Let $X \setminus \{F_1,F_2,F_3\} = \{F,F_5,F_6\}$. We note that $F_5$ (respectively, $F_6$) has degree~$2$ or~$3$ in $M_X$, and is adjacent to at least one of $e_1$ and $e_2$. Hence, in this case when $|S|=4$, there are two vertices of degree less than~$4$ adjacent to $e_1$ or $e_2$. Recall that in $H$, the edges $e_1$ and $e_2$ intersects $F$ in the same vertex, implying that in $H$, two edges in $\{e_3,e_4,e_5,e_6\}$ intersect $F$ in the same vertex. Renaming vertices of $S$, if necessary, we may assume that in $H$, the edges $e_3$ and $e_4$ both intersect $F$ in the same vertex, and therefore form a $E^*(X)$-pair associated with $F$.
If we had considered the $E^*(X)$-pair $(e_3,e_4)$ instead of the $E^*(X)$-pair $(e_1,e_2)$, then we note that every neighbor of $e_3$ and $e_4$ has degree at least~$4$ in $M_X$. This contradicts our choice of the $E^*(X)$-pair, $(e_1,e_2)$. Therefore, $|S|=4$ is not possible.

We now return to the case when $|S|=5$. Given a $E^*(X)$-pair, $(e_1,e_2)$, we showed that there exists a vertex $F(e_1,e_2)$ in $X \setminus \{F\}$ that has degree~$2$ in $M_X$, and has $e_1$ and $e_2$ as its neighbors. Since $|S| \le 4$ is not possible, this implies that for every $E^*(X)$-pair, $(e,e')$, associated with some vertex $F' \in X$, there exists a vertex $F(e,e')$ in $X \setminus \{F'\}$ that has degree~$2$ in $M_X$, and has $e$ and $e'$ as its neighbors. We note that if a vertex $F^* \in X$ has degree~$4+i$ in $M_X$, then it is associated with at least $i$ $E^*(X)$-pairs. Thus, since the vertex $F(e_1,e_2)$ has degree~$2$ in $M_X$, and since $|E(M_X)|=28$, there are at least six $E^*(X)$-pairs, which are all distinct. To each one we associated a vertex of degree~$2$ in $M_X$ that belongs to $X$. However this implies that there are only $12$ edges in $M_X$, a contradiction to $|E(M_X)|=28$. This completes the proof of Claim~H.20.~\smallqed

\medskip
\ClaimX{H.21}
The case $|X| = 5 $ and $|\partial(X)| = 0$ cannot occur.

\ProofClaimX{H.21}
Suppose, to the contrary, that $|X| = 5$ and $|\partial(X)| = 0$. By Claim~H.18, $|X_4| \ge 3$ and $|E(M_X)| = 4|X| + 4 = 24$.
Also, by Claim~H.18, $M_X$ does not contain triple edges. By Claim~H.19 and Observation~\ref{property:special}(k), no $F \in X$ has double edges to two distinct vertices of $E^*(X)$ in $M_X$. By Claim~H.18(c), all vertices in $E^*(X)$ have degree~$4$ in $M_X$.

Since $|E(M_X)|=24$ and $|X|=5$, there is a vertex $F$ of $X$ of degree at least~$5$ in $M_X$. By Observation~\ref{property:special}
this implies that there exists a $E^*(X)$-pair, say $(e_1,e_2)$, associated with $F$.
We now consider the bipartite multigraph $M_X' = M_X - \{e_1,e_2,F\}$ with partite sets $X' = X \setminus \{F\}$ and $E'_X = E^*(X) \setminus \{e_1,e_2\}$. If there exists a perfect matching in $M_X'$, then we can find a $\tau(X)$-transversal covering $E^*(X)$, contradicting Claim~B. Therefore, there is no perfect matching in $M_X'$. By Hall's Theorem, there is a nonempty subset $S \subseteq E'_X$ such that in $M_X'$, $|N(S)|<|S|$. We note that $1 \le |S| \le 4$.

If $|S|=1$, then the vertex in $S$ has degree at most~$2$ in $M_X$, a contradiction. If $|S|=2$, then the two vertices in $S$ must have double edges to $F$ (and to the vertex in $N_{M_X'}(S)$), contradicting the fact that no vertex in $X$ has double edges to two distinct vertices of $E^*(X)$ in $M_X$. Therefore, $|S| \ge 3$. Suppose that $|S|=3$, and let $S=\{e_3,e_4,e_5\}$. For each $i \in [3]$, the vertex $e_i \in S$ has degree~$4$ in $M_X$ and $N_{M_X}(e_i) \subseteq \{F\} \cup N_{M_X'}(S)$, implying that $|N_{M_X'}(S)|=2$ and that $e_i \in S$ is adjacent to a double edge. Let $N_{M_X'}(S) = \{F_1,F_2\}$. Since  no vertex in $X$ is adjacent to two double edges,  we may assume, renaming vertices if necessary, that $e_3$ has a double edge to $F$, $e_4$ has a double edge to $F_1$, and $e_5$ has a double edge to $F_2$. However,  this implies that none of $F$, $F_1$ or $F_2$ belong to $X_4$, contradicting the fact that $|X_4| \ge 3$. Therefore, $|S| \ge 4$.

We have shown that $|S|=4$. Let $S=\{e_3,e_4,e_5,e_6\}$. If $|N_{M_X'}(S)| \le 2$, then let $X^* = X \setminus (N_{M_X}(S) \cup \{F\})$ and note
that $|E^*(X^*)| \le |\{e_1,e_2\}| \le |X^*|$, a contradiction to Claim~D. Therefore, $|N_{M_X'}(S)| = 3$. Let $N_{M_X'}(S) = \{F_1,F_2,F_3\}$ and let $X \setminus N_{M_X'}(S) = \{F,F(e_1,e_2)\}$. We note  that $N_{M_X}(F(e_1,e_2)) \subseteq \{e_1,e_2\}$. By Claim~E, $|N_{M_X}(F(e_1,e_2))| \ge 2$, implying that  $N_{M_X}(F(e_1,e_2)) = \{e_1,e_2\}$. Since no vertex of $X$ has double edges in $M_X$ to two distinct vertices in $E^*(X)$, we note that $d_{M_X}(F(e_1,e_2)) \le 3$. Thus, for each $E^*(X)$-pair, $(e',e'')$ say, there is a vertex $F' \in X$ such that $N_{M_X}(F(e',e'')) = \{e',e''\}$ and $d_{M_X}(F(e',e'')) \le 3$.

We show next that if $F' \in X$ has degree $4+i$ in $M_X$, then it is associated with at least $i$ $E^*(X)$-pairs. We may assume that $F' \notin X_4$, for otherwise this property clearly holds if $F' \in X_4$. If there are three neighbors $\{f_1,f_2,f_3\}$ of $F'$ in $M_X$, no two of which form a $E^*(X)$-pair associated with $F'$, then every vertex in $N_{M_X}(F') \setminus \{f_1,f_2,f_3\}$ forms a $E^*(X)$-pair with one of the vertices in $\{f_1,f_2,f_3\}$, by Observation~\ref{property:special}(j).
This gives us at least $|N_{M_X}(F')|-3$ distinct $E^*(X)$-pairs associated with $F'$. Since no vertex in $X$ has double edges to two distinct vertices in $E^*(X)$,
we note that $|N_{M_X}(F')| \ge d_{M_X}(F')-1 = 3 + i$ and the result follows. Therefore, we may assume that every set of three neighbors of $F'$ in $M_X$ yields at least one $E^*(X)$-pair associated with $F'$. If every two neighbors of $F'$ in $M_X$ form a $E^*(X)$-pair associated with $F'$, then we clearly have enough $E^*(X)$-pairs, so we assume that there are two neighbors, $f_1$ and $f_2$ say, of $F'$ in $M_X$ which do not form a $E^*(X)$-pair associated with $F'$.
However, since every set of three neighbors of $F'$ in $M_X$ contains a $E^*(X)$-pair, for every $f' \in N_{M_X}(F') \setminus \{f_1,f_2\}$ at least one of $(f',f_1)$ and $(f',f_2)$ is a $E^*(X)$-pair associated with $F'$, yielding at least $d_{M_X}(F') - 2 = 2 + i$ $E^*(X)$-pairs associated with $F'$. Therefore, we have shown that if $F' \in X$ has degree $4+i$ in $M_X$, then it is associated with at least $i$ $E^*(X)$-pairs. Recall that $d_{M_X}(F(e_1,e_2)) \le 3$ and  $|E(M_X)|=24$. Thus, there are at least five $E^*(X)$-pairs.

If $(e,e')$ is an arbitrary $E^*(X)$-pair, then, by the linearity of $H$, $(e,e')$ can be associated with at most one vertex from $X_4$. Hence, since $|X \setminus X_4| \le 2$, the $E^*(X)$-pair can be associated with at most three vertices of $X$. Since there are at least five $E^*(X)$-pairs, there are therefore at least $\lceil \frac{5}{3} \rceil = 2$ distinct $E^*(X)$-pairs associated with distinct vertices of $X$. Let $(e',e'')$ and $(f',f'')$ be two such $E^*(X)$-pairs. Thus, there exist two vertices in $X$, say $F'$ and $F''$, of degree at most~$3$ in $M_X$, where $F'$ has neighborhood $\{e',e''\}$ and $F''$ has neighborhood $\{f',f''\}$, in $M_X$.

Suppose that both $F'$ and $F''$ have degree~$3$ in $M_X$. Then, both $F'$ and $F''$ are incident with double edges in $M_X$ and therefore do not belong to $X_4$. Since $|X_4| \ge 3$, this implies that $X \setminus \{F',F''\} = X_4$. Hence, no vertex in $X \setminus \{F',F''\}$ is incident with a double edge in $M_X$, and therefore has degree at most~$|E^*(X)| = 6$ in $M_X$. Since $|E(M_X)|=24$, the degree sequence of the vertices in $X$ must be $(3,3,6,6,6)$ in $M_X$. However in this case, each of the three vertices of $X$ of degree~$6$ belongs to $X_4$, and is therefore intersected by six edges in $H$, which gives rise to at least two associated $E^*(X)$-pairs. Further, by the linearity of $H$, the $E^*(X)$-pairs associated with two distinct vertices of $X$ that belong to $X_4$ are distinct. Hence, there are at least six distinct $E^*(X)$-pairs associated with distinct vertices of $X$. However, as observed earlier, for each $E^*(X)$-pair, $(e,f)$, there is a vertex in $X$ whose neighborhood is precisely $\{e,f\}$, implying that $|X| \ge 6$, a contradiction. Therefore, at least one of $F'$ and $F''$ has degree~$2$ in $M_X$.

Since $|E(M_X)|=24$, and since the sum of the degrees of $F'$ and $F''$ in $M_X$ is at most~$5$, the sum of the degrees of the three vertices of $X \setminus \{F',F''\}$ in $M_X$ is at least~$19$. As observed earlier,  if $F' \in X$ has degree $4+i$ in $M_X$, then it is associated with at least $i$ $E^*(X)$-pairs. Hence there are at least seven $E^*(X)$-pairs. Since each $E^*(X)$-pair can be associated with at most three vertices of $X$, there are therefore at least $\lceil \frac{7}{3} \rceil = 3$ distinct $E^*(X)$-pairs associated with distinct vertices of $X$. Thus, there exist at least three vertices in $X$ of degree at most~$3$ in $M_X$. Since $|E(M_X)|=24$, the sum of the remaining two vertices of $X$ in $M_X$ is at least~$15$, implying that $X$ has a vertex of degree at least~$8$ in $M_X$. Such a vertex has triple edges to a vertex of $E^*(X)$ or double edges to two distinct vertices in $E^*(X)$, a contradiction. This completes the proof of Claim~H.21.~\smallqed


\medskip
\ClaimX{H.22}
The case $|X| = 4 $ and $|\partial(X)| = 0$ cannot occur.

\ProofClaimX{H.22}
Suppose, to the contrary, that $|X| = 4$ and $|\partial(X)| = 0$. By Claim~H.18, $|E(M_X)| = 4|X| + 4 = 20$. Also, by Claim~H.18, $M_X$ does not contain triple edges. By Claim~H.19 and Observation~\ref{property:special}(k), no $F \in X$ has double edges to two distinct vertices of $E^*(X)$ in $M_X$. By Claim~H.18(c), all vertices in $E^*(X)$ have degree~$4$ in $M_X$.

Since $|E(M_X)|=20$ and $|X|=4$, there is a vertex $F$ of $X$ of degree at least~$5$ in $M_X$. By Observation~\ref{property:special}
this implies that there exists a $E^*(X)$-pair associated with $F$. Among all such $E^*(X)$-pairs, we choose the pair $(e_1,e_2)$ to maximize the number of edges between $\{e_1,e_2\}$ and $F$ in $M_X$. Since $F$ has double edges to at most one vertex of $E^*(X)$ in $M_X$, we note that there are at most three edges between $\{e_1,e_2\}$ and $F$ in $M_X$.

We now consider the bipartite multigraph $M_X' = M_X - \{e_1,e_2,F\}$ with partite sets $X' = X \setminus \{F\}$ and $E'_X = E^*(X) \setminus \{e_1,e_2\}$. If there exists a perfect matching in $M_X'$, then we can find a $\tau(X)$-transversal covering $E^*(X)$, contradicting Claim~B. Therefore, there is no perfect matching in $M_X'$. By Hall's Theorem, there is a nonempty subset $S \subseteq E'_X$ such that in $M_X'$, $|N(S)| < |S|$. We note that $1 \le |S| \le 3$.

If $|S|=1$, then the vertex in $S$ has degree at most~$2$ in $M_X$, a contradiction. If $|S|=2$, then the two vertices in $S$ must have double edges to $F$ (and to the vertex in $N_{M_X'}(S)$), contradicting the fact that no vertex in $X$ has double edges to two distinct vertices of $E^*(X)$ in $M_X$. Therefore, $|S| = 3$.

Let $S=\{e_3,e_4,e_5\}$. For each $i \in [3]$, the vertex $e_i \in S$ has degree~$4$ in $M_X$ and $N_{M_X}(e_i) \subseteq \{F\} \cup N_{M_X'}(S)$, implying that $|N_{M_X'}(S)|=2$ and that $e_i \in S$ is adjacent to a double edge. Let $N_{M_X'}(S) = \{F_1,F_2\}$. Since  no vertex in $X$ is adjacent to two double edges,  we may assume, renaming vertices if necessary, that $e_3$ has a double edge to $F$, $e_4$ has a double edge to $F_1$, and $e_5$ has a double edge to $F_2$. Furthermore, $e_4$ and $e_5$ also have (single) edges to $F$. Thus, $F$ has degree~$6$ in $M_X$, with a single edge to each of $e_1, e_2, e_4, e_5$ and a double edge to $e_3$. By Observation~\ref{property:special}
this implies that there exists a $E^*(X)$-pair, $(e_3,e)$, associated with $F$, for some $e \in \{e_1,e_2,e_4,e_5\}$. The number of edges between $\{e_1,e_2\}$ and $F$ in $M_X$ is~$2$, while the number of edges between $\{e_3,e\}$ and $F$ in $M_X$ is~$3$. This contradicts our choice of the $E^*(X)$-pair $\{e_1,e_2\}$, and completes the proof of Claim~H.22.~\smallqed

\medskip
\ClaimX{H.23}
The case $|\partial(X)| = 0$ cannot occur.

\ProofClaimX{H.23}
Suppose, to the contrary, that $|\partial(X)| = 0$. By Claim~H.20, Claim~H.21 and Claim~H.22 we must have $|X| \le 3$.
By Claim~H.18, $M_X$ does not contain triple edges. By Claim~H.19 and Observation~\ref{property:special}(k), no $F \in X$ has double edges to two distinct vertices of $E^*(X)$ in $M_X$. Since every vertex of $E^*(X)$ has degree~$4$ in $M_X$,  we note that $|X| \le 2$ is therefore impossible. Therefore, $|X|=3$. Let $E^*(X) = \{e_1,e_2,e_3,e_4\}$. Since every vertex of $E^*(X)=\{e_1,e_2,e_3,e_4\}$ has a double edge to a vertex of $X$ in $M_X$, by the Pigeonhole Principle, some vertex in $X$ has double edges to two distinct vertices in $E^*(X)$, a contradiction. This completes the proof of Claim~H.23.~\smallqed

\medskip
\ClaimX{H.24}
$|X| =  1$.

\ProofClaimX{H.24}  By Claim~H.10, Claim~H.14, Claim~H.15, and Claim~H.23, the case $|X| \ge 4$ cannot occur. By Claim~H.11, Claim~H.16, Claim~H.17, and Claim~H.23, the case $|X| = 3$ cannot occur. By Claim~H.6, Claim~H.8, Claim~H.12, Claim~H.17, and Claim~H.23, the case $|X| = 2$ cannot occur. Therefore, $|X| =  1$.~\smallqed

\medskip
Throughout the remaining subclaims of Claim~H, we implicitly use the fact that $|X|=1$, and we let $e_1$ and $e_2$ be the two edges of $E^*(X)$ that intersect $X$ in $H$.

\medskip
\ClaimX{H.25}
$|\partial(X)| \ge 4$.

\ProofClaimX{H.25}
Suppose, to the contrary, that $|\partial(X)| \le 3 $.
Recall that $e_1$ and $e_2$ were the edges intersecting $X$ in $H$.
If $|V(e_1) \cap V(X)| =1$, then $\partial(X) \subseteq V(e_1)$ and therefore $|V(e_2) \cap V(X)| \ge 3$ as $H$ is linear. Therefore, by Observation~\ref{property:special}($\ell$), we can cover $E^*(X)$ with a $\tau(X)$-transversal, a contradiction to Claim~B.  If $|V(e_1) \cap V(X)| \ge 3$, we analogously get a contradiction using  Observation~\ref{property:special}($\ell$) and Claim~B. If $|V(e_1) \cap V(X)| = 2$, then $|V(e_2) \cap V(X)| \ge 2$, as $|\partial(X)| \le 3$.
We now obtain a contradiction to Claim~B using  Observation~\ref{property:special}(k).
Therefore, $|\partial(X)| \ge 3 $.~\smallqed

\medskip
\ClaimX{H.26}
If $H'$ is linear, then
one of the following hold. \\
\hspace*{1cm}
{\rm (a)} $\defic(H') \in \{8,10\}$.
\\
\hspace*{1cm}
{\rm (b)} If $\defic(H')=8$, then $X=X_4$.

\ProofClaimX{H.26}
By Claim~H.4(a), we note that
$8|X_4| + 5|X_{14}| + 4|X_{11}| + |X_{21}|  > 13|X| - 6|X'| - \defic(H')$. As $|X|=1$ and $|X'|=0$, we note that $\defic(H') \ge 6$. By Claim~H.3, $\defic(H') = \defic_{H'}(Y)$, where $|Y| = 1$ and $e$ is an edge of the hypergraph in the special $H'$-set $Y$. Thus, either $Y$ is an $H_4$-component of $H'$, in which case $\defic(H')=8$, or $Y$ is an $H_{10}$-component of $H'$, in which case $\defic(H')=10$. This proves Part~(a).

Suppose that $\defic(H')=8$, and so $Y$ is an $H_4$-component of $H'$ consisting of the edge~$e$. We now let $H^* = H - V(X) - V(e)$. Equivalently, $H^*$ is obtained from $H'$ by deleting the $H_4$-component $Y$. If $\defic(H^*) > 0$ and $Y^*$ is a special $H^*$-set such that $\defic(H^*) = \defic_{H^*}(Y^*) > 0$, then $\defic(H') \ge \defic_{H'}(Y \cup Y^*) > 8$, a contradiction. Hence, $\defic(H^*) = 0$. By Observation~\ref{property:special}(e), we can choose a $\tau(X)$-set, $T_X$, to contain a vertex of $e_1$. Every $\tau(H^*)$-set can be extended to a transversal of $H$ by adding to the set $T_X$ and an arbitrary vertex of $e_2$, implying that $\tau(H) \le \tau(H^*) + \tau(X) + 1$. We note that $n(H^*) = n(H) -  n(X) - 4$ and $m(H^*) = m(H) - m(H) - 2$. We show that $X = X_4$. Suppose, to the contrary, $X \ne X_4$. Thus, $X = X_i$ for some $i \in \{11,14,21\}$. If $X = X_{11}$, then $45\tau(X) = 6n(X) + 13m(X) + 4$. If $X = X_{14}$, then $45\tau(X) = 6n(X) + 13m(X) + 5$. If $X = X_{21}$, then $45\tau(X) = 6n(X) + 13m(X) + 1$. In all cases, $45\tau(X) \le 6n(X) + 13m(X) + 5$. Thus,

\[
\begin{array}{lcl}
\Phi(H^*) & = & \xi(H^*) - \xi(H)
\2 \\
& = & 45\tau(H^*) - 6n(H^*) - 13m(H^*) - \defic(H^*)   \\
& & \hspace*{1cm} - 45\tau(H) + 6n(H) + 13m(H) + \defic(H) \2 \\
& \ge & 45(\tau(H) - \tau(X) - 1) - 6(n(H) -  n(X) - 4) \\
& & \hspace*{1cm} - 13(m(H) - m(H) - 2)  \\
& & \hspace*{1.5cm} - 45\tau(H) + 6n(H) + 13m(H) \2 \\
& = & -45\tau(X) + 6n(X) + 13m(X) - 45 + 24 + 26 \\
& \ge & -5 - 45 + 24 + 26 \\
& = & 0,
\end{array}
\]
a contradiction to Claim~G. This proves Part~(b), and completes the proof of Claim~H.26.~\smallqed


\medskip

\ClaimX{H.27}
$\defic(H') \le 8$.

\ProofClaimX{H.27}
By Claim~H.3, either $\defic(H') = 0$ or $\defic(H') = \defic_{H'}(Y)$ where $|Y| = 1$ and $e$ is an edge of the hypergraph in the special $H'$-set $Y$. For the sake of contradiction suppose that $\defic(H') > 8$, which implies that $\defic(H') = \defic_{H'}(Y)$ where $|Y| = 1$ and $Y$ is an $H_{10}$-component, $F$, of $H'$ and $e$ is an edge of $F$. By Claim~H.25, $|\partial(X)| \ge 4$.

Suppose that $|\partial(X)| = 4$. If $|V(e_1) \cap V(X)| =1$ and $|V(e_2) \cap V(X)| =1$, then the linearity of $H$ implies that $|\partial(X)| \ge 5$, a contradiction. Renaming the edges $e_1$ and $e_2$, we may assume that  $|V(e_2) \cap V(X)| \ge 2$. By Claim~B, there is no $\tau(X)$-transversal covering both $e_1$ and $e_2$. Therefore, by Observation~\ref{property:special}(k) and~\ref{property:special}($\ell$), we note  that  $|V(e_1) \cap V(X)| = 1$ and $|V(e_2) \cap V(X)| = 2$. Thus, by  Observation~\ref{property:special}(m), we have $X = X_{11}$. Let $T_1 = V(e_1) \cap V(X)$ and let $T_2 = V(e_2) \cap V(X)$, and so $T_1$ and $T_2$ are vertex-disjoint subsets of $X$ such that $|T_1| = 1$ and $|T_2| = 2$ where $T_2$ contains two vertices that are not adjacent in $X$. Further by Observation~\ref{property:special}(m), one degree-$1$ vertex in $X = X_{11}$ belongs to $T_1$ and the other degree-$1$ vertex to $T_2$, and the second vertex of $T_2$ is adjacent to the vertex of $T_1$. As $|\partial(X)| = 4$, the edges $e_1$ and $e_2$ intersect in a vertex in $\partial(X)$. As every edge in $H_{10}$ is equivalent, the above observations imply that we get that $H = H_{21,6}$, contradicting Claim~C. Hence, $|\partial(X)| \ge 5$.

Since $e$ is an edge of $F$, we note that $|\partial(X) \cap V(F)| \ge 4$. Suppose that $|\partial(X) \cap V(F)| \ge 5$. Let $z \in (\partial(X) \cap V(F))\setminus V(e)$ be arbitrary. Removing from $F$ the edge $e$ and all edges incident to $z$, we are left with only two edges. Further, these two remaining edges intersect in a vertex, say $z'$. Let $H''$ be obtained from $H'$ by removing the two vertices $z$ and $z'$ (and all edges incident with $z$ and $z'$) and removing the edge $e$; that is, $H'' = H' - \{z,z'\} - e$. Equivalently, $H''$ is obtained from $H'$ by deleting the $H_{10}$-component $F$. If $\defic(H'') > 0$ and $Y''$ is a special $H''$-set such that $\defic(H'') = \defic_{H''}(Y'') > 0$, then $\defic(H') \ge \defic_{H'}(Y \cup Y'') > \defic(H')$, a contradiction. Hence, $\defic(H'') = 0$. Every $\tau(H'')$-transversal can be extended to a transversal of $H$ by adding to a $\tau(X)$-transversal the two vertices $\{z,z'\}$, implying that $\tau(H) \le \tau(H'') + \tau(X) + 2$. We note that $45\tau(X) = 6n(X) + 13m(X) + \defic(X) \le 6n(X) + 13m(X) + 10$. Thus,
\[
\begin{array}{lcl}
\Phi(H'') & = & \xi(H'') - \xi(H)
\2 \\
& = & 45\tau(H'') - 6n(H'') - 13m(H'') - \defic(H'')   \\
& & \hspace*{1cm} - 45\tau(H) + 6n(H) + 13m(H) + \defic(H) \2 \\
& \ge & 45(\tau(H) - \tau(X) - 2) - 6(n(H) - n(X) - 10) \\
& & \hspace*{1cm} - 13(m(H) - m(X) - 6)  \\
& & \hspace*{1.5cm} - 45\tau(H) + 6n(H) + 13m(H) \2 \\
& = & -45\tau(X) + 6n(X) + 13m(X) - 90 + 60 + 78 \\
& \ge & -10 - 90 + 60 + 78 \\
& > & 0,
\end{array}
\]
a contradiction to Claim~G. Hence, $|\partial(X) \cap V(F)| = 4$, implying that $\partial(X) \cap V(F) = V(e)$. Let $w \in \partial(X) \setminus V(F)$ be arbitrary. Assume we had picked the edge $e$ to contain three vertices from $V(F)$ and the vertex $w$.  In this case, $H'$ would be linear. Further, $\defic(H') = \defic_{H'}(Y')$, where $|Y'| = 1$ and $e$ is an edge of the hypergraph in the special $H'$-set $Y'$. However, $Y$ is neither an $H_4$-component nor an $H_{10}$-component of $H'$, implying that $\defic(H') \notin \{8,10\}$, contradicting Claim~H.26.
This completes the proof of Claim~H.27.~\smallqed



\ClaimX{H.28}
$H'$ is not linear.

\ProofClaimX{H.28}
Suppose, to the contrary, that $H'$ is linear.
By Claim~H.26 and Claim~H.27, $X=X_4$ and $\defic(H')=8$. As $X=X_4$ and $H$ is linear, we note that $|\partial(X)| \ge 5$. By Claim~H.3, the edge $e$ is an isolated edge in $H'$ and therefore there are at least four isolated vertices in $\partial(X)$ in $H'-e$. Suppose that there is a non-isolated vertex in $\partial(X)$ in $H - V(X)$. In this case, choosing the edge $e$ to contain this vertex together with three isolated vertices of $H - V(X)$ that belong to $\partial(X)$, yields a new linear hypergraph $H'$ in which the newly chosen edge $e$ does not belong to a $H_4$- or $H_{10}$-component, contradicting Claim~H.3 and Claim~H.26. Therefore, all vertices in $\partial(X)$ are isolated vertices in $H - V(X)$, implying that $|V(H)| \in \{9,10\}$, $|E(H)|=3$ and $\tau(H)=2$, which implies that the theorem holds for $H$, a contradiction. Thus, $\xi(H) = 45\tau(H) - 6n(H) - 13m(H) - \defic(H) \le 90 - 54 - 39 = -3$, contradicting the fact that $\xi(H) > 0$.~\smallqed



\ClaimX{H.29}
The following holds. \\
\hspace*{1cm}
{\rm (a)} $X=X_4$.
\\
\hspace*{1cm}
{\rm (b)} The edge $e$ contains no vertex of degree~$1$ in $H'$.
\\
\hspace*{1cm}
{\rm (c)} Every edge in $H' - e$ intersects $\partial(X)$ in at most two vertices.

\ProofClaimX{H.29}
By Claim~H.28, $H'$ is not linear, implying that the edge $e$ overlaps some other edge in $H'$. If the edge $e$ contains a degree-$1$ vertex in $H'$, then, by Claim~H.5(a), $8 \ge 8|X_4| + 5|X_{14}| + 4|X_{11}| + |X_{21}|  \ge 13|X| - 6|X'| - 3 = 10$, a contradiction.  Therefore, the edge $e$ does not contain a degree-$1$ vertex in $H'$. Hence, by Claim~H.5(b), $8|X_4| + 5|X_{14}| + 4|X_{11}| + |X_{21}|  \ge  13|X| - 6|X'| - 7 = 6$, which implies that $X=X_4$. This completes the proof of Part~(a) and Part~(b). To prove Part~(c), suppose to the contrary that some edge, say $f$, in $H' - e$ intersects $\partial(X)$ in at least three vertices. In this case, the edge $f$ must intersect $e_1$ or $e_2$ in at least two vertices, a contradiction to $H$ being linear.~\smallqed

\medskip


\ClaimX{H.30}
There exists two distinct edges in $H' - e$ that both intersect $\partial(X)$ in exactly two vertices.



\ProofClaimX{H.30}
By Claim~H.29, $X = X_4$. The linearity of $H$ implies that $|\partial(X)| \ge 5$. Suppose that at most one edge in $H'-e$ intersects $\partial(X)$ in more than one vertex. By Claim~H.29, such an edge intersects $\partial(X)$ in exactly two vertices, implying that the edge $e$ could have be chosen so that $H'$ is linear, contradicting Claim~H.28. Therefore there are at least two distinct edges in $H'-e$ that intersects $\partial(X)$ in more than one vertex. By Claim~H.29, they each intersect  $\partial(X)$ in exactly two vertices, as claimed.~\smallqed

By Claim~H.30, there exists two distinct edges in $H' - e$ that both intersect $\partial(X)$ in exactly two vertices. Let $f_1$ and $f_2$ be two such edges in $H' - e$. Thus, $f_1$ and $f_2$ are edges in $H - V(X)$ and $|V(f_1) \cap \partial(X)| = |V(f_2) \cap \partial(X)| = 2$.



\ClaimX{H.31}
The following holds. \\
\hspace*{1cm}
{\rm (a)} Every vertex in $\partial(X) \cap (V(f_1) \cup V(f_2))$ has degree~$1$ in $H'-e$.
\\
\hspace*{1cm}
{\rm (b)} Every vertex in $\partial(X)$ belongs to at most one of $f_1$ and $f_2$.


\ProofClaimX{H.31}
Suppose, to the contrary, that $ x \in \partial(X) \cap (V(f_1) \cup V(f_2))$ and $d_{H'-e}(x) \ne  1$. If $d_{H'-e}(x) = 0$, then we could have chosen the edge $e$ to contain the vertex $x$, implying that $e$ would contain a vertex of degree~$1$ in $H' - e$, contradicting Claim~H.29(b). Therefore, $d_{H'-e}(x) = 2$ noting that $\Delta(H) \le 3$ and $x \in \partial(X)$. Renaming $f_1$ and $f_2$, if necessary, we may assume that $x \in V(f_1)$. Let $y$ be the vertex in $V(f_1) \cap \partial(X)$ different from $x$.
Let $H^*$ be obtained from $H' - e - f_1$ by adding two new vertices, say $s_1$ and $s_2$, and the edge $s=\{x,y,s_1,s_2\}$.

We show that $\tau(H) \le \tau(H^*) + 1$. Let $T^*$ be a minimum transversal in $H^*$. As $T^*$ covers the edge $\{x,y,s_1,s_2\}$ and $s_1$ and $s_2$ have degree~$1$ in $H^*$, we can choose $T^*$ to contain  $x$ or $y$, implying that $T^*$ covers the edge $f_1$. Further since $\{x,y\} \subset \partial(X)$, the set $T^*$ covers at least one of $e_1$ and $e_2$ in $H$, say $e_1$. We can therefore add to $T^*$ a vertex from $X \cap V(e_2)$ in order to cover the edge in $X$ and the edge $e_2$, thereby getting a transversal for $H$. This proves that $\tau(H) \le |T^*| + 1 =  \tau(H^*) + 1$. We note that $H^*$ is linear, and so $H^* \in \cH_4$. Further, $n(H) = n(H^*) + 2$ and $m(H) = m(H^*) + 3$. By Claim~G, we have $\Phi(H^*)<0$. Thus,
\[
\begin{array}{lcl}
0 > \Phi(H^*) & = & \xi(H^*) - \xi(H)
\2 \\
& = & 45\tau(H^*) - 6n(H^*) - 13m(H^*) - \defic(H^*)   \\
& & \hspace*{1cm} - 45\tau(H) + 6n(H) + 13m(H) + \defic(H) \2 \\
& \ge & 45(\tau(H) - 1) - 6(n(H) - 2) - 13(m(H) - 3) \\
& & \hspace*{1cm} - \defic(H^*) - 45\tau(H) + 6n(H) + 13m(H)  \2 \\
& = & -45 + 12 + 39 - \defic(H^*) \\
& = & 6 - \defic(H^*),
\end{array}
\]
and so $\defic(H^*) \ge 7$. Let $Y \subseteq H^*$ be a special $H^*$-set such that $\defic_{H^*}(Y) = \defic(H^*) \ge 7$. If $|E^*_{H^*}(Y)| \ge |Y|$, then $\defic_{H^*}(Y) \le 10|Y| - 13|E^*_{H^*}(Y)| \le 10|Y| - 13|Y| < 0$, a contradiction. Hence, $|E^*_{H^*}(Y)| \le |Y| - 1$. If $s \notin E(Y)$, then
\[
|E^*_H(X \cup Y)| \le |E^*_{H^*}(Y)| + |\{e_1,e_2,f_1\}| \le (|Y|-1)+3 = |X \cup Y| + 1,
\]
contradicting the maximality of $|X|$. Therefore, $s \in E(Y)$. Thus, since the edge $s$ contains two degree-$1$ vertices, namely $s_1$ and $s_2$, and at least one vertex of degree~$2$, namely $x$, in $H^*$,  we note that $|Y| \ge 2$ and that $s$ is a $H_4$-component in $Y$. If $|E^*_{H^*}(Y)| = |Y| - 1$, then $\defic_{H^*}(Y) \le 10(|Y| - 1) + 8 - 13|E^*_{H^*}(Y)| = 10|Y| - 2 - 13(|Y|-1) = -3|Y| + 11 \le 5$, a contradiction. Hence, $|E^*_{H^*}(Y)| \le |Y| - 2$. Therefore,
\[
|E_H^*(X \cup (Y \setminus \{s\})| \le |E^*_{H^*}(Y)| + |\{e_1,e_2,f_1\}| \le (|Y|-2)+3 = |X \cup (Y \setminus \{s\})| + 1,
\]
contradicting the maximality of $|X|$. This completes the proof of Part~(a). Part~(b) follows directly from Part~(a).~\smallqed

\medskip
Recall that $|V(f_1) \cap \partial(X)| = |V(f_2) \cap \partial(X)| = 2$. By Claim~H.31, every vertex in $\partial(X)$ belongs to at most one of $f_1$ and $f_2$. We now choose the edge $e$ to contain the four vertices in $\partial(X) \cap (V(f_1) \cup V(f_2))$. %
We next define a new hypergraph $H_f$ as follows. Let $H''$ be constructed from $H$ by removing the four vertices in $X$ and the four vertices in $e$ and removing the edge in $E(X)$ and the four edges $e_1,e_2,f_1,f_2$. We now define the edge $f$ as follows. If $f_1$ and $f_2$ have no common neighbor in $V(H') \setminus \partial(X)$, then $|(V(f_1) \cup V(f_2)) \setminus V(e)| = 4$ and, in this case, we let $f$ contain these four vertices. If $f_1$ and $f_2$ have a common neighbor in $V(H') \setminus \partial(X)$, then $|(V(f_1) \cup V(f_2)) \setminus V(e)| = 3$ and, in this case, we let $f$ contain these three vertices as well as a vertex from $\partial(X) \setminus V(e)$. Let $H_f$ be obtained from $H''$ by adding to it the edge $f$; that is, $H_f = H'' \cup f$.


\medskip


\ClaimX{H.31}
$\tau(H) \le \tau(H_f) + 2$.


\ProofClaimX{H.31}
Let $T_f$ be a minimum transversal in $H_f$. In order to cover the edge $f$, there is a vertex $z$ in $T_f$ that belongs to $f$. Suppose that $z \in V(f_1) \cup V(f_2)$. Thus, at least one of the edges in $\{f_1,f_2\}$, say $f_2$, is covered by $T_f$. Let $w_1$ be any vertex in $V(e) \cap V(f_2)$. We note that $w_1$ covers the edge $f_2$ and at least one of the vertices in $\{e_1,e_2\}$, say $e_1$. Let $w_2$ be the vertex in $V(X) \cap V(e_2)$. We note that $w_2$ covers the edge in $X$ and the edge $e_2$.
Hence, $T_f \cup \{w_1,w_2\}$ is a transversal in $H$, which implies that $\tau(H) \le \tau(H_f) + 2$. Suppose that $z \notin V(f_1) \cup V(f_2)$. This implies that $|(V(f_1) \cup V(f_2)) \setminus V(e)| = 3$ and $z \in \partial(X) \setminus V(e)$. We note that $z$ covers at least one of the vertices in $\{e_1,e_2\}$, say $e_1$. Let $z_2$ be the vertex in $V(X) \cap V(e_2)$. We note that $w_2$ covers the edge in $X$ and the edge $e_2$. Let $z_1$ be the vertex common to $f_1$ and $f_2$. Hence, $T_f \cup \{z_1,z_2\}$ is a transversal in $H$, which implies that $\tau(H) \le \tau(H_f) + 2$.~\smallqed


\medskip


\ClaimX{H.32}
If $Y$ is a special $H_f$-set and $f \notin E(Y)$, then $|E_{H_f}^*(Y)| \ge |Y|+1$.

\ProofClaimX{H.32}
Suppose, to the contrary, that $|E_{H_f}^*(Y)| \le |Y|$. We now prove the following claims.


\ClaimX{H.32.1}
$f \in E_{H_f}^*(Y)$.

\ProofClaimX{H.32.1}
If $f \notin E_{H_f}^*(Y)$, then
$|E_H^*(X \cup Y)| \le |E_{H_f}^*(Y)| + |\{e_1,e_2\}| \le |Y| + 2 = |X \cup Y| +1$, contradicting the maximality of $|X|$.~\smallqed

\ClaimX{H.32.2}
$|E_{H_f}^*(Y)| = |Y|$.

\ProofClaimX{H.32.2}
If $|E_{H_f}^*(Y)| < |Y|$, then $|E_H^*(X \cup Y)| \le |E_{H_f}^*(Y)| - |\{f\}| + |\{f_1,f_2,e_1,e_2\}| \le (|Y|-1) -1 + 4 = |X \cup Y| +1$, contradicting the maximality of $|X|$.~\smallqed

\medskip
%
We construct a bipartite graph $G_Y$, with partite sets $Y \cup Y_{10}$ (that is, there are two copies of every element in $Y_{10}$) and $E_{H_f}^*(Y)$, where an edge joins $g \in E_{H_f}^*(Y)$ and $F \in Y \cup Y_{10}$ in $G_Y$ if and only if the edge $g$ intersects the subhypergraph $F$ of $Y$ in~$H_f$.

\ClaimX{H.32.3}
There is no matching in $G_Y$ saturating every vertex in $E_{H_f}^*(Y)$.

\ProofClaimX{H.32.3}
Suppose, to the contrary, that there is a matching, $M$, in $G_Y$ saturating every vertex in $E_{H_f}^*(Y)$. In this case, by Observation~\ref{property:special}(g) and~\ref{property:special}(h), using the matching $M$ we can find a $\tau(Y)$-transversal in $H_f$ covering all edges in $E(Y) \cup E_{H_f}^*(Y)$. Let $H^* = H_f - V(Y)$; that is, $H^*$ is obtained from $H_f$ by removing all vertices in $Y$ and all edges in $E(Y) \cup E_{H_f}^*(Y)$.

We show that $\defic(H^*) = 0$. Suppose, to the contrary, that $\defic(H^*) > 0$. Let $Y^*$ be a special $H^*$-set with $\defic_{H^*}(Y^*) > 0$. This implies that $|E_{H^*}^*(Y^*)| < |Y^*|$. By Claim~H.32.2, $|E_{H_f}^*(Y)| = |Y|$. Thus, $|E_{H_f}^*(Y \cup Y^*)| = |E_{H_f}^*(Y)| + |E_{H^*}^*(Y^*)| < |Y| + |Y^*|$. Hence we have shown the existence of a special $H_f$-set, $Y \cup Y^*$, such that $f \notin E(Y \cup Y^*)$ and $|E_{H_f}^*(Y \cup Y^*)| < |Y| + |Y^*|$, contradicting Claim~H.32.2. Hence, $\defic(H^*) = 0$.

As $f \in E_{H_f}^*(Y)$, we note that $f \notin E(H^*)$. Therefore, $H^*$ is linear. Recall that for any given special hypergraph, $F$, we have $\defic(F) = 45 \tau(F) - 6 |V(F)| - 13 |E(F)|$. Thus, $\defic(H_{10})=10$, $\defic(H_4)=8$, $\defic(H_{14})=5$, $\defic(H_{11})=4$, and $\defic(H_{21})=1$. In particular, $\defic(F) \le 10$ for any given special hypergraph, $F$. By Claim~H.31, $\tau(H) \le \tau(H_f) + 2$. Thus, since $\tau(H_f) \le \tau(H^*) + \tau(Y)$, we have $\tau(H) \le \tau(H^*) + \tau(Y) + 2$. We note that $n(H) = n(H^*) + n(Y) + 8$ and $m(H) = m(H^*) + m(Y) + |E_{H_f}^*(Y)| + 4$. Thus,

\[
\begin{array}{lcl}
\Phi(H^*) & = & \xi(H^*) - \xi(H)
\2 \\
& = & 45\tau(H^*) - 6n(H^*) - 13m(H^*) - \defic(H^*)   \\
& & \hspace*{1cm} - 45\tau(H) + 6n(H) + 13m(H) + \defic(H) \2 \\
& \ge & 45(\tau(H) - \tau(Y) - 2) - 6(n(H) - n(Y) - 8) \\
& & \hspace*{1cm} - 13(m(H) - m(Y) - |E_{H_f}^*(Y)| - 4)  \\
& & \hspace*{1.5cm} - 45\tau(H) + 6n(H) + 13m(H) \2 \\
& = & 10 + 13|E_{H_f}^*(Y)| - 45\tau(Y) + 6n(Y) + 13m(Y)  \1 \\
& = & \displaystyle{ 10 + 13|Y| - \sum_{F \in Y} \defic(F)  } \1 \\
& = & \displaystyle{ 10 + \sum_{F \in Y} (13 - \defic(F) )  } \1 \\
& > & 10,
\end{array}
\]
a contradiction to Claim~G. This completes the proof of Claim~H.32.3.~\smallqed

\medskip
By Claim~H.32.3, there is no matching in $G_Y$ saturating every vertex in $E_{H_f}^*(Y)$. By Hall's Theorem, there is a nonempty subset $S \subseteq E_{H_f}^*(Y)$ such that $|N_{G_Y}(S)|<|S|$.
We now consider the bipartite graph $G_Y'$, with partite sets $Y$ and $E_{H_f}^*(Y)$, where an edge joins $g \in E_{H_f}^*(Y)$ and $F \in Y$ in $G_Y$ if and only if the edge $g$ intersects the subhypergraph $F$ of $Y$ in~$H_f$. Thus, $G_Y'$ is obtained from $G_Y$ by deleting all duplicated vertices associated with copies of $H_{10}$ in $Y$. We note that $|N_{G_{Y'}}(S)| \le |N_{G_{Y}}(S)| < |S|$.
We now consider the special $H$-set, $Y' = Y \setminus N_{G_{Y'}}(S)$. Recall that by Claim~H.32.2, $|E_{H_f}^*(Y)| = |Y|$. Thus,
$|E_{H_f}^*(Y')| = |E_{H_f}^*(Y)| - |S| = |Y| - |S| = |Y'| + |N_{G_{Y'}}(S)| - |S| < |Y'|$,
contradicting Claim~H.32.2. This completes the proof of Claim~H.32.~\smallqed

\medskip

\ClaimX{H.33}
Let $Y$ be a special $H_f$-set. If $f \notin E(Y)$ and $f \in E_{H_f}^*(Y)$ and $|Y_{10}|>0$, then $|E_{H_f}^*(Y)| \ge |Y|+2$.

\ProofClaimX{H.33}
Suppose, to the contrary, that $|E_{H_f}^*(Y)| \le |Y|+1$. By Claim~H.32, $|E_{H_f}^*(Y)| \ge |Y| + 1$. Consequently,  $|E_{H_f}^*(Y)| = |Y| + 1$. Let $G_Y$ be the bipartite graph as defined in the proof of Claim~H.32.

\ClaimX{H.33.1}
There is no matching in $G_Y$ saturating every vertex in $E_{H_f}^*(Y)$.

\ProofClaimX{H.33.1}
Suppose, to the contrary, that there is a matching, $M$, in $G_Y$ saturating every vertex in $E_{H_f}^*(Y)$. We proceed now analogously as in the proof of Claim~H.32.3. If $Y^*$ is a special $H^*$-set with $\defic_{H^*}(Y^*) > 0$, then $|E_{H^*}^*(Y^*)| < |Y^*|$ and $|E_{H_f}^*(Y \cup Y^*)| = |E_{H_f}^*(Y)| + |E_{H^*}^*(Y^*)| \le (|Y|+1) + (|Y^*|-1) = |Y| + |Y^*|$, contradicting Claim~H.32. Hence, $\defic(H^*) = 0$. Proceeding now exactly as in the proof of Claim~H.32.3, we show that $\Phi(H^*) > 0$, contradicting Claim~G.~\smallqed

\medskip
By Claim~H.32.3, there is no matching in $G_Y$ saturating every vertex in $E_{H_f}^*(Y)$. By Hall's Theorem, there is a nonempty subset $S \subseteq E_{H_f}^*(Y)$ such that $|N_{G_Y}(S)|<|S|$. If $S = E_{H_f}^*(Y)$, then $N_{G_Y}(S) = Y \cup Y_{10}$. However, by assumption  $|Y_{10}|>0$, and so $|N_{G_Y}(S)| \ge |Y| + 1 = |E_{H_f}^*(Y)| = |S|$, a contradiction. Hence, $S$ is a proper subset of $E_{H_f}^*(Y)$. Let $G_Y'$ be the bipartite graph defined in the proof of Claim~H.32. We note that $|N_{G_{Y'}}(S)| \le |N_{G_{Y}}(S)| < |S|$.
We now consider the special $H$-set, $Y' = Y \setminus N_{G_{Y'}}(S)$. If $Y' = \emptyset$, then, $|N_{G_{Y'}}(S)| = |Y| = |E_{H_f}^*(Y)| - 1 \ge |S|$, a contradiction. Hence, $Y' \ne \emptyset$. Further,
$|E_{H_f}^*(Y')| = |E_{H_f}^*(Y)| - |S| = (|Y| + 1) - |S| = |Y'| + 1 + |N_{G_{Y'}}(S)| - |S| < |Y'| + 1$,
and so $|E_{H_f}^*(Y')| \le |Y'|$, contradicting Claim~H.32. This completes the proof of Claim~H.33.~\smallqed

\medskip


\ClaimX{H.34}
If $Y$ is a special $H_f$-set and $f \in E(Y)$, then $|E_{H_f}^*(Y)| \ge |Y|-1$.

\ProofClaimX{H.34}
Suppose, to the contrary, that $|E_{H_f}^*(Y)| < |Y|-1$. This implies that $|Y| \ge 2$.  Let $F \in Y$ be the special subhypergraph containing $f$ and let $Y' = Y \setminus  F$. In this case, $|E_{H_f}^*(Y')| \le |E_{H_f}^*(Y)|  < |Y|-1 = |Y'|$. Therefore, $Y'$ is a special $H_f$-set and $f \notin Y'$ satisfying   $|E_{H_f}^*(Y')| \le  |Y'|$, contradicting Claim~H.32.~\smallqed

\medskip


\ClaimX{H.35}
$H_f$ is not linear.

\ProofClaimX{H.35}
Suppose, to the contrary, that $H_f$ is linear. Thus, by Claim~G and Claim~H.31, we have
\[
\begin{array}{lcl}
0 > \Phi(H_f) & = & \xi(H_f) - \xi(H)
\2 \\
& = & 45\tau(H_f) - 6n(H_f) - 13m(H_f) - \defic(H_f)   \\
& & \hspace*{1cm} - 45\tau(H) + 6n(H) + 13m(H) + \defic(H) \2 \\
& \ge & 45(\tau(H) - 2) - 6(n(H) - 8) - 13(m(H) - 4)  \\
& & \hspace*{1.5cm} - 45\tau(H) + 6n(H) + 13m(H) \2 \\
& = & -45 \cdot 2 + 6 \cdot 8 + 13 \cdot 4 - \defic(H_f)   \1 \\
& = & 10 - \defic(H_f),
\end{array}
\]
and so, $\defic(H_f) \ge 11$. Let $Y$ be a special $H_f$-set such that $\defic(H_f) = \defic_{H_f}(Y)$.
As $\defic_{H_f}(Y) \ge 11$, we note that $|E_{H_f}^*(Y)| \le |Y| -2$. This is a contradiction to Claim~H.32, Claim~H.33 and Claim~H.34, no matter whether $f \notin E(Y)$ or $f \in E(Y)$.~\smallqed

\medskip
We now return to the proof of Claim~H one final time.  By Claim~H.35, there exists an edge $g$ in $H_f$ that overlaps $f$. Let $\{u,v\} \subseteq V(f) \cap V(g)$. Let $H_f^*$ be the hypergraph obtained from $H_f$ by removing the edges $f$ and $g$, adding two new vertices $z_1$ and $z_2$, and adding a new edge $e_{fg} = \{u,v,z_1,z_2\}$. Since $H_f - f$ is linear, we note that $H_f^*$ is linear. Let $T_f^*$ be a minimum transversal of $H_f^*$. If $z_1$ or $z_2$ belong to $T_f^*$, we can replace it with $u$ or $v$, implying that $T_f^*$ is a transversal of $H_f$. Thus, by Claim~H.31, $\tau(T_f^*) \ge \tau(H_f) \ge \tau(H) - 2$. Thus, by Claim~G, we have

\[
\begin{array}{lcl}
0 > \Phi(H_f^*) & = & \xi(H_f^*) - \xi(H)
\2 \\
& = & 45\tau(H_f^*) - 6n(H_f^*) - 13m(H_f^*) - \defic(H_f^*)   \\
& & \hspace*{1cm} - 45\tau(H) + 6n(H) + 13m(H) + \defic(H) \2 \\
& \ge & 45(\tau(H) - 2) - 6(n(H) - 6) - 13(m(H) - 5)  \\
& & \hspace*{1.5cm} - 45\tau(H) + 6n(H) + 13m(H) \2 \\
& = & -45 \cdot 2 + 6 \cdot 6 + 13 \cdot 5 - \defic(H_f^*)   \1 \\
& = & 11 - \defic(H_f^*),
\end{array}
\]

\noindent
and so, $\defic(H_f^*) \ge 12$.  Let $Y$ be a special $H_f^*$-set such that $\defic(H_f^*) = \defic_{H_f^*}(Y)$. As $\defic_{H_f^*}(Y) \ge 12$ we note that $|E_{H_f^*}^*(Y)| \le |Y| -2$. If $e_{fg} \notin E(Y)$, then $|E_{H_f}^*(Y)| \le |Y| -2 + |\{f,g\}| = |Y|$, contradicting Claim~H.32. Therefore, $e_{fg} \in E(Y)$. Let $F_{fg} \in Y$ be the special subhypergraph containing the edge $e_{fg}$. Since $e_{fg}$ contains at least two degree-$1$ vertices in $H_f^*$, we note that $F_{fg} \in Y_4$. Let $Y_{fg} = Y \setminus  \{F_{fg}\}$, and note that $f \notin E(Y_{fg})$. Further,
\[
|E_{H_f}^*(Y_{fg})| \le |E_{H_f^*}^*(Y)|  + |\{f,g\}| = |E_{H_f^*}^*(Y)|  + 2.
\]

As observed earlier, $|E_{H_f^*}^*(Y)| \le |Y| -2$. Suppose that $|E_{H_f^*}^*(Y)|  = |Y| - 2$. In this case, since $\defic_{H_f^*}(Y) \ge 12$, we note that $|Y_{10}| > 0$. Since $F_{fg} \in Y_4$, this implies that $|(Y_{fg})_{10}|>0$. If $f \in E_{H_f}^*(Y_{fg})$, then $|E_{H_f}^*(Y_{fg})| \le |E_{H_f^*}^*(Y)|  + 2 = |Y| = |Y_{fg}|+1$, contradicting Claim~H.33. If $f \notin E_{H_f}^*(Y_{fg})$, then $|E_{H_f}^*(Y_{fg})| \le |E_{H_f^*}^*(Y)|  + |\{g\}| = |E_{H_f^*}^*(Y)|  + 1 = |Y| - 1 = |Y_{fg}|$, contradicting Claim~H.32. Therefore, $|E_{H_f^*}^*(Y)|  \le |Y| - 3$, implying that $|E_{H_f}^*(Y_{fg})| \le |E_{H_f^*}^*(Y)|  + 2 \le |Y| - 1 = |Y_{fg}|$, once again contradicting Claim~H.32. This completes the proof of Claim~H.~\smallqed

\medskip


\ClaimX{I}
If $Y$ be a special $H$-set, then the following holds. \\
\hspace*{1cm}
{\rm (a)} $|E^*(Y)| \ge |Y|+2$. \\
\hspace*{1cm}
{\rm (b)} If $|Y_{10}| \ge 2$, then $|E^*(Y)| \ge |Y|+3$.

\ProofClaimX{I}
Part~(a) follows immediately from our choice of the $H$-special set $X$ and by Claim~H. To prove Part~(b), suppose, to the contrary, that $|Y_{10}| \ge 2$ and $|E^*(Y)| \le |Y|+2$. By Part~(a), $|E_H^*(Y)| \ge |Y|+2$. Consequently, $|E^*(Y)| = |Y|+2$. We construct a bipartite graph $G_Y'$, with partite sets $Y \cup Y_{10}$ (that is, there are two copies of every element in $Y_{10}$) and $E^*(Y)$, where an edge joins $e \in E^*(Y)$ and $F \in Y \cup Y_{10}$ in $G_Y'$ if and only if the edge $e$ intersects the subhypergraph $F$ of $Y$ in~$H$.

Suppose that there is a matching, $M$, in $G_Y'$ saturating every vertex in $E^*(Y)$. By Observation~\ref{property:special}(g) and~\ref{property:special}(h), using the matching $M$ we can find a $\tau(Y)$-transversal in $H$ covering all edges in $E(Y) \cup E^*(Y)$, contradicting Claim~B. Hence, there is no matching in $G_Y'$ saturating every vertex in $E^*(Y)$.

By Hall's Theorem, there is a nonempty subset $S \subseteq E^*(Y)$ such that $|N_{G_Y'}(S)|<|S|$. If $S = E^*(Y)$, then $N_{G_Y'}(S) = Y \cup Y_{10}$. However, by assumption  $|Y_{10}| \ge 2$, and so $|N_{G_Y'}(S)| = |Y| + |Y_{10}| \ge |Y| + 2 = |E^*(Y)| = |S|$, a contradiction. Hence, $S$ is a proper subset of $E^*(Y)$. We now consider the special $H$-set, $Y'$, obtained from $Y$ by deleting every element in $Y$ that is a neighbor of $S$ in $H$.  We note that at most $|N_{G_Y'}(S)|$ such elements have been deleted from $Y$, possibly fewer since some elements $F \in Y_{10}$ may get deleted twice, and so $|Y'| \ge |Y| - |N_{G_Y'}(S)|$. Therefore, $|E^*(Y')| \le |E^*(Y)| - |S| < (|Y|+2) - |N_{G_Y'}(S)| \le |Y'| + 2$. Thus, $|E^*(Y')| - |Y'| < 2$.
We therefore have a contradiction to the   choice of the special $H$-set $X$ which was chosen so that $|E^*(X)| - |X|$ is minimum.~\smallqed

\medskip
In particular, Claim~I implies that if $|X_{10}| \ge 2$, then $|E^*(X)| \ge |X|+3$.

\medskip

\ClaimX{J}
There is no $H_{10}$-subhypergraph in $H$.

\ProofClaimX{J}
Suppose, to the contrary, that $F$ is a $H_{10}$-subhypergraph in $H$. By Claim~C, $F$ is not a component of $H$, and so there exists an edge $e \in E^*(F)$. We choose such an edge $e$ so that $|V(e) \cap V(F)|$ is maximum possible. Since $H$ is linear, we note that $|V(e) \cap V(F)| \le 2$. Let $x \in V(e) \cap V(F)$ be arbitrary. Let $E(F) = \{e_1,e_2,e_3,e_4,e_5\}$, where $e_1$ and $e_2$ are the two edge of $F$ that contain~$x$. We note that $d_H(x) = 3$, and $e,e_1,e_2$ are the three edges of $H$ that contains~$x$. Let $H' = H - x$. Let $E' = \{e_3,e_4,e_5\}$ and note that $E' \subseteq E(H')$ and every pair of edges in $E'$ intersect in $H'$. Every $\tau(H')$-transversal can be extended to a transversal of $H$ by adding to it the vertex~$x$, and so $\tau(H) \le \tau(H') + 1$. By Claim~G,
\[
\begin{array}{lcl}
0 > \Phi(H') & = & \xi(H') - \xi(H)
\2 \\
& = & 45\tau(H') - 6n(H') - 13m(H') - \defic(H')   \\
& & \hspace*{1cm} - 45\tau(H) + 6n(H) + 13m(H) + \defic(H) \2 \\
& \ge & 45(\tau(H) - 1) - 6(n(H) - 1) - 13(m(H) - 3)  \\
& & \hspace*{1.5cm} - 45\tau(H) + 6n(H) + 13m(H) \2 \\ & = & - 45 + 6 \cdot 1  + 13 \cdot 3 - \defic(H')   \1 \\
& = & - \defic(H'),
\end{array}
\]
and so, $\defic(H') > 0$. Let $Y$ be a special $H'$-set such that $\defic(H') = \defic_{H'}(Y)$. By Claim~I(a), $|E_H^*(Y)| \ge |Y|+2$.

\medskip
\ClaimX{J.1}
$|Y|=2$, $|E_{H'}^*(Y)|=1$, and $0 < \defic_{H'}(Y) \le 5$.

\ProofClaimX{J.1}
If $|E_{H'}^*(Y)| \le |Y| - 2$, then $|E_{H}^*(Y)| \le |E_{H'}^*(Y)| + |\{e,e_1,e_2\}| \le |Y|+1$, a contradiction. Therefore, $|E_{H'}^*(Y)| \ge |Y| - 1$. As $\defic_{H'}(Y) \ge 1$, we note that $|E_{H'}^*(Y)| \le |Y| -1$. Consequently, $|E_{H'}^*(Y)| = |Y| - 1$.

If $|Y| \ge 5$, then $\defic(H') \le 10|Y| - 13|E_{H'}^*(Y)| = 10|Y| - 13(|Y| - 1) = -3|Y| + 13 < 0$, a contradiction. Hence, $|Y| \le 4$. Suppose that $|Y| \ge 3$. If $|Y_{10}| \le 1$, then $\defic(H') \le 10 + 8(|Y|-1) - 13|E_{H'}^*(Y)| = 8|Y| + 2 - 13(|Y| - 1) = -5|Y| + 15 \le 0$, a contradiction. Hence, $|Y_{10}| \ge 2$. However, $|E_{H}^*(Y)| \le |E_{H'}^*(Y)| + |\{e,e_1,e_2\}| = |Y| + 2$, contradicting Claim~I(b). Hence, $|Y| \le 2$.

Suppose that $|Y|=1$, and let $Y = \{F_1\}$. By Claim~I(a), $|Y| + 2 \le |E_{H}^*(Y)| \le |E_{H'}^*(Y)| + |\{e,e_1,e_2\}| = |Y| + 2$, implying that $|E_H^*(Y)| = |Y|+2$ and $E^*(F_1) = \{e,e_1,e_2\}$.
Therefore, $E' \subseteq E(F_1)$, which implies that $V(F) \setminus \{x\} \subseteq V(F_1)$. Therefore, the edges $e_1$ and $e_2$ intersect $F_1$ in three vertices, while the edge $e$ intersects $F_1$ in at least one vertex. If $e$ intersects $F_1$ in at least two vertices, then by Observation~\ref{property:special}(n) we can cover $E^*(F_1)$ with a $\tau(F_1)$-transversal in $H$, a contradiction to Claim~B. Therefore, we may assume that the edge $e$ intersects $F_1$ in only one vertex.

By Observation~\ref{property:special}($\ell{}$), there exists a $\tau(F_1)$-transversal, $T_1$, of $F_1$ covering both $e_1$ and $e_2$. Let $H_1=H - T_1$.  We will now show that  $\defic(H_1) \le 8$. Suppose, to the contrary, that $\defic(H_1) \ge 9$ and let $Y_1$ be a special $H_1$-set such that $\defic(H_1) = \defic_{H_1}(Y_1)$.
If the edge $e$ does not belong to any special subhypergraph in $Y_1$, then $\defic_{H'}(Y_1 \cup F_1) >  \defic_{H'}(F_1) = \defic(H')$, a contradiction.
Therefore, $e \in E(F_e)$ for some special subhypergraph in $Y_1$.  As the edge $e$ has at least two vertices of degree~$1$ in $H_1$ (namely, the vertex $x$ and the vertex in $V(F_1) \cup V(e)$), the subhypergraph $F_e$ is an $H_4$-component. However since $\defic_{H_1}(Y_1) \ge 9$, the special set $\{F_1\} \cup (Y_1 \setminus \{F_e\})$ has higher deficiency than $F_1$ in $H'$, a contradiction. Therefore, $\defic(H_1) \le 8$, as desired.

As there is an edge, $e_1$ (and $e_2$), that intersects $F_1$ in three vertices we note that the deficiency of $F_1$ is at most $5$. Note that we have removed the edges $e_1$ and $e_2$ as well as all edges in $F_1$ to get from $H$ to $H_1$ as well as all vertices in $F_1$ except
the vertex in $V(e) \cap V(F_1)$. Therefore, $n(H_1) = n(H) - (n(F_1)-1)$, $m(H_1) = m(H) - m(F_1) - |\{e_1,e_2\}|$, and $\tau(H) \le \tau(H_1) + \tau(F_1)$.
As observed earlier, $\defic(F_1) \le 5$. Thus, $45\tau(F_1) \le 6n(F_1)+13m(F_1) +5$, and therefore
\[
\begin{array}{lcl}
0 > \Phi(H_1) & = & \xi(H_1) - \xi(H) \2 \\
& = & 45\tau(H_1) - 6n(H_1) - 13m(H_1) - \defic(H_1)   \\
& & \hspace*{1cm} - 45\tau(H) + 6n(H) + 13m(H) + \defic(H) \2 \\
& \ge & 45(\tau(H) - \tau(F_1)) - 6(n(H) - n(F_1)+1) - 13(m(H) - m(F_1)-2) - 8 \\
& & \hspace*{1.5cm} - 45\tau(H) + 6n(H) + 13m(H) \2 \\
& = & - 45\tau(F_1) + 6 \cdot (n(F_1)-1)  + 13 \cdot (m(F_1)+2) - 8   \1 \\
& \ge  & - 5 - 6 \cdot 1  + 13 \cdot 2 - 8  \1 \\
& > & 0,
\end{array}
\]
a contradiction. Therefore, $|Y|=2$ and $|E_{H'}^*(Y)| = |Y| - 1 = 1$. As observed earlier, if $|Y_{10}| \ge 2$, then $|E_{H}^*(Y)| \le |Y| + 2$, contradicting Claim~I(b). Hence, $|Y_{10}| \le 1$, implying that $\defic(H') = \defic_{H'}(Y) \le 10 + 8 - 13|E_{H'}^*(Y)| = 18 - 13 = 5$. This completes the proof of Claim~J.1.~\smallqed

\medskip
By Claim~J.1, we have $|Y|=2$, $|E_{H'}^*(Y)|=1$, and $0 < \defic_{H'}(Y) \le 5$. Let $Y=\{F_1,F_2\}$. Since $F_1$ and $F_2$ are vertex disjoint and every pair of edges in $E'$ intersect in $H'$, we note that $E(F_1) \cap E' = \emptyset$ or $E(F_2) \cap E' = \emptyset$. Renaming $F_1$ and $F_2$ if necessary, we may assume that $E(F_2) \cap E' = \emptyset$. If neither $e_1$ nor $e_2$ intersects $F_2$ in $H$, then $E_H^*(F_2) \subseteq \{e\} \cup E_{H'}^*(Y)$, implying that $|E_H^*(F_2)| \le 1 + |E_{H'}^*(Y)| = 2$. However, by Claim~I(a), $|E_H^*(F_2)| \ge |F_2| + 2 = 3$, a contradiction. Therefore, $e_1$ or $e_2$ intersects $F_2$ in $H$. Renaming $e_1$ and $e_2$ if necessary, we may assume that $e_1$ intersects $F_2$ in $H$.


\medskip

\ClaimX{J.2}
$F_2 \notin Y_{10}$.

\ProofClaimX{J.2}
Suppose, to the contrary, that $F_2 \in Y_{10}$. Let $z \in V(e_1) \cap V(F_2)$ be arbitrary. The vertex $z$ is incident with two edges from $F$ and two edges from $F_2$ in $H$ and these four edges are distinct, contradicting the fact that the maximum degree in $H$ is three. Therefore, $F_2 \notin Y_{10}$.~\smallqed

\medskip
\ClaimX{J.3}
$F_2 \notin Y_4$.

\ProofClaimX{J.3}
Suppose, to the contrary, that $F_2 \in Y_4$. Since $H$ is linear and $F_2 \in Y_4$, every edge in $E_H^*(F_2)$ intersects $F_2$ in exactly one vertex. Let $V(e_1) \cap V(F_2) = \{z_1\}$. The vertex $z_1$ is incident with two edges from $F$ and the one edge from $F_2$ in $H$, and so $d_H(z)=3$.

Suppose that both $e$ and $e_2$ intersects $F_2$ in $H$. Let $V(e_2) \cap V(F_2) = \{z_2\}$. Since $x$ is the only vertex common to both $e_1$ and $e_2$, we note that $z_1 \ne z_2$. Let $f$ be the edge of $F_2$. By the linearity of $H$ and since $F$ is a $H_{10}$-subhypergraph in $H$, we note that $V(f) \cap V(F) = \{z_1,z_2\}$. Thus, the edge $f \in E^*_H(F)$ and $|V(f) \cap V(F)| = 2$, implying by our choice of the edge $e$ that $|V(e) \cap V(F)| = 2$. This in turn implies by the linearity of $H$ and the structure of $F$, that the edge $e$ does not intersect $F_2$, a contradiction.

As observed above, at least one of $e$ and $e_2$ does not intersect $F_2$ in $H$. Thus, $|E_H^*(F_2)| \le 2 + |E_{H'}^*(Y)| = 3$, implying that there exists a vertex $q$ in $F_2$ that has degree~$1$ in $H$. If we had chosen to use the vertex $z_1$ instead of $x$ when constructing $H'$, the vertex $q$ would become isolated and therefore deleted from $H'$, implying that $n(H') \le n(H) - 2$. Thus, $0 > \Phi(H') \ge - 45 + 6 \cdot 2  + 13 \cdot 3 - \defic(H') = 6 - \defic(H')$. Hence, $\defic(H') > 6$, which is a contradiction to Claim~J.1. Therefore, $F_2 \notin F_4$.~\smallqed

\medskip
\ClaimX{J.4}
$F_1 \notin Y_{10}$.

\ProofClaimX{J.4}
Suppose, to the contrary, that $F_1 \in Y_{10}$.

Suppose that $|E(F_1) \cap E'| = 0$. By Claim~I(a), at least two of the edges in $\{e,e_1,e_2\}$ intersect $F_1$ in $H$, implying that there is a vertex $q \in V(F_2) \cap V(F)$. The vertex $q$ is incident with two edges from $F$ and two edges from $F_2$ in $H$ and these four edges are distinct, contradicting the fact that the maximum degree in $H$ is three. Therefore, $|E(F_1) \cap E'| \ge 1$. If $|E(F_1) \cap E'| = 1$, then $|E_{H'}^*(Y)| \ge 2$, as all edges in $E'$ intersect each other, contradicting Claim~J.1. Therefore, $|E(F_1) \cap E'| \ge 2$. If $|E(F_1) \cap E'| = 3$ (and so, $E' \subset E(F_1)$), then, since $F_1$ and $F_2$ are vertex disjoint, the edges $e_1$ and $e_2$ do not intersect $F_2$, implying that $E_H^*(F_2) \subseteq \{e\} \cup E_{H'}^*(Y)$ and hence that $|E_H^*(F_2)| \le 1 + |E_{H'}^*(Y)| = 2$, contradicting Claim~I(a). Therefore, $|E(F_1) \cap E'| = 2$.

Renaming the edges in $E'$ if necessary, we may assume that $e_3,e_4 \in E(F_1)$, which implies that $E_{H'}^*(Y) = \{e_5\}$. As $e_5$ intersects both $e_3$ and $e_4$, the edge $e_5$ intersects $F_1$ in at least two vertices; that is, $e_5 \in E^*_H(F_1)$ and $|V(e_5) \cap V(F_1)| \ge 2$ where recall that $F_1$ is a $H_{10}$-subhypergraph in $H$. By the maximality of $|V(e) \cap V(F)|$, we may therefore assume that $|V(e) \cap V(F)| \ge 2$. Let $\{x,y\} \subseteq V(e) \cap V(F)$, where recall that $H' = H - x$. By the linearity of $H$ and since $F$ is a $H_{10}$-subhypergraph in $H$, we note that $V(e) \cap V(F) = \{x,y\}$ and $y \notin V(e_1) \cup V(e_2)$, which implies that $y$ belongs to two of the edges in $\{e_3,e_4,e_5\}$. Therefore, $y \in V(F_1)$ and all four edges $e,e_1,e_2,e_5$ intersect $F_1$.

Recall that $E_{H'}^*(Y) = \{e_5\}$. If the edge $e_5$ does not intersect $F_2$, then neither do the edges $e_1$ and $e_2$, implying that $|E_H^*(F_2)| \le |\{e\}| = 1$, contradicting Claim~I(a). Therefore, $V(e_5) \cap V(F_2) \ne \emptyset$. This implies that one of the edges in $\{e_1,e_2,e_3,e_4\}$ intersects $F_2$ as every vertex of $e_5$ intersects one of the edges
in $\{e_1,e_2,e_3,e_4\}$. Since $F_1$ and $F_2$ are vertex disjoint and $e_3,e_4 \in E(F_1)$, the edges $e_3$ and $e_4$ do not intersect $F_2$. Thus, the common vertex of $e_1$ and $e_5$ belongs to $F_2$ or the common vertex of $e_2$ and $e_5$ belongs to $F_2$. Renaming $e_1$ and $e_2$ if necessary, we may assume that the common vertex of $e_1$ and $e_5$ belongs to $F_2$, and therefore we can cover $e_1$ and $e_5$ with a $\tau(F_2)$-transversal. As $F_1$ is a $H_{10}$ subhypergraph in $H$ we can cover $e$ and $e_2$ with a $\tau(F_1)$-transversal by Observation~\ref{property:special}(h). Therefore, we can cover $E^*(Y)$ in $H$ with a
$\tau(Y)$-transversal in $H$, contradicting Claim~B. This completes the proof of Claim~J.4.~\smallqed

\medskip
We now return to the proof of Claim~J. Recall that $Y=\{F_1,F_2\}$ and that by Claim~J.1, $|E_{H'}^*(Y)|=1$. By Claims~J.2 and J.3, $F_2 \notin Y_4 \cup Y_{10}$. By Claim~J.4, $F_1 \notin F_{10}$. Thus, $\defic(H') = \defic_{H'}(Y) \le 8 + 5 - 13 = 0$, a contradiction to $\defic(Y) > 0$. This completes the proof of Claim~J.~\smallqed

\medskip
Recall that the boundary of a set $Z$ of vertices in a hypergraph $H$ is the set $N_H(Z) \setminus Z$, denoted $\partial_H(Z)$ or simply $\partial(Z)$ if $H$ is clear from context.

\medskip
\ClaimX{K}
Let $Z \subseteq V(H)$ be an arbitrary nonempty set of vertices and let $H' = H - Z$. Then either $|E_{H'}^*(Y)| \ge |Y|$ for all special $H'$-sets $Y$ in $H'$ (and therefore $\defic(H')=0$) or there exists a transversal $T'$ in $H'$, such that
\[
45 |T'| \le 6n(H') + 13m(H') + \defic(H') \hspace*{0.5cm} \mbox{and}  \hspace*{0.5cm} T' \cap \partial(Z) \ne \emptyset.
\]

\ProofClaimX{K}
Suppose that $|E_{H'}^*(Y)| < |Y|$ for some special $H'$-set, $Y$, in $H'$, as otherwise the claim holds. Among all such special $H'$-sets $Y$, let $Y$ be chosen so that

\hspace*{1cm} (1) $|Y| - |E_{H'}^*(Y)|$ is maximum.\\
\indent \hspace*{1cm} (2) Subject to (1), $|Y|$ is minimum.

If $\partial(Z) \cap V(Y) = \emptyset$, then  $|E_{H}^*(Y)| < |Y|$, a contradiction to Claim~I(a).
Therefore, $\partial(Z) \cap V(Y) \ne \emptyset$ and let $q$ be any vertex in $\partial(Z) \cap V(Y)$.
 We now consider the bipartite graph, $G_Y$, with partite sets $Y$ and $E_{H'}^*(Y) \cup \{q\}$, where an edge joins $e \in E_{H'}^*(Y)$ and $F \in Y$ in $G_Y$
if and only if the edge $e$ intersects the subhypergraph $F$ of $Y$ in $H'$ and $q$ is joined to $F \in Y$ in $G_Y$ if and only if $q \in V(F)$.

\ClaimX{K.1}
If there is a matching in $G_Y$ saturating every vertex in $E_{H'}^*(Y) \cup \{q\}$, then Claim~K holds.

\ProofClaimX{K.1}
Suppose that there exists a matching in $G_Y$ that saturates every vertex in $E_{H'}^*(Y) \cup \{q\}$.
This implies that there exists a $\tau(Y)$-transversal, $T_Y$, in $H'$ covering $E_{H'}^*(Y)$ and with $q \in T_Y$.
Let $H_Y$ be obtained from $H'$ by deleting the vertices $V(Y)$ and edges $E(Y) \cup E_{H'}^*(Y)$; that is, $H_Y = H' - T_Y$.
Suppose that $\defic(H_Y) > 0$, and let $Y'$ be a special $H_Y$-set in $H_Y$ such that $\defic_{H_Y}(Y') = \defic(H_Y) > 0$.
Since $|E_{H_Y}^*(Y')| \le |Y'|-1$, we note that $|E_{H'}^*(Y \cup Y')| = |E_{H'}^*(Y)| + |E_{H_Y}^*(Y')| < |Y| + |Y'|$ and
\[
|Y \cup Y'| - |E_{H'}^*(Y \cup Y')| = |Y| + |Y'| - (|E_{H'}^*(Y)|+|E_{H_Y}^*(Y')|) \ge |Y| - |E_{H'}^*(Y)| + 1,
\]
contradicting the choice of the special $H'$-set $Y$. Therefore, $\defic(H_Y)=0$.
Since $H$ is a counterexample with minimum value of $n(H)+m(H)$, and since $n(H_Y) + m(H_Y) < n(H) + m(H)$,  we note that $45 \tau(H_Y) \le 6n(H_Y) + 13m(H_Y) + \defic(H_Y) = 6n(H_Y) + 13m(H_Y)$.
As $45|T_Y| = 6n(Y) + 13m(Y) + 13|E_{H'}^*(Y)| + \defic_{H'}(Y)$ we get the following,
\[
\begin{array}{rcl}
45|T'| & \le & 45|T_Y| + 45\tau(H_Y) \\
 & \le & 6n(Y) + 13m(Y) + 13|E_{H'}^*(Y)| + \defic_{H'}(Y) + 6n(H_Y) + 13m(H_Y) \\
 & \le & 6n(H') + 13m(H') +  \defic_{H'}(Y).
\end{array}
\]
As $q$ belonged to $T_Y$ this proves Claim~K in this case.~\smallqed

\medskip

By Claim~K.1, we may consider the case when there is no perfect matching in $G_Y$ saturating every vertex in $E_{H'}^*(Y) \cup \{q\}$.
By Hall's Theorem, there is therefore a nonempty subset $S \subseteq E_{H'}^*(Y) \cup \{q\}$ such that $|N_{G_Y}(S)| < |S|$.
We now consider the special $H$-set $Y' = Y \setminus N_{G_Y}(S)$.
Since $|Y'| = |Y| - |N_{G_Y}(S)|$ and $|E_{H'}^*(Y')| \le |(E_{H'}^*(Y) \cup \{q\}) \setminus S| = |E_{H'}^*(Y)| + 1 - |S| \le |E_{H'}^*(Y)| - |N_{G_Y}(S)|$,
we note that $|Y'| - |E_{H'}^*(Y')| \ge |Y| - |E_{H'}^*(Y)|$ and $|Y'|<|Y|$, contradicting our choice of $Y$.
This completes the proof of Claim~K.~\smallqed

\medskip


Now let $f$ be the function defined in Table~1.
\begin{center}
\begin{tabular}{|c||c|c|c|c|c|} \hline
$i$ & $1$ & $2$ & $3$ & $4$ & $\ge 5$ \\ \hline
$f(i)$ & $39$ & $33$ & $27$ & $23$ & $22$ \\ \hline
\end{tabular}
\end{center}
\begin{center}
\vskip -0.25cm
\textbf{Table~1.} The function $f$.
\end{center}


\medskip
\ClaimX{L}
Let $Z \subseteq V(H)$ be an arbitrary nonempty set of vertices that intersects at least two edges of $H$, and let $H' = H - Z$. If $\defic(H') \le 21$ and $|\partial(Z)| \ge 1$, then there exists a transversal, $T'$, in $H'$, such that $T' \cap \partial(Z) \ne \emptyset$ and the following holds. \\[-25pt]
\begin{enumerate}
\item $45 |T'| \le 6n(H') + 13m(H') + f(|\partial(Z)|)$.
\item If $|\partial(Z)| \ge 5$ and $H'$ does not contain two intersecting edges $e$ and $f$, such that \vskip 0.1cm \hspace*{0.75cm} (i) $\partial(Z) \subseteq (V(e) \cup V(f)) \setminus(V(e) \cap V(f))$, \\ \hspace*{0.75cm} (ii) $e$ contains three degree-$1$ vertices, and \\ \hspace*{0.75cm} (iii) $|\partial(Z) \cap V(e)|,|\partial(Z) \cap V(f)| \ge 2$,
\vskip 0.1cm
then $45 |T'| \le 6n(H') + 13m(H') + f(|\partial(Z)|)-1 = 6n(H') + 13m(H') + 21$.
\end{enumerate}

\ProofClaimX{L}
Suppose, to the contrary, that there exists a set $Z \subseteq  V(H)$, such that $H' = H - Z$ with $\defic(H') \le 21$ and $|\partial(Z)| \ge 1$, but there exists no transversal $T'$ in $H'$ satisfying (a) or (b) in the statement of the claim. Among all such sets, let $Z$ be chosen so that $|Z|$ is as large as possible. Let $T'$ be a smallest possible transversal in $H'$ containing a vertex from $\partial(Z)$.

\ClaimX{L.1}
For any special $H'$-set, $Y$, in $H'$ we have $|E_{H'}^*(Y)| \ge |Y|$. Furthermore $\defic(H')=0$.

\ProofClaimX{L.1}
If $|E_{H'}^*(Y)| < |Y|$ for some special $H'$-set, $Y$, in $H'$, then by Claim~K there exists a transversal $T'$ in $H'$, such that  $45 |T'| \le 6n(H') + 13m(H') + \defic(H')$ and $T' \cap \partial(Z) \ne \emptyset$. Since $\defic(H') \le 21$, this implies that $45 |T'| \le 6n(H') + 13m(H') + 21 < 6n(H') + 13m(H') + f(|\partial(Z)|)$, contradicting our definition of $Z$. Therefore  $|E_{H'}^*(Y)| \ge |Y|$ for all $H'$-sets $Y$, which furthermore implies that $\defic(H')=0$.~\smallqed

\medskip
By supposition, $|\partial(Z)| \ge 1$. Let $Z'$ be a set of new vertices (not in $H$), where $|Z'| = \max(0,4 - |\partial(Z)|)$. We note that $|Z' \cup \partial(Z)| \ge 4$ and $0 \le |Z'| \le 3$. Further if $|\partial(Z)| \ge 4$, then $Z' = \emptyset$. Let $H_e$ be the hypergraph obtained from $H'$ by adding the set $Z'$ of new vertices and adding a $4$-edge, $e$, containing four vertices in $Z' \cup \partial(Z)$. Note that $H_e$ may not be linear as the edge $e$ may overlap other edges in $H_e$. Every $\tau(H_e)$-transversal is a transversal of $H'$. Further, if $Z' = \emptyset$, then in order to cover the edge $e$ every $\tau(H_e)$-transversal contains a vertex of $\partial(Z)$, while if $Z' \ne \emptyset$, then the edge $e$ contains at least one vertex of $\partial(Z)$ and we can choose a $\tau(H_e)$-transversal to contain such a vertex of $\partial(Z)$. Therefore, $|T'| \le \tau(H_e)$, where recall that $T'$ is a smallest possible transversal in $H'$ containing a vertex from $\partial(Z)$.

\medskip


\ClaimX{L.2}
$H_e$ is not linear and $|\partial(Z)| \ge 2$.

\ProofClaimX{L.2}
Suppose, to the contrary, that $H_e$ is linear. By construction, $n(H_e) \le n(H)$ and $m(H_e) = m(H') + 1 < m(H)$. Therefore, $n(H_e) + m(H_e) < n(H) + m(H)$. Since $H$ is a counterexample with minimum value of $n(H)+m(H)$, we note that $45 \tau(H_e) \le 6n(H_e) + 13m(H_e) + \defic(H_e)$.

\ClaimX{L.2.1}
$\defic(H_e) \ge 9$.

\ProofClaimX{L.2.1}
Suppose, to the contrary, that $\defic(H_e) \le 8$. Recall that $m(H_e) = m(H') + 1$. If $Z' = \emptyset$, then $n(H_e) = n(H')$ and either $|\partial(Z)| = 4$ and $f(|\partial(Z)|) = 23$ or $|\partial(Z)| \ge 5$ and $f(|\partial(Z)|) = 22$. In this case,

\[
\begin{array}{rcl}
45 |T'| & \le & 45 \tau(H_e) \\
& \le & 6n(H_e) + 13m(H_e) + \defic(H_e) \\
& = & 6n(H') + 13m(H')| + 13 + 8 \\
& \le & 6n(H') + 13m(H') + f(|\partial(Z)|) - 1,
\end{array}
\]

\noindent
a contradiction. If $|Z'| = 1$, then $n(H_e) = n(H') + 1$, $|\partial(Z)| = 3$ and $f(|\partial(Z)|) = 27$, implying that $45 |T'| \le 6n(H') + 13m(H')| + 6 + 13 + 8 = 6n(H') + 13m(H') + f(|\partial(Z)|)$, a contradiction. If $|Z'| = 2$, then $n(H_e) = n(H') + 2$, $|\partial(Z)| = 2$ and $f(|\partial(Z)|) = 33$, implying that $45 |T'| \le 6n(H') + 13m(H')| + 12 + 13 + 8 = 6n(H') + 13m(H') + f(|\partial(Z)|)$, a contradiction. If $|Z'| = 3$, then $n(H_e) = n(H') + 3$, $|\partial(Z)| = 1$ and $f(|\partial(Z)|) = 39$, implying that $45 |T'| \le 6n(H') + 13m(H')| + 18 + 13 + 8 = 6n(H') + 13m(H') + f(|\partial(Z)|)$, a contradiction. This completes the proof of Claim~L.2.1.~\smallqed

\medskip
By Claim~L.2.1, $\defic(H_e) \ge 9$. Let $Y_e$ be a special $H_e$-set in $H_e$ such that $\defic_{H_e}(Y_e) = \defic(Y_e) \ge 9$. Suppose that $|Y_e| \ge 2$. In this case $|E_{H_e}^*(Y_e)| \le |Y_e|-2$, for otherwise $\defic(H_e) < 9$. Let $Y_e'$ be equal to $Y_e$, except we remove any special subhypergraph from $Y_e$ if it contains the edge $e$. We note that $|E_{H'}(Y_e')| \le |E_{H_e}^*(Y_e)| \le |Y_e|-2 \le |Y_e'|-1$, contradicting Claim~L.1. Therefore, $|Y_e|=1$. Since $\defic(H_e) \ge 9$ and $|Y_e|=1$, the  special $H_e$-set consists of one  $H_{10}$-component in $H_e$ and $\defic(H_e) = 10$. By Claim~J, this $H_{10}$-component, $Y_e$, in $H_e$ contains the edge $e$.

If $|\partial(Z)| \le 3$, then $Z' \ne \emptyset$ and the edge $e$ would contain a vertex of degree~$1$ in $H_e$, implying that $Y_e$ would contain a vertex of degree~$1$. However, $H_{10}$ is $2$-regular, a contradiction. Hence, $|\partial(Z)| \ge 4$. Thus, $Z' = \emptyset$ and $n(H_e) = n(H')$. If $|\partial(Z)| = 4$, then $f(|\partial(Z)|) = 23$, implying that $45 |T'| \le 45 \tau(H_e) \le 6n(H_e) + 13m(H_e) + \defic(H_e) = 6n(H') + 13m(H')| + 13 + 10 = 6n(H') + 13m(H') + f(|\partial(Z)|)$,
a contradiction. Hence, $|\partial(Z)| \ge 5$, and so $f(|\partial(Z)|) = 22$.

If there is a vertex $v \in \partial(Z) \setminus V(Y_e)$, then we can change $e$ by removing from it an arbitrary vertex and adding to it the vertex $v$ instead. With this new choice of the added edge $e$, we note that $H_e$ is once again linear but now the component in $H_e$ containing the edge $e$ contains at least~$11$ vertices, implying that there would be no $H_{10}$-component in $H_e$ containing the edge $e$, contradicting our earlier arguments. Therefore, $\partial(Z) \subseteq V(Y_e)$. Since $H$ is connected, this implies that $H' = Y_e - e$. Let $u \in \partial(Z) \setminus V(e)$ be arbitrary. We note that $u$ exists since $|\partial(Z)| \ge 5$ and $|V(e)| = 4$. We note further that there are exactly two edges in $H' - u$ and these two edges intersect in a vertex, say $w$. Hence, $\{u,w\}$ is a transversal in $H'$ containing a vertex from $\partial(Z)$, implying that $45 |T'| \le 90 < 60 + 52 + 21 = 6n(H') + 13m(H') + f(|\partial(Z)|) - 1$, a contradiction. Therefore, $H_e$ cannot be linear. This implies that $|\partial(Z)| \ge 2$, which
completes the proof of Claim~L.2.~\smallqed

\medskip

\ClaimX{L.3}
For any non-empty special $H'$-set, $Y$, in $H'$ we have $|E_{H'}^*(Y)| \ge |Y|+1$ or $|\partial(Z)\setminus V(Y)| \le 1$.

\ProofClaimX{L.3}
Suppose, to the contrary, that there is a non-empty special  $H'$-set, $Y$ in $H'$ where $|E_{H'}^*(Y)| \le |Y|$ and $|\partial(Z)\setminus V(Y)| \ge 2$.
Assume that $|Y|$ is minimum possible with the above property. By Claim~L.1, $|E_{H'}^*(Y)| = |Y|$.
By Claim~E, we note that $\partial(Z) \cap V(Y) \ne \emptyset$. Let $Z^* = Z \cup V(Y)$ and let $H^*=H - Z^*$. Further, let $Q^* = \partial(Z^*)$, and so $Q^* = N_H(Z^*) \setminus Z^*$. Since $|\partial(Z) \setminus V(Y)| \ge 2$, we have $|Q^*| \ge 2$.

\ClaimX{L.3.1}
$\defic(H^*) \le 21$.

\ProofClaimX{L.3.1}
Suppose, to the contrary, that $\defic(H^*) \ge 22$. Let $Y^*$ be a special $H^*$-set with $\defic_{H^*}(Y^*) \ge 22$.  By Claim~J, there is no $H_{10}$-subhypergraph in $H$, implying that $\defic_{H^*}(Y^*) \le 8|Y^*| - 13|E_{H^*}^*(Y^*)|$. If $|E_{H^*}^*(Y^*)| \ge |Y^*| - 2$, then $22 \le \defic_{H^*}(Y^*) \le 8|Y^*| - 13|Y^*| + 26 = -5|Y^*| + 26$, implying that $|Y^*| = 0$, a contradiction. Thus, $|E_{H^*}^*(Y^*)| \le |Y^*| - 3$. Hence, $|E_{H'}^*(Y \cup Y^*)| \le |E_{H'}^*(Y)| + |E_{H^*}^*(Y^*)| \le |Y| + |Y^*|-3 < |Y \cup Y^*|$, contradicting Claim~L.1.~\smallqed

Let $T^*$ be a minimum transversal in $H^*$ that contains a vertex from $Q^*$. By the maximality of $|Z|$, we note that $45|T^*| \le 6n(H^*) + 13m(H^*) + f(|Q^*|)$.

\ClaimX{L.3.2}
$|T'| \le |T^*| + \tau(Y)$.

\ProofClaimX{L.3.2}
By construction of $T^*$, either $\partial(Z) \cap T^* \ne \emptyset$ or $T^*$ intersects an edge in $E_{H'}^*(Y)$. We now consider the bipartite graph, $G_Y$, with partite sets $Y$ and $E_{H'}^*(Y) \cup \{q\}$, where an edge joins $e \in E_{H'}^*(Y)$ and $F \in Y$ in $G_Y$ if and only if the edge $e$ intersects the subhypergraph $F$ of $Y$ in $H'$ and $q$ is joined to $F \in Y$ in $G_Y$ if and only if $\partial(Z) \cap V(F) \ne \emptyset$. We now construct $G_Y'$ from $G_Y$ by removing either $q$ if $\partial(Z) \cap T^* \ne \emptyset$ or by removing some $e' \in E_{H'}^*(Y)$ where  $T^*$ intersects the edge $e'$. We note that since $|E_{H'}^*(Y)| = |Y|$, the partite sets in $G_Y'$ have the same size.

Suppose there is no perfect matching in $G_Y'$. By Hall's Theorem, there is a nonempty subset $S \subseteq V(G_Y')\setminus Y$ such that $|N_{G_Y'}(S)|<|S|$. We now consider the special $H'$-set, $Y' = Y \setminus N_{G_Y'}(S)$. Then, $|Y'| = |Y| - |N_{G_Y'}(S)|$. Since possibly $q \in S$, we therefore have $|E_{H'}^*(Y')| \le |E_{H'}^*(Y)| - (|S|-1) \le |Y| - |N_{G_Y'}(S)| = |Y'|$. Thus, $Y'$ is a non-empty special $H'$-set such that $|E_{H'}^*(Y')| \le |Y'|$ and $|Y'|<|Y|$. Further, since $Y' \subset Y$ and $|\partial(Z) \setminus V(Y)| \ge 2$, we note that $|\partial(Z) \setminus V(Y')| \ge 2$. This contradicts our choice of $Y$. Therefore, there is a perfect matching in $G_Y'$, implying that we can find a $\tau(Y)$-transversal that together with $T^*$ intersects all the edges in $E_{H'}^*(Y)$ and intersects $\partial(Z)$. This implies that $|T'| \le |T^*| + \tau(Y)$, as desired.~\qed

\medskip
As observed earlier, $45|T^*| \le 6n(H^*) + 13m(H^*) + f(|Q^*|)$. Since $Y$ is a non-empty special $H'$-set, we note that $45 \tau(Y) = 6n(Y) + 13m(Y) + \defic_{H'}(Y) \le 6n(Y) + 13m(Y) + \defic_{H'}(Y) + 13|E_{H'}^*(Y)|$. By Claim~L.3.2, $|T'| \le |T^*| + \tau(Y)$. Therefore,
\[
\begin{array}{rcl}
45 |T'| & \le & 45 |T^*| + 45 \tau(Y) \1 \\
        & \le & 6n(H^*) + 13m(H^*) + f(|Q^*|) + 6n(Y) + 13m(Y) + 13|E_{H'}^*(Y)| + \defic_{H'}(Y) \1 \\
        & = & 6n(H') + 13m(H') + f(|Q^*|) + \defic_{H'}(Y). \\
\end{array}
\]

As observed earlier, $|E_{H'}^*(Y)| = |Y|$, implying that $\defic_{H'}(Y) < 0$. In fact, by Claim~J, $\defic_{H'}(Y) = 8|Y_4| + 5|Y_{14}| + 4|Y_{11}| + |Y_{21}| - 13|Y| = -5|Y_4| - 8|Y_{14}| -9|Y_{11}| -12|Y_{21}|$. Since $|Q^*| \ge 2$, we note that $f(|Q^*|) \le 33$.

Suppose that $|Y| \ge 2$. In this case, $\defic_{H'}(Y) \le -10$, implying that $f(|Q^*|) + \defic_{H'}(Y) \le 33 - 10 = 22 \le f(|\partial(Z)|)$. By supposition, the transversal $T'$ in $H'$ does not satisfy (a) or (b) in the statement of the claim. Hence we immediately obtain a contradiction unless $f(|Q^*|) = 33$, $\defic_{H'}(Y) = -10$ and $|\partial(Z)| \ge 5$. This in turn implies that $|Y|=2$, $|Y_4|=2$ and $|Q^*|=2$. However if $|Y|=|Y_4|=2$ and $|E_{H'}(Y)|=2$, then by the linearity of $H$ the two edges intersecting $Y$ both have two vertices not in $Y$, implying that $|Q^*| \ge 3$, a contradiction to $|Q^*|=2$. Therefore, $|Y|=1$.
Let $E_{H'}^*(Y) = \{f\}$. Recall that by supposition, $|\partial(Z) \setminus V(Y)| \ge 2$.

Suppose that $|V(f) \cap V(Y)| \ge 2$. By the linearity of $H$ we note that $Y_4 = \emptyset$. Thus, $\defic_{H'}(Y) = - 8|Y_{14}| -9|Y_{11}| - 12|Y_{21}| \le - 8$, and so $45 \tau(Y) = 6n(Y) + 13m(Y) + \defic_{H'}(Y) \le 6n(Y) + 13m(Y) - 8$. By Observation~\ref{property:special}(k) we can find a $\tau(Y)$-transversal, $T_Y$, that contains a vertex in $\partial(Z)$ and a vertex in $f$. Let $H''$ be the hypergraph obtained from $H'$ by removing the vertices $V(Y)$ and edges $E(Y) \cup E_{H'}^*(Y)$; that is, $H'' = H' - T_Y$. If $\defic(H'')>0$, then letting $Y''$ be a special $H''$-set with $\defic_{H''}(Y'') > 0$ we note that $|E_{H'}^*(Y \cup Y'')| \le |E_{H'}^*(Y)| + |E_{H''}^*(Y'')| \le 1 + (|Y''| -1) < |Y'' \cup Y|$, contradicting Claim~L.1. Therefore, $\defic(H'')=0$, and so $45\tau(H'') \le 6n(H'') + 13m(H'') + \defic(H'') = 6n(H'') + 13m(H'')$. Hence,
\[
\begin{array}{rcl}
45 |T'| & \le & 45|T_Y| + 45\tau(H'') \1 \\
       & \le & (6n(Y) + 13m(Y) - 8) + (6n(H'') + 13m(H'')) \1 \\
       & = & 6n(H') + 13(m(H') - 1) - 8  \1 \\
        & < & 6n(H') + 13m(H') - 8, \\
\end{array}
\]

\noindent
contradicting the supposition that the transversal $T'$ in $H'$ does not satisfy (a) or (b) in the statement of the claim. Therefore, $|V(f) \cap V(Y)| = 1$, implying that $|Q^*| \ge 3$. Recall that
\[
45 |T'| \le 6n(H') + 13m(H') + f(|Q^*|) + \defic_{H'}(Y).
\]

If $Y_4 = \emptyset$, then $\defic_{H'}(Y) \le -8$, and so $f(|Q^*|) + \defic_{H'}(Y) \le 27 - 8 = 19$. If $Y_4 \ne \emptyset$ and $|Q^*| \ge 4$, then $f(|Q^*|) + \defic_{H'}(Y) \le 23 - 5 = 18$. In both cases, $45 |T'| \le 6n(H') + 13m(H') + f(|Q^*|) + \defic_{H'}(Y) \le 6n(H') + 13m(H') + f(|\partial(Z)|) - 1$, a contradiction. Therefore, $|Y|=|Y_4|=1$ and $|Q^*|=3$, implying that $f(|Q^*|) + \defic_{H'}(Y) = 27 - 5 = 22$ and $45 |T'| \le 6n(H') + 13m(H') + 22 \le 6n(H') + 13m(H') + f(|\partial(Z)|)$. This proves part~(a) of Claim~L.

In order to prove part~(b) of Claim~L, we consider next the case when $|\partial(Z)| \ge 5$. Since $|Q^*| = 3$, we note that $\partial(Z) \subseteq V(Y) \cup V(f)$. Let $E(Y) = \{e'\}$. By our earlier observations, $e'$ and $f$ intersect in exactly one vertex, say $v'$. Further, the edge $e'$ contains three degree-$1$ vertices in $H'$. If $v'$ belongs to $\partial(Z)$, then we can find a $\tau(Y)$-transversal that contains a vertex in $\partial(Z)$ and a vertex in $f$, which analogously to our previous arguments gives us a contradiction. Therefore, $V(Y) \cap V(f) \cap \partial(Z) = \emptyset$, implying that $\partial(Z) \subseteq (V(e') \cup V(f)) \setminus \{v'\}$. Since $|\partial(Z)| \ge 5$, this implies that $|\partial(Z) \cap V(e')| \ge 2$ and $|\partial(Z) \cap V(f)| \ge 2$. Hence, $|\partial(Z)| \ge 5$ and $H$ contains two intersecting edges $e$ and $f$ such that (i), (ii) and (iii) hold in the statement of part~(b) in Claim~L. Hence, from our earlier observations, part~(b) of Claim~L holds. This completes the proof of Claim~L.3.~\smallqed

\medskip
We now return to the proof of Claim~L.  By Claim~L.2 there exists an edge $f$ in $H'$ that overlaps the edge $e$ in $H_e$. Let $\{u,v\} \subseteq V(e) \cap V(f)$ be arbitrary. Let $H_{ef}$ be obtained from $H-e-f$ by adding two new vertices $z_1$ and $z_2$ and a new edge $g=\{u,v,z_1,z_2\}$. Since $H$ is linear, so too is $H_{ef}$ is linear.
We show that $\defic(H_{ef}) \le 9$. Suppose, to the contrary, that $\defic(H_{ef}) > 9$ and let $Y_{ef}$ be a special $H_{ef}$-set with $\defic_{H_{ef}}(Y_{ef}) > 9$. By Claim~J, $(Y_{ef})_{10} = \emptyset$, implying that $|E_{H_{ef}}^*(Y_{ef})| \le |Y_{ef}|-2$. If $g \notin E(Y_{ef})$, then
\[
|E_{H'}(Y_{ef})| \le |E_{H_{ef}}(Y_{ef})| + |\{f\}| \le (|Y_{ef}|-2) + 1 < |Y_{ef}|,
\]
contradicting Claim~L.1. Therefore, $g \in E(Y_{ef})$. Since $g$ has two degree-$1$ vertices in $H_{ef}$, we note that $g$ is the edge of a $H_4$-hypergraph, $R_g$, in $Y_{ef}$. In this case,
\[
|E_{H'}^*(Y_{ef} \setminus \{R_g\})| \le |E_{H_{ef}}^*(Y_{ef})| + |\{f\}| \le (|Y_{ef}|-2) + 1 = |Y_{ef} \setminus \{R_g\}|.
\]

As $\{u,v\} \subseteq \partial(Z)$ and $u$ and $v$ do not belong to $Y_{ef} \setminus \{R_g\}$ (as they belong to $V(R_g)$) we obtain a contradiction to Claim~L.3. Therefore, $\defic(H_{ef}) \le 9$. Thus, by the minimality of $n(H) + m(H)$,
\[
\begin{array}{rcl}
 45|T'| \, \le \, 45 \tau(H_{ef}) & \le & 6n(H_{ef}) + 13m(H_{ef}) + \defic(H_{ef}) \\
 & \le & 6(n(H')+2) + 13m(H') + \defic(H_{ef}) \\
& = & 6n(H') + 13m(H') + 12 + \defic(H_{ef}) \\
& \le & 6n(H') + 13m(H') +  21 \\
& \le & 6n(H') + 13m(H') + f(|\partial(Z)|) - 1,
\end{array}
\]
\noindent
a contradiction. This completes the proof of Claim~L.~\smallqed

\medskip
We call a component of a $4$-uniform, linear hypergraph that contains two vertex disjoint copies of $H_4$ that are both intersected by a common edge and such that each copy of $H_4$ has three vertices of degree~$1$ and one vertex of degree~$2$ a \emph{double}-\emph{$H_4$}-\emph{component}. We call these two copies of $H_4$ the $H_4$-\emph{pair} of the double-$H_4$-component, and the edge that intersects them the \emph{linking edge}. We note that a double-$H_4$-component contains at least ten vertices, namely eight vertices from the $H_4$-pair and at least two additional vertices that belong to the linking edge.

\medskip
\ClaimX{M}
If $x$ is an arbitrary vertex of $H$ of degree~$3$, then one of the following holds. \\[-27pt]
\begin{enumerate}
\item  $\defic(H-x) = 8$ and the hypergraph $H - x$ contains an $H_4$-component that is intersected by all three edges incident with~$x$.
\item $\defic(H-x) = 3$ and the hypergraph $H - x$ contains a double-$H_4$-component. Further, the $H_4$-pair in this component is intersected by all three edges incident with~$x$.
\end{enumerate}

\ProofClaimX{M}
Let $x$ be an arbitrary vertex of $H$ of degree~$3$ and let $e_1$, $e_2$ and $e_3$ be the three edges incident with $x$. By the linearity of $H$, the vertex $x$ is the only common vertex in $e_i$ and $e_j$ for $1 \le i < j \le 3$. Let $H' = H - x$. We note that $n(H') = n(H) - 1$ and $m(H') = m(H) - 3$. Every transversal in $H'$ can be extended to a transversal in $H$ by adding to it the vertex~$x$. Hence, applying the inductive hypothesis to $H'$, we have that
$45\tau(H) \le 45(\tau(H') + 1) \le 6n(H') + 13m(H') + \defic(H') + 45 \le 6n(H) + 13m(H) + \defic(H')$. If $\defic(H') = 0$, then $45\tau(H) \le 6n(H) + 13m(H)$, implying that $\xi(H) \le 0$, contradicting the fact that $H$ is a counterexample to the theorem. Hence, $\defic(H') > 0$. Let $X$ be a special $H'$-set satisfying $\defic(H') = \defic_{H'}(X)$.

\ClaimX{M.1}
$|E_{H'}^*(X)| = |X| - 1$. Further, all three edges incident with the vertex $x$ intersect $X$.

\ProofClaimX{M.1} Since $\defic(H') > 0$, we note that $|E_{H'}^*(X)| \le |X| - 1$. By Claim~I and since $E_{H}^*(X) \setminus  E_{H'}^*(X) \subseteq \{e_1,e_2,e_3\}$, we note that $|X| + 2 \le |E_{H}^*(X)| \le |E_{H'}^*(X)| + 3$, implying that $|E_{H'}^*(X)| \ge |X| - 1$. Consequently,  $|E_{H'}^*(X)| = |X| - 1$ and all three edges incident with the vertex $x$ intersect $X$.~\smallqed

\ClaimX{M.2}
$|X| = 1$ or $|X| = 2$ and $X = X_4$.

\ProofClaimX{M.2} By Claim~M.1, $|E_{H'}^*(X)| = |X| - 1$ and all three edges incident with the vertex $x$ intersect $X$. By Claim~J, there is no $H_{10}$-subhypergraph in $H$. Thus,
\[
\begin{array}{lcl}
\defic_{H'}(X) & = & 8|X_4| + 5|X_{14}| + 4|X_{11}| + |X_{21}| - 13|E_{H'}^*(X)| \\
& = & 8|X_4| + 5|X_{14}| + 4|X_{11}| + |X_{21}| - 13(|X| - 1) \\
& = & 13 - 5|X_4| - 8|X_{14}| - 9|X_{11}| - 12|X_{21}|. \\
\end{array}
\]
\indent
Since $\defic_{H'}(X) > 0$, either $|X| = 1$ or $|X| = 2$ and $X = X_4$.~\smallqed

\ClaimX{M.3}
If $|X| = 2$, then $H - x$ contains a double-$H_4$-component and $\defic(H') = 3$.

\ProofClaimX{M.3} Suppose that $|X| = 2$. By Claim~M.2, $X = X_4$, and so $\defic_{H'}(X) = 3$. By Claim~M.1, $|E_{H'}^*(X)| = 1$. Let $e$ be the edge in $E_{H'}^*(X)$. If $e$ intersects only one of the copies of $H_4$ in $X$, then the other copy of $H_4$ is an $H_4$-component in $H - x$ implying that $\defic(H') \ge 8$, contradicting the fact that $\defic(H') = \defic_{H'}(X) = 3$. Therefore, $e$ intersects both copies of $H_4$ in $X$. By the linearity of $H$ the edge $e$ intersects each copy in one vertex, implying that $H - x$ contains a double-$H_4$-component.~\smallqed

\ClaimX{M.4}
If $|X| = 1$, then $H - x$ contains an $H_4$-component and $\defic(H') = 8$.

\ProofClaimX{M.4} Suppose that $|X| = 1$, and let $R$ be the special subhypergraph in $X$. By Claim~M.1 we note that $|E_{H'}^*(X)| = 0$, and so $R$ is a component of $H - x$. Suppose, to the contrary, that $|X_4| = 0$. Among all degree-$3$ vertices, we choose the vertex $x$ so that $X = X_j$ where $j$ is a maximum; that is, $X$ is chosen to contain a special subhypergraph of maximum possible size~$j$. By supposition, $j \ne 4$, and so $j \in \{11,14,21\}$.

Suppose that each edge incident with $x$ intersects $R$ in at most two vertices. By Claim~M.1, all three edges incident with $x$ intersect $R$ in at least one vertex, implying that in this case at least three neighbors of $x$ in $H$ belong to $R$. Since every special subhypergraph different from $H_4$ contains at most two vertices of degree~$1$, we can choose a neighbor $y$ of $x$ in $R$ such that $d_{H'}(y) \ge 2$. Since every neighbor of $x$ in $R$ has degree at most~$2$ in $R$ since $H$ has maximum degree~$3$, this implies that $d_{H'}(y) = 2$,  and so $d_H(y) = 3$. This implies by Observation~\ref{property:special}(p) that $H' - y$ is connected or is disconnected with exactly two components, one of which consists of an isolated vertex. In both cases, since all three edges incident with~$x$ intersect $R$ we note that either $H - y$ is connected and has cardinality at least~$n(R) + 3$ or $H - y$ is disconnected with two components, one of which consists of an isolated vertex with the other component of cardinality at least~$n(R) + 2$.  Analogously as with the vertex~$x$, there is a special $(H - y)$-set, $Y$, with $\defic(H-y) = \defic_{H-y}(Y)$ satisfying $|Y| = 1$ or $|Y| = 2$ and $Y = Y_4$.

If $|Y| = 1$, then by our earlier observations, $Y$ consists of a special subhypergraph of cardinality at least~$n(R) + 2$, contradicting the maximality of $R$. Hence, $|Y| = 2$ and $Y = Y_4$. Analogously as with the vertex $x$ (see Claim~M.3), $H - y$ contains a double-$H_4$-component. Let $R_1$ and $R_2$ be the $H_4$-pair of this double-$H_4$-component, and let $h^*$ be the linking edge that intersects them. We note that $Y = \{R_1,R_2\}$. By Observation~\ref{property:special}(p), $R - y$ does not contain a double-$H_4$-component. Further, $R-y$ does not contain a component of order~$4$.

Suppose that the edge of $R_1$ belongs to $R$. If the edge $h^*$ does not belong to $R$, then $R_1$ is a component in $R-y$ (of order~$4$), a contradiction. Hence, the edge $h^*$ also belongs to $R$. If the edge in $R_2$ belongs to $R$, then $R-y$ would contain a double-$H_4$-component, a contradiction. Hence, the edge in $R_2$ does not belong to $R$ and must therefore contain the vertex~$x$. However in this case $R_2$ is intersected by at least two edges in $R - y$, namely $h^*$ and the edge incident with $x$ that intersects $R_2$. Therefore, the pair $R_1$ and $R_2$ is not the $H_4$-pair of a double-$H_4$-component in $H - y$, a contradiction. Therefore, the edge in $R_1$ does not belong to $R$. Analogously, the edge in $R_2$ does not belong to $R$.

Since the $H_4$-pair, $R_1$ and $R_2$, in this double-$H_4$-component is intersected by all three edges incident with~$y$, we note that both $R_1$ and $R_2$ intersect $V(R)$. However as they do not belong to $R$, they must both contain the vertex~$x$, which
implies that they are not vertex disjoint, a contradiction.  Therefore, at least one edge incident with $x$ intersects $R$ in three vertices.

Renaming edges if necessary, we may assume that the edge $e_1$ intersects $R$ in three vertices. From the structure of special subhypergraphs, this implies that we can choose a vertex $y \in V(e_1) \cap V(R)$ so that $R - y$ is connected (of cardinality~$n(R) - 1$). If at least one neighbor of $x$ does not belong to $R$, then analogously as before $H - y$ is a connected hypergraph of cardinality at least~$n(R) + 1$, contradicting the maximality of $R$. Hence, letting $T_i = V(e_i) \cap V(R)$, we note that $|T_i| = 3$ for all $i \in [3]$. Further, by the linearity of $H$, the sets $T_1$, $T_2$ and $T_3$ are vertex disjoint independent sets in $R$. By Observation~\ref{property:special}(n) there exists a $\tau(R)$-transversal, $T^*$ say, that contains a vertex from each of the sets $T_i$ for $i \in [3]$. Since $H$ is connected and all neighbors of $x$ belong to $R$, we note that $H' = R$, $V(H) = V(R) \cup \{x\}$ and $E(H) = E(R) \cup \{e_1,e_2,e_3\}$. Thus, $T^*$ is a transversal of $H$, implying that $45\tau(H) \le 45|T^*| = 45\tau(H') \le 6n(H') + 13m(H') + \defic(H') = 6n(H) + 13m(H) + \defic(H') - 45$. By Claim~M.2, $\defic(H') \le 8$, and so $45\tau(H) < 6n(H) + 13m(H)$, a contradiction. Therefore, $X = X_4$, implying that $R = H_4$ is an $H_4$-component of $H - x$ and $\defic(H') = 13 - 5|X_4| = 8$.~\smallqed

\medskip
By Claims~M.2, M.3 and M.4, $H - x$ contains an $H_4$-component or a double-$H_4$-component. This completes the proof of Claim~M.~\smallqed


\ClaimX{N}
No edge in $H$ contains two degree-$1$ vertices of $H$.

\ProofClaimX{N}
Suppose, to the contrary, that $e = \{v_1,v_2,v_3,v_4\}$ is an edge in $H$, where $1 = d_H(v_1) = d_H(v_2) \le d_H(v_3) \le d_H(v_4)$. Suppose that $d_H(v_4) = 3$. Let $H'$ be obtained from $H - v_4$ by deleting all resulting isolated vertices (including $v_1$ and $v_2$). By Claim~M, $\defic(H') \le 8$. We note that $n(H') \le n(H) - 3$ and $m(H') = m(H) - 3$. Applying the inductive hypothesis to $H'$, we have that $45\tau(H) \le 45(\tau(H') + 1) \le 6n(H') + 13m(H') + \defic(H') + 45 \le 6n(H) - 3 \times 6 + 13m(H) - 3 \times 13 + \defic(H') + 45 \le  6n(H) + 13m(H) - 18 - 39 + 8 + 45  < 6n(H) + 13m(H)$, a contradiction. Therefore, $d_H(v_4) \le 2$, implying that at most two edges of $H$ intersect the edge $e$. Taking $Y$ to the special $H$-set consisting only of the special subhypergraph $H_4$ with $e$ as its edge, we have $|E^*(Y)| \le 2 = |Y|+1$, contradicting Claim~I(a).~\smallqed

\medskip
\ClaimX{O}
Let $H'$ be a $4$-uniform, linear hypergraph with no $H_{10}$-subhypergraph satisfying $\defic(H') > 0$. If $Y$ is a special $H'$-set satisfying $\defic(H') = \defic_{H'}(Y)$, then $|E_{H'}^*(Y)| = |Y| - i$ for some $i \ge 1$ and
\[
\defic(H') \le 13i - 5|Y_4| - 8|Y_{14}| - 9|Y_{11}| - 12|Y_{21}|.
\]
In particular, $\defic(H') \le 13i - 5|Y|$. Further, if $\defic(H') > 8j$ for some $j \ge 0$, then $|E_{H'}^*(Y)| \le |Y| - (j+1)$.

\ProofClaimX{O}
Since $\defic_{H'}(Y) > 0$, we note that $|E_{H'}^*(Y)| = |Y| - i$ for some $i \ge 1$. Thus, $\defic_{H'}(Y) = 8|Y_4| + 5|Y_{14}| + 4|Y_{11}| + |Y_{21}| - 13|E_{H'}^*(Y)| = 8|Y_4| + 5|Y_{14}| + 4|Y_{11}| + |Y_{21}| - 13(|Y| - i) = 13i - 5|Y_4| - 8|Y_{14}| - 9|Y_{11}| - 12|Y_{21}|$, and so $\defic_{H'}(Y) \le 13i - 5|Y|$.

Further, let $\defic(H') > 8j$ for some $j \ge 0$ and suppose, to the contrary, that $|E_{H'}^*(Y)| \ge |Y| - j$. Hence,  $|E_{H'}^*(Y)| = |Y| - k$ where $k \le j$. Therefore, $\defic(H') \le 13k - 5|Y| \le 13k - 5k = 8k \le 8j$, a contradiction. Therefore, $|E_{H'}^*(Y)| \le |Y| - (j+1)$.~\smallqed

\medskip
We show next that the removal of any specified vertex of degree~$3$ from $H$ produces an $H_4$-component. In what follows, we use the following notation for simplicity. If $e$ and $f$ are intersecting edges of $H$, then we denote the vertex in this intersection by $(ef)$; that is, $e \cap f = \{ (ef) \}$.

\medskip
\ClaimX{P}
If $x$ is an arbitrary vertex of $H$ of degree~$3$, then $H - x$ contains an $H_4$-component.

\ProofClaimX{P}
Let $x$ be an arbitrary vertex of $H$ of degree~$3$. Let $e_1$, $e_2$ and $e_3$ be the three edges incident with $x$ and let $H' = H - x$. By the linearity of $H$, the vertex $x$ is the only common vertex in $e_i$ and $e_j$ for $1 \le i < j \le 3$. Suppose, to the contrary, that $H'$ does not contain an $H_4$-component. By Claim~M, $H'$ therefore contains a double-$H_4$-component, say $C'$, and $\defic(H') = 3$. Let $f_1$ and $f_2$ be the two edges belonging to the $H_4$-pair in $C'$ and let $h$ be the linking edge in $H'$ that intersects $f_1$ and $f_2$. By definition of a double-$H_4$-component, we note that $f_i$ has three vertices of degree~$1$ and one vertex of degree~$2$ in $H'$ for $i \in [2]$. We proceed further with the following series of subclaims.

\ClaimX{P.1}
$2 \le |f_i \cap N_H(x)| \le 3$ for $i \in [2]$.

\ProofClaimX{P.1}
By the linearity of $H$, $|f_i \cap e_j| \le 1$ for every $j \in [3]$, and so $|f_i \cap N_H(x)| \le 3$. If $|f_i \cap N_H(x)| \le 1$, then the edge $f_i$ contains at least two vertices of degree~$1$ in $H$, contradicting Claim~N.~\smallqed

\ClaimX{P.2}
Any two edges amongst $e_1$, $e_2$ and $e_3$ can be covered by two vertices one of which belongs to $f_1$ and the other to $f_2$.

\ProofClaimX{P.2}
Renaming edges if necessary, it suffices to consider the two edges $e_1$ and $e_2$. By Claim~P.1, the edge $f_1$ intersects at least one of $e_1$ and $e_2$, say $e_1$. If $f_2$ intersects $e_2$, then we can cover $e_1$ by the vertex $(e_1f_1)$ and we can cover $e_2$ by the vertex $(e_2f_2)$. If $f_2$ does not intersect $e_2$, then by Claim~M and Claim~P.1, the edge $f_1$ intersects $e_2$ and the edge $f_2$ intersects $e_1$, and we can therefore cover $e_1$ by the vertex $(e_1f_2)$ and we can cover $e_2$ by the vertex $(e_2f_1)$.~\smallqed

\medskip
Let $h \setminus (f_1 \cup f_2) = \{a,b\}$, and let $h \cap f_1 = \{c\}$ and $h \cap f_2 = \{d\}$. Thus, $h = \{a,b,c,d\}$. By definition of a double-$H_4$-component, we note that $d_H(c) = d_{H'}(c) = 2$ and $d_H(d) = d_{H'}(d) = 2$.

\ClaimX{P.3}
$d_{H'}(a) \ge 2$ and $d_{H'}(b) \ge 2$.

\ProofClaimX{P.3}
Suppose, to the contrary, that $d_{H'}(a) = 1$ or $d_{H'}(b) = 1$. Interchanging the names of $a$ and $b$ if necessary, we may assume that $d_{H'}(a) = 1$. Let $H''$ be obtained from $H - \{x,c,d\}$ by deleting all resulting isolated vertices (including the vertex $a$, as well as three vertices of degree~$1$ in each of $f_1$ and $f_2$). We note that $n(H'') \le n(H') - 10$. Further, the edges $e_1$, $e_2$, $e_3$, $f_1$, $f_2$ and $h$ are deleted from $H$ when constructing $H''$, and so $m(H'') = m(H) - 6$. Applying the inductive hypothesis to $H''$, we have that $45\tau(H) \le 45(\tau(H'') + 3) \le 6n(H'') + 13m(H'') + \defic(H'') + 135 \le 6n(H) - 6 \times 10 + 13m(H) - 13 \times 6 + \defic(H'') + 135 = 6n(H) + 13m(H) + \defic(H'') - 3$. If $\defic(H'') \le 3$, then $45\tau(H) \le 6n(H) + 13m(H)$, a contradiction. Hence, $\defic(H'') \ge 4$. Let $Y''$ be a special $H''$-set satisfying $\defic_{H''}(Y) = \defic(H'') > 0$. By Claim~O, $|E_{H''}^*(Y'')| \le |Y''| - 1$. Let $Y_i$ be the $H_4$-subhypergraph of $H$ with $E(Y_i) = \{f_i\}$ for $i \in [2]$, and consider the special $H$-set $Y = Y'' \cup \{Y_1,Y_2\}$. We note that $|E_{H}^*(Y)| \le |E_{H''}^*(Y'')| + |\{e_1,e_2,e_3,h\}| \le (|Y''| - 1) + 4 = |Y''| + 3 = |Y| + 1$, contradicting Claim~I(a).~\smallqed

\medskip
As observed earlier, the edge $f_i$ intersects the edge $e_j$ in at most one vertex for every $i \in [2]$ and $j \in [3]$. For $j \in [3]$, let $s_j$ be a vertex different from $x$ that belongs to $e_j$ but not to $f_1 \cup f_2$.

\ClaimX{P.4}
The following holds. \\[-29pt]
\begin{enumerate}
\item $d_H(s_i) \ge 2$ for all $i \in [3]$.
\item If $d_H(s_i) = 3$ for some $i \in [3]$, then $e_i \cap f_1 = \emptyset$ or $e_i \cap f_2 = \emptyset$.
\item $d_H(s_i) = 2$ for at least one $i \in [3]$.
\end{enumerate}


\ProofClaimX{P.4}
(a) Suppose, to the contrary, that the vertex $s_i$ has degree~$1$ in $H$ for some $i \in [3]$, and so $e_i$ is the only edge in $H$ containing $s_i$. Letting $H'' = H - \{x,s_i\}$, we note that $n(H'') = n(H) - 2$, $m(H') = m(H) - 3$ and $\defic(H'') = \defic(H') = 3$. Every transversal in $H''$ can be extended to a transversal in $H$ by adding to it the vertex~$x$. Hence, applying the inductive hypothesis to $H''$, we have that $45\tau(H) \le 45(\tau(H'') + 1) \le 6n(H'') + 13m(H'') + \defic(H'') + 45 \le 6n(H) + 3m(H) - 6 \cdot 2 - 13 \cdot 3 + 3 + 45 = 6n(H) + 13m(H) - 3$, a contradiction. Therefore, $d_H(s_i) \ge 2$ for all $i \in [3]$.

(b) Let $d_H(s_i) = 3$ for some $i \in [3]$, and suppose, to the contrary, both edges $f_1$ and $f_2$ intersect $e_i$. Renaming indices if necessary, we may assume $i = 1$. By Claim~M, $H - s_1$ either contains an $H_4$-component that is intersected by the edge $e_1$ or a double-$H_4$-component in which the $H_4$-pair is intersected by the edge $e_1$. However such a component of $H - s_1$ would contain the vertex~$x$ and the four edges $e_2,e_3,f_1,f_2$, which is not possible.

(c) Suppose that $d_H(s_1) = d_H(s_2) = d_H(s_3) = 3$. By Claim~P.1, $e_i \cap f_1 \ne \emptyset$ and $e_i \cap f_2 \ne \emptyset$ for some $i \in [3]$. This, however, contradicts part~(b) above.~\smallqed

\medskip
We show next that the edge $h$ contains no neighbor of $x$.

\ClaimX{P.5}
$N_H(x) \cap V(h) = \emptyset$.

\ProofClaimX{P.5}
Suppose, to the contrary, that $y \in N_H(x) \cap V(h)$.

\ClaimX{P.5.1}
$y \notin f_1 \cup f_2$.

\ProofClaimX{P.5}
Suppose, to the contrary, that $y \in f_1 \cup f_2$. Renaming edges if necessary, we may assume that $y \in f_1 \cap e_3$. Since $y$ belongs to the three edges $e_3, f_1, h$, we note that $d_H(y) = 3$. Further since $f_1 \cap h = \{y\}$, we note that the edge $f_1$ contains a vertex of degree~$1$ in $H$ that is not a neighbor of $x$. By Claim~N, no edge contains two degree-$1$ vertices, implying that in this case $f_1$ intersects all three edges $e_1$, $e_2$ and $e_3$. Let $v_i = f_1 \cap e_i$ for $i \in [3]$. By Claim~M, $H - y$ contains an $H_4$-component or a double-$H_4$-component.

If $H - y$ contains an $H_4$-component $R$, then by Claim~M the component $R$ intersects all three edges incident with $y$. In particular, $R$ intersects the edge $f_1$ and therefore contains the vertex $v_1$ or $v_2$. In both cases, the component $R$ contains the vertex $x$ and both edges $e_1$ and $e_2$, a contradiction.

Hence, $H - y$ contains a double-$H_4$-component $R'$. By Claim~M, all three edges incident with $y$ intersect the $H_4$-pair in $R'$. In particular, the edge $f_1$ intersects the $H_4$-pair in $R'$, implying that $e_1$ or $e_2$ belongs to the $H_4$-pair in $R'$. In both cases, the component $R'$ contains the vertex $x$ and both edges $e_1$ and $e_2$. Renaming $e_1$ and $e_2$ if necessary, we may assume that $e_1$ belongs to the $H_4$-pair in $R'$, and so $e_1$ contains three vertices of degree~$1$ and one vertex, namely $x$, of degree~$2$ in $H - y$. Let $e_1 = \{x,v_1,u,v\}$. We note that neither $u$ nor $v$ belong to $f_1$ or $e_3$, and at most one of $u$ and $v$ belong to $h$. Thus, at least one of $u$ and $v$, say $v$, has degree~$1$ in $H$. Thus, $v$ is isolated in $H' = H - x$. Hence, letting $H'' = H - \{x,y\}$, we note that $n(H'') = n(H) - 2$, $m(H') = m(H) - 3$ and $\defic(H'') = \defic(H') = 3$. Every transversal in $H''$ can be extended to a transversal in $H$ by adding to it the vertex~$x$. Hence, applying the inductive hypothesis to $H''$, we have that $45\tau(H) \le 45(\tau(H'') + 1) \le 6n(H'') + 13m(H'') + \defic(H'') + 45 \le 6n(H) - 6 \cdot 2 + 13m(H) - 13 \cdot 3 + \defic(H') + 45 = 6n(H) + 13m(H) - 3$, a contradiction.~\smallqed

\medskip
By Claim~P.5.1, $y \notin f_1 \cup f_2$. Thus, $y \in \{a,d\}$. Recall that $y \in N_H(x)$, and so $y$ is incident with $e_1$, $e_2$ or $e_3$. Thus, by Claim~P.3, $d_H(y) = d_{H'}(y) + 1 \ge 3$. Consequently, $d_H(y) = 3$. Renaming edges if necessary, we may assume that $y$ is incident with $e_1$. Let $g$ be the third edge (different from $e_1$ and $h$) incident with $y$.  By Claim~P.2, the edges $e_2$ and $e_3$ can be covered by two vertices one of which belongs to $f_1$ and the other to $f_2$. Thus, $e_2 \cap f_i \ne \emptyset$ and $e_3 \cap f_{3-i} \ne \emptyset$ for some $i \in [2]$. Let $Z = f_1 \cup f_2 \cup \{x,y\}$ and let $H^* = H - Z$. We note that $n(H^*) = n(H) - 10$. The edges $e_1,e_2,e_3,f_1,f_2,g,h$ are deleted from $H$ when constructing $H^*$, and so $m(H^*) = m(H) - 7$.

Every transversal in $H^*$ can be extended to a transversal in $H$ by adding to it the vertices $y$, $(e_2f_i)$ and $(e_3f_{3-i})$. Hence, applying the inductive hypothesis to $H^*$, we have that $45\tau(H) \le 45(\tau(H^*) + 3) \le 6n(H^*) + 13m(H^*) + \defic(H^*) + 45 \le 6n(H) + 13m(H) - 6 \cdot 10 - 13 \cdot 7 + \defic(H^*) + 45 \times 3 = 6n(H) + 13m(H) + \defic(H^*) + 16$. If $\defic(H^*) \le 16$, then $45\tau(H) \le 6n(H) + 13m(H)$, a contradiction. Hence, $\defic(H^*) \ge 17$. Recall that $H' = H - x$ and note that $H^*$ is obtained from $H'$ by deleting the edges $f_1, f_2, g, h$ and the resulting isolated vertex $y$. We now consider the hypergraph $H''$ obtained from $H'$ by deleting the edges $f_1, f_2, h$; that is, $H'' = H' - f_1 - f_2 - h$. We note that $\defic(H'') \ge \defic(H^*) - 13 \ge 4$ since the edge $g$, which contributes at most~$13$ to $\defic(H^*)$ has been added back to $H^*$ (along with the vertex~$y$). Let $Y''$ be a special $H''$-set satisfying $\defic_{H''}(Y'') = \defic(H'')$. Let $Y_i$ be the $H_4$-subhypergraph of $H$ with $E(Y_i) = \{f_i\}$ for $i \in [2]$, and consider the special $H'$-set $Y = Y'' \cup \{Y_1,Y_2\}$. We note that each of $Y_1$ and $Y_2$ contribute~$8$ to $\defic_{H'}(Y)$, and the edge $h$ contributes~$-13$ to $\defic_{H'}(Y)$, implying that $\defic_{H'}(Y) = \defic_{H''}(Y'') + 8 + 8 - 13 \ge 4 + 8 + 8 - 13 = 7$, contradicting the fact that $\defic(H') = 3$. This completes the proof of Claim~P.5.~\smallqed

\medskip
For $i \in [3]$, let $s_i$ be a vertex in $e_i \setminus \{x\}$ that does not belong to $f_1 \cup f_2$. By Claim~P.5, the edge $h$ contains no neighbor of $x$ in $H$. In particular, $h \cap \{s_1,s_2,s_3\} = \emptyset$. Further, we note by Claim~P.5 that every vertex that belongs to the edge $h$ has the same degree in $H$ and in $H'$, that is, $d_H(v) = d_{H'}(v)$ for all $v \in \{a,b,c,d\} = V(h)$. Recall that $d_{H}(c) = d_{H}(d) = 2$ and by Claim~P.3, $d_{H}(a) \ge 2$ and $d_{H}(b) \ge 2$.

\ClaimX{P.6}
$d_{H}(a) = 2$ and $d_{H}(b) = 2$.

\ProofClaimX{P.6}
Suppose, to the contrary, that $d_{H}(a) = 3$ or $d_{H}(b) = 3$. Interchanging the names of $a$ and $b$ if necessary, we may assume that $d_{H}(a) = 3$. Let $h$, $h_1$ and $h_2$ be the three edges incident with $a$. By Claim~M, $\defic(H-a) = 3$ or $\defic(H-a) = 8$. We consider the two possibilities in turn, and show that neither case can occur.

\ClaimX{P.6.1}
The case $\defic(H-a) = 3$ cannot occur.

\ProofClaimX{P.6.1}
Suppose that $\defic(H-a) = 3$. By Claim~M, $H - a$ therefore contains a double-$H_4$-component, say $C^*$. Let $g_1$ and $g_2$ be the two edges belonging to the $H_4$-pair in $C^*$ and let $h^*$ be the linking edge in $H-a$ that intersects $g_1$ and $g_2$. By definition of a double-$H_4$-component, we note that $g_i$ has three vertices of degree~$1$ and one vertex of degree~$2$ in $H-a$ for $i \in [2]$. By Claim~M, the edges $g_1$ and $g_2$ combined are intersected by all three edges $h, h_1,h_2$ incident with $a$. Analogously as in the proof of Claim~P.1, $2 \le |g_i \cap N_H(a)| \le 3$ for $i \in [2]$.
By Claim~P.1 and P.5, we note that the edges $g_1$ and $g_2$ are distinct from the edges $f_1$ and $f_2$. Further, we note that the vertices $c$ and $d$ both have degree~$1$ in $H-a$ and belong to the edges $f_1$ and $f_2$, respectively, in $H - a$. Thus, neither $g_1$ nor $g_2$ contain the vertex $c$ or $d$, while at least one of $g_1$ and $g_2$ contains the vertex $b$. Renaming $g_1$ and $g_2$ if necessary, we may assume that the vertex $b$ belongs to the edge $g_1$. Thus, $b$ has degree~$1$ in $H - a$.

Let $Z = f_1 \cup f_2 \cup g_1 \cup g_2 \cup \{x,a\}$ and let $H^* = H - Z$. We note that $n(H^*) = n(H) - 18$. The edges $e_1,e_2,e_3,f_1,f_2,g_1,g_2,h,h_1,h_2,h^*$ are deleted from $H$ when constructing $H^*$, and so $m(H^*) = m(H) - 11$.

We show firstly that $\defic(H^*) \le 21$. Suppose, to the contrary, that $\defic(H^*) \ge 22$. Let $Y^*$ be a special $H^*$-set satisfying $\defic_{H^*}(Y^*) = \defic(H^*)$. By Claim~O, $|E_{H^*}^*(Y^*)| \le |Y^*| - 3$. Let $Y_1,Y_2,Y_3,Y_4$ be the $H_4$-subhypergraph of $H$ with edge set $f_1,f_2,g_1,g_2$, respectively, and consider the special $H$-set $Y = Y^* \cup \{Y_1,Y_2,Y_3,Y_4\}$. We note that $|E_{H}^*(Y)| \le |E_{H^*}^*(Y^*)| + |\{e_1,e_2,e_3,h,h_1,h_2,h^*\}| \le (|Y^*| - 3) + 7 = |Y^*| + 4 = |Y|$, contradicting Claim~I(a). Therefore, $\defic(H^*) \le 21$.

We show next that every transversal in $H^*$ that contains a vertex in $\partial(Z)$ can be extended to a transversal in $H$ by adding to it five vertices. We note that if $\partial(Z)$ contains a vertex from some edge $e$ of $H$, then $e \in \{e_1,e_2,e_3,h_1,h_2,h^*\}$. Let $T^*$ be a transversal in $H^*$ that contains a vertex in $\partial(Z)$.
Suppose that $T^*$ contains a vertex in $e_1$, $e_2$ or $e_3$. Renaming edges if necessary, we may assume $T^*$ contains a vertex in $e_1$. By Claim~P.2, the edges $e_2$ and $e_3$ can be covered by two vertices one of which belongs to $f_1$ and the other to $f_2$. Thus, $e_2 \cap f_i \ne \emptyset$ and $e_3 \cap f_{3-i} \ne \emptyset$ for some $i \in [2]$. In this case, $T^* \cup \{(e_2f_i),(e_3f_{3-i}),a,(g_1h^*),(g_2h^*)\}$ is a transversal in $H$ of size~$|T^*| + 5$.
Suppose that $T^*$ contains a vertex in $h_1$ or $h_2$, say in $h_1$. We note that $h_2 \cap g_j \ne \emptyset$ for some $j \in [2]$. In this case, $T^* \cup \{x,c,d,(g_jh_2),(g_{3-j}h^*)\}$ is a transversal in $H$ of size~$|T^*| + 5$.
Suppose that $T^*$ contains a vertex in $h^*$. We note that $h_1 \cap g_i \ne \emptyset$ and $h_2 \cap g_{3-i} \ne \emptyset$ for some $i \in [2]$. In this case, $T^* \cup \{x,c,d,(g_ih_1),(g_{3-i}h_2)\}$ is a transversal in $H$ of size~$|T^*| + 5$.
In all cases, $T^*$ can be extended to a transversal in $H$ by adding to it five vertices.

For $i \in [2]$, let $\theta_i$ be a vertex in $h_i$ that does not belong to $g_1 \cup g_2$. Let $\theta_3$ and $\theta_4$ be the two vertices in $h^*$ that do not belong to $g_1$ or $g_2$. Analogously as in Claim~P.5, we note that $N_H(a) \cap h^* = \emptyset$. In particular, the vertices $\theta_1$, $\theta_2$, $\theta_3$ and $\theta_4$ are distinct, implying that $|\partial(Z)| \ge 4$. Recall that $\defic(H^*) \le 21$. By Claim~L, there exists a transversal, $T^*$, in $H^*$, that contains a vertex in $\partial(Z)$ such that $45 |T^*| \le 6n(H^*) + 13m(H^*) + f(|\partial(Z)|) \le 6n(H^*) + 13m(H^*) + 23$. As observed earlier, $T^*$ can be extended to a transversal in $H$ by adding to it five vertices. Hence, $45\tau(H) \le 45(|T^*| + 5) \le (6n(H^*) + 13m(H^*) + 23) + 225 = 6n(H) + 13m(H) - 6 \cdot 18 - 13 \cdot 11 + 23 + 225 = 6n(H) + 13m(H) - 3$, a contradiction. Hence the case $\defic(H-a) = 3$ cannot occur.~\smallqed

\ClaimX{P.6.2}
The case $\defic(H-a) = 8$ cannot occur.

\ProofClaimX{P.6.2}
Suppose that $\defic(H-a) = 8$. By Claim~M, $H - a$ therefore contains an $H_4$-component. Let $g$ be the edge belonging in this $H_4$-component. By Claim~M, the edge $g$ is intersected by all three edges $h, h_1,h_2$ incident with the vertex~$a$. We note that the edge $g$ is distinct from $f_1$ and $f_2$, and contains the vertex~$b$ which therefore has degree~$1$ in $H - a$. For $i \in [2]$, let $z_{i1}$ and $z_{i2}$ be the two vertices in $h_i \setminus \{a\}$ that do not belong to the edge $g$, and let $z_{i3}$ be the vertex $(h_ig)$ that is common to $h_i$ and $g$. Thus, $h_i = \{a,z_{i1},z_{i2},z_{i3}\}$.

Let $Z = f_1 \cup f_2 \cup g \cup \{x,a\}$ and let $H^* = H - Z$. We note that $n(H^*) = n(H) - 14$. The edges $e_1,e_2,e_3,f_1,f_2,g,h,h_1,h_2$ are deleted from $H$ when constructing $H^*$, and so $m(H^*) = m(H) - 9$.

We show firstly that $\defic(H^*) \le 21$. Suppose, to the contrary, that $\defic(H^*) \ge 22$. Let $Y^*$ be a special $H^*$-set satisfying $\defic_{H^*}(Y^*) = \defic(H^*)$. By Claim~O, $|E_{H^*}^*(Y^*)| \le |Y^*| - 3$. Let $R_1,R_2,R_3$ be the $H_4$-subhypergraphs of $H$ with edge set $f_1,f_2,g$, respectively, and consider the special $H$-set $Y = Y^* \cup \{R_1,R_2,R_3\}$. We note that $|E_{H}^*(Y)| \le |E_{H^*}^*(Y^*)| + |\{e_1,e_2,e_3,h_1,h_2\}| \le (|Y^*| - 3) + 6 = |Y^*| + 3 = |Y|$, contradicting Claim~I(a). Therefore, $\defic(H^*) \le 21$.

We show next that every transversal in $H^*$ that contains a vertex in $\partial(Z)$ can be extended to a transversal in $H$ by adding to it four vertices. We note that if $\partial(Z)$ contains a vertex from some edge $e$ of $H$, then $e \in \{e_1,e_2,e_3,h_1,h_2\}$. Let $T^*$ be a transversal in $H^*$ that contains a vertex in $\partial(Z)$.
Suppose that $T^*$ contains a vertex in $e_1$, $e_2$ or $e_3$. Renaming edges if necessary, we may assume $T^*$ contains a vertex in $e_1$. By Claim~P.2, $e_2 \cap f_i \ne \emptyset$ and $e_3 \cap f_{3-i} \ne \emptyset$ for some $i \in [2]$, implying that $T^* \cup \{(e_2f_i),(e_3f_{3-i}),a,b\}$ is a transversal in $H$ of size~$|T^*| + 4$.
If $T^*$ contains a vertex in $h_1$ or $h_2$, say in $h_1$, then  $T^* \cup \{(gh_2),x,c,d\}$ is a transversal in $H$ of size~$|T^*| + 4$.
In all cases, $T^*$ can be extended to a transversal in $H$ by adding to it four vertices.

Recall that $s_i$ is a vertex in $e_i \setminus \{x\}$ that does not belong to $f_1$ or $f_2$. Thus, $s_1,s_2,s_3$ are distinct vertices in $\partial(Z)$. Recall that $\defic(H^*) \le 21$. As observed earlier, $T^*$ can be extended to a transversal in $H$ by adding to it four vertices, and so $\tau(H) \le |T^*| + 4$. If $H^*$ contains two intersecting edges $e$ and $f$, such that $\partial(Z) \subseteq (V(e) \cup V(f)) \setminus(V(e) \cap V(f))$, then $|\partial(Z)| \le 6$. Hence if $|\partial(Z)| \ge 7$, then by Claim~L(b) there exists a transversal, $T^*$, in $H^*$, that contains a vertex in $\partial(Z)$ such that $45 |T^*| \le 6n(H^*) + 13m(H^*) + 21$. Thus, $45\tau(H) \le 45(|T^*| + 4) \le (6n(H^*) + 13m(H^*) + 21) + 180 = 6n(H) + 13m(H) - 6 \cdot 14 - 13 \cdot 9 + 21 + 180 = 6n(H) + 13m(H)$, a contradiction. Hence, $|\partial(Z)| \le 6$. Renaming vertices if necessary, we may assume that $s_1 = z_{11}$.

If $d_H(s_1) = 2$, then the vertex $s_1$ is isolated in $H^*$. In this case, adding $s_1$ to the set $Z$ we note that $n(H^*) = n(H) - 15$. Further $z_{12},z_{21},z_{22}$ are distinct vertices in $\partial(Z)$, and so $|\partial(Z)| \ge 3$ and $f(|\partial(Z)|) \le 27$. By Claim~L, there exists a transversal, $T^*$, in $H^*$, that contains a vertex in $\partial(Z)$ such that $45 |T^*| \le 6n(H^*) + 13m(H^*) + f(|\partial(Z)|) \le 6n(H) + 13m(H) + 27 - 6 \cdot 15 - 13 \cdot 9 = 6n(H) + 13m(H) - 180$, implying that $45\tau(H) \le 45(|T^*| + 4) \le 6n(H) + 13m(H)$, a contradiction. Hence, $d_H(s_1) = 3$.

Since $z_{11},z_{12},z_{21},z_{22}$ are distinct vertices in $\partial(Z)$, we note that $|\partial(Z)| \ge 4$. If $|\partial(Z)| = 4$, then $\{s_1,s_2,s_3\} \subset \{ z_{11},z_{12},z_{21},z_{22} \}$, implying analogously as above that $d_H(s_i) = 3$ for all $i \in [3]$, contradicting Claim~P.5. Therefore, $|\partial(Z)| \ge 5$, and so $f(|\partial(Z)|) = 22$. If the statement of Claim~L(b) holds, then $45 |T^*| \le 6n(H^*) + 13m(H^*) + 21$, implying as before that $45\tau(H) \le 45(|T^*| + 4) \le 6n(H) + 13m(H)$, a contradiction. Hence, $H^*$ contains two intersecting edges $e$ and $f$ having properties (i), (ii) and (iii) in Claim~L(b). Recall that $(ef)$ denotes the vertex in the intersection of $e$ and $f$. Thus, $\partial(Z) \subseteq (V(e) \cup V(f)) \setminus \{(ef)\}$, $e$ contains three degree-$1$ vertices, and $|\partial(Z) \cap V(e)|,|\partial(Z) \cap V(f)| \ge 2$.

Renaming vertices if necessary, we may assume that $s_1 = z_{11}$. If $z_{11} \in e$, let $z \in \partial(Z) \cap V(f)$, while if $z_{11} \in f$, let $z \in \partial(Z) \cap V(e)$. Let $Z^{\bullet} = Z \cup (e \setminus \{(ef)\}) \cup \{z\}$ and let $H^{\bullet} = H - Z^{\bullet}$. We note that $n(H^{\bullet}) = n(H) - 14 - 4 = n(H) - 18$. The edges $e,e_1,e_2,e_3,f,f_1,f_2,g,h,h_1,h_2$ are deleted from $H$ when constructing $H^{\bullet}$, and so $m(H^{\bullet}) = m(H) - 11$. By Claim~P.2, $e_2 \cap f_i \ne \emptyset$ and $e_3 \cap f_{3-i} \ne \emptyset$ for some $i \in [2]$. Every transversal in $H^{\bullet}$ can be extended to a transversal in $H$ by adding to it the five vertices in the set $\{b,(e_2f_i),(e_3f_{3-i}),s_1,z\}$, and so $\tau(H) \le \tau(H^{\bullet}) + 5$. Applying the inductive hypothesis to $H^{\bullet}$, we have that $45\tau(H) \le 45(\tau(H^{\bullet}) + 5) \le 6n(H^{\bullet}) + 13m(H^{\bullet}) + \defic(H^{\bullet}) + 225 \le 6n(H) + 13m(H) + \defic(H^{\bullet}) - 6 \times 18 - 13 \times 11 + 225 \le  6n(H) + 13m(H) + \defic(H^{\bullet}) - 26$. If $\defic(H^{\bullet}) \le 26$, then $45\tau(H) \le 6n(H) + 13m(H)$, a contradiction. Hence, $\defic(H^{\bullet}) \ge 27$. Let $Y^{\bullet}$ be a special $H^{\bullet}$-set satisfying $\defic(H^{\bullet}) = \defic_{H^{\bullet}}(Y^{\bullet})$. By Claim~O, $|E_{H^{\bullet}}^*(Y^{\bullet})| \le |Y^*| - 4$. Let $Y_1,Y_2,Y_3$ be the $H_4$-subhypergraphs of $H$ with edge set $f_1,f_2,g$, respectively, and consider the special $H$-set $Y = Y^{\bullet} \cup \{Y_1,Y_2,Y_3\}$. We note that $|E_{H}^*(Y)| \le |E_{H^{\bullet}}^*(Y^{\bullet})| + |\{e,e_1,e_2,e_3,f,h,h_1,h_2\}| \le (|Y^{\bullet}| - 4) + 8 = |Y^{\bullet}| + 4 = |Y|+1$, contradicting Claim~I(a). Hence the case $\defic(H-a) = 8$ cannot occur.~\smallqed

Since neither the case $\defic(H-a) = 3$ nor the case $\defic(H-a) = 8$ can occur, this completes the proof of Claim~P.6.~\smallqed

\medskip
Recall that $s_i$ is a vertex in $V(e_i) \setminus \{x\}$ that does not belong to $f_1$ or to $f_2$, and that by Claim~P.4, $d_H(s_i) \ge 2$.

\ClaimX{P.7}
If $d_H(s_i) = 3$ and $d_H(s_j) = 2$ for some $i$ and $j$ where $1 \le i,j \le 3$, then the vertices $s_i$ and $s_j$ do not belong to a common edge.

\ProofClaimX{P.7}
Renaming edges $e_1,e_2,e_3$ if necessary, let $d_H(s_1) = 3$ and $d_H(s_2) = 2$ and suppose, to the contrary, that there is an edge $g$ containing both $s_1$ and $s_2$. Let $g^{\bullet}$ be the third edge containing $s_1$ different from $e_1$ and $g$. By Claim~M, $\defic(H-s_1) = 3$ or $\defic(H-s_1) = 8$. We consider the two possibilities in turn, and show that neither case can occur.

\ClaimX{P.7.1}
The case $\defic(H-s_1) = 3$ cannot occur.

\ProofClaimX{P.7.1}
Suppose that $\defic(H-s_1) = 3$. By Claim~M, $H - s_1$ therefore contains a double-$H_4$-component, say $F$. Let $q_1$ and $q_2$ be the two edges belonging to the $H_4$-pair in $F$ and let $q^{\bullet}$ be the linking edge in $H-s_1$ that intersects $q_1$ and $q_2$. Recall that $V(h) = \{a,b,c,d\}$, where the edges $f_1$ and $h$ intersect in the vertex $c$ and where the edges $f_2$ and $h$ intersect in the vertex $d$. By Claim~P.6 and our earlier observations, $d_H(v) = 2$ for all $v \in V(h)$. Let $q^{\bullet} \cap q_1 = \{c_q\}$ and $q^{\bullet} \cap q_2 = \{d_q\}$.

Let $Z = (f_1 \cup f_2 \cup q_1 \cup q_2 \cup \{x,s_1,s_2\}) \setminus \{c,d,c_q,d_q\}$ and let $H^* = H - Z$. We note that $n(H^*) = n(H) - 15$. The edges $e_1,e_2,e_3,f_1,f_2,g,g^{\bullet},q_1,q_2$ are deleted from $H$ when constructing $H^*$, and so $m(H^*) = m(H) - 9$.

\ClaimX{P.7.1.1}
$\defic(H^*) \le 21$.

\ProofClaimX{P.7.1.1}
Suppose, to the contrary, that $\defic(H^*) \ge 22$. Let $Y^*$ be a special $H^*$-set satisfying $\defic_{H^*}(Y^*) = \defic(H^*)$. By Claim~O, $|E_{H^*}^*(Y^*)| \le |Y^*| - 3$. Let $F_1,F_2,F_3,F_4$ be the $H_4$-subhypergraphs of $H$ with edge set $f_1,f_2,q_1,q_2$, respectively. Let $Y^{**}$ be the special $H^*$-set obtained from $Y^*$ by removing all special subhypergraphs in $Y^*$, if any, that contain the edge $h$ or the edge $q^{\bullet}$, and consider the special $H$-set $Y = Y^{**} \cup \{F_1,F_2,F_3,F_4\}$.

Since we remove at most two  special subhypergraphs from $Y^*$ when constructing $Y$, we note that $|Y| \ge |Y^*| - 2$.  Further, we note that $|E_{H^*}^*(Y^{**})| \le |E_{H^*}^*(Y^*)| \le |Y^*| - 3$, and $|E_{H}^*(Y)| \le |E_{H^*}^*(Y^{**})| + |\{e_1,e_2,e_3,g,g^{\bullet},h,q^{\bullet}\}| \le (|Y^*| - 3) + 7 = |Y^*| + 4$. If $|Y| \ge |Y^*| - 1$ or if $|E_{H^*}^*(Y^*)| \le |Y^*| - 4$, then $|E_{H}^*(Y)| \le |Y| + 1$, contradicting Claim~I(a). Therefore, $|Y| = |Y^*| - 2$ and $|E_{H^*}^*(Y^*)| = |Y^*| - 3$. In particular, since $|E_{H^*}^*(Y^*)| \ge 0$, we note that $|Y^*| \ge 3$.
By Claim~O, $22 \le \defic(H^*) \le 13 \cdot 3 - 5|Y_4^*| - 8|Y_{14}^*| - 9|Y_{11}^*| - 12|Y_{21}^*|$, implying that $|Y^*| = |Y_4^*| = 3$ and  $|E_{H^*}^*(Y^*)| = 0$. Let $Y^* = \{R_1,R_2,R_3\}$, where $R_1$ and $R_2$ are the special subhypergraphs in $Y^*$ containing the edges $h$ and $q^{\bullet}$, respectively. We note that all three subhypergraph $R_1$, $R_2$ and $R_3$ are $H_4$-components in $H^*$. In particular, $E(R_1) = \{h\}$ and $E(R_2) = \{q^{\bullet}\}$.  Let $E(R_3) = \{\theta\}$.

If $q_1$ or $q_2$ intersects the edge $h$, then $h = q^{\bullet}$, a contradiction. Hence, neither $q_1$ nor $q_2$ intersect $h$. By Claim~P.6, $d_{H}(a) = 2$ and $d_{H}(b) = 2$, implying by the linearity of $H$ that both edges $g$ and $g^{\bullet}$ intersect $h$. Analogously, both edges $e_2$ and $e_3$ intersect the edge $q^{\bullet}$. Thus, no vertex in $N_H[x]$ belongs to the edge $\theta$, and no vertex in $N_H[s_1]$ belongs to the edge $\theta$. Thus, no edge in $H$ intersects the edge $\theta$, implying that $R_3$ is an $H_4$-component in $H$, a contradiction. This completes the proof of Claim~P.7.1.1.~\smallqed

\ClaimX{P.7.1.2}
Every transversal in $H^*$ that contains a vertex in $\partial(Z)$ can be extended to a transversal in $H$ by adding to it four vertices.

\ProofClaimX{P.7.1.2}
We note that if $\partial(Z)$ contains a vertex from some edge $e$ of $H$, then $e \in \{e_2,e_3,f_1,f_2,g,g^{\bullet},q_1,q_2\}$. Let $T^*$ be a transversal in $H^*$ that contains a vertex in $\partial(Z)$.
By Claim~P.2, any two edges amongst $e_1$, $e_2$ and $e_3$ (incident with $x$) can be covered by two vertices one of which belongs to $f_1$ and the other to $f_2$. Analogously, any two edges amongst $e_1$, $g$ and $g^{\bullet}$ (incident with $s_1$) can be covered by two vertices one of which belongs to $q_1$ and the other to $q_2$.
If $T^*$ contains a vertex in $f_1$ or $f_2$, say in $f_1$, then $T^* \cup \{(f_2h),(gq_i),(g^{\bullet}q_{3-i}),x\}$ is a transversal in $H$ of size~$|T^*| + 4$ for some $i \in [2]$.
If $T^*$ contains a vertex in $q_1$ or $q_2$, say in $q_1$, then $T^* \cup \{(e_2f_i),(e_3f_{3-i}),s_1,v_2\}$ is a transversal in $H$ of size~$|T^*| + 4$ for some $i \in [2]$, where $v_2$ is an arbitrary vertex in $q_2$.
If $T^*$ contains a vertex in $e_2$ or $e_3$, say in $e_2$, then $T^* \cup  \{(e_1f_i),(e_3f_{3-i}),(gq_j),(g^{\bullet}q_{3-j})\}$ is a transversal in $H$ of size~$|T^*| + 4$ for some $i \in [2]$ and $j \in [2]$.
If $T^*$ contains a vertex in $g$ or $g^{\bullet}$, say in $g$, then $T^* \cup \{(g^{\bullet}q_i), (e_1q_{3-i}),(gq_j), (e_2f_j),(e_3f_{3-j})\}$ is a transversal in $H$ of size~$|T^*| + 4$ for some $i \in [2]$ and $j \in [2]$.
In all cases, $T^*$ can be extended to a transversal in $H$ by adding to it four vertices. This completes the proof of Claim~P.7.1.2.~\smallqed

By Claim~P.5, the vertex $s_3$ does not belong to the edge $h$, implying that $\{a,b,s_3\} \subseteq \partial(Z)$, and so $|\partial(Z)| \ge 3$ and $f(|\partial(Z)|) \le 27$. By Claim~P.7.1.1, $\defic(H^*) \le 21$. By Claim~L, there exists a transversal, $T^*$, in $H^*$, that contains a vertex in $\partial(Z)$ such that $45 |T^*| \le 6n(H^*) + 13m(H^*) + f(|\partial(Z)|) \le 6n(H^*) + 13m(H^*) + 27$. By Claim~P.7.1.2, $T^*$ can be extended to a transversal in $H$ by adding to it four vertices. Hence, $45\tau(H) \le 45(|T^*| + 4) \le (6n(H^*) + 13m(H^*) + 27) + 180 = 6n(H) + 13m(H) - 6 \cdot 15 - 13 \cdot 9 + 27 + 180 = 6n(H) + 13m(H)$, a contradiction. Hence the case $\defic(H-s_1) = 3$ cannot occur.~\smallqed

\ClaimX{P.7.2}
The case $\defic(H-s_1) = 8$ cannot occur.

\ProofClaimX{P.7.2}
Suppose that $\defic(H-s_1) = 8$. By Claim~M, $H - a$ therefore contains an $H_4$-component. Let $q$ be the edge belonging in this $H_4$-component. By Claim~M, the edge $q$ is intersected by all three edges $e_1, g,g^{\bullet}$ incident with the vertex~$s_1$. Let $Z = f_1 \cup f_2 \cup q \cup \{x,s_1,s_2\}$ and let $H^* = H - Z$. We note that $n(H^*) = n(H) - 15$. The edges $e_1,e_2,e_3,f_1,f_2,g,g^{\bullet},h,q$ are deleted from $H$ when constructing $H^*$, and so $m(H^*) = m(H) - 9$.

We show firstly that $\defic(H^*) \le 21$. Suppose, to the contrary, that $\defic(H^*) \ge 22$. Let $Y^*$ be a special $H^*$-set satisfying $\defic_{H^*}(Y^*) = \defic(H^*)$. By Claim~O, $|E_{H^*}^*(Y^*)| \le |Y^*| - 3$. Let $F_1,F_2,F_3$ be the $H_4$-subhypergraphs of $H$ with edge set $f_1,f_2,q$, respectively, and consider the special $H$-set $Y = Y^* \cup \{F_1,F_2,F_3\}$. We note that $|E_{H}^*(Y)| \le |E_{H^*}^*(Y^*)| + |\{e_1,e_2,e_3,g,g^{\bullet},h\}| \le (|Y^*| - 3) + 6 = |Y^*| + 3 = |Y|$, contradicting Claim~I(a). Therefore, $\defic(H^*) \le 21$.

We show next that every transversal in $H^*$ that contains a vertex in $\partial(Z)$ can be extended to a transversal in $H$ by adding to it four vertices. We note that if $\partial(Z)$ contains a vertex from some edge $e$ of $H$, then $e \in \{e_2,e_3,g,g^{\bullet},h\}$. Let $T^*$ be a transversal in $H^*$ that contains a vertex in $\partial(Z)$.
If $T^*$ contains a vertex in $e_2$ or $e_3$, say in $e_2$, then $T^* \cup \{(e_3f_i),(hf_{3-i}),s_1,w\}$ is a transversal in $H$ of size~$|T^*| + 4$ for some $i \in [2]$ where $w$ is an arbitrary vertex in the edge $q$. %
If $T^*$ contains a vertex in $h$, then $T^* \cup \{(e_2f_i),(e_3f_{3-i}),s_1,w\}$ is a transversal in $H$ of size~$|T^*| + 4$ for some $i \in [2]$ where $w$ is an arbitrary vertex in the edge $q$.
If $T^*$ contains a vertex in $g$ or $g^{\bullet}$, say in $g$, then $T^* \cup \{(g^{\bullet}q),c,d,x\}$ is a transversal in $H$ of size~$|T^*| + 4$.
In all cases, $T^*$ can be extended to a transversal in $H$ by adding to it four vertices.

By Claim~P.5, the vertex $s_3$ does not belong to the edge $h$, implying that $\{a,b,s_3\} \subseteq \partial(Z)$, and so $|\partial(Z)| \ge 3$ and $f(|\partial(Z)|) \le 27$. An identical argument as in the last paragraph of the proof of Claim~P.7.1 yields the contradiction $45\tau(H) \le 6n(H) + 13m(H)$. Hence the case $\defic(H-s_1) = 8$ cannot occur.~\smallqed

Since neither the case $\defic(H-s_1) = 3$ not the case $\defic(H-s_1) = 8$ can occur, this completes the proof of Claim~P.7.~\smallqed

\ClaimX{P.8}
There is no edge containing all of $s_1$, $s_2$ and $s_3$.

\ProofClaimX{P.8}
Suppose, to the contrary, that there is an edge, $g$ say, in $H$ containing all of $s_1$, $s_2$ and $s_3$. By Claim~P.4, $d_H(s_i) \ge 2$ for all $i \in [3]$, and $d_H(s_i) = 2$ for at least one $i \in [3]$. Thus, by Claim~P.7, $d_H(s_1) = d_H(s_2) = d_H(s_3) = 2$. Let $Z = f_1 \cup f_2 \cup \{x,s_1,s_2,s_3\}$ and let $H^* = H - Z$. We note that $n(H^*) = n(H) - 12$. The edges $e_1,e_2,e_3,f_1,f_2,g,h$ are deleted from $H$ when constructing $H^*$, and so $m(H^*) = m(H) - 7$.

If $\defic(H^*) \ge 22$, then let $Y^*$ be a special $H^*$-set satisfying $\defic_{H^*}(Y^*) = \defic(H^*)$. By Claim~O, $|E_{H^*}^*(Y^*)| \le |Y^*| - 3$. Let $F_1,F_2$ be the $H_4$-subhypergraphs of $H$ with edge set $f_1,f_2$, respectively, and consider the special $H$-set $Y = Y^* \cup \{F_1,F_2\}$. We note that $|E_{H}^*(Y)| \le |E_{H^*}^*(Y^*)| + |\{e_1,e_2,e_3,g,h\}| \le (|Y^*| - 3) + 5 = |Y^*| + 2 = |Y|$, contradicting Claim~I(a). Therefore, $\defic(H^*) \le 21$.

We show next that every transversal in $H^*$ that contains a vertex in $\partial(Z)$ can be extended to a transversal in $H$ by adding to it three vertices. We note that if $\partial(Z)$ contains a vertex from some edge $e$ of $H$, then $e \in \{e_1,e_2,e_3,g,h\}$. Let $T^*$ be a transversal in $H^*$ that contains a vertex in $\partial(Z)$.
If $T^*$ contains a vertex in $e_1$, $e_2$ or $e_3$, say in $e_1$, then $T^* \cup \{s_2,(e_3f_i),(hf_{3-i})\}$ is a transversal in $H$ of size~$|T^*| + 3$ for some $i \in [2]$.
If $T^*$ contains a vertex in $h$, then $T^* \cup \{s_1,(e_2f_i),(e_3f_{3-i})\}$ is a transversal in $H$ of size~$|T^*| + 3$ for some $i \in [2]$.
If $T^*$ contains a vertex in $g$, then $T^* \cup \{c,d,x\}$ is a transversal in $H$ of size~$|T^*| + 3$.
In all cases, $T^*$ can be extended to a transversal in $H$ by adding to it three vertices.

As observed earlier, $\defic(H^*) \le 21$. By Claim~L, there exists a transversal, $T^*$, in $H^*$, that contains a vertex in $\partial(Z)$ such that $45 |T^*| \le 6n(H^*) + 13m(H^*) + f(|\partial(Z)|)$. As observed earlier, $T^*$ can be extended to a transversal in $H$ by adding to it three vertices. Hence, $45\tau(H) \le 45(|T^*| + 3) \le (6n(H^*) + 13m(H^*) + f(|\partial(Z)|)) + 135 = 6n(H) + 13m(H) + f(|\partial(Z)|) - 6 \cdot 12 - 13 \cdot 7 + 135 = 6n(H) + 13m(H) + f(|\partial(Z)|)) - 28$. If $|\partial(Z)| \ge 3$, then $f(|\partial(Z)|) \le 27$ and $45\tau(H) \le 6n(H) + 13m(H)$, a contradiction. Hence, $|\partial(Z)| \le 2$. Since $\{a,b\} \subseteq \partial(Z)$, we note that $|\partial(Z)| \ge 2$. Consequently, $|\partial(Z)| = 2$ and $\partial(Z) = \{a,b\}$. This implies that $f_i$ intersects each of the edges $e_1,e_2,e_3$ for $i \in [2]$ and that either $a$ or $b$ belong to the edge $g$. Renaming the vertices $a$ and $b$ if necessary, we may assume that $g = \{a,s_1,s_2,s_3\}$. This implies that the only edges in $H$ that intersect the edges $f_1$, $f_2$ and $g$ are $e_1,e_2,e_3,h$. We now let $F_1,F_2,F_3$ be the $H_4$-subhypergraphs of $H$ with edge set $f_1,f_2,g$, respectively, and consider the special $H$-set $Y = \{F_1,F_2,F_3\}$. We note that $|E_{H}^*(Y)| \le |\{e_1,e_2,e_3,h\}| = 4 = |Y|+ 1$, contradicting Claim~I(a). This completes the proof of Claim~P.8.~\smallqed

\medskip
\ClaimX{P.9}
If $d_H(s_i) = d_H(s_j) = 3$ for some $i$ and $j$ where $1 \le i,j \le 3$ and $i \ne j$, then the vertices $s_i$ and $s_j$ do not belong to a common edge.

\ProofClaimX{P.9}
Renaming edges $e_1,e_2,e_3$ if necessary, let $d_H(s_1) = d_H(s_2) = 3$ and suppose, to the contrary, that there is an edge $g$ containing both $s_1$ and $s_2$. By Claim~P.8, the edge $g$ does not contain the vertex~$s_3$. By Claim~P.4, we may assume, renaming the edges $f_1$ and $f_2$ if necessary, that $e_1 \cap f_1 \ne \emptyset$, $e_1 \cap f_2 = \emptyset$, and that $e_2 \cap f_1 = \emptyset$ and $e_2 \cap f_2 \ne \emptyset$. Thus, $e_1$ intersects $f_1$ but not $f_2$, and $e_2$ intersects $f_2$ but not $f_1$. By Claim~P.1, this implies that $e_3$ intersects both $f_1$ and $f_2$; that is, $e_3 \cap f_1 \ne \emptyset$ and $e_3 \cap f_2 \ne \emptyset$. For $i \in [2]$, let $w_i$ be the vertex in $e_i$ different from $(e_if_i)$, $x$ and $s_i$, and so $V(e_i) = \{x,s_i,w_i,(e_if_i)\}$. We note that the vertex $(e_if_i)$ has degree~$2$ in $H$ for $i \in [2]$.

\ClaimX{P.9.1}
There is no edge containing both $w_1$ and $w_2$.

\ProofClaimX{P.9.1}
Suppose, to the contrary, that there is an edge $f'$ that contains both $w_1$ and $w_2$. We now consider the hypergraph $H - s_1$. By Claim~M, either $f'$ belongs to an $H_4$-component or to a double-$H_4$-component. Since $f'$ is intersected by the edge $e_2$ in $H - s_1$, the edge $f'$ cannot belong to an $H_4$-component. Hence, $f'$ belongs to a double-$H_4$-component. Let $f''$ be the second edge belonging to the $H_4$-pair in this double-$H_4$-component. We note that $e_2$ is the linking edge in $H - s_1$ that intersects $f'$ and $f''$. Analogously as in Claim~P.5, $N_H(s_1) \cap V(e_2) = \emptyset$. However, $s_2 \in N_H(s_1) \cap V(e_2)$, a contradiction.~\smallqed

\ClaimX{P.9.2}
There is no edge containing $s_3$ and one of $w_1$ or $w_2$.

\ProofClaimX{P.9.2}
Suppose, to the contrary, that there is an edge $f'$ that contains $s_3$ and one of $w_1$ or $w_2$. By symmetry, we may assume that $w_1 \in V(f')$. We now consider the hypergraph $H - s_1$. By Claim~M, either $f'$ belongs to an $H_4$-component or to a double-$H_4$-component. Since $f'$ is intersected by the edge $e_3$ in $H - s_1$, the edge $f'$ cannot belong to an $H_4$-component. Hence, $f'$ belongs to a double-$H_4$-component. Let $f''$ be the second edge belonging to the $H_4$-pair in this double-$H_4$-component. We note that $e_3$ is the linking edge in $H - s_1$ that intersects $f'$ and $f''$. Analogously as in Claim~P.5, $N_H(s_1) \cap V(e_3) = \emptyset$. However, $x \in N_H(s_1) \cap V(e_3)$, a contradiction.~\smallqed

\medskip
By Claim~P.9.1 and~P.9.2, $\{w_1,w_2,s_3\}$ is an independent set. Recall that $h = \{a,b,c,d\}$, where the edges $f_1$ and $h$ intersect in the vertex $c$ and where the edges $f_2$ and $h$ intersect in the vertex $d$. Further, $d_H(c) = d_{H'}(c) = 2$ and $d_H(d) = d_{H'}(d) = 2$. In particular, we note that $\{d,w_1,w_2,s_3\}$ is an independent set. Let $Z = (f_1 \cup f_2 \cup \{x\}) \setminus \{c,d\}$ and let $H^{\bullet}$ be obtained from $H - Z$ by adding to it the edge $f^{\bullet} = \{d,w_1,w_2,s_3\}$. We note that $n(H^{\bullet}) = n(H) - 7$ and $m(H^{\bullet}) = m(H) - 4$. Further since $\{d,w_1,w_2,s_3\}$ is an independent set in $H$, the hypergraph $H^{\bullet}$ is linear.

We show next that $\tau(H) \le \tau(H^{\bullet}) + 2$. Let $T^{\bullet}$ be a $\tau(H^{\bullet})$-transversal. In order to cover the edge $f^{\bullet}$, the transversal $T^{\bullet}$ contains at least one vertex in $f^{\bullet}$. If $w_1 \in T^{\bullet}$, let  $T = T^{\bullet} \cup \{(e_2f_2),(e_3f_1)\}$. If $w_2 \in T^{\bullet}$, let $T = T^{\bullet} \cup \{(e_1f_1),(e_3f_{2})\}$. If $s_3 \in T^{\bullet}$, let $T = T^{\bullet} \cup \{(e_1f_1),(e_2f_2)\}$. If $d \in T^{\bullet}$, let $T = T^{\bullet} \cup \{x,c\}$. In all cases, $T$ is a transversal in $H$ of size~$|T^{\bullet}| + 2$, and so $\tau(H) \le |T^{\bullet}| + 2 = \tau(H^{\bullet}) + 2$.

Applying the inductive hypothesis to $H^{\bullet}$, we have that $45\tau(H) \le 45(\tau(H^{\bullet}) + 2) \le 6n(H^{\bullet}) + 13m(H^{\bullet}) + \defic(H^{\bullet}) + 90 \le 6n(H) + 13m(H) + \defic(H^{\bullet}) - 7 \times 6 - 4 \times 13 + 90 =  6n(H) + 13m(H) + \defic(H^{\bullet}) - 4$. If $\defic(H^{\bullet}) \le 4$, then $45\tau(H) \le 6n(H) + 13m(H)$, a contradiction. Hence, $\defic(H^{\bullet}) \ge 5$. Let $Y^{\bullet}$ be a special $H^{\bullet}$-set satisfying $\defic_{H^{\bullet}}(Y^{\bullet}) = \defic(H^{\bullet})$. By Claim~O, $|E_{H^{\bullet}}^*(Y^{\bullet})| \le |Y^{\bullet}| - 1$.

\ClaimX{P.9.3}
$f^{\bullet} \in Y^{\bullet}$.

\ProofClaimX{P.9.3}
Suppose, to the contrary, that $f^{\bullet} \notin Y^{\bullet}$. If $f^{\bullet} \notin E_{H^{\bullet}}^*(Y^{\bullet})$, then $h \notin Y^{\bullet}$, implying that $|E_{H}^*(Y^{\bullet})| \le |E_{H^{\bullet}}^*(Y^{\bullet})| + |\{e_1,e_2\}| \le (|Y^{\bullet}| - 1) + 2 = |Y^{\bullet}| + 1$, contradicting Claim~I. Hence, $f^{\bullet} \in E_{H^{\bullet}}^*(Y^{\bullet})$. Since $H - Z = H^{\bullet} - f^{\bullet}$, this implies that $\defic(H-Z) = \defic(H^{\bullet} - f^{\bullet}) = \defic(H^{\bullet}) + 13 \ge 18$. Let $Y$ be a special $(H-Z)$-set satisfying $\defic_{H-Z}(Y) = \defic(H-Z)$.
If $h \notin Y$, then neither $f_1$ nor $f_2$ intersect any special subhypergraphs in $Y$, implying that $\defic_{H-x}(Y) = \defic_{H-Z}(Y) \ge 18$, contradicting Claim~M. Hence, $h \in Y$. Let $R_1$ be the special subhypergraph in $Y$ that contains the edge $h$. We note that both $c$ and $d$ have degree~$1$ in $H - Z$, implying that the subhypergraph $R_1$ is an $H_4$-component in $H-Z$ and therefore contributes~$8$ to the deficiency of $Y$ in $H - Z$. Since neither $f_1$ nor $f_2$ intersect any special subhypergraph in $Y$ different from $R_1$, this implies that $\defic_{H-x}(Y \setminus \{R_1\}) \ge 18 - 8 = 10$, once again contradicting Claim~M.~\smallqed

\medskip
By Claim~P.9.3, $f^{\bullet} \in Y^{\bullet}$. Let $R_1$ be the special subhypergraph in $Y^{\bullet}$ that contains the edge $f^{\bullet}$. Suppose that $|Y^{\bullet}| \ge 2$. Let $|E_{H^{\bullet}}^*(Y^{\bullet})| = |Y^{\bullet}| - i$ for some $i \ge 1$. By Claim~O and by our earlier observations, $5 \le \defic(H^{\bullet}) \le 13i - 5|Y^{\bullet}| \le 13i - 10$, implying that $i \ge 2$; that is, $|E_{H^{\bullet}}^*(Y^{\bullet})| \le |Y^{\bullet}| - 2$. Let $R_2$ and $R_3$ be the $H_4$-subhypergraphs of $H$ with edge set $f_1$ and $f_2$, respectively. Let $Y^* = (Y^{\bullet} \setminus \{R_1\}) \cup \{R_2,R_3\}$, and so $|Y^*| = |Y^{\bullet}| + 1$. Thus, $|E_{H}^*(Y^*)| \le  |E_{H^{\bullet}}^*(Y^{\bullet})| + |\{e_1,e_2,e_3,h\}| \le (|Y^{\bullet}| - 2) + 4 = |Y^*| + 1$, contradicting Claim~I. Hence, $|Y^{\bullet}| = 1$; that is, $Y^{\bullet} = \{R_1\}$.

Recall that $\defic(H^{\bullet}) \ge 5$ and $|E_{H^{\bullet}}^*(Y^{\bullet})| \le |Y^{\bullet}| - 1$. Thus since $|Y^{\bullet}| = 1$, we note that $|E_{H^{\bullet}}^*(Y^{\bullet})| = 0$, implying that $R_1$ is a component in $H^{\bullet}$. Since the vertex $s_3$ belongs to $R_1$ and is contained in both edges $f^{\bullet}$ and $h$, the component $R_1$ is not an $H_4$-component in $H^{\bullet}$, implying that $R_1$ is an $H_{14}$-component in $H^{\bullet}$ and $\defic(H^{\bullet}) = 5$.

We now consider the hypergraph $H - s_1$. By Claim~M, either $H - s_1$ contains an $H_4$-component or a double-$H_4$-component. We note that in either case, the vertex $w_1$ belongs to such a component. If $H - s_1$ contains an $H_4$-component, then this component contains a vertex of degree~$1$ in $H$ different from $w_1$. Suppose that $H - s_1$ contains a double-$H_4$-component. In this case, let $q_1$ and $q_2$ be the two edges belonging to the $H_4$-pair in this component and let $q^{\bullet}$ be the linking edge in $H - s_1$ that intersects $q_1$ and $q_2$. Renaming $q_1$ and $q_2$ if necessary, we may assume that $w_1 \in V(q_1)$. Thus, $q_1 \cap V(e_1) = \emptyset$, implying that $q_2$ contains a vertex of degree~$1$ in $H$. In both cases, there is a vertex of degree~$1$ in $H$ that belongs to the component of $H - s_1$ that contains $w_1$, implying that the component $R_1$ of $H^{\bullet}$ contains a vertex of degree~$1$ in $H$. Further, the component $R_1$ of $H^{\bullet}$ contains the vertex $c$ which has degree~$1$ in $H^{\bullet}$ and degree~$2$ in $H$. Thus, $R_1$ has at least two vertices of degree~$1$ in $H^{\bullet}$, implying that $R_1$ is an $H_{14,4}$-component. We now consider the edge $h$ in $R_1$ that contains the vertex~$c$ of degree~$1$ in $H^{\bullet}$. All three edges intersecting the edge $h$ in $R_1 \cong H_{14,4}$ have at least one vertex of degree~$3$. Since $f^{\bullet}$ is one of the three edges that intersect $h$ in $R_1$, this implies in particular that the edge $f^{\bullet}$ contains at least one vertex of degree~$3$ in $R_1$. However, no vertex of $f^{\bullet}$ has degree~$3$ in $H^{\bullet}$, a contradiction. This completes the proof of Claim~P.9.~\smallqed

\ClaimX{P.10}
Every edge contains at most one vertex from the set $\{s_1,s_2,s_3\}$.

\ProofClaimX{P.10}
Suppose, to the contrary, that there is an edge, $g$ say, in $H$ containing $s_i$ and $s_j$ where $1 \le i < j \le 3$. Renaming the edges $e_1, e_2, e_3$ if necessary, we may assume that $\{s_1,s_2\} \subset V(g)$. Let $g = \{s_1,s_2,u,v\}$. By Claim~P.8, the edge $g$ does not contain the vertex $s_3$, and so $s_3 \notin \{u,v\}$. By Claim~P.4, Claim~P.7 and Claim~P.9, $d_H(s_1) = d_H(s_2) = 2$. Let $Z = f_1 \cup f_2 \cup \{s_1,s_2,x\}$ and let $H^* = H - Z$. We note that $n(H^*) = n(H) - 11$. The edges $e_1,e_2,e_3,f_1,f_2,g,h$ are deleted from $H$ when constructing $H^*$, and so $m(H^*) = m(H) - 7$.

We show firstly that $\defic(H^*) \le 21$. Suppose, to the contrary, that $\defic(H^*) \ge 22$. Let $Y^*$ be a special $H^*$-set satisfying $\defic_{H^*}(Y^*) = \defic(H^*)$. By Claim~O, $|E_{H^*}^*(Y^*)| \le |Y^*| - 3$. Let $R_1$ and $R_2$ be the $H_4$-subhypergraphs of $H$ with edges $f_1$ and $f_2$, respectively, and consider the special $H$-set $Y = Y^* \cup \{R_1,R_2\}$. We note that $|E_{H}^*(Y)| \le |E_{H^*}^*(Y^*)| + |\{e_1,e_2,e_3,g,h\}| \le (|Y^*| - 3) + 5 = |Y^*| + 2 = |Y|$, contradicting Claim~I(a). Therefore, $\defic(H^*) \le 21$.

We show next that every transversal in $H^*$ that contains a vertex in $\partial(Z)$ can be extended to a transversal in $H$ by adding to it three vertices. We note that if $\partial(Z)$ contains a vertex from some edge $e$ of $H$, then $e \in \{e_1,e_2,e_3,g,h\}$. Let $T^*$ be a transversal in $H^*$ that contains a vertex in $\partial(Z)$.
If $T^*$ contains a vertex in $e_1$ or in $e_2$, say in $e_1$ by symmetry, then $T^* \cup \{(e_2g),(e_3f_{i}),(hf_{3-i})\}$ is a transversal in $H$ for some $i \in [2]$.
If $T^*$ contains a vertex in $g$, then $T^* \cup \{c,d,x\}$ is a transversal in $H$.
If $T^*$ contains a vertex in $h$, then $T^* \cup \{(e_1g),(e_2f_{i}),(e_3f_{3-i})\}$ is a transversal in $H$ for some $i \in [2]$.
In all cases, $T^*$ can be extended to a transversal in $H$ by adding to it three vertices.

By Claim~L, there exists a transversal, $T^*$, in $H^*$, that contains a vertex in $\partial(Z)$ such that $45 |T^*| \le 6n(H^*) + 13m(H^*) + f(|\partial(Z)|) = 6n(H) + 13m(H) + f(|\partial(Z)|) - 11 \cdot 6 - 7 \cdot 13 = 6n(H) + 13m(H) + f(|\partial(Z)|) - 157$. As observed earlier, $T^*$ can be extended to a transversal in $H$ by adding to it three vertices. Hence, $45\tau(H) \le 45(|T^*| + 3) \le 6n(H) + 13m(H) + f(|\partial(Z)|) - 22$. If $|\partial(Z)| \ge 5$, then $f(|\partial(Z)|) = 22$, implying that $45\tau(H) \le 6n(H) + 13m(H)$, a contradiction. Hence, $|\partial(Z)| \le 4$.

We note that $\{a,b\} \subseteq \partial(Z)$, $\{u,v\} \subseteq \partial(Z)$, and $s_3 \in \partial(Z)$. By Claim~P.5, the edge $h$ contains no neighbor of $x$. In particular, $s_3 \notin \{a,b\}$, and so $\{a,b,s_3\} \subseteq \partial(Z)$. As observed earlier, $s_3 \notin \{u,v\}$. By the linearity of $H$, we note that $\{a,b\} \ne \{u,v\}$. Since $|\partial(Z)| \le 4$, this implies that $|\{a,b\} \cap \{u,v\}| = 1$ and $|\partial(Z)| = 4$. Renaming vertices if necessary, we may assume that $b = v$. By Claim~P.6, $d_H(b) = 2$. Since $b = v$, the vertex $b$ is therefore incident only with the edge $g$ and $h$ in $H$, implying that $d_{H^*}(b) = 0$.

We now add the vertex $b$ to the set $Z$. With this addition of $b$ to $Z$, we note that $n(H^*) = n(H) - 12$, $m(H^*) = m(H) - 7$, and $45 |T^*| \le 6n(H) + 13m(H) + f(|\partial(Z)|) - 157 - 6 = 6n(H) + 13m(H) + f(|\partial(Z)|) - 163$. Further since $\{a,s_3,u\} \subseteq \partial(Z)$, we note that $|\partial(Z)| \ge 3$ and $f(|\partial(Z)|) \le 27$. Thus, $45\tau(H) \le 45(|T^*| + 3) \le 6n(H) + 13m(H) + f(|\partial(Z)|) - 28 < 6n(H) + 13m(H)$, a contradiction. This completes the proof of Claim~P.10.~\smallqed

\medskip
We now continue with our proof of Claim~P. By Claim~P.10, $\{s_1,s_2,s_3\}$ is an independent set. Recall that $h = \{a,b,c,d\}$, where the edges $f_1$ and $h$ intersect in the vertex $c$ and where the edges $f_2$ and $h$ intersect in the vertex $d$. Further, $d_H(c) = d_{H'}(c) = 2$ and $d_H(d) = d_{H'}(d) = 2$. In particular, we note that $\{d,s_1,s_2,s_3\}$ is an independent set. Let $Z = (f_1 \cup f_2 \cup \{x\}) \setminus \{c,d\}$ and let $H^{\bullet}$ be obtained from $H - Z$ by adding to it the edge $f^{\bullet} = \{d,s_1,s_2,s_3\}$. We note that $n(H^{\bullet}) = n(H) - 7$ and $m(H^{\bullet}) = m(H) - 4$. Further since $\{d,s_1,s_2,s_3\}$ is an independent set in $H$, the hypergraph $H^{\bullet}$ is linear. By Claim~P.1, at least one of the edges $e_1,e_2,e_3$ intersects both $f_1$ and $f_2$. Renaming edges if necessary, we may assume that $e_3$ intersects both $f_1$ and $f_2$; that is, $e_3 \cap f_1 \ne \emptyset$ and $e_3 \cap f_2 \ne \emptyset$. Let $Y^{\bullet}$ be a special $H^{\bullet}$-set satisfying $\defic_{H^{\bullet}}(Y^{\bullet}) = \defic(H^{\bullet})$.

\ClaimX{P.11}
The following properties hold in the hypergraph $H^{\bullet}$. \\
\hspace*{0.75cm} {\rm (a)} $\tau(H) \le \tau(H^{\bullet}) + 2$. \\
\hspace*{0.75cm}  {\rm (b)} $\defic(H^{\bullet}) \ge 5$ and $|E_{H^{\bullet}}^*(Y^{\bullet})| \le |Y^{\bullet}| - 1$. \\
\hspace*{0.75cm}  {\rm (c)} There is no $H_{10}$-subhypergraph in $H^{\bullet}$. \\
\hspace*{0.75cm}  {\rm (d)} $f^{\bullet} \in Y^{\bullet}$.
\\
\hspace*{0.75cm}  {\rm (e)} $|Y^{\bullet}| = 1$.

\ProofClaimX{P.11}
(a) Let $T^{\bullet}$ be a $\tau(H^{\bullet})$-transversal. In order to cover the edge $f^{\bullet}$, the transversal $T^{\bullet}$ contains at least one vertex in $f^{\bullet}$. Suppose that $s_j \in T^{\bullet}$ for some $j \in [3]$. Renaming edges $e_1,e_2,e_3$ if necessary, we may assume that $s_1 \in T^{\bullet}$. In this case, $T^{\bullet} \cup \{(e_2f_i),(e_3f_{3-1})\}$ is a transversal of $H$ for some $i \in [2]$. If $d \in T^{\bullet}$, then $T^{\bullet} \cup \{x,c\}$ is a transversal of $H$. In both cases, we produce a transversal in $H$ of size~$|T^{\bullet}| + 2$, and so $\tau(H) \le |T^{\bullet}| + 2 = \tau(H^{\bullet}) + 2$.

(b) Applying the inductive hypothesis to $H^{\bullet}$, we have that $45\tau(H) \le 45(\tau(H^{\bullet}) + 2) \le 6n(H^{\bullet}) + 13m(H^{\bullet}) + \defic(H^{\bullet}) + 90 \le 6n(H) + 13m(H) + \defic(H^{\bullet}) - 7 \times 6 - 4 \times 13 + 90 =  6n(H) + 13m(H) + \defic(H^{\bullet}) - 4$. If $\defic(H^{\bullet}) \le 4$, then $45\tau(H) \le 6n(H) + 13m(H)$, a contradiction. Hence, $\defic(H^{\bullet}) \ge 5$. By Claim~O, $|E_{H^{\bullet}}^*(Y^{\bullet})| \le |Y^{\bullet}| - 1$.

(c) Suppose, to the contrary, that there is a $H_{10}$-subhypergraph, say $F$, in $H^{\bullet}$. By Claim~J, $F$ is not a subhypergraph of $H$, implying that the added edge $f^{\bullet}$ is an edge of $F$ and therefore the vertex $d$ belongs to $F$. We note that every vertex of $F$ has degree~$2$ in $F$ and therefore degree at least~$2$ in $H^{\bullet}$. Since the vertex $d$ has degree~$2$ in $H^{\bullet}$ and is contained in the edges $f^{\bullet}$ and $h$ in $H^{\bullet}$, the edge $h$ must belong to $F$. This implies that the vertex $c$ belongs to $F$. However the vertex $c$ has degree~$1$ in $H^{\bullet}$, and therefore degree~$1$ in $F$, a contradiction.

(d) Suppose, to the contrary, that $f^{\bullet} \notin Y^{\bullet}$. Recall that the edge $e_3$ intersects both $f_1$ and $f_2$. If $f^{\bullet} \notin E_{H^{\bullet}}^*(Y^{\bullet})$, then $h \notin Y^{\bullet}$, implying that $|E_{H}^*(Y^{\bullet})| \le |E_{H^{\bullet}}^*(Y^{\bullet})| + |\{e_1,e_2\}| \le (|Y^{\bullet}| - 1) + 2 = |Y^{\bullet}| + 1$, contradicting Claim~I. Hence, $f^{\bullet} \in E_{H^{\bullet}}^*(Y^{\bullet})$.
If $h \notin Y$, then neither $f_1$ nor $f_2$ intersect any special subhypergraphs in $^{\bullet}$, implying that $|E_{H}^*(Y^{\bullet})| \le |E_{H^{\bullet}}^*(Y^{\bullet})| + |\{e_1,e_2,e_3\}| - |\{f^{\bullet}| \le (|Y^{\bullet}| - 1) + 3 - 1 = |Y^{\bullet}| + 1$, contradicting Claim~I. Hence, $h \in Y$.
Let $R_1$ be the special subhypergraph in $Y$ that contains the edge $h$. We note that both $c$ and $d$ have degree~$1$ in $R_1$, implying that the subhypergraph $R_1$ is an $H_4$-component in $H^{\bullet}$. By Claim~P.5 and Claim~P.6, both vertices $a$ and $b$ have degree~$2$ in $H^{\bullet}$. By the linearity of $H$, the edge different from $h$ that contains $a$ and the edge different from $h$ that contains $b$ are distinct. These two edges, together with the edge $f^{\bullet}$, all intersect the $H_4$-component $R_1$, implying that $|E_{H}^*(Y^{\bullet})| \ge 3$. Thus, by Claim~O, $|E_{H}^*(Y^{\bullet})| \le |Y^{\bullet}| - 3$. Hence, $|E_{H}^*(Y^{\bullet})| \le |E_{H^{\bullet}}^*(Y^{\bullet})| + |\{e_1,e_2,e_3,f_1,f_2\}| - |\{f^{\bullet}| \le (|Y^{\bullet}| - 3) + 5 - 1 = |Y^{\bullet}| + 1$, contradicting Claim~I.

(e) By Part~(d), $f^{\bullet} \in Y^{\bullet}$. Let $R_1$ be the special subhypergraph in $Y^{\bullet}$ that contains the edge $f^{\bullet}$. If $|Y^{\bullet}| \ge 3$, then by Claim~O and by our earlier observations, $|E_{H^{\bullet}}^*(Y^{\bullet})| \le |Y^{\bullet}| - 2$.  Let $R_2$ and $R_3$ be the $H_4$-subhypergraphs of $H$ with edge set $f_1$ and $f_2$, respectively. Let $Y^* = (Y^{\bullet} \setminus \{R_1\}) \cup \{R_2,R_3\}$, and so $|Y^*| = |Y^{\bullet}| + 1$. Thus, $|E_{H}^*(Y^*)| \le  |E_{H^{\bullet}}^*(Y^{\bullet})| + |\{e_1,e_2,e_3,h\}| \le (|Y^{\bullet}| - 2) + 4 = |Y^*| + 1$, contradicting Claim~I. Hence, $|Y^{\bullet}| \le 2$.
Suppose that $|Y^{\bullet}| = 2$. Let $R_2$ be the special subhypergraph in $Y^{\bullet}$ different from $R_1$, where recall that $f^{\bullet} \in E(R_1)$.  As observed earlier, $\defic(H^{\bullet}) \ge 5$ and there is no $H_{10}$-subhypergraph in $H^{\bullet}$. Thus, $5 \le \defic(H^{\bullet}) \le 8|Y^{\bullet}| - 13|E_{H^{\bullet}}^*(Y^{\bullet})| = 16 - 13|E_{H^{\bullet}}^*(Y^{\bullet})|$, implying that $E_{H^{\bullet}}^*(Y^{\bullet}) = \emptyset$. Letting $Y^* = \{R_2\}$, this in turn implies that $E_{H}^*(Y^*) = \emptyset$ and therefore that $\defic(H) > 0$, a contradiction. Hence, $|Y^{\bullet}| = 1$. This completes the proof of Claim~P.11.~\smallqed

\medskip
By Claim~P.11(d), $f^{\bullet} \in Y^{\bullet}$. By Claim~P.11(e), $|Y^{\bullet}| = 1$. Let $R_1$ be the special subhypergraph in $Y^{\bullet}$ that contains the edge $f^{\bullet}$; that is, $Y^{\bullet} = \{R_1\}$. By Claim~P.11, we note that $5 \le \defic(H^{\bullet}) \le 8|Y^{\bullet}| - 13|E_{H^{\bullet}}^*(Y^{\bullet})| = 8 - 13|E_{H^{\bullet}}^*(Y^{\bullet})|$, implying that $E_{H^{\bullet}}^*(Y^{\bullet}) = \emptyset$; that is, $R_1$ is a component in $H^{\bullet}$. By our way in which $H^{\bullet}$ is constructed, this implies that $R_1 =  H^{\bullet}$. Thus, $n(H) = n(H^{\bullet}) + 7 = 14 + 7 = 21$. Since the vertex $d$ belongs to $R_1$ and is contained in both edges $f^{\bullet}$ and $h$, the component $R_1$ is not an $H_4$-component in $H^{\bullet}$, implying that $R_1$ is an $H_{14}$-component in $H^{\bullet}$ and $\defic(H^{\bullet}) = 5$. Since $R_1$ contains the edge $h$, the component $R_1$ has at least one vertex of degree~$1$, namely the vertex $c$ which has degree~$1$ in $H^{\bullet}$. This implies that $R_1 \notin \{H_{14,5},H_{14,6}\}$; that is, $R_1 \cong H_{14,i}$ for some $i \in [4]$.

Suppose that $R_1 = H_{14,1}$. In this case, $h$ is the edge in $R_1$ that contains the (unique) vertex of degree~$1$ in $R_1$. Further, the edge $h$ is intersected by three edges in $R_1$, one of which is the edge $f^{\bullet}$. The structure of $H_{14,1}$ implies that we can choose a  $\tau(H^{\bullet})$-transversal, $T^{\bullet}$, to contain three vertices of $f^{\bullet}$, one of which is the vertex $d$ (that belongs to both $f^{\bullet}$ and $h$). Renaming the edges $e_1,e_2,e_3$ if necessary, we may assume that $\{e_2,e_3\} \subset T^{\bullet}$. If the edges $e_1$ and $f_1$ intersect, let $T = T^{\bullet} \cup \{(e_1f_1)\}$. If the edges $e_1$ and $f_1$ do not intersect, then the edges $e_1$ and $f_2$ intersect and we let $T = (T^{\bullet} \setminus \{d\}) \cup \{c,(e_1f_2)\}$. In both cases, the set $T$ is a transversal in $H$ of size~$|T^{\bullet}| + 1$, and so $\tau(H) \le \tau(H^{\bullet}) + 1$. Thus, $45\tau(H) \le 45(\tau(H^{\bullet}) + 1) \le 6n(H^{\bullet}) + 13m(H^{\bullet}) + \defic(H^{\bullet}) + 45 \le 6n(H) + 13m(H) + 5 + 45 - 7 \times 6 - 4 \times 13 <  6n(H) + 13m(H)$, a contradiction.

Suppose that $R_1 = H_{14,3}$. In this case, $h$ is the edge in $R_1$ that contains the (unique) vertex of degree~$1$ in $R_1$. Further, the edge $h$ is intersected by three edges in $R_1$, one of which is the edge $f^{\bullet}$. Using analogous arguments as in the previous case, due to the structure of $H_{14,3}$ we can choose a $\tau(H^{\bullet})$-transversal, $T^{\bullet}$, to contain three vertices of $f^{\bullet}$, one of which is the vertex $d$, implying as before that $\tau(H) \le \tau(H^{\bullet}) + 1$, producing a contradiction.

Suppose that $R_1 = H_{14,4}$. In this case, $R_1$ contains two vertices of degree~$1$, and $h$ is one of the two edges in $R_1$ that contain a vertex of degree~$1$ in $R_1$. Further, the edge $h$ is intersected by three edges in $R_1$, one of which is the edge $f^{\bullet}$. Using analogous arguments as in the previous two cases, due to the structure of $H_{14,4}$ we can choose a $\tau(H^{\bullet})$-transversal, $T^{\bullet}$, to contain at least three vertices of $f^{\bullet}$, one of which is the vertex $d$, implying as before that $\tau(H) \le \tau(H^{\bullet}) + 1$, producing a contradiction.

By the above, we must have $R_1 = H_{14,2}$. In this case, $h$ is the edge in $R_1$ that contains the (unique) vertex of degree~$1$ in $R_1$. Further, the edge $h$ is intersected by three edges in $R_1$, one of which is the edge $f^{\bullet}$. We note that in this case,
each vertex $s_i$ has degree~$2$ in $R_1$ and therefore degree~$2$ in $H$.

Suppose that $f_1$ or $f_2$ (or both $f_1$ and $f_2$) does not intersect one of the edges $e_1,e_2,e_3$. Renaming vertices and edges, if necessary, we may assume that $f_1$ does not intersect the edge $e_2$. By Claim~P.1, the edge $f_1$ therefore intersects both $e_1$ and $e_3$. Let $x_2$ be the vertex in the edge $e_2$ different from $x$, $s_2$ and $(e_2f_2)$. We note that we could have chosen $s_2$ to be the vertex $x_2$, implying by our earlier observations that $d_H(x_2) = 2$ and $\{s_1,x_2,s_3\}$ is an independent set.
If the vertex $x_2$ belongs to the set $V(R_1)$, then $x_2$ would be adjacent to $s_1$ or $s_3$ in $H$ or $x_2$ and $s_2$ would be adjacent in $H - e_2$, a contradiction. Hence, $x_2 \notin V(R_1)$. Interchanging the roles of $x_2$ and $s_2$ in our earlier arguments and letting $f^{\bullet} = \{d,s_1,x_2,s_3\}$, the special subhypergraph, say  $R_1^{\bullet}$, in $Y^{\bullet}$ that contains the edge $f^{\bullet}$ is a $H_{14,2}$-component of $H^{\bullet}$ that does not contain the vertex $s_2$. This, however, is a contradiction since the vertex of degree~$3$ in $R_1$ is also the degree~$3$ vertex in $R_1^{\bullet}$, implying that $s_2 \in V(R_1^{\bullet})$, a contradiction.

Therefore, the edge $f_i$ intersects all three edges $e_1,e_2,e_3$ for $i \in [2]$. The hypergraph $H$ is now determined and $H = H_{21,5}$, contradicting Claim~C.~\smallqed


\medskip
\ClaimX{Q}
No edge in $H$ contains two degree-$3$ vertices.

\ProofClaimX{Q}
We first need the following subclaim.

\ClaimX{Q.1}
There does not exist three degree-$3$ vertices in $H$ that are pairwise adjacent.

\ProofClaimX{Q.1}
For the sake of contradiction, suppose to the contrary that $x_1$, $x_2$ and $x_3$ are pairwise adjacent degree-$3$ vertices in $H$. First consider the case when some edge, $e$, in $H$ contains all three vertices. Let $V(e)=\{x_1,x_2,x_3,x_4\}$.
By Claim~P, the hypergraph $H-x_i$ contains an $H_4$-component, $R_i$, for all $i \in [3]$. Further, let $f_i$ be the edge in $R_i$. By Claim~M, the edge $f_i$ is intersected by all three edges incident with~$x_i$. Since every vertex in $f_i$ has degree at most~$2$ in $H$ and $d(x_1)=d(x_2)=d(x_3)=3$, we note that $x_4 \in V(f_i)$ for all $i \in [3]$, which implies that $f_1=f_2=f_3$, which is impossible as the two edges different from $e$ containing $x_1$ intersect $f_1$
and therefore $R_1$ is not an $H_4$-component in $H-x_2$.

Therefore there exists three distinct edges $e_{12}$, $e_{13}$ and $e_{23}$, such that $e_{ij}$ contains $x_i$ and $x_j$ for all $1 \le i < j \le 3$. Let $R_i$ be the $H_4$-component in $H-x_i$ for $i \in [3]$, and let $h_i$ be the edge in $R_i$. Note that $V(h_i) \cap \{x_1,x_2,x_3\} = \emptyset$ for all $i \in [3]$ as no vertex in $h_i$ has degree $3$. If $h = h_1 = h_2$, then the edge $h$ intersects all edges containing $x_1$ (as $h=h_1$) and all edges intersecting $x_2$ (as $h=h_2)$, but then $h$ does not belong to an $H_4$ component in $H-x_1$ (or in $H-x_2$), a contradiction. Hence, $h_1$, $h_2$ and $h_3$ are distinct edges. Furthermore, if $V(h_1) \cap V(h_2) \ne \emptyset$ and $y \in V(h_1) \cap V(h_2)$, then $y$ is contained in the three edges $e_{12}$, $h_1$ and $h_2$, and so $d(y) = 3$, a contradiction to all vertices in $h_1$ having degree at most two. Therefore, $h_1$, $h_2$ and $h_3$ are non-intersecting edges.

For each $i \in [3]$, let $f_i$ be the edge incident with $x_i$, different from $e_{12}$, $e_{13}$ and $e_{23}$, which exists as $d(x_i)=3$ and two of the three edges $e_{12}$, $e_{13}$ and $e_{23}$ contain $x_i$. By the above definitions, we note that $h_1$ intersects the three edges $e_{12}$, $e_{13}$ and $f_1$, which all contain $x_1$. Also, $h_2$ intersects the three edges $e_{12}$, $e_{23}$ and $f_2$, which all contain $x_2$. And finally,
$h_3$ intersects the three edges $e_{13}$, $e_{23}$ and $f_3$, which all contain $x_3$.

Let $Z=V(e_{12}) \cup V(e_{13}) \cup V(e_{23}) \cup V(h_1) \cup V(h_2) \cup V(h_3)$ and note that $|Z|=15$, as $|V(e_{12}) \cup V(e_{13}) \cup V(e_{23})|=9$ and each $h_i$ contains two vertices not in $V(e_{12}) \cup V(e_{13}) \cup V(e_{23})$.
Also note that the only edges in $H$ intersecting $Z$ are $e_{12},e_{13},e_{23},h_1,h_2,h_3,f_1,f_2,f_3$.
Let $H'=H \setminus Z$ and note that $n(H')=n(H)-15$ and $m(H')=m(H)-9$. Further, we note that $|V(f_i) \setminus Z|=2$ for $i \in [3]$. Moreover since $H$ is linear, $V(f_1) \setminus Z$ and $V(f_2) \setminus Z$ can have at most one vertex in common, implying that $|\partial(Z)| \ge 3$.

We wish to use Claim~L, so first we need to show that $\defic(H') \le 21$. Suppose to the contrary that $\defic(H') \ge 22$. Let $Y'$ be a special $H'$-sets in $H'$ with $\defic_{H'}(Y') = \defic(H') \ge 22$. By Claim~O, $|E_{H'}^*(Y')| \le |Y'| - 3$. This implies that $|E_{H}^*(Y')| \le |E_{H'}^*(Y')| + |\{f_1,f_2,f_3\}| \le (|Y'|-3)+3 = |Y'|$, a contradiction to Claim~I(a). Therefore, $\defic(H') \le 21$.

As observed earlier, $|\partial(Z)| \ge 3$, and so $f(|\partial(Z)|) \le 27$. By Claim~L, there exists a transversal, $T'$, in $H'$, that contains a vertex in $\partial(Z)$ such that $45|T'|  \le  6n(H') + 13m(H') + f(|\partial(Z)|) \le 6(n(H)-15) + 13(m(H)-9) + 27 = 6n(H) + 13m(H) - 180$.  Since $T'$ contains a vertex in $\partial(Z)$, the transversal $T'$ contains a vertex from $V(f_1)$, $V(f_2)$ or $V(f_3)$. Without loss of generality, $T'$ contains a vertex from $V(f_1)$. With this assumption, the set $T = T' \cup \{x_3,(e_{12}h_1),(f_2h_2),w\}$ where $w$ is an arbitrary vertex from $V(h_3)$, is a transversal in $H$ of size~$|T'| + 4$. Hence, $45\tau(H) \le 45|T| = 45(|T'| + 4) \le 6n(H) + 13m(H) - 180 + 4 \cdot 45 = 6n(H) + 13m(H)$, a contradiction. This proves  Claim~Q.1~\smallqed

\medskip
We now return to the proof of Claim~Q. Suppose, to the contrary, that there exists an edge $e = \{x_1,x_2,x_3,x_4\}$, such that $d(x_1)=d(x_2)=3$ in $H$. Let $e$, $f_1$ and $f_2$ be the three edges in $H$ containing $x_1$ and let $e$, $g_1$ and $g_2$ be the three edges in $H$ containing $x_2$. Let $R_i$ be the $H_4$-component in $H-x_i$ and let $h_i$ be the edge of $R_i$ for $i \in [2]$, and so $E(R_i) = \{h_i\}$. Renaming vertices if necessary, we may assume that $x_3 \in V(h_1)$ and $x_4 \in V(h_2)$, as $h_i$ intersects all three edges incident with $x_i$ for $i \in [2]$ and $h_1$ and $h_2$ are non-intersecting edges. Let $Z' = V(h_1) \cup V(h_2) \cup \{x_1,x_2\}$ and let
\vskip 0.1cm \hspace*{0.35cm}
$Q_1 = V(g_1) \setminus \{(g_1h_2),x_2\}$, \hspace*{1cm}  $Q_2 = V(g_2) \setminus  \{(g_2h_2),x_2\}$,
\\ \hspace*{0.75cm}
$Q_3 = V(f_1) \setminus  \{(f_1h_1),x_1\}$,
\hspace*{1cm}  $Q_4 = V(f_2) \setminus  \{(f_2h_1),x_1\}$.

Let $Q' = Q_1 \cup Q_2 \cup Q_3 \cup Q_4$ and note that $Q' = \partial(Z') = N_H(Z') \setminus Z'$. Let $H'=H-Z'$. We note that $n(H')=n(H)-10$ and  $m(H') = m(H)- |\{e,f_1,f_2,g_1,g_2,h_1,h_2\}| = m(H) - 7$. Let $R_1 = Q_1 \cup Q_2$ and $R_2 = Q_3 \cup Q_4$.
As $|R_1|=|R_2|=4$ and $Q' = R_1 \cup R_2$ we note that $4 \le |Q'| \le 8$.

\ClaimX{Q.2}
If there is a transversal, $T'$, in $H'$ containing a vertex from $Q'$, then $\tau(H) \le |T'|+3$. Furthermore, if $T'$ contains a vertex from three of the sets $Q_1,Q_2,Q_3,Q_4$, then $\tau(H) \le |T'|+2$.

\ProofClaimX{Q.2}
Let $T'$ be a transversal in $H'$ containing a vertex from $Q'$. Renaming vertices if necessary, we may assume that $T'$ intersects $Q_1$, and therefore covers the edge $g_1$. With this assumption, the set $T = T' \cup \{x_1,x_3,(g_2h_2)\}$ is a transversal in $H$ of size~$|T'| + 3$, where we note that the vertex $x_3$ covers $h_1$, the vertex $x_1$ covers $e$, $f_1$ and $f_2$, and the vertex $(g_2h_2)$ covers $g_2$ and $h_2$. Hence, $\tau(H) \le |T| = |T'|+3$. Furthermore, suppose that the transversal $T'$ of $H'$ contains a vertex from three of the sets $Q_1$, $Q_2$, $Q_3$, $Q_4$. Without loss of generality, we may assume that $T'$ intersects $Q_1$, $Q_2$ and $Q_3$ and therefore covers the edges $g_1$, $g_2$ and $f_1$. In this case, the set $T = T' \cup \{(eh_2),(f_2h_1)\}$ is a transversal in $H$ of size~$|T'| + 2$, implying that $\tau(H) \le |T| = |T'|+2$. This completes the proof of Claim~Q.2~\smallqed

\ClaimX{Q.3}
Every vertex in $R_1 \cap R_2$ has degree~$2$ in $H$.

\ProofClaimX{Q.3}
Every vertex in $R_1 \cap R_2$ is intersected by one of $g_1$ or $g_2$ and one of $f_1$ or $f_2$ and therefore has degree at least~$2$. By Claim~Q.1, such a vertex in $R_1 \cap R_2$ cannot have degree~$3$, which completes the proof of Claim~Q.3~\smallqed

\ClaimX{Q.4}
$|Q'| \ge 5$.

\ProofClaimX{Q.4}
For the sake of contradiction, suppose that $|Q'| \le 4$. As observed earlier, $|Q'| \ge 4$. Consequently, $|Q'|=4$ and $Q' = R_1 = R_2$, implying by Claim~Q.3 that every vertex in $Q$ has degree~$2$. The hypergraph $H$ is now determined and by linearity we note that $H$ is isomorphic to $H_{14,4}$, a contradiction. This proves Claim~Q.4~\smallqed

\ClaimX{Q.5}
$|E_{H'}^*(Y)| \ge |Y|$ for all special $H'$-sets $Y$ in $H'$ and therefore $\defic(H') = 0$.

\ProofClaimX{Q.5}
We first prove that $\defic(H') \le 16$. Suppose to the contrary that $\defic(H') \ge 17$. Let $Y$ be a special $H'$-sets with $\defic_{H'}(Y) =  \defic(H') \ge 17$. By Claim~O, $|E^*_{H'}(Y)| \le |Y|-3$. Let $R_1$ and $R_2$ be the $H_4$-subhypergraphs of $H$ with edges $h_1$ and $h_2$, respectively. Let $Y^* = Y \cup \{R_1,R_2\}$, and so $|Y^*| = |Y| + 2$. Thus, $|E^*_{H}(Y^*)| \le |E^*_{H'}(Y)| + |\{e,f_1,f_2,g_1,g_2\}| \le |Y|-3+5 = |Y^*|$, contradicting Claim~I(a). Therefore  $\defic(H') \le 16$.

For the sake of contradiction suppose that $|E_{H'}^*(Y)| < |Y|$ for some special $H'$-sets $Y$ in $H'$. By Claim~K, there exists a transversal $T'$ in $H'$, such that $45 |T'| \le 6n(H') + 13m(H') + \defic(H')$ and $T' \cap Q' \ne \emptyset$. By Claim~Q.2 we note that $\tau(H) \le |T'|+3$. Recall that $n(H')=n(H)-10$ and $m(H') = m(H) - 7$. Thus, $45\tau(H) \le 45(|T'| + 3) \le 6n(H') + 13m(H') + \defic(H') + 135 \le 6(n(H)-10) + 13(m(H)-7) + 16 + 135 = 6n(H) + 13m(H)$, a contradiction. Therefore,  $|E_{H'}^*(Y)| \ge |Y|$ for all special $H'$-sets $Y$ in $H'$, and so $\defic(H') = 0$ which proves Claim~Q.5~\smallqed

\ClaimX{Q.6}
$|Q'| = 8$.

\ProofClaimX{Q.6}
For the sake of contradiction suppose that $|Q'| \ne 8$. By Claim~Q.4 and the fact that $|Q'| \le 8$ we obtain $5 \le |Q'| \le 7$. Therefore, $R = R_1 \cap R_2$ is non-empty and by Claim~Q.3 every vertex in $R$ has degree~$2$ in $H$ and therefore degree zero in $H'$.

Let $Z'' = Z' \cup R$, let $Q'' = Q' \setminus R$ and note that $Q'' = \partial (Z'')$. Let $H'' = H - Z''$. By Claim~Q.5 we note that $\defic(H') = 0$ and therefore $\defic(H'') = 0$ also holds (as we only removed isolated vertices). Note that $n(H'') = n(H) - 10 - |R|$ and $m(H'') = m(H') = m(H) - 7$. Further, we note that $|Q'| = 8-|R|$ and $|Q''| = |Q'|-|R| = 8-2|R|$. By Claim~L, there exists a transversal, $T''$, in $H''$ containing a vertex from $Q''$, such that
\[
\begin{array}{rcl}
45|T''| & \le & 6n(H'')+m(H'')+f(|Q''|) \\
       & \le & 6 (n(H)-10-|R|) +13(m(H)-7) + f(8-2|R|) \\
       & = & 6 n(H)+13m(H) - 3\cdot45 - 16 - 6|R| + f(8-2|R|). \\
\end{array}
\]

If $|R|=1$, then $f(8-2|R|) - 16 - 6|R| = 22-22=0$.
If $|R|=2$, then $f(8-2|R|) - 16 - 6|R| = 23-28=-6$.
If $|R|=3$, then $f(8-2|R|) - 16 - 6|R| = 33-34=-1$.
So in all cases, $f(8-2|R|) - 16 - 6|R| \le 0$, implying that $45|T''| \le 6 n(H) + 13m(H) - 3\cdot45$. Thus by Claim~Q.2, $45\tau(H) \le 45(|T'|+3) \le  6 n(H)+13m(H)$, a contradiction. This completes the proof of Claim~Q.6~\smallqed

By Claim~Q.6 we note that $Q_1$, $Q_2$, $Q_3$ and $Q_4$ are vertex disjoint. We will now define the following four different hypergraphs.

Let $H_1^*$ be the hypergraph obtained from $H'$ by adding three vertices $x_{23}$, $x_{24}$ and $x_{34}$
and three hyperedges $h^*_2 = Q_2 \cup \{x_{23},x_{24}\}$,
 $h^*_3 = Q_3 \cup \{x_{23},x_{34}\}$ and
 $h^*_4 = Q_4 \cup \{x_{24},x_{34}\}$.

Analogously let $H_2^*$ be the obtained from $H'$ by adding three vertices $y_{13}$, $y_{14}$ and $y_{34}$
and three hyperedges $e^*_1 = Q_1 \cup \{y_{13},y_{14}\}$,
 $e^*_3 = Q_3 \cup \{y_{13},y_{34}\}$ and
 $e^*_4 = Q_4 \cup \{y_{14},y_{34}\}$.

Let $H_3^*$ be the obtained from $H'$ by adding three vertices $z_{13}$, $z_{14}$ and $z_{34}$
and three hyperedges $f^*_1 = Q_1 \cup \{z_{12},z_{14}\}$,
 $f^*_2 = Q_3 \cup \{z_{12},z_{24}\}$ and
 $f^*_4 = Q_4 \cup \{z_{14},z_{24}\}$.

Finally let $H_4^*$ be the obtained from $H'$ by adding three vertices $w_{12}$, $w_{13}$ and $w_{23}$ and three hyperedges $g^*_1 = Q_1 \cup \{w_{12},w_{13}\}$, $g^*_2 = Q_2 \cup \{w_{12},w_{23}\}$ and $g^*_3 = Q_3 \cup \{w_{13},w_{23}\}$.

\ClaimX{Q.7}
$\tau(H) \le \tau(H_i^*) + 2$, for $i \in [4]$.

\ProofClaimX{Q.7}
Let $T^*$ be a minimum transversal in $H_1^*$.  If $T^* \cap \{x_{23}, x_{24}, x_{34}\} = \emptyset$, then
$T^*$ covers $Q_2$, $Q_3$ and $Q_4$, and by Claim~Q.2, $\tau(H) \le |T^*|+2$. We may therefore, without loss of generality, assume that $x_{23} \in T^*$. We may also assume that $T^* \cap \{x_{24}, x_{34}\} = \emptyset$, since otherwise we could have
picked a vertex from $Q_4$ instead of this vertex (or vertices). Therefore, $Q_4$ is covered by $T^* \setminus \{x_{23}\}$. By Claim~Q.2, $\tau(H) \le |T^* \setminus \{x_{23}\} |+3= |T^*|+2$, which completes the proof for $H_1^*$. Analogously, the claim holds for  $H_i^*$, with $i=2,3,4$.~\smallqed

\ClaimX{Q.8}
$\defic(H_i^*) \ge 5$, for $i \in [4]$.

\ProofClaimX{Q.8}
For the sake of contradiction, suppose to the contrary that $\defic(H_1^*) \le 4$ and let $T^*$ be a minimum transversal in $H^*_1$. Applying the inductive hypothesis to $H''$, we have by Claim~Q.7 that $45 \tau(H) \le 45(|T^*|+2) \le 6 n(H_1^*) + 13m(H_1^*) + \defic(H_1^*)  + 90 \le 6(n(H) - 7) + 13(m(H)-4) + 4 + 90 = 6 n(H)+13m(H)$, a contradiction. Hence, $\defic(H_1^*) \ge 5$. Analogously, $\defic(H_i^*) \ge 5$, for $i=2,3,4$.~\smallqed

Let $Y_i^*$ be a special $H_i^*$-set satisfying $\defic(H_i^*) = \defic_{H_i^*}(Y_i^*)$ for $i \in [4]$.

\ClaimX{Q.9}
$|Y_1^*|=1$ and $Y_1^*$ is an $H_{14}$-component in $H_1^*$ containing the edges $h_2^*,h_3^*,h_4^*$.

\ProofClaimX{Q.9}
If none of $h_2^*$, $h_3^*$ or $h_4^*$ belong to $E(Y_1^*)$, then in this case $|E^*(Y_1^*)| \le |Y_1^*|-1$ in $H_1^*$, which implies that $|E^*(Y_1^*)| \le |Y_1^*|-1$ in $H'$, contradicting Claim~Q.5. Therefore, without loss of generality, we may assume that $h_2^* \in E(Y_1^*)$.

Suppose that there is an $H_4$-component, $R^*$, in $Y_1^*$ that contains the edge $h_2^*$. Thus $E(R^*) = \{h_2^*\}$. In this case, we note that $h_3^*, h_4^* \in E^*(Y_1^*)$ in $H_1^*$. As $\defic(H_1^*) \ge 5$ and $|E^*(Y_1^*)| \ge 2$, we have $|E^*(Y_1^*)| \le |Y_1^*|-2$ in $H_1^*$. Let $Y_1' = Y_1^* \setminus \{R_1\}$. Since $h_3^*, h_4^* \in E^*(Y_1^*)$, we note that $|E^*(Y_1')| \le |Y_1^*|-2-2 = |Y_1'|-3$ in $H'$, contradicting Claim~Q.5.  Therefore, $R^*$ is not an $H_4$-component in $Y_1^*$. Analogously, there is no $H_4$-component in $Y_1^*$ that contains the edge $h_3^*$ or $h_4^*$.

If two of the three edges $h_2^*,h_3^*,h_4^*$ belong to $E(Y_1^*)$ in $H_1^*$, then these two edges, which intersect,  both contain a degree-$1$ vertex in $Y_1^*$. However, there is no special hypergraph with two intersecting edges both containing a degree-$1$ vertex, a contradiction. Therefore, we must have all three edges $h_2^*,h_3^*,h_4^*$ belonging to $E(Y_1^*)$.

Now consider the case when $|Y_1^*| \ge 2$. Let $Y_1^*=\{R_1^*,R_2^*,\ldots,R_\ell^*\}$, where $\ell \ge 2$. Renaming the elements of $Y_1^*$, we may assume that $\{h_2^*,h_3^*,h_4^*\} \subseteq E(R_1^*)$. As $\defic(H_1^*) \ge 5$, we must have $|E^*(Y_1^*)| \le |Y_1^*|-2$ in $H_1^*$. This implies that in $H$ we have  $|E^*(Y_1^* \setminus \{R_1^*\})| \le |Y_1^*|-2 + |\{g_1\}| = |Y_1^* \setminus \{R_1^*\}|$, contradicting Claim~I(a). Therefore, $|Y_1^*|=1$ and $E^*(Y_1^*) = \emptyset$ in $H^*$.
As  all three edges $h_2^*,h_3^*,h_4^*$ belong to $E(Y_1^*)$, we note that $Y_1^*$ is not an $H_4$-component.  As  $\defic(H_1^*) \ge 5$ we must therefore have that $Y_1^*$ is an $H_{14}$-component, completing the proof of Claim~Q.9.~\smallqed

We note that an analogous claim to Claim~Q.9 also holds for $H_2^*$, $H_3^*$ and $H_4^*$.

\ClaimX{Q.10}
Both vertices in at least two of the sets $Q_1$, $Q_2$, $Q_3$, $Q_4$ have degree~$1$ in~$H'$.

\ProofClaimX{Q.10}
For the sake of contradiction, suppose to the contrary that $Q_2$, $Q_3$ and $Q_4$ contain a vertex which does not have degree~$1$ in $H'$. If it is an isolated vertex, then without loss of generality assume that $y \in Q_2$ is isolated in $H'$. Consider the hypergraph $H'' = H - y$ and let $Z'' = Z' \cup \{y\}$. By Claim~Q.5 we note that $\defic(H'') = \defic(H') = 0$. Furthermore, $H'' = H - Z''$. Let $Q'' = N(Z'') \setminus Z''$, and so $Q'' = \partial(Z'')$. By Claim~Q.6, we note that $|Q''| = |Q' \setminus \{y\}|=7$. By Claim~L, there exists a transversal, $T''$, in $H''$ containing a vertex from $Q''$, such that the following holds.
\[
\begin{array}{rcl}
45|T''| & \le & 6n(H'')+m(H'')+f(|Q''|) \\
       & \le & 6 (n(H)-11) +13(m(H)-7) + f(7) \\
       & = & 6 n(H)+13m(H) - 66 - 91 + 22 \\
       & = & 6 n(H)+13m(H) - 3 \times 45. \\
\end{array}
\]

By Claim~Q.2, this implies that $45\tau(H) \le 45(|T''|+3) \le 6 n(H)+13m(H)$, a contradiction.
Therefore, $Q_2$, $Q_3$ and $Q_4$ do not contain isolated vertices in $H'$ and they therefore
all contain a vertex of degree at least~$2$ (in fact, equal to~$2$), which implies that they have degree~$3$ in $H_1^*$. As all these vertices belong to $Y_1^*$, we get a contradiction to the fact that no $H_{14}$-component contains three degree-$3$ vertices.
This completes the proof of Claim~Q.10~\smallqed

By Claim~Q.10, we assume without loss of generality that both vertices in $Q_1$ and in $Q_2$ have degree~$1$ in $H'$, which implies that the vertices in $Q_1$ have degree~$1$ in $H_1^*$ and all vertices in the edge $h_2^*$ (which contains $Q_2$) have degree~$2$ in $H_1^*$.

\ClaimX{Q.11}
$Q' \subseteq V(Y_1^*)$.

\ProofClaimX{Q.11}
By Claim~Q.9 we note that $Q_i \subseteq V(Y_1^*)$ for $i=2,3,4$. For the sake of contradiction, suppose to the contrary that $Q_1 \not\subseteq V(Y_1^*)$. As $Y_1^*$ is an $H_{14}$-component and all vertices in $h_2^*$ have degree~$2$ in $H_1^*$ and therefore also in $Y_1^*$, we note that $Y_1^* - h_2^*$ is a component containing fourteen vertices. This implies that the component containing $\{x_{23},x_{24},x_{34}\}$ in $H_2^*$ has more than fourteen vertices as it contains all vertices of $V(Y_1^*)$ as well as $Q_1$. This contradiction to Claim~Q.9 completes the proof of Claim~Q.11.~\smallqed

By Claim~Q.11 we note that $Y_1^*$ contains two degree-$1$ vertices (namely, the two vertices in $Q_1$) and an edge consisting of only degree-$2$ vertices (namely, the edge $h_2^*$).  However no $H_{14}$ component has these properties, a contradiction. This completes the proof of Claim~Q.~\smallqed

\medskip
By Claim~M and Claim~P, if $x$ is an arbitrary vertex of $H$ of degree~$3$, then $H - x$ contains an $H_4$-component that is intersected by all three edges incident with~$x$, and $\defic(H-x) = 8$. By Claim~Q, no edge in $H$ contains two degree-$3$ vertices.

\medskip
\ClaimX{R}
Every degree-$3$ vertex in $H$ has at most one neighbor of degree~$1$.

\ProofClaimX{R}
Suppose, to the contrary, that there exists a vertex $x$ of degree~$3$ in $H$ with at least two degree-$1$ neighbors, say $y$ and $z$. Let $H' = H - \{x,y,z\}$, and note that $n(H') = n(H) - 3$ and $m(H') = m(H) - 3$. Let $e$ be the edge in the $H_4$-component of $H-x$. By Claim~M and Claim~P, the edge $e$ is intersected by all three edges incident with $x$, and contains a vertex (of degree~$1$ in $H$) not adjacent to $x$. Further by Claim~M, $\defic(H - x) = 8$. We note that $\defic(H') = \defic(H - x)$ since $H'$ is obtained from $H - x$ by deleting the two isolated vertices $y$ and $z$ in $H - x$. Every transversal in $H'$ can be extended to a transversal in $H$ by adding to it the vertex $x$, and so $\tau(H) \le \tau(H') + 1$. Applying the inductive hypothesis to $H'$, we have that $45\tau(H) \le 45(\tau(H') + 1) \le 6n(H') + 13m(H') + \defic(H') + 45 = 6(n(H)-3) + 13(m(H)-3) + 8 + 45 = 6n(H) + 13m(H) - 4 < 6n(H) + 13m(H)$, a contradiction.~\smallqed

\medskip
By Claim~R, every degree-$3$ vertex in $H$ has at most one neighbor of degree~$1$. We now define the operation of \emph{duplicating} a degree-$3$ vertex $x$ as follows.

Let $e_1$, $e_2$ and $e_3$ be the three edges incident with $x$. By Claim~Q and Claim~R, every neighbor of $x$ has degree~$2$, except possibly for one vertex which has degree~$1$. Renaming edges if necessary, we may assume that the edge $e_1$ contains no vertex of degree~$1$, and therefore every vertex in $e_1$ different from $x$ has degree~$2$. We now delete the edge $e_1$ from $H$, and add a new vertex $x'$ and a new edge $e_1' = (V(e_1) \setminus \{x\}) \cup \{x'\}$ to $H$. We note that in the resulting hypergraph the vertex $x$ now has degree~$2$ (and is incident with the edges $e_2$ and $e_3$) and the new vertex $x'$ has degree~$1$ with all its three neighbors of degree~$2$. We call $x'$ the \emph{vertex duplicated copy} of $x$.

Let $H'$ be obtained from $H$ by duplicating every degree-$3$ vertex as described above. By construction, $H'$ is a linear $4$-uniform connected hypergraph with minimum degree $\delta(H') \ge 1$ and maximum degree $\Delta(H') \le 2$. For $i \in [2]$, let $n_i(H')$ be the number of vertices of degree~$i$ in $H'$. Then, $n(H') = n_1(H') + n_2(H')$ and $4m(H') = 2n_2(H') + n_1(H')$. We proceed further with the following properties of the hypergraph $H'$.

\medskip
\ClaimX{S}
The following properties hold in the hypergraph $H'$. \\
\hspace*{0.75cm} {\rm (a)} $\tau(H) = \tau(H')$. \\
\hspace*{0.75cm}  {\rm (b)} $\defic(H') = 0$.

\ProofClaimX{S}
In order to show that $\tau(H) = \tau(H')$ we show that the operation of duplicating a degree-$3$ vertex $x$ leaves the transversal number unchanged. Let $x'$, $e_1$, $e_2$ and $e_3$ be defined as in the description of duplication.  Let $e_x$ be the edge in the $H_4$-component in $H-x$ and let $H_x$ be the hypergraph obtained from $H$ by duplicating $x$. Let $x_1^*$ be the vertex $(e_1e_x)$ in $V(e_x) \cap V(e_1)$.

Let $T_x$ be a transversal in $H_x$. As $d_{H'}(x')=1$, we may assume that $x' \notin T_x$. Now we note that $T_x$ is also a transversal in $H$ and therefore $\tau(H) \le \tau(H_x)$.
Conversely assume that $T$ is a transversal in $H$. If $x \notin T$, let $T' = T$. If $x \in T$, then since $T$ is a transversal in $H$, there exists a vertex $y \in V(e_x) \cap T$. In this case, we let $T' = (T \setminus \{y\}) \cup \{x_1^*\}$. In both cases, $|T'| = |T|$ and $T'$ is a transversal in $H_x$, implying that $\tau(H_x) \le |T'| = |T| = \tau(H)$. Consequently, $\tau(H) = \tau(H_x)$, which implies that $\tau(H) = \tau(H')$ and (a) holds.

In order to prove part~(b) we, for the sake of contradiction, suppose that $\defic(H') > 0$. Let $Y$ be a special $H'$-set satisfying $\defic_{H'}(Y) = \defic(H') > 0$. By Claim~N, no edge in $H$ contains two vertices of degree~$1$.  We note that this is still the case after every duplication and therefore no edge in $H'$ contains two vertices of degree~$1$. This implies that every $H_4$-component in $Y_4$ is intersected by at least three edges of $E^*(Y)$.

If $Y \setminus Y_4 \ne \emptyset$, let $R \in Y \setminus Y_4$ be arbitrary. Since $\Delta(H') \le 2$, we note that $R \in \{H_{10},H_{11},H_{14,5},H_{14,6}\}$. Let $x$ be a vertex of degree~$3$ in $H$ and let $x'$ be the duplicated copy of $x$ in $H'$. Furthermore let $e_1'$ be the edge containing $x'$ in $H'$ and let $e_x$ be the edge in the $H_4$-component in $H-x$. Finally let $e_2$ and $e_3$ be the edges containing $x$ in $H'$.  Since $e_x$ and $e_1'$ both contain a degree-$1$ vertex and intersect each other  we note that they both cannot belong to $R$. If $e_1' \in E(R)$, then it would therefore contain two degree-$1$ vertices in $R$ (namely $x'$ and the vertex $(e_xe_1')$ in $V(e_x) \cap V(e_1')$), a contradiction.  Therefore, $e_1' \notin E(R)$. Analogously $e_x \notin E(R)$. The edges $e_2$ and $e_3$ cannot both belong to $R$ as they intersect and would both contain degree-$1$ vertices in $R$ (namely the vertices $(e_2e_x)$ and $(e_3e_x)$ in $V(e_2) \cap V(e_x)$ and $V(e_3) \cap V(e_x)$, respectively). However if $e_2 \in E(R)$, then $e_2$ would contain two degree-$1$ vertices in $R$ (namely $x$ and the vertex $(e_2e_x)$ in $V(e_2) \cap V(e_x)$), a contradiction.  Therefore, $e_2 \notin R$ and analogously $e_3 \notin R$.  This implies that $R$ was also a special hypergraph in $H$. By Claim~J, $R$ is not isomorphic to $H_{10}$, and so $R \in \{H_{11},H_{14,5},H_{14,6}\}$. Furthermore if $R$ is intersected by $k$ edges in $H'$, then it is also intersected by $k$ edges in $H$, and so by Claim~I(a) we note that at least three edges intersect $R$ in $H$ and therefore also in $H'$.
Therefore, all components in $Y$ are intersected by at least three edges from $|E^*(Y)|$. This implies that
\[
4|E^*(Y)| \ge \sum_{e \in E^*(Y)} |V(e) \cap V(Y)| \ge 3|Y|,
\]
and so $13|E^*(Y)| \ge \frac{39}{4} |Y| > 8|Y|$.  As observed earlier, $Y_{10} = \emptyset$. Hence, $\defic(Y) \le 8|Y|-13|E^*(Y)| < 0$, a contradiction. Therefore
$\defic(H') = 0$. This completes the proof of Claim~S.~\smallqed


\medskip
We now consider the multigraph $G$ whose vertices are the edges of $H'$ and whose edges correspond to the $n_2(H')$ vertices of degree~$2$ in $H'$: if a vertex of $H'$ is contained in the edges $e$ and $f$ of $H'$, then the corresponding edge of the multigraph $G$ joins vertices $e$ and $f$ of $G$. By the linearity of $H'$, the multigraph $G$ is in fact a graph, called the \emph{dual} of $H'$. We shall need the following properties about the dual $G$ of the hypergraph $H'$.

\ClaimX{T}
The following properties hold in the dual, $G$, of the hypergraph $H'$. \\
\indent {\rm (a)} $G$ is connected, $n(G) = m(H')$ and $m(G) = n_2(H')$. \\
\indent {\rm (b)} $\Delta(G) \le 4$, and so $m(G) \le 2n(G)$ and $8n(G) + 6m(G) \le 20n(G)$. \\
\indent {\rm (c)} $\tau(H') = m(H') - \alpha'(G)$.

\ProofClaimX{T} Since $H'$ is connected, so too is $G$ by construction. Further by construction, $n(G) = m(H')$ and $m(G) = n_2(H')$. Since $H'$ is $4$-uniform and $\Delta(H') \le 2$, we see that $\Delta(G) \le 4$, implying that $m(G) \le 2n(G)$ and $8n(G) + 6m(G) \le 20n(G)$. This establishes Part~(a) and Part~(b). Part~(c) is well-known (see, for example,~\cite{DoHe14,HeYe17}), so we omit the details of this property.~\smallqed

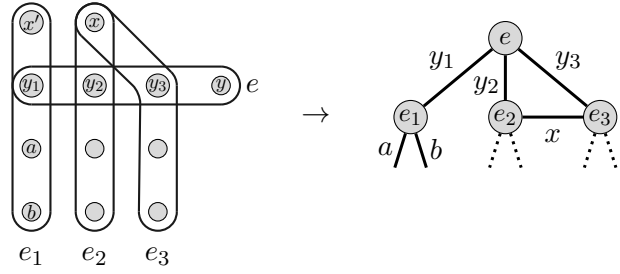
\begin{figure}[htb]
\begin{center}
\tikzstyle{vertexX}=[circle,draw, fill=gray!30, minimum size=10pt, scale=0.7, inner sep=0.1pt]
\tikzstyle{vertexY}=[circle,draw, fill=gray!30, minimum size=14pt, scale=0.9, inner sep=0.21pt]
\begin{tikzpicture}[scale=0.42]
\node at (1, -0.4) {$e_1$};
\node at (3, -0.4) {$e_2$};
\node at (5, -0.4) {$e_3$};
\node at (8, 5) {$e$};

\node (a0) at (3.0,7.0) [vertexX] {$x$};
\node (a1) at (1.0,5.0) [vertexX] {$y_1$};
\node (a2) at (3.0,5.0) [vertexX] {$y_2$};
\node (a3) at (5.0,5.0) [vertexX] {$y_3$};
\node (a4) at (1.0,3.0) [vertexX] {$a$};
\node (a5) at (3.0,3.0) [vertexX] {};
\node (a6) at (5.0,3.0) [vertexX] {};
\node (a7) at (1.0,1.0) [vertexX] {$b$};
\node (a8) at (3.0,1.0) [vertexX] {};
\node (a9) at (5.0,1.0) [vertexX] {};
\node (bb) at (7.0,5.0) [vertexX] {$y$};
\node (a10) at (1.0,7.0) [vertexX] {$x'$};
\draw[color=black!90, thick,rounded corners=4pt] (0.40000003248683635,0.999802555824115) arc (180.01885453028638:359.93733542082634:0.6); 
\draw[color=black!90, thick,rounded corners=4pt] (1.599980696647063,7.004812863067338) arc (0.4595994973493447:179.98114899285022:0.6); 
\draw[color=black!90, thick,rounded corners=4pt] (0.40000003248683635,0.999802555824115) -- (0.4000000324746966,7.0001974072816875);
\draw[color=black!90, thick,rounded corners=4pt] (1.599980696647063,7.004812863067338) -- (1.5999962965560313,4.997891892069423) -- (1.5999996411440085,0.999343778192253);
\draw[color=black!90, thick,rounded corners=4pt] (3.5999115296852535,6.989696769890026) arc (-0.9839343621107446:180.01636223573303:0.6); 
\draw[color=black!90, thick,rounded corners=4pt] (2.4000000324868362,0.999802555824115) arc (180.01885453028638:359.93733542082634:0.6); 
\draw[color=black!90, thick,rounded corners=4pt] (3.5999115296852535,6.989696769890026) -- (3.5999962965560313,4.997891892069423) -- (3.5999996411440085,0.999343778192253);
\draw[color=black!90, thick,rounded corners=4pt] (2.4000000324868362,0.999802555824115) -- (2.4000000244659043,6.999828655070412);
\draw[color=black!90, thick,rounded corners=4pt] (3.4229899529306573,7.4255343696104) arc (45.17180874636795:225.7817582678152:0.6); 
\draw[color=black!90, thick,rounded corners=4pt] (5.599999642053717,5.000655389511118) arc (0.06258510065257994:44.82819125363204:0.6); 
\draw[color=black!90, thick,rounded corners=4pt] (4.400035491598568,1.0065259986690296) arc (179.37680074391378:359.93741489934746:0.6); 
\draw[color=black!90, thick,rounded corners=4pt] (3.4229899529306573,7.4255343696104) -- (5.4255343696104,5.422989952930657);
\draw[color=black!90, thick,rounded corners=4pt] (5.599999642053717,5.000655389511118) -- (5.599999642053717,0.9993446104888825);
\draw[color=black!90, thick,rounded corners=4pt] (4.400035491598568,1.0065259986690296) -- (4.44589310811197,4.769857538984682) -- (2.5815640102268604,6.569986834547391);
\draw[color=black!90, thick,rounded corners=4pt] (1.0103032301099744,5.5999115296852535) arc (89.01606563788926:270.016362235733:0.6); 
\draw[color=black!90, thick,rounded corners=4pt] (7.0001974441758845,4.400000032486837) arc (-89.98114546971368:89.93733542082629:0.6); 
\draw[color=black!90, thick,rounded corners=4pt] (1.0103032301099744,5.5999115296852535) -- (3.002108107930577,5.599996296556031) -- (7.0006562218077475,5.5999996411440085);
\draw[color=black!90, thick,rounded corners=4pt] (7.0001974441758845,4.400000032486837) -- (1.000171344929588,4.400000024465904);

\node at (10, 4) {$\rightarrow$};

\node (e) at (16.0,6.5) [vertexY] {$e$};
\node (e1) at (13.0,4.0) [vertexY] {$e_1$};
\node (e2) at (16.0,4.0) [vertexY] {$e_2$};
\node (e3) at (19.0,4.0) [vertexY] {$e_3$};

\draw [line width=0.04cm] (e) -- (e1);
\draw [line width=0.04cm] (e) -- (e2);
\draw [line width=0.04cm] (e) -- (e3);
\draw [line width=0.04cm] (e2) -- (e3);
\draw [line width=0.04cm] (e1) -- (12.5,2.4);
\draw [line width=0.04cm] (e1) -- (13.5,2.4);
\draw [dotted, line width=0.04cm] (e2) -- (15.5,2.4);
\draw [dotted, line width=0.04cm] (e2) -- (16.5,2.4);
\draw [dotted, line width=0.04cm] (e3) -- (18.5,2.4);
\draw [dotted, line width=0.04cm] (e3) -- (19.5,2.4);

\node at (14, 5.8) {$y_1$};
\node at (15.4, 5.0) {$y_2$};
\node at (18, 5.8) {$y_3$};
\node at (17.5, 3.4) {$x$};

\node at (12.2, 3.0) {$a$};
\node at (13.8, 3.0) {$b$};
 \end{tikzpicture}
\end{center}
\vskip -0.75 cm
\caption{The transformation creating a quadruple.} \label{quadruple}
\end{figure}

\medskip
Suppose that $x$ is a vertex of degree~$3$ in $H$, and let $x'$ be the vertex duplicated from $x$ when constructing $H$. Adopting our earlier notation, let $e_1$, $e_2$ and $e_3$ be the three edges incident with $x$. Further, let $e$ be the edge in the $H_4$-component of $H-x$, and let $y$ be the vertex of degree~$1$ in $e$. Let $y_i$ be the vertex common to $e$ and $e_i$ for $i \in [3]$, and so $y_i = (ee_i)$ and $V(e) = \{y,y_1,y_2,y_3\}$. We note that in the graph $G$, which is the dual of $H'$, the vertex $e$ has degree~$3$ and is adjacent to the vertices $e_1$, $e_2$ and $e_3$. Further, we note that in the graph $G$, the vertex $e_1$ has degree~$3$, while the vertices $e_2$ and $e_3$ are adjacent and have degree at most~$4$. Further, the edge $y_i$ in $G$ is the edge $ee_i$, while the edge $x$ in $G$ is the edge $e_2e_3$. The vertex $e_1$ is adjacent in $G$ to neither $e_2$ nor $e_3$. This set of four vertices $\{e,e_1,e_2,e_3\}$ in the graph $G$ we call a \emph{quadruple} in $G$. We illustrate this quadruple in $G$ in Figure~\ref{quadruple}. We denote the set of (vertex-disjoint) quadruples in $G$ by $Q$.

We shall need the following additional property about the dual $G$ of the hypergraph $H'$.

\ClaimX{U}
If $G$ is the dual of the hypergraph $H'$, then
\[
\begin{array}{cl}
& \, 45\tau(H') \le 6n(H') + 13m(H') - 6|Q| \\ \Leftrightarrow & \, 45\alpha'(G) \ge 8n(G) + 6m(G) + 6|Q|.
\end{array}
\]

\ProofClaimX{U} Recall that $n(H') = n_1(H') + n_2(H')$ and $4m(H') = 2n_2(H') + n_1(H')$. The following holds by Claim~T:
\[
\begin{array}{crcl}
& 45\tau(H') & \le & 6n(H') + 13m(H') - 6|Q| \\
\Leftrightarrow & 45(m(H') - \alpha'(G)) & \le & 6n(H') + 13m(H') - 6|Q|  \\
\Leftrightarrow & 45 \alpha'(G) & \ge & 32m(H') - 6n(H') + 6|Q|  \\
\Leftrightarrow & 45 \alpha'(G) & \ge & (16n_2(H') + 8n_1(H')) - (6n_2(H') + 6n_1(H')) + 6|Q| \\
\Leftrightarrow & 45 \alpha'(G) & \ge & 10n_2(H') + 2n_1(H') + 6|Q| \\
\Leftrightarrow & 45 \alpha'(G) & \ge & 8m(H') + 6n_2(H') + 6|Q| \\
\Leftrightarrow & 45 \alpha'(G) & \ge & 8n(G) + 6m(G) + 6|Q|. \\
\end{array}
\]
This completes the proof of Claim~U.~\smallqed

\medskip
Let $S$ be a set of vertices in $G$ such that $(n(G) +|S|-\oc(G-S))/2$ is minimum. By the Tutte-Berge Formula,
\begin{equation}
\alpha'(G) = \frac{1}{2} \left( n(G) +|S|-\oc(G-S) \right).
 \label{Eq2}
\end{equation}
We now consider two cases, depending on whether $S = \emptyset$ or $S \ne \emptyset$.

\medskip
\ClaimX{V} If $S \ne \emptyset$, then $45\tau(H') \le 6n(H') + 13m(H') - 6|Q|$.

\ProofClaimX{V} Suppose that $S \ne \emptyset$. For $i \ge 1$, let $n_i(G - S)$ denote the number of components on $G - S$ of order~$i$. Let $n_5^1(G - S)$ be the number of components of $G - S$ isomorphic to $K_5 - e$ and let $n_5^2(G - S)$ denote all remaining components of $G - S$ on five vertices (with at most eight edges), and so $n_5(G-S) =  n_5^1(G - S) + n_5^2(G - S)$. For notational convenience, let $n = n(G)$, $m = m(G)$, $n_5^1 = n_5^1(G - S)$, $n_5^2 = n_5^2(G - S)$, and $n_i = n_i(G-S)$ for $i \ge 1$. Let $\ZZ^+$ denote the set of all positive integers, and let $\ZZ^+_{\even}$ and $\ZZ^+_{\odd}$ denote the set of all even and odd integers, respectively, in $\ZZ^+$. Further for a fixed $j \in \ZZ^+$, let $\ZZ_{\ge j} = \{ i \in \ZZ \mid i \ge j\}$, $\ZZ_{\even}^j = \{ i \in \ZZ_{\ge j} \mid i \,\, \mbox{even} \}$, and $\ZZ_{\odd}^j = \{ i \in \ZZ_{\ge j} \mid i \,\, \mbox{odd} \}$. We note that
\begin{equation}
n = |S| + \sum_{i \in \ZZ^+} i \cdot n_i. \label{Eq2a}
\end{equation}
By Equation~(\ref{Eq2}) and Equation~(\ref{Eq2a}), and since
\[
\oc(G-S) = \sum_{i \in \ZZ^+_{\odd}} n_i,
\]
we have the equation
\begin{equation}
45\alpha'(G) = 45|S| + \frac{45}{2} \left( (\sum_{i \in \ZZ_{\odd}^3} (i-1) \cdot n_i)  + \sum_{i \in \ZZ_{\even}^2} i \cdot n_i  \right). \label{Eq3}
\end{equation}

\ClaimX{V.1} $\displaystyle{m \le 4|S| + n_2 + 3n_3 + 6n_4 + 9n_5^1 + 8n_5^2 + \sum_{i \in \ZZ^6} (2i-1) n_i - |Q|}$.

\ProofClaimX{V.1} Since $G$ is connected and $\Delta(G) \le 4$, we note that if $F$ is a component of $G-S$ of order~$i$, then $m(F) \le 2i - 1$. Further, every component of $G - S$ of order~$5$ is either isomorphic to $K_5 - e$ or contains at most eight edges, while every component of $G - S$ of order~$2$, $3$ and $4$ contains at most~$1$,~$3$ and~$6$ edges, respectively. The above observations imply that  \[
m \le 4|S| + n_2 + 3n_3 + 6n_4 + 9n_5^1 + 8n_5^2 + \sum_{i \in \ZZ^6} (2i-1) n_i.
\]

We show next that each quadruple in the graph $G$ decreases the count on the right hand side expression of the above inequality by at least~$1$. Adopting our earlier notation, consider a quadruple $\{e,e_1,e_2,e_3\}$. Recall that the vertices $e$ and $e_1$ both have degree~$3$ in $G$, and there is no vertex in $G$ that is adjacent to both $e$ and $e_1$. Further, recall that the vertices $e_2$ and $e_3$ are adjacent in $G$. If $e$ or $e_1$ or if both $e_2$ and $e_3$ belong to the set $S$, then the quadruple decreases the count $4|S|$ by at least~$1$.  Hence, we may assume that $e$, $e_1$ and $e_2$ all belong to a component, $C$ say, of $G - S$. In particular, we note that $C$ has order at least~$3$. Abusing notation, we say that the component $C$ contains the quadruple $\{e,e_1,e_2,e_3\}$, although possibly the vertex $e_3$ may belong to $S$. Since no vertex in $G$ is adjacent to both $e$ and $e_1$, we note that if the component $C$ has order~$3$,~$4$ or $5$, then it contains at most~$2$,~$4$ and~$7$ edges, respectively. Further, since we define a component to contain a quadruple if it contains at least three of the four vertices in the quadruple, we note  in this case when the component $C$ has order at most~$5$ that it contains exactly one quadruple. Further, this quadruple decreases the count $3n_3 + 6n_4 + 9n_5^1 + 8n_5^2$ by at least~$1$.

It remains for us to consider a component $F$ of $G-S$ of order~$i \ge 6$ that contains $q$ quadruples, and to show that these $q$ quadruples decrease the count $2i - 1$ by at least~$q$. We note that each quadruple contains a pair of adjacent vertices of degree~$3$ in $G$. Further, at least one vertex $v$ in $F$ is joined to at least one vertex of $S$ in $G$, implying that $d_F(v) < d_G(v)$. These observations imply that $2m(F) = \sum_{v \in V(F)} d_F(v) \le 4n(F) - 2q - 1 = 4i - 2q - 1$, and therefore that $m(F) \le 2i - q - 1$. Hence, these $q$ quadruples contained in $F$ combined decrease the count $2i - 1$ by at least~$q$. This completes the proof of Claim~V.1.~\smallqed

By Claim~V.1 and by Equation~(\ref{Eq2a}), we have
\begin{eqnarray}
8n + 6m + 6|Q| & \le & 32|S| + 8n_1 + 22n_2 + 42n_3 + 68n_4 \nonumber \1 \\
& & \, + \, 88n^2_5 + 94n^1_5 + \sum_{i \in \ZZ^6} (20i-6) n_i.
 \label{Eq4}
\end{eqnarray}
Let
\[
\Sigma_{\even} = \sum_{i \in \ZZ_{\even}^6} \left( \frac{5}{2} i + 6 \right) \cdot n_i \hspace*{0.7cm} \mbox{and} \hspace*{0.7cm} \Sigma_{\odd} = \sum_{i \in \ZZ_{\odd}^7} \frac{1}{2} \left(  5i - 33 \right) \cdot n_i.
\]

We note that every (odd) component in $G$ isomorphic to $K_1$ corresponds to a subhypergraph $H_4$ in $H'$, while every (odd) component in $G$ isomorphic to $K_5 - e$ corresponds to a subhypergraph $H_{11}$ in $H'$. Hence the odd components of $G$ isomorphic to $K_1$ or isomorphic to $K_5 - e$ correspond to a special $H'$-set, $X$ say, where $|X| = |X_4| + |X_{11}|$, $|X_4| = n_1$ and $|X_{11}| = n_5^1$. Further, the set $S$ of vertices in $G$ correspond to the set $E^*(X)$ of edges in $H'$, and so $|E^*(X)| \le |S|$. Thus,
\[
\defic_{H'}(X) = 8|X_4| + 4|X_{11}| - 13|E^*(X)| \ge 8n_1 + 4n_5^1 - 13|S|.
\]
By Claim~S(b), $\defic_{H'}(X) \le \defic(H') = 0$, and therefore we have that
\begin{equation}
13|S| \ge 8n_1 + 4n_5^1.
 \label{Eq5}
\end{equation}
By Equation~(\ref{Eq3}), and by Inequalities~(\ref{Eq4}) and~(\ref{Eq5}), and noting that $\Sigma_{\even} \ge 0$ and $\Sigma_{\odd} \ge 0$, the following now holds.
\[
\begin{array}{lcl}
45\alpha'(G) & \stackrel{(\ref{Eq3})}{=} & \displaystyle{ 32|S| + 8n_1 + 22n_2 + 42n_3 + 68n_4 + 88n^2_5 + 94n^1_5 + \sum_{i \in \ZZ^6} (20i-6) n_i } \\
& & \displaystyle{ + 13|S| - 8n_1 + 23n_2 + 3n_3 + 22n_4 + 2n^2_5 - 4n_5^1 +  \Sigma_{\even} + \Sigma_{\odd} } \2 \\
& \stackrel{(\ref{Eq4})}{\ge} & (8n + 6m + 6|Q|) + (13|S| - 8n_1 - 4n_5^1) \1 \\
& \stackrel{(\ref{Eq5})}{\ge} & 8n + 6m + 6|Q|.
\end{array}
\]
Claim~V now follows from Claim~U.~\smallqed


\medskip
\ClaimX{W} If $S = \emptyset$, then $45\tau(H) \le 6n(H) + 13m(H) - 6|Q|$.

\ProofClaimX{W} Suppose that $S = \emptyset$. Then, $\alpha'(G) = (n(G) - \oc(G))/2$.
Since $G$ is connected by Claim~T, we have the following.

\[
\alpha'(G) = \left\{ \begin{array}{lcl}
\frac{1}{2}n(G) & \hspace{0.5cm} & \mbox{if $n(G)$ is even} \2 \\
\frac{1}{2}(n(G) -1) & & \mbox{if $n(G)$ is odd}. \end{array} \right.
\]

By Claim~T(b), $\Delta(G) \le 4$. As every quadruple in $G$ contains two vertices of degree~$3$,
\[
2m(G) = \sum_{v \in G} d_G(v) \le 4n(G) - 2|Q|,
\]
implying that $12n(G) \ge 6m(G) + 6|Q|$. If $n(G)$ is even, then $\alpha'(G) = n(G)/2$, and so
\[
45 \alpha'(G) = \frac{45}{2}n(G) > 8n(G) + 12n(G) \ge 8n(G) + 6m(G) + 6|Q|.
\]
This completes the case when $n(G)$ is even by Claim~U. Suppose next that $n(G)$ is odd. In this case $45\alpha'(G) = 45(n(G)-1)/2$, and so
\[
90\alpha'(G) = 45n(G)-45 = 21n(G) + 24n(G) - 45 \ge 21n(G) + 12m(G) + 12|Q| -45.
\]

If $5n(G) \ge 45$, then $90\alpha'(G) \ge 16n(G)+12m(G)+12|Q|$, which completes the proof by Claim~U. We may therefore assume that $5n(G)<45$, implying that $n(G) \in \{1,3,5,7\}$. Since every quadruple contains four vertices, we must therefore have $|Q| \le 1$.

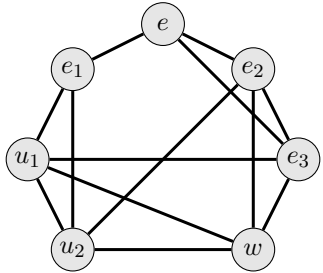
\begin{figure}[htb]
\begin{center}
\tikzstyle{vertexX}=[circle,draw, fill=black!10, minimum size=18pt, scale=0.9, inner sep=0.1pt]
\begin{tikzpicture}[scale=0.6]
 \draw (0,0) node {\mbox{ }};
\draw (4,3) node {\mbox{ }};
\node (e) at (4.0,6.0) [vertexX] {$e$};
\node (e1) at (2.0,5.0) [vertexX] {$e_1$};
\node (e2) at (6.0,5.0) [vertexX] {$e_2$};
\node (u1) at (1.0,3.0) [vertexX] {$u_1$};
\node (e3) at (7.0,3.0) [vertexX] {$e_3$};
\node (u2) at (2.0,1.0) [vertexX] {$u_2$};
\node (w) at (6.0,1.0) [vertexX] {$w$};

\draw [line width=0.04cm] (e) -- (e1);
\draw [line width=0.04cm] (e) -- (e2);
\draw [line width=0.04cm] (e) -- (e3);
\draw [line width=0.04cm] (e1) -- (u1);
\draw [line width=0.04cm] (e1) -- (u2);
\draw [line width=0.04cm] (e2) -- (e3);
\draw [line width=0.04cm] (e2) -- (w);
\draw [line width=0.04cm] (e2) -- (u2);
\draw [line width=0.04cm] (u1) -- (u2);
\draw [line width=0.04cm] (u1) -- (e3);
\draw [line width=0.04cm] (u1) -- (w);
\draw [line width=0.04cm] (e3) -- (w);
\draw [line width=0.04cm] (u2) -- (w);
 \end{tikzpicture}
\end{center}
\vskip -1 cm
\caption{The graph $G$ if $|Q|=1$.} \label{FigG}
\end{figure}

We first consider the case when $|Q|=1$.  In this case $n(G) \ge 6$, as the quadruple contains four vertices and one vertex (called $e_1$ in the definition of a quadruple) has two neighbours outside the quadruple. Therefore, $n(G)=7$. Since two vertices in the quadruple have degree at most~$3$, $2m(G) = \sum_{v \in V(G)} d(v) \le 4n(G)-2 = 26$, and so $m(G) \le 13$.  If $m(G) \le 12$, then $45\alpha'(G) = 45 \cdot 3 > 8\cdot 7 + 6\cdot 12 + 6 \ge  8n(G) +  6 m(G) + 6|Q|$, and the desired result follows from Claim~U.  Therefore, we may assume that $m(G)=13$, for otherwise the case is complete. In this case, all vertices in $G$ have degree~$4$ except
for two vertices in the quadruple which have degree~$3$. Let $\{e,e_1,e_2,e_3\}$ be the vertices in the quadruple in $G$, such that $d(e)=d(e_1)=3$ and $e$, $e_2$ and $e_3$ form a $3$-cycle in $G$. Define $u_1$ and $u_2$ such that $N(e_1) = \{e,u_1,u_2\}$ and define $w$, such that $V(G)=\{e,e_1,e_2,e_3,u_1,u_2,w\}$.
As $d(w)=4$ and $w$ is not adjacent to $e$ or $e_1$ we have $N(w) = \{e_2,e_3,u_1,u_2\}$. Therefore, $e_2$ is adjacent to $e$, $e_3$ and $w$.  Its fourth neighbour is either $u_1$ or $u_2$. Renaming $u_1$ and $u_2$ if necessary, we may assume that $N(e_2)=\{e,e_3,w,u_2\}$.  This implies that $u_1$ must be adjacent to $u_2$ and $e_3$ and $G$ is the graph shown in Figure~\ref{FigG}.

If we draw the corresponding hypergraph whose dual is the graph $G$, we note that it is obtained by duplicating the degree-$3$ vertex in $H_{14,3}$. However, $H$ is not equal to $H_{14,3}$ by Claim~C, implying that $G$ cannot be the graph in Figure~\ref{FigG}, a contradiction. This completes the case when $|Q|=1$. Therefore, we may assume that $|Q|=0$.

If $n(G) = 1$, then $H = H_4$, contradicting Claim~C. Hence, $n(G) \in \{3,5,7\}$.
Suppose that $n(G) = 3$.
Then, $\alpha'(G) =  1$ and $m(G) \le 3$.
In this case, $8n(G) + 6m(G) \le 8 \cdot 3 + 6 \cdot 3 = 42 < 45 = 45\alpha'(G)$, which by Claim~U completes the proof.
Hence we may assume that $n(G) = 5$ or $n(G) = 7$.

Suppose that $n(G) = 5$. Then, $\alpha'(G) =  2$ and by Claim~T(b), $m(G) \le 10$.
If $m(G) = 10$, then $G = K_5$. In this case, $H$ is a $4$-uniform $2$-regular linear intersecting hypergraph.
However, $H_{10}$ is the unique such hypergraph as shown, for example, in~\cite{DoHe14,HeYe17}.
Thus if $m(G) = 10$, then $H = H_{10}$, contradicting Claim~C.
Hence, $m(G) \le 9$.
If $m(G) = 9$, then $G = K_5 - e$, where $e$ denotes an arbitrary edge in $K_5$. In this case, $H = H_{11}$, contradicting Claim~C.
Hence, $m(G) \le 8$. Thus, $8n(G) + 6m(G) \le 8 \cdot 5 + 6 \cdot 8 = 88 < 90 = 45\alpha'(G)$, which by Claim~U completes the proof in this case.

Finally suppose that $n(G) = 7$. Then, $\alpha'(G) =  3$ and by Claim~T(b), $m(G) \le 14$. Suppose that $m(G) = 14$. Then, $G$ is a $4$-regular graph of order~$7$. Equivalently, the complement, $\barG$, of $G$ is a $2$-regular graph of order~$7$.
If $\barG = C_3 \cup C_4$, then $H = H_{14,2}$. If $\barG = C_7$, then $H = H_{14,4}$. Both cases contradict Claim~C.
Hence, $m(G) \le 13$. Thus, $8n(G) + 6m(G) \le 8 \cdot 7 + 6 \cdot 13 = 134 < 135 = 45\alpha'(G)$,  which by Claim~U completes the proof.~\smallqed

\medskip
Recall that $n(H') = n(H) + |Q|$ and $m(H') = m(H)$. By Claim~S(a), Claim~V and Claim~W, $45\tau(H) = 45\tau(H') \le 6n(H') + 13m(H') - 6|Q| = 6n(H) + 13m(H)$, a contradiction. This completes the proof of Theorem~\ref{t:thm_key}.~\qed

\subsection{Proof of Theorem~\ref{t:mainZ}}
\label{S:mainZ}

We are finally in a position to present a proof of our main result, namely Theorem~\ref{t:mainZ}. Recall its statement.

\noindent \textbf{Theorem~\ref{t:mainZ}}. \emph{If $H$ is a $4$-uniform, linear hypergraph on $n$ vertices with $m$ edges, then \\ $\tau(H) \le \frac{1}{5}(n + m)$.
}

\noindent
\textbf{Proof of Theorem~\ref{t:mainZ}.} Let $H$ be a $4$-uniform, linear hypergraph on $n$ vertices with $m$ edges. We show that $\tau(H) \le (n+m)/5$. We proceed by induction on~$n$. If $n = 4$, then $H$ consists of a single edge, and $\tau(H) = 1 = (n+m)/5$. Let $n \ge 5$ and suppose that the result holds for all $4$-uniform, linear hypergraphs on fewer than $n$ vertices. Let $H$ be a $4$-uniform, linear hypergraph on $n$ vertices with $m$ edges.

Suppose that $\Delta(H) \ge 4$. Let $v$ be a vertex of maximum degree in $H$, and consider the $4$-uniform, linear hypergraph $H' = H - v$ on $n'$ vertices with $m'$ edges. We note that $n' = n - 1$ and $m' = m - \Delta(H) \le m - 4$. Every transversal in $H'$ can be extended to a transversal in $H$ by adding to it the vertex~$v$. Hence, applying the inductive hypothesis to $H'$, we have that
$\tau(H) \le \tau(H') + 1 \le (n'+m')/5 + 1 \le (n+m-5)/5 + 1  = (n+m)/5$. Hence, we may assume that $\Delta(H) \le 3$, for otherwise the desired result follows. With this assumption, we note that $4m \le 3n$. Applying Theorem~\ref{t:thm_key} to the hypergraph $H$, we have
\[
45\tau(H) \le 6n(H) + 13m(H) + \defic(H).
\]

If $\defic(H) = 0$, then $45\tau(H) \le 6n(H) + 13m(H) = (9n + 9m) + (4m - 3n) \le 9(n+m)$, and so $\tau(H) \le (n+m)/5$. Hence, we may assume that $\defic(H) > 0$, for otherwise the desired result follows. Among all special non-empty $H$-sets, let $X$ be chosen so that $|E^*(X)| - |X|$ is minimum. We note that since $\defic(H) > 0$, $|E^*(X)| - |X| < 0$. As in Section~\ref{S:defic}, we associate with the set $X$ a bipartite graph, $G_X$, with partite sets $X$ and $E^*(X)$, where an edge joins $e \in E^*(X)$ and $H' \in X$ in $G_X$ if and only if the edge $e$ intersects the subhypergraph $H'$ of $X$ in~$H$. Suppose that there is no matching in $G_X$ that matches $E^*(X)$ to a subset of $X$.  By Hall's Theorem, there is a nonempty subset $S \subseteq E^*(X)$ such that $|N_{G_X}(S)|<|S|$. We now consider the special $H$-set, $X' = X \setminus N_{G_X}(S)$, and note that $|X'| = |X| - |N_{G_X}(S)| > |E^*(X)| - |S| \ge 0$ and $|E^*(X')| = |E^*(X)| - |S|$. Thus, $X'$ is a special non-empty $H$-set satisfying

\[
\begin{array}{lcl}
|E^*(X')| - |X'| & =  & (|E^*(X)| - |S|) - (|X| - |N_{G_X}(S)|) \\
& = & (|E^*(X)| - |X|) + (|N_{G_X}(S)| - |S|) \\
& < & |E^*(X)| - |X|,
\end{array}
\]
contradicting our choice of the special $H$-set $X$. Hence, there exists a matching in $G_X$ that matches $E^*(X)$ to a subset of $X$. By Observation~\ref{property:special}(g), there exists a minimum $X$-transversal, $T_X$, that intersects every edge in $E^*(X)$. By Observation~\ref{property:special}, every special hypergraph $F$ satisfies $\tau(F) \le (n(F) + m(F))/5$. Hence, letting
\[
n(X) = \sum_{F \in X} n(F) \hspace*{0.5cm} \mbox{and} \hspace*{0.5cm} m(X) = \sum_{F \in X} m(F),
\]
we note that
\[
|T_X| = \sum_{F \in X} \tau(F) \le \sum_{F \in X} \frac{n(F) + m(F)}{5} = \frac{n(X) + m(X)}{5}.
\]

We now consider the $4$-uniform, linear hypergraph $H' = H - V(X)$ on $n'$ vertices with $m'$ edges. We note that $n' = n - n(X)$ and $m' = m - m(X) - |E^*(X)| \le m - m(X)$. Every transversal in $H'$ can be extended to a transversal in $H$ by adding to it the set $T_X$. Hence, applying the inductive hypothesis to $H'$, we have that

\[
\begin{array}{lcl}
\tau(H) & \le  & \tau(H') + |T_X| \1 \\
& \le & \frac{1}{5}(n'+m') + |T_X| \1 \\
& \le & \frac{1}{5}(n+m) - \frac{1}{5}(n(X)+m(X)) +  |T_X| \1 \\
& \le & \frac{1}{5}(n+m).
\end{array}
\]
This completes the proof of Theorem~\ref{t:mainZ}.~\qed

\section{Proof of Theorem~\ref{t:linearDeg2}}
\label{S:linearDeg2}

In this section, we present a proof of  Theorem~\ref{t:linearDeg2}. We first prove the following lemma. Let $c(H)$ denote the number of components of a hypergraph $H$. Recall that if $X$ is a special $H$-set, we write $E^*(X)$ to denote the set  $E_H^*(X)$ if the hypergraph $H$ is clear from context.

\begin{lem}
If $H$ is a $4$-uniform, linear hypergraph and $X$ is a special $H$-set, then $3|E_H^*(X)| \ge |X| - c(H)$. \label{l:specialHset}
\end{lem}
\proof We proceed by induction on $|E_H^*(X)| = k \ge 0$. If $k = 0$, then $c(H) \ge |X|$, and so $3|E^*(X)| = 0 \ge |X| - c(H)$. This establishes the base case. Suppose $k \ge 1$ and the result holds for special $H$-sets, $X$, such that $|E_H^*(X)| < k$. Let $X$ be a special $H$-set satisfying $|E_H^*(X)| = k$. Let $e \in E_H^*(X)$ and consider the hypergraph $H' = H - e$. We note that $H'$ is a $4$-uniform, linear hypergraph, and that $c(H') \le c(H) + 3$. Further, the set $X$ is a special $H'$-set satisfying $|E_{H'}^*(X)| = k - 1$. Applying the inductive hypothesis to the hypergraph $H' \in \cH_4$ and to the special $H'$-set, $X$, we have $3(|E_H^*(X)| - 1) = 3|E_{H'}^*(X)| \ge |X| - c(H') \ge |X| - (c(H) + 3)$, implying that $3|E_H^*(X)| \ge |X| - c(H)$.~\qed

\medskip
We are now in a position to present a proof of  Theorem~\ref{t:linearDeg2}. Recall its statement, where $\cF$ is the family defined in Section~\ref{S:Applic3}.

\noindent \textbf{Theorem~\ref{t:linearDeg2}}. \emph{Let $H \ne H_{10}$ be a $4$-uniform, connected, linear hypergraph with maximum degree~$\Delta(H) \le 2$ on $n$ vertices with $m$ edges. Then, $\tau(H) \le \frac{3}{16}(n + m) + \frac{1}{16}$, with equality if and only if $H \in \{H_{14,5}, H_{14,6}\}$ or $H \in \cF$.
}

\noindent
\textbf{Proof of Theorem~\ref{t:linearDeg2}.} Let $H \ne H_{10}$ be a $4$-uniform, connected, linear hypergraph with maximum degree~$\Delta(H) = 2$ on $n$ vertices with $m$ edges. Suppose firstly that $H$ is a special hypergraph. By assumption, $H \ne H_{10}$. Since $\Delta(H) \le 2$, we note that $H \in \{H_4,H_{11},H_{14,5}, H_{14,6}\}$. If $H = H_{11}$, then by Observation~\ref{property:special}(c),
$\tau(H) = \frac{3}{16}(n + m)$. If $H \in \{H_4,H_{14,5}, H_{14,6}\}$, then by Observation~\ref{property:special}(a) and~\ref{property:special}(d), $\tau(H) = \frac{3}{16}(n + m) + \frac{1}{16}$. Hence, we may assume that $H$ is not a special hypergraph, for otherwise the desired result holds, noting that $H_4 \in \cF$.

Since $\Delta(H) \le 2$, we observe that $m \le \frac{1}{2}n$. By Theorem~\ref{t:thm_key}, $45\tau(H) \le 6n + 13m + \defic(H)$. Suppose that $\defic(H) = 0$. Then, since $0 \le \frac{1}{2}n - m$, we have
\[
\begin{array}{lcl}
45\tau(H) & \le & 6n + 13m \\
& \le & 6n + 13m + \frac{14}{3}(\frac{1}{2}n - m) \\
& = & (6 + \frac{14}{6})n + (13 - \frac{14}{3})m \\
& = & \frac{25}{3}(n+m),
\end{array}
\]
or, equivalently, $\tau(H) \le \frac{5}{27}(n+m) <  \frac{3}{16}(n + m)$. Hence, we may assume that $\defic(H) > 0$, for otherwise the desired result follows. Let $X$ be a special $H$-set such that $\defic(H) = \defic_H(X)$.
If $H_{10}$ belongs to $X$, then, since $\Delta(H) \le 2$ and $H$ is connected, $H = H_{10}$, a contradiction. If $H_{14,5}$ or $H_{14,6}$ belong to $X$, then, analogously, $H \in \{H_{14,5}, H_{14,6}\}$, contradicting our assumption that $H$ is not a special hypergraph.  Thus, if $F \in X$, then $F \in \{H_4,H_{11}\}$, noting that $\Delta(H) \le 2$. Recall that if $F$ is a hypergraph, we denote by $n_1(F)$ the number of vertices of degree~$1$ in $H$. Let
\[
n_1(X) = \sum_{F \in X} n_1(F),
\]
and note that $n_1(H) \ge n_1(X) - 4|E^*(X)|$, since every edge in $E^*(X)$ contains at most four vertices whose degree is~$1$ in some subhypergraph $F \in X$. Since $\Delta(H) \le 2$ and $H$ is $4$-uniform, we note that $4m = 2n - n_1(H)$, or, equivalently, $n_1(H) = 2n - 4m$.
Let $\beta = \frac{73}{64}$. If $F = H_4$, then $n_1(F) = 4$ and $\defic(F) = 8 = (8 - 4\beta) + 4\beta = (8 - 4\beta) + n_1(F) \cdot \beta$. If $F = H_{11}$, then $n_1(F) = 2$ and  $\defic(F) = 4 < 5.71875 = 8 - 2\beta =  (8 - 4\beta) + n_1(F) \cdot \beta$. Hence, if $F \in X$, then  $\defic(F) \le (8 - 4\beta) + n_1(F) \cdot \beta$, with strict inequality if $F = H_{11}$. Therefore,
\[
\begin{array}{lcl}
\defic(H) & = &  8|X_4| + 4|X_{11}| - 13|E^*(X)| \1 \\
& = & \displaystyle{ \left( \sum_{F \in X} \defic(F) \right) - 13|E^*(X)| } \1 \\
& \le & (8 - 4\beta)|X| + n_1(X) \cdot \beta - 13|E^*(X)|,
\end{array}
\]
with strict inequality if $X \ne X_4$. By Lemma~\ref{l:specialHset}, $|E^*(X)| \ge \frac{1}{3}(|X| - 1)$. We also note that $n \ge 4|X_4| + 11|X_{11}| \ge 4|X|$, and so $|X| \le \frac{n}{4}$. Thus by our previous observations, \[
\begin{array}{lcl}
45\tau(H) & \le & 6n + 13m + \defic(H) \5 \\
& \le &  6n + 13m + (8 - 4\beta)|X| + n_1(X) \cdot \beta - 13|E^*(X)| \5 \\
& \le &  6n + 13m + (8 - 4\beta)|X| + (n_1(H) + 4|E^*(X)|) \cdot \beta - 13|E^*(X)| \5 \\
& = &  6n + 13m + (8 - 4\beta)|X| + n_1(H)\cdot \beta - |E^*(X)|(13 - 4\beta) \5 \\
& \le & 6n + 13m + (8 - 4\beta)|X| + n_1(H)\cdot \beta - \frac{1}{3}(|X| - 1)(13 - 4\beta) \5 \\
& = & 6n + 13m + (8 - 4\beta - \frac{1}{3}(13 - 4\beta))|X| + (2n - 4m) \cdot \beta + \frac{1}{3}(13 - 4\beta) \5 \\
& \le & 6n + 13m + (8 - 4\beta - \frac{1}{3}(13 - 4\beta))\frac{n}{4} + (2n - 4m) \cdot \beta + \frac{1}{3}(13 - 4\beta) \5 \\
& = & (\frac{83}{12} + \frac{4}{3}\beta) n + (13 - 4 \beta)m  + \frac{1}{3}(13 - 4\beta) \5 \\
& = &  (13 - 4 \beta)(n+m)  + \frac{1}{3}(13 - 4\beta) \5 \\
& = &  \frac{135}{16}(n+m)  + \frac{135}{48}, \5 \\
\end{array}
\]
or, equivalently, $\tau(H) \le \frac{3}{16}(n + m) + \frac{1}{16}$. This establishes the desired upper bound.

Recall that $H$ is a $4$-uniform, connected, linear hypergraph with maximum degree at most~$2$. Suppose that $\tau(H) = \frac{3}{16}(n + m) + \frac{1}{16}$. Then we must have equality throughout the above inequality chain. This implies that $X = X_4$, $V(H) = V(X)$, $n = 4|X|$, $E(H) = E(X) \cup E^*(X)$, and $|E^*(X)| = \frac{1}{3}(|X| - 1)$. We show by induction on $n \ge 4$ that these conditions imply that $H \in \cF$. When $n = 4$, $|X| = 1$ and $|E^*(X)| = 0$, and so $H = H_4 \in \cF$. This establishes the base case. Suppose that $n > 4$. Thus, $|X| \ge 2$ and $n = 4|X| \ge 8$. We now consider the bipartite graph, $G_X$, with partite sets $X$ and $E^*(X)$, where an edge joins $e \in E^*(X)$ and $H' \in X$ in $G_X$ if and only if the edge $e$ intersects the subhypergraph $H'$ of $X$ in $H$. Since $H$ is $4$-uniform and linear, each vertex in $E^*(X)$ has degree~$4$ in $G_X$. Let $n_1 = n_1(G)$, and so $n_1$ is the number of vertices of degree~$1$ in $G$. Counting the edges in $G$, we note that $\frac{4}{3}(|X| - 1) = 4|E^*(X)| = m(G) \ge n_1 + 2(|X| - n_1)$, implying that $n_1 \ge \frac{1}{3}(2|X| + 4)$. By the Pigeonhole Principle, there is a vertex of $E^*(X)$ adjacent in $G$ to at least
\[
\frac{n_1}{|E^*(X)|} = \frac{ \left(\frac{2|X|+4}{3} \right) }{ \left(\frac{|X|-1}{3} \right) } = \frac{2|X| + 4}{|X| - 1} = 2 + \frac{6}{|X|-1}
\]

\noindent
vertices of degree~$1$ (that belong to $X$). Thus, since $|X| \ge 2$ here, some vertex $e \in E^*(X)$ in $G_X$ is adjacent to three vertex of degree~$1$, say $x_1$, $x_2$ and $x_3$. Let $x_4$ be the remaining neighbor of $e$ in $G_X$. We now consider the hypergraph $H'$ obtained from $H$ by deleting the 12 vertices from the three special $H_4$-subhypergraphs, say $F_1$, $F_2$, and $F_3$, corresponding to $x_1$, $x_2$ and $x_3$, respectively, and deleting the hyperedge corresponding to $e$. Since $H$ is connected and linear, so too is $H'$. Let $X' = X \setminus \{F_1,F_2,F_3\}$, and so $|X'| = |X| - 3$.  We note that $|E_{H'}^*(X')| = |E_{H}^*(X)| - 1 = \frac{1}{3}(|X| - 1) - 1 = \frac{1}{3}(|X'| - 1)$. Further, $X' = X_4$, $V(H') = V(X')$, $n' = 4|X'|$, and $E(H') = E(X') \cup E_{H'}^*(X')$. Applying the inductive hypothesis to $H'$, we deduce that $H' \in \cF$. The original hypergraph $H$ can now be reconstructed from $H'$ by adding back the three deleted edges and $12$ deleted vertices in $F_1 \cup F_2 \cup F_3$, and adding back the deleted edge $e$ that contains the vertex $x_4 \in V(H')$ and contains one vertex from each edge in $F_1$, $F_2$ and $F_3$. Thus, $H \in \cF$. This completes the proof of Theorem~\ref{t:linearDeg2}.~\qed

\section{Closing Conjecture}

In this paper, we have shown that Conjecture~\ref{ConjF} fails for large $k$ but holds for small $k$. More precisely, we prove (see Theorem~\ref{t:main4}) that
$5 \le k_{\min} \le 166$, where recall that in Problem~\ref{prob1} we define $k_{\min}$ as the smallest value of~$k$ for which Conjecture~\ref{ConjF} fails. We believe that Conjecture~\ref{ConjF} holds when $k = 5$, and state this new conjecture formally as follows.

\begin{conj}
 \label{t:conj5unif}
$k_{\min} \ge 6$.
\end{conj}



\end{document}